\documentclass[12pt, reqno]{amsart}
\usepackage{ifthen,amsfonts,amsmath,amssymb,epic,eepic}
\usepackage{epsfig,color}

\makeatother

\theoremstyle{definition}

\def\fnum{equation}
\newtheorem{Thm}[\fnum]{Theorem}
\newtheorem{Cor}[\fnum]{Corollary}

\newtheorem{Lem}[\fnum]{Lemma}

\newtheorem{Def}[\fnum]{Definition}
\newtheorem{DefLem}[\fnum]{Definition/Lemma}

\newtheorem{Rem}[\fnum]{Remark}
\newtheorem{Pro}[\fnum]{Proposition}

\newcommand{\cali}{{\it{i}}}

\newcommand{\nn}{{\bf{n}}}

\newcommand{\Hol}{{\text{Hol}}}

\newcommand{\diam}{{\text {diam}}}

\newcommand{\dist}{{\text {dist}}}

\newcommand{\comp}{{\text{comp}}}

\def\ZZ{{\bold Z}}
\def\RR{{\bold R}}

\def\CC{{\bold C }}
\newcommand{\dv}{{\text {div}}}
\newcommand{\e}{{\text {e}}}
\newcommand{\Area}{{\text {Area}}}
\newcommand{\Length}{{\text {Length}}}

\newcommand{\cF}{{\mathcal{F}}}

\newcommand{\cI}{{\mathcal{I}}}

\newcommand{\cL}{{\mathcal{L}}}
\newcommand{\cB}{{\mathcal{B}}}

\newcommand{\cP}{{\mathcal{P}}}
\newcommand{\cS}{{\mathcal{S}}}
\newcommand{\cSt}{{\mathcal{S}_{neck}}}
\newcommand{\cSta}{{\mathcal{S}_{neck}^1}}
\newcommand{\cStb}{{\mathcal{S}_{neck}^2}}

\newcommand{\cSu}{{\mathcal{S}_{ulsc}}}

\newcommand{\barga}{{\Gamma_{Clos}}}
\newcommand{\bargaprime}{{\Gamma'_{Clos}}}

\newcommand{\eqr}[1]{(\ref{#1})}
\newcommand{\an}{{\text{An}}}
\newcommand{\dd}{{\text{d}}}

\newcommand{\cone}{{\bf{C}}}

\begin{document}

\title[Fixed genus]
{The space of embedded minimal surfaces of fixed genus in a
$3$-manifold V; Fixed genus}

\author{Tobias H. Colding}%
\address{Courant Institute of Mathematical Sciences\\
251 Mercer Street\\ New York, NY 10012}
\author{William P. Minicozzi II}%
\address{Department of Mathematics\\
Johns Hopkins University\\
3400 N. Charles St.\\
Baltimore, MD 21218}
\thanks{The authors were partially supported by NSF Grants DMS
0104453 and DMS 0104187}

%%\date{\today}

\email{colding@cims.nyu.edu and minicozz@math.jhu.edu}

\begin{abstract}
This paper is the fifth and final in a series on embedded minimal
surfaces.  Following our earlier papers on  disks, we prove here
two main structure theorems for {\it{non-simply connected}}
embedded minimal surfaces of any given fixed genus.

The first of these asserts that any such surface
{\underline{without}} small necks can be obtained by gluing
together two oppositely--oriented double spiral staircases; see
Figure \ref{f:01a}.

The second gives a pair of pants decomposition of any such surface
when there {\underline{are}} small necks, cutting the surface
along a collection of short curves; see Figure \ref{f:01b}. After
the cutting, we are left with graphical pieces that are defined
over a disk with either one or two sub--disks removed (a
topological disk with two sub--disks removed is called a pair of
pants).

Both of these structures occur as different extremes in the
two-parameter family of minimal surfaces known as the Riemann
examples.

 The results of \cite{CM3}--\cite{CM6} have already been  used by
many authors; see, e.g.,  the surveys  \cite{MeP}, \cite{P},
\cite{Ro} and the introduction in \cite{CM6} for some of these
applications. There is much current research on  minimal surfaces
with infinite topology.  Some of the results of the present paper
were announced previously and  have already been widely used to
study infinite topology minimal surfaces; see, e.g., \cite{MeP},
\cite{MePRs1}, \cite{MePRs2}, \cite{MePRs3}, and \cite{P}.
\end{abstract}

\maketitle

\numberwithin{equation}{section}

\begin{figure}[htbp]
    \center{The two main structure theorems for non-simply connected surfaces:}
        \vskip2mm
    \setlength{\captionindent}{20pt}
    \begin{minipage}[t]{0.5\textwidth}
    \centering\input{gen01a.pstex_t}
    \caption{Absence of necks:  The surface can be obtained by gluing together two
oppositely--oriented double spiral staircases.}
    \label{f:01a}
    \end{minipage}\begin{minipage}[t]{0.5\textwidth}
    \centering\input{gen01b.pstex_t}
    \caption{Presence of necks:  The surface can be decomposed into a collection of pair of pants
    by cutting along short curves.}
    \label{f:01b}
    \end{minipage}
\end{figure}

%%%%%%%%%%%%%
{\small \tableofcontents}
%%%%%%%%%%%%%

\section{Introduction} \label{s:s0}

This paper is the fifth and final in a series where we describe
the space of all properly embedded minimal surfaces of fixed genus
in a fixed (but arbitrary) closed $3$-manifold.  We will see that
the key is to understand the structure of an embedded minimal
planar domain in a ball in $\RR^3$.  Since the case of disks was
considered in the first four papers, the focus here is on
non-simply connected planar domains.

{\bf{We will first restrict to the case of planar domains, i.e.,
when the surfaces have genus zero. In particular, the main
theorems will first be stated and proved for planar domains.}} We
will
 see that the general case of fixed genus requires only minor
changes. The necessary changes to the main theorems and the
modifications needed for their proofs will be given in Part
\ref{s:highergenus}.

Sequences of planar domains which are not simply connected are,
after passing to a subsequence, naturally divided into two
separate cases depending on whether or not the topology is
concentrating at points.  To distinguish between these cases, we
will say that a sequence of surfaces $\Sigma_i^2\subset \RR^3$ is
{\it{uniformly locally simply connected}} (or ULSC)  if for each
compact subset $K$ of $\RR^3$, there exists a constant $r_0 > 0$
(depending on $K$) so that for every $x \in K$, all $r \leq
r_0${\footnote{If each component of the  intersection of a minimal
surface with a ball of radius $r_0$ is a disk, then so are the
intersections with all sub--balls by  the convex hull property
(see, e.g., lemma C.1 in \cite{CM6}). Therefore, it would be
enough that \eqr{e:ulsc2} holds for $r=r_0$.}}, and every surface
$\Sigma_i$
\begin{equation}    \label{e:ulsc2}
 {\text{each connected component of }} B_{r}(x) \cap \Sigma_i {\text{ is
 a disk.}}
\end{equation}
  For instance, a sequence of rescaled catenoids
   where the necks shrink to zero is not ULSC, whereas a
sequence of rescaled helicoids is.

Another way of locally distinguishing sequences where the topology
does not concentrate from sequences where it does comes from
analyzing the singular set. The singular set $\cS$ is defined to
be the set of points where the curvature is blowing up.  That is,
a point $y$ in $\RR^3$ is in $\cS$ for a sequence $\Sigma_i$ if
\begin{equation}
    \sup_{B_{r}(y)\cap
\Sigma_i}|A|^2\to\infty {\text{ as $i \to \infty$ for all $r>0$}}
.
\end{equation}
We will show that for embedded minimal surfaces $\cS$ consists of
two types of points. The first type is roughly modelled on
rescaled helicoids
 and the second on
rescaled catenoids:
\begin{itemize}
\item
A point $y$ in $\RR^3$ is in $\cSu$ if the curvature for the
sequence $\Sigma_i$ blows up at $y$ and the sequence is ULSC in a
neighborhood of $y$.
\item
A point $y$ in $\RR^3$ is in $\cSt$ if the sequence is not ULSC in
any neighborhood of $y$. In this case, a sequence of closed
non-contractible curves $\gamma_i \subset \Sigma_i$ converges to
$y$.
\end{itemize}
The sets $\cSt$ and $\cSu$ are obviously disjoint and
  the curvature blows up at both, so   $\cSt \cup \cSu \subset \cS$.   An easy
argument will later show that, after passing to a subsequence, we
can assume that
\begin{equation}    \label{e:}
    \cS=  \cSt \cup \cSu \, .
\end{equation}
Note that $\cSt = \emptyset$ is equivalent to that the sequence is
ULSC as is the case for sequences of rescaled helicoids.  On the
other hand, $\cSu = \emptyset$ for sequences of rescaled
catenoids.  These definitions of $\cSu$ and $\cSt$ are specific to
the genus zero case that we are focusing on now;  the slightly
different definitions in the higher genus case  can be found
around equation \eqr{e:genusstable}.

We will show that every sequence $\Sigma_i$ has a subsequence that
is either    ULSC or for which $\cSu$ is empty. This is the
 next ``no mixing'' theorem.  We will see later that these two different cases give
 two very different structures.

\begin{Thm}     \label{c:main}
If $\Sigma_i \subset B_{R_i}=B_{R_i}(0)\subset \RR^3$ is a
sequence
 of compact embedded minimal planar domains{\footnote{The theorem holds also for sequences with fixed genus;
 see Part \ref{s:highergenus}.}}
 %\label{fn:repeat}}}  %here
  with $\partial
\Sigma_i\subset
\partial B_{R_i}$ where $R_i\to \infty$, then there is a subsequence with either
$\cSu = \emptyset$ or $\cSt = \emptyset$.
\end{Thm}

In view of Theorem \ref{c:main} and the earlier results for disks,
it is natural to first analyze sequences that are ULSC, so where
$\cSt=\emptyset$, and second analyze sequences where $\cSu$ is
empty.  We will do this next.

\vskip2mm
 As already mentioned, our main theorems deal
with sequences $\Sigma_i \subset B_{R_i}=B_{R_i}(0)\subset \RR^3$
 of compact embedded minimal planar domains with $\partial
\Sigma_i\subset
\partial B_{R_i}$ where $R_i\to \infty$.  We will assume
here that these planar domains are not disks  (recall that the
case of disks was dealt with in \cite{CM3}--\cite{CM6}). In
particular, we will assume that for each $i$, there exists some
$y_i \in \RR^3$ and $s_i>0$ so that
\begin{equation}    \label{e:notulsc}
{\text{ some component of }} B_{s_i}(y_i) \cap \Sigma_i {\text{
 is not a disk.}}
\end{equation}
Moreover,  if the non-simply connected balls $B_{s_i}(y_i)$ ``run
off to infinity'' (i.e., if each connected component of
$B_{R_i'}(0) \cap \Sigma_i$ is a disk for some $R_i' \to \infty$),
then the results of \cite{CM3}--\cite{CM6} apply.  Therefore,
 after passing to a subsequence, we can assume that the
surfaces are uniformly not  disks, namely, that  there exists some
$R>0$ so that  \eqr{e:notulsc} holds with $s_i=R$ and $y_i=0$  for
all $i$.

 In general, we will allow our sequence of surfaces to
have bounded genus. Recall that for a surface $\Sigma$ with
boundary $\partial \Sigma$, the {\it genus}  of $\Sigma$   is the
genus of the closed surface $\hat{\Sigma}$ obtained by adding a
disk to each boundary circle. The genus of a union of disjoint
surfaces is the sum of the genuses. Therefore, a surface with
boundary has nonnegative genus; the genus is zero if and only if
it is a planar domain. For example, the disk and the annulus both
have genus zero; on the other hand, a closed surface of genus $g$
with any number of disks removed has genus $g$.

Common for both the ULSC case and the case where $\cSu$ is empty
is that the limits are always  laminations by flat parallel planes
and  the singular sets are always  closed subsets contained in the
union of the planes.  This is the content of the next theorem:

\begin{Thm} \label{t:tab}
 Let $\Sigma_i \subset B_{R_i}=B_{R_i}(0)\subset \RR^3$
be a sequence of compact embedded minimal planar
domains{\footnote{The theorem holds also for sequences with fixed genus;
 see Part \ref{s:highergenus}.}}
%\footnotemark[\ref{fn:repeat}]    %%here
with $\partial
\Sigma_i\subset \partial B_{R_i}$ where $R_i\to \infty$. If
\begin{equation}
\sup_{B_1\cap \Sigma_i}|A|^2\to \infty \, ,
 \end{equation}
   then there exists a subsequence $\Sigma_j$,
     a lamination $\cL=\{x_3=t\}_{ \{ t \in \cI \} }$ of $\RR^3$
 by parallel
 planes (where $\cI \subset \RR$ is a closed set), and a closed nonempty set
 $\cS$ in the union of the leaves of $\cL$ such that
 after a rotation of $\RR^3$:
 \begin{enumerate}
\item[(A)] For each $1>\alpha>0$, $\Sigma_j\setminus \cS$
converges in the $C^{\alpha}$-topology to the lamination $\cL
\setminus \cS$.
\item[(B)]  $\sup_{B_{r}(x)\cap \Sigma_j}|A|^2\to\infty$ as $j \to \infty$ for
all $r>0$ and $x\in \cS$.  (The curvatures blow up along $\cS$.)
\end{enumerate}
\end{Thm}

\vskip2mm Before discussing the general ULSC case, it is useful to
recall the case of disks.  One consequence of
\cite{CM3}--\cite{CM6} is that there are only two
{\underline{local models}} for ULSC sequences of embedded minimal
surfaces.  That is, locally in a ball in $\RR^3$, one of following
holds:
\begin{itemize}
\item  The curvatures are bounded and the surfaces are locally
{\underline{graphs}} over a plane.
\item
The curvatures blow up  and the surfaces are locally
{\underline{double spiral staircases}}.
\end{itemize}
Both of these cases are illustrated by taking a sequence of
rescalings of the helicoid; the first case occurs away from the
axis, while the second case occurs on the axis. Namely, recall
that the helicoid is the minimal surface $\Sigma$ in $\RR^3$
parametrized by
\begin{equation}
(s\cos t,s\sin t,t)\text{ where }s,\,t\in \RR\, .
\end{equation}
If we take a sequence $\Sigma_i = a_i \, \Sigma$ of rescaled
helicoids where $a_i \to 0$, then the curvature blows up along the
vertical axis but is bounded away from this axis.  Thus, we get
that
\begin{itemize} \item The intersection of the rescaled helicoids
with a ball {\underline{away from}} the vertical axis gives a
collection of graphs over the plane $\{ x_3 = 0 \}$.
\item The intersection of the
rescaled helicoids with a ball {\underline{centered on}} the
vertical axis gives a double spiral staircase.
\end{itemize}

\begin{figure}[htbp]
    \setlength{\captionindent}{20pt}
    %\begin{minipage}[t]{0.5\textwidth}
    \centering\input{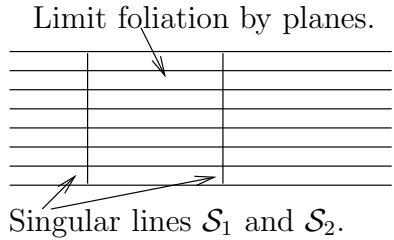}
    \caption{Theorem \ref{t:t5.1}: Limits of sequences of non-simply
connected, yet ULSC, surfaces with
    curvature blowing up.  The singular set consists of two lines
$\cS_1$ and $\cS_2$ and the limit  is a
    foliation by flat parallel planes.}
    \label{f:0}
    %\end{minipage}
\end{figure}

Loosely speaking, our next result  shows that when the sequence is
ULSC (but not simply connected), a subsequence converges to a
foliation by parallel planes away from two lines $\cS_1$ and
$\cS_2$; see Figure \ref{f:0}. The lines $\cS_1$ and $\cS_2$ are
disjoint and orthogonal to the leaves of the foliation and the two
lines are precisely the points where the curvature is blowing up.
This is similar to the case of disks, except that we get two
singular curves for non-disks as opposed to just one singular
curve for disks (the precise statement for disks is recalled in
Part \ref{p:disks}).

\begin{Thm} \label{t:t5.1}
 Let  a sequence $\Sigma_i$, limit lamination $\cL$, and singular
 set $\cS$
be  as in Theorem \ref{t:tab}.{\footnote{The theorem holds also
for sequences with fixed genus with one minor change in the
conclusion and one in the hypothesis.  The change in the
hypothesis is that we do not assume \eqr{e:notulsc}. The change in
the conclusion is that there might be either one or two singular
curves. The hypothesis
 \eqr{e:notulsc} is used in the genus zero case to
 show that there cannot be just
one singular curve.  The reason that we will not assume
\eqr{e:notulsc} in the fixed genus case is that there can be
either one or two singular curves in this case regardless;
 see Part \ref{s:highergenus}.}} Suppose that each
$\Sigma_i$ satisfies
 \eqr{e:notulsc}  with
$s_i=R > 1$ and $y_i=0$. If every $\Sigma_i$ is ULSC and
\begin{equation}
\sup_{B_1\cap \Sigma_i}|A|^2\to \infty \, ,
 \end{equation}
 then the limit lamination $\cL$ is the foliation $\cF = \{ x_3 = t\}_t$
  and the singular set $\cS$ is the union of two  disjoint
lines $\cS_1$ and $\cS_2$
 such that:
 \begin{enumerate}
\item[($C_{ulsc}$)]
Away from $\cS_1 \cup \cS_2$, each $\Sigma_j$ consists of exactly
two multi-valued graphs spiraling together.
 Near $\cS_1$ and $\cS_2$, the pair of multi-valued graphs form double
spiral staircases
with opposite orientations at $\cS_1$ and $\cS_2$. Thus, circling
only $\cS_1$ or only $\cS_2$  results in going either up or down,
while a path circling both $\cS_1$ and $\cS_2$ closes up (see
Figure \ref{f:mg2}).
\item[($D_{ulsc}$)] $\cS_1$ and $\cS_2$
are  orthogonal to the leaves of the foliation.
\end{enumerate}
\end{Thm}

Notice that Theorem \ref{t:t5.1} shows that if the fixed genus
ULSC surfaces $\Sigma_j$ have curvature blowing up, then they
essentially have genus zero.  More precisely, given an arbitrarily
large ball $B_R \subset \RR^3$, then $B_R \cap \Sigma_j$ has genus
zero for $j$ sufficiently large.  To see this, combine the double
spiral staircase structure near the two singular curves that holds
for ULSC sequences (cf. Figure \ref{f:mg2}) with the smooth
convergence elsewhere.

Despite the similarity  of Theorem \ref{t:t5.1}
 to the case of disks, it is worth noting that the results for
 disks do not alone give this theorem.  Namely, even though the
 ULSC sequence consists {\underline{locally}} of disks, the compactness
result for disks
 was in the {\underline{global}} case where the radii go to
 infinity.  One might wrongly think that Theorem \ref{t:t5.1} could be proven
 using the results for disks and a blow up argument. However, local
 examples constructed in \cite{CM15} shows the difficulty with such an
 argument.{\footnote{In  \cite{CM15}, we constructed
  a sequence of embedded minimal disks $\Sigma_i$ in
the unit ball $B_1$ with $\partial \Sigma_i \subset
\partial B_1$ where
 the curvatures blow up only at $0$.  This sequence converges to a lamination of
 $B_1 \setminus \{ 0\}$ that
 cannot be extended smoothly to a lamination of $B_1$; that is to say, $0$ is not a removable
 singularity.  This should be contrasted with Theorem \ref{t:t5.1}
 where every singular point is a removable singularity
 for the limit foliation by parallel planes.}}
 We shall explain this further later together with what
 else is needed for the proof.

\vskip2mm
When the sequence is no longer ULSC, then there are other local
models for the surfaces.  The simplest example is a sequence of
rescaled catenoids; the catenoid is the minimal surface in $\RR^3$
parametrized by
\begin{equation}    \label{e:cat}
 (\cosh s\, \cos t,\cosh s\, \sin t,s)  \text{ where }s,\,t\in \RR\,
 .
\end{equation}
 A sequence of rescaled  catenoids converges with
multiplicity two to the flat plane.  The convergence is in the
$C^{\infty}$ topology except at $0$ where $|A|^2 \to \infty$. This
sequence of rescaled catenoids is not ULSC because the simple
closed geodesic on the catenoid -- i.e., the unit circle in the
$\{ x_3 = 0 \}$ plane -- is non-contractible and the rescalings
shrink it down to the origin.

One can get other types of curvature blow-up by considering the
family   of embedded minimal planar domains known as the Riemann
examples.{\footnote{See
http://www.msri.org/publications/sgp/jim/geom/minimal/library/riemann/index.html
for a description, as well as computer graphics, of these
surfaces.}}
 Modulo translations and rotations, this is a
two-parameter family of periodic minimal surfaces, where the
parameters can be thought of as the size of the necks and the
angle  from one fundamental domain to the next. By choosing the
two parameters appropriately, one can produce sequences of Riemann
examples that illustrate both of the two structure theorems (cf.
Figures \ref{f:01a} and \ref{f:01b}):
\begin{enumerate}
\item
If we take a sequence of Riemann examples where the neck size is
fixed and the angles go to $\frac{\pi}{2}$, then the surfaces with
angle near $\frac{\pi}{2}$ can be obtained by gluing together two
oppositely--oriented double spiral staircases. Each double spiral
staircase looks like a helicoid. This sequence of Riemann examples
converges to a foliation by parallel planes.  The convergence is
smooth away from the axes of the two helicoids (these two axes are
the singular set $\cS$ where the curvature blows up). The sequence
is ULSC since the size of the necks is fixed and  thus illustrates
the first structure theorem, Theorem \ref{t:t5.1}.
\item
If we take a sequence of examples where the neck sizes go to zero,
then we get a sequence  that is {\it{not}} ULSC.  However, the
surfaces can be cut along   short curves into collections of
graphical pairs of pants.  The short curves converge to points and
the graphical pieces converge to flat planes except at these
points, illustrating the second structure theorem, Theorem
\ref{t:t5.2a} below.
\end{enumerate}

With these examples in mind, we are now ready to state our second
main structure theorem describing the case where $\cSu$ is empty.

\begin{Thm} \label{t:t5.2a}
 Let  a sequence $\Sigma_i$, limit lamination $\cL$, and singular
 set $\cS$
be  as in Theorem \ref{t:tab}.{\footnote{The theorem holds also
for sequences with fixed genus with one small change in
($C_{neck}$). Namely, the number of points in ($C_{neck}$) is
bounded by two plus the bound for the genus;
 see Part \ref{s:highergenus}.\label{fn:repeat3}}}  If
$\cSu = \emptyset$ and
\begin{equation}
\sup_{B_1\cap \Sigma_i}|A|^2\to \infty \, ,
 \end{equation}
   then $\cS=\cSt$ by \eqr{e:} and
 \begin{enumerate}
\item[($C_{neck}$)]
Each point $y$ in $\cS$ comes with a sequence of
{\underline{graphs}} in $\Sigma_j$ that converge  to the plane $\{
x_3 = x_3 (y) \}$. The convergence is in the $C^{\infty}$ topology
away from the point $y$ and possibly also one other point  in $\{
x_3 = x_3 (y) \} \cap \cS$.  If the convergence is   away from one
point, then these graphs are defined over annuli; if the
convergence is away from two points, then the graphs are defined
over disks with two subdisks removed.
\end{enumerate}
\end{Thm}

\vskip2mm Theorem \ref{t:t5.2a}, as well as Theorem \ref{c:main},
are proven by first analyzing sequences of minimal surfaces
without any assumptions on the sets $\cSu$ and $\cSt$.  In this
general case, we show that a subsequence
  converges to a lamination $\cL'$ divided into regions where
  Theorem \ref{t:t5.1} holds and  regions
  where Theorem \ref{t:t5.2a} holds.  This
convergence is in the smooth topology away from the singular set
$\cS$ where the curvature blows up.  Moreover, each point of $\cS$
comes with a plane and these planes are essentially contained in
$\cL'$; see (P) below. The set of heights of the planes is a
closed subset $\cI \subset \RR$ but may not be all of $\RR$ as it
was in Theorem \ref{t:t5.1} and may not even be connected. The
behavior of the sequence is different at the two types of singular
points in $\cS$ - the set $\cSt$ of ``catenoid points'' and the
set $\cSu$ of ULSC singular points. We will see that $\cSu$
consists of a union of Lipschitz curves transverse to the
lamination $\cL$. This structure of $\cSu$ implies that the set of
heights in $\cI$ which intersect $\cSu$ is a union of intervals;
thus this part of the lamination is foliated. In contrast, we will
not get any structure of the set of ``catenoid points'' $\cSt$;
see (D) below. Given a point $y$ in $\cSt$, we will get a sequence
of graphs in $\Sigma_j$ converging to a plane through $y$; see
(C1) below. This convergence will be in the smooth topology away
from either one or two singular points, one of which is $y$.
Moreover, this limit plane through $y$ will be a leaf of the
lamination $\cL$.

 The precise statement of the compactness theorem for sequences
 that are neither necessarily ULSC nor with $\cSu = \emptyset$
 is the following (see Figure \ref{f:0b}):

\begin{Thm} \label{t:t5.2}
 Let $\Sigma_i \subset B_{R_i}=B_{R_i}(0)\subset \RR^3$
be a sequence of compact embedded minimal planar
domains{\footnote{The theorem holds also for sequences with fixed genus;
 see Part \ref{s:highergenus}.}}
%\footnotemark[\ref{fn:repeat3}] %here
with $\partial
\Sigma_i\subset \partial B_{R_i}$ where $R_i\to \infty$. If
\begin{equation}
\sup_{B_1\cap \Sigma_i}|A|^2\to \infty \, ,
 \end{equation}
   then there is a subsequence $\Sigma_j$, a closed set $\cS$,
    and a lamination $\cL'$ of
$\RR^3 \setminus \cS$ so that:
\begin{enumerate}
\item[(A)] For each $1>\alpha>0$, $\Sigma_j\setminus \cS$
converges in the $C^{\alpha}$-topology to the lamination $\cL'$.
\item[(B)]  $\sup_{B_{r}(x)\cap \Sigma_j}|A|^2\to\infty$ as $j \to \infty$ for
all $r>0$ and $x\in \cS$.  (The curvatures blow up along $\cS$.)
\item[(C1)]
($C_{neck}$) from Theorem \ref{t:t5.2a} holds for each point
 $y$ in $\cSt$.
\item[(C2)]
($C_{ulsc}$) from Theorem \ref{t:t5.1} holds locally near $\cSu$.
 More precisely, each point $y$ in $\cSu$ comes with a sequence
of {\underline{multi-valued graphs}} in $\Sigma_j$ that converge
to the plane $\{ x_3 = x_3 (y) \}$. The convergence is in the
$C^{\infty}$ topology away from the point $y$ and possibly also
one other point   in $\{ x_3 = x_3 (y) \} \cap \cSu$.  These two
possibilities correspond to the two types of multi-valued graphs
defined in Section \ref{s:s1}.
\item[(D)]     The set $\cSu$ is a union of Lipschitz curves
transverse to the lamination.  The leaves intersecting $\cSu$ are
planes foliating an open subset of $\RR^3$ that does not intersect
$\cSt$. For the set $\cSt$, we make no claim about the structure.
\item[(P)]
Together (C1) and (C2) give a sequence of graphs or multi-valued
graphs converging to a plane through each point of $\cS$. If $P$
is one of these planes, then each leaf of $\cL'$ is either
disjoint from $P$ or is contained in $P$.
\end{enumerate}
\end{Thm}

\begin{figure}[htbp]
    \setlength{\captionindent}{20pt}
    %\begin{minipage}[t]{0.5\textwidth}
    \centering\input{gen0b.pstex_t}
    \caption{Theorem \ref{t:t5.2}: Limits of sequences of {\underline{non-ULSC}} surfaces with
    curvature blowing up.
     The limit  is a
    {\underline{lamination}} of $\RR^3 \setminus \cS$.
    The singular set $\cS$ consists of {\underline{two}}
     types of points - the ones in $\cSt$ and the ones in $\cSu$.
    Note that the set $\cSt$ is automatically closed, while the set
    $\cSu$ is not.
    The set $\cSu$ is a union of Lipschitz curves; the injectivity radius goes to
    zero at the endpoints of these curves, so these endpoints are
    in  $\cSt$.  Finally,  the part of the lamination
    containing $\cSu$ is foliated by planar leaves.}
    \label{f:0b}
    %\end{minipage}
\end{figure}

 Note that Theorem \ref{t:t5.2}
 is a technical tool that will be used to
 prove the main compactness theorem in the non-ULSC case, Theorem
 \ref{t:t5.2a}. In particular, Theorem \ref{t:t5.2} itself   will be superseded
 by the stronger compactness theorems in the ULSC and non-ULSC
 cases, Theorem \ref{t:t5.1} and Theorem
 \ref{t:t5.2a}.  This is because eventually we will know by the no mixing
 theorem that either $\cSt =
\emptyset$ or $\cSu = \emptyset$, so that these cover all possible
cases.  Moreover, the assertions in Theorem \ref{t:t5.1} and
Theorem
 \ref{t:t5.2a} are stronger than those in Theorem
 \ref{t:t5.2}.

After proving  Theorem \ref{t:t5.2} in Part \ref{p:prove2}, we
will be ready in Part \ref{p:nomix} to prove
 the
no mixing theorem,   Theorem \ref{c:main}.

In Part \ref{p:prove3}, we will then
 complete the
 proof of Theorem \ref{t:t5.2a}.  The main point left, which is not part of Theorem
 \ref{t:t5.2},
 is to prove that {\underline{every}} leaf of
  the lamination
  $\cL$ in Theorem \ref{t:t5.2a} is a plane.  In contrast,
Theorem \ref{t:t5.2} gives a plane through each point of $\cSt$,
but does not claim that the leaves of $\cL'$ are planar.

 Finally, since
the no mixing theorem implies that Theorem \ref{t:t5.1} and
Theorem
 \ref{t:t5.2a} cover all cases,
   Theorem \ref{t:tab} will be a corollary of
 these two theorems.

\vskip2mm  We refer to the introduction of \cite{CM6} and the
surveys \cite{MeP}, \cite{P}, and \cite{Ro} for related results,
including applications of the results of \cite{CM3}--\cite{CM6} as
well as the results of this paper.

\subsection{Brief outline of the paper and overview of the proofs}

In Section \ref{s:s1}, we will define the two notions of
multi-valued graphs which will be needed to explain and prove the
two main theorems.

 Part \ref{p:disks} is devoted to recalling some of the earlier results
for disks given in \cite{CM3}--\cite{CM6} and \cite{CM10}.  The
first of these shows that embedded minimal disks are either graphs
or are part of a double spiral staircase.  The second result that
we recall is the one-sided curvature estimate.  Finally, we will
recall the chord-arc bound for embedded minimal disks proven in
\cite{CM10}.

In Part \ref{p:partsi}, we will first define the singular set
$\cS$ and prove the convergence to the lamination $\cL'$ away from
$\cS$.  The rest of the part focuses on  describing  a
neighborhood of each point in the ULSC singular set $\cSu$ and the
leaves of $\cL'$ whose closure intersects $\cSu$.  A key point
will be  that the results of \cite{CM3}--\cite{CM6} for disks will
give a sequence of multi-valued graphs in the $\Sigma_j$'s
 near each point
$x \in \cSu$.  Moreover, these multi-valued graphs close up in the
limit to give a leaf of $\cL'$ which extends smoothly across $x$.
Such a leaf is said to be collapsed; in a neighborhood of $x$, the
leaf can be thought of as a limit of double-valued graphs where
the upper sheet collapses onto the lower. We will show that every
collapsed leaf is stable, has at most two points of $\cSu$ in its
closure, and these points are removable singularities. These
results on collapsed leaves will be applied first in the USLC case
in the next part and then later to get the structure of the ULSC
regions of the limit in general, i.e., (C2) and (D) in Theorem
\ref{t:t5.2}.

In Part \ref{p:prove1}, we prove Theorem \ref{t:t5.1} that gives
the convergence of a ULSC sequence to a foliation by parallel
planes away from two singular curves.  Roughly speaking, there are
two main steps to the proof:
\begin{enumerate}
\item
Show that each collapsed leaf is in fact a plane punctured at two
points of $\cS$ and,  moreover, the sequence has the structure of
a double spiral staircase near both of these points, with opposite
orientations at the two points.
\item
Show that leaves which are nearby a collapsed leaf of $\cL'$ are
also planes punctured at two points of $\cS$.  (We call this
``properness''.)
\end{enumerate}

In Part \ref{p:prove2} we consider general sequences of minimal
surfaces that are neither necessarily ULSC nor with $\cSu =
\emptyset$ and we prove the general compactness theorem, Theorem
 \ref{t:t5.2}.  Recall that this theorem
asserts that the limit lamination $\cL'$ can be divided into two
disjoint sub-laminations.  One of which is the support of a region
where (a subsequence of) the surfaces are ULSC and all of the
results about ULSC sequences from Part \ref{p:prove1} hold, such
as the structure of the singular set and the multi-valued graphs
structure.  In the other region, curvature blow up comes
exclusively from neck pinching and, thus, in this region there are
no helicoid like points.  The key steps for proving the general
structure theorem are the following:
\begin{enumerate}
\item Finding a {\underline{stable}} plane through each point of
$\cSt$.  This plane will be a limit of a sequence of stable
graphical annuli that lie in the complement of the surfaces.
\item
Finding graphs in $\Sigma_j$ that converge to a plane through each
point of $\cSt$.  To do this, we look in regions between
consecutive necks and show that in any such region the surfaces
are ULSC.  The one-sided curvature estimate will then allow us to
show that these regions are graphical.
\item
Using (1) and (2) we then analyze the ULSC regions of a limit.
That is, we show that if the closure of a leaf in $\cL'$
intersects $\cSu$, then it has a neighborhood that is ULSC.  This
will allow us to use the argument for the proof of Theorem
\ref{t:t5.1} to get the same structure for such a neighborhood as
we did in case where the entire surfaces where ULSC.
\end{enumerate}

In Part \ref{p:nomix}, we will use the structure obtained in
Theorem \ref{t:t5.2} to show the no mixing theorem; Theorem
\ref{c:main}. The key here is to show that if $\cSu$ is non-empty,
then $\cSu$ cannot stop.

In Part \ref{p:prove3}, we will complete the proofs of Theorem
\ref{t:tab} and Theorem \ref{t:t5.2a}. The only thing that remains
to be proven is that {\underline{every}} leaf of the lamination
$\cL'$  is contained in a plane. We have already proven that the
leaves of $\cL'$ are planes when the sequence is ULSC; thus, by
the no mixing theorem, the only remaining case is when $\cS = \cSt
\ne \emptyset$.
  We will divide
the proof that the leaves of $\cL'$ are contained in planes into
two cases, depending on whether or not the leaf is complete.
 In both cases, we will use a flux argument to rule out a
   non-flat leaf of $\cL'$.

In Part \ref{s:highergenus}, we describe the necessary changes to the main theorems and the
modifications needed for their proofs when the sequence has positive genus.

\section{Multi-valued graphs}   \label{s:s1}

To explain the theorems stated in the introduction and their
proofs, we will need two notions of  multi-valued graphs - namely,
the one used in \cite{CM3}--\cite{CM6} and a generalization.

\begin{figure}[htbp]
    \setlength{\captionindent}{20pt}
    %\begin{minipage}[t]{0.5\textwidth}
    \centering\input{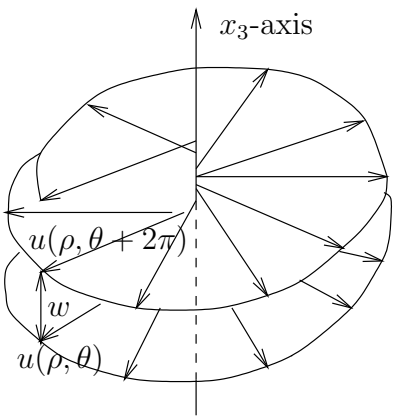}
    \caption{A multi-valued graph over the singly-punctured plane.}
    \label{f:mg1}
    %\end{minipage}
\end{figure}

In  \cite{CM3}--\cite{CM6}, we defined multi-valued graphs as
multi-sheeted covers of the punctured plane. To be precise, let
$D_r$ be the disk in the plane centered at the origin and of
radius $r$ and let $\cP$ be the universal cover of the punctured
plane $\CC\setminus \{0\}$ with global polar coordinates $(\rho,
\theta)$ so $\rho>0$ and $\theta\in \RR$. Given $0 \leq r\leq s$
and $\theta_1 \leq \theta_2$, define the ``rectangle''
$S_{r,s}^{\theta_1 , \theta_2} \subset \cP$ by
\begin{equation}
   S_{r,s}^{\theta_1 , \theta_2} = \{(\rho,\theta)\,|\,r\leq \rho\leq s\, ,\, \theta_1 \leq \theta\leq
\theta_2 \} \, .
\end{equation}
An $N$-valued graph of a function $u$ on the annulus $D_s\setminus
D_r$ is a single valued graph over  (see Figure \ref{f:mg1})
\begin{equation}
   S_{r,s}^{-N\pi , N\pi} = \{(\rho,\theta)\,|\,r\leq \rho\leq s\, ,\, |\theta|\leq
N\,\pi\} \, .
\end{equation}
($\Sigma_{r,s}^{\theta_1 , \theta_2}$ will denote the subgraph of
$\Sigma$ over the smaller rectangle $S_{r,s}^{\theta_1 ,
\theta_2}$).
  As in the earlier papers in the series, the multi-valued graphs that we will
consider will never close up; in fact they will all be embedded.
Note that embedded corresponds to that the separation never
vanishes.  Here the separation $w$ is the difference in height
between consecutive sheets and is  therefore given by
\begin{equation}        \label{e:sepw}
w(\rho,\theta)=u(\rho,\theta+2\pi)-u(\rho,\theta)\, .
\end{equation}
In the case where $\Sigma$ is the helicoid [i.e., $\Sigma$ can be
parametrized by $(s\,\cos t,s\,\sin t,t)$ where $s,\,t\in \RR$],
then
\begin{equation}
    \Sigma\setminus x_3-\text{axis}=\Sigma_1\cup \Sigma_2 \, ,
\end{equation}
 where $\Sigma_1$, $\Sigma_2$ are $\infty$-valued graphs.
$\Sigma_1$ is the graph of the function $u_1(\rho,\theta)=\theta$
and $\Sigma_2$ is the graph of the function
$u_2(\rho,\theta)=\theta+\pi$. In either case the separation
$w=2\,\pi$.

 Locally, the above multi-valued graphs give the complete picture
 for a ULSC sequence.   However, the global picture can consist of
 several different multi-valued graphs glued together.  To
 allow for this, we are forced to consider multi-valued graphs
 defined over the universal cover of $\CC \setminus P$ where $P$
 is a discrete subset of the complex plane $\CC$ (see Figure \ref{f:mg2}).
 We will see that the bound on the
 genus implies that $P$ consists of at most two points.   The
 basic example of such a multi-valued graph
 comes from the family of minimal surfaces known as
 the Riemann examples.

\begin{figure}[htbp]
    \setlength{\captionindent}{20pt}
    \centering\input{gen1b.pstex_t}
    \caption{A multi-valued graph over the doubly-punctured plane.
    The spiral staircases near each puncture are oppositely--oriented.}
    \label{f:mg2}
\end{figure}

%numbering for equations

\setcounter{part}{0}
\numberwithin{section}{part} %number sections within parts
\renewcommand{\rm}{\normalshape} %makes rm mean roman
\renewcommand{\thepart}{\Roman{part}}
\setcounter{section}{1}

\part{Results for disks from \cite{CM3}--\cite{CM6}}
\label{p:disks}

The results for non-simply connected minimal surfaces that are
proven in this paper rely on the earlier results for disks given
in \cite{CM3}--\cite{CM6}. For completeness and easy reference, we
start by recalling those.

\section{The lamination theorem and one-sided curvature estimate}

The first theorem that we recall shows that embedded minimal disks
are either graphs or are part of    double spiral staircases;
moreover, a sequence of such disks with curvature blowing up
converges to a foliation by parallel planes away from a singular
curve $\cS$.  This theorem is modelled on rescalings of the
helicoid and the precise statement is as follows (we state the
version for extrinsic balls; it was extended to intrinsic balls
in \cite{CM10}):

\begin{Thm} \label{t:t0.1}
(Theorem 0.1 in \cite{CM6}.)
 Let $\Sigma_i \subset B_{R_i}=B_{R_i}(0)\subset \RR^3$
be a sequence of embedded minimal {\underline{disks}} with
$\partial \Sigma_i\subset \partial B_{R_i}$ where $R_i\to \infty$.
If \begin{equation}
    \sup_{B_1\cap \Sigma_i}|A|^2\to \infty \, ,
    \end{equation}
 then there exists a subsequence, $\Sigma_j$, and
a Lipschitz curve $\cS:\RR\to \RR^3$ such that after a rotation of
$\RR^3$:
\begin{enumerate}
\item[\underline{1.}] $x_3(\cS(t))=t$.  (That is, $\cS$ is a graph
over the $x_3$-axis.)
\item[\underline{2.}]  Each $\Sigma_j$
consists of exactly two multi-valued graphs away from $\cS$ (which
spiral together).
\item[\underline{3.}] For each $1>\alpha>0$,
$\Sigma_j\setminus \cS$ converges in the $C^{\alpha}$-topology to
the foliation, $\cF=\{x_3=t\}_t$, of $\RR^3$.
\item[\underline{4.}] $\sup_{B_{r}(\cS (t))\cap
\Sigma_j}|A|^2\to\infty$ for all $r>0$, $t\in \RR$.  (The
curvatures blow up along $\cS$.)
\end{enumerate}
\end{Thm}

The second theorem that we need to recall asserts that every
embedded minimal disk lying above a plane, and coming close to the
plane near the origin, is a graph.   Precisely this is the {\it
intrinsic one-sided curvature estimate} which follows by combining
\cite{CM6} and \cite{CM10}:

\begin{Thm}  \label{t:t2}
There exists $\epsilon>0$, so that if
\begin{equation}
    \Sigma \subset \{x_3>0\} \subset \RR^3
\end{equation}
is an embedded minimal {\underline{disk}} with $\cB_{2R} (x)
\subset \Sigma \setminus
\partial \Sigma$ and $|x|<\epsilon\,R$,  then
\begin{equation}        \label{e:graph}
\sup_{ \cB_{R}(x) } |A_{\Sigma}|^2 \leq R^{-2} \, .
\end{equation}
\end{Thm}

Theorem \ref{t:t2} is in part used to prove the regularity of the
singular set where the curvature is blowing up.

\vskip2mm Note that the assumption in Theorem \ref{t:t0.1} that
the surfaces are disks is crucial and cannot even be replaced by
assuming that the sequence is ULSC. To see this, observe that  one
can choose a one-parameter family of Riemann examples which is
ULSC but where the singular set $\cS$ is given by a
{\underline{pair}} of vertical lines. Likewise, the assumption in
Theorem \ref{t:t2} that $\Sigma$ is simply connected is crucial as
can be seen from the example of a rescaled catenoid, see
\eqr{e:cat}.
 Under rescalings the catenoid converges (with
multiplicity two) to the flat plane. Thus a neighborhood of the
neck can be scaled arbitrarily close to a plane but the curvature
along the neck becomes unbounded as it gets closer to the plane.
Likewise, by considering the universal cover of the catenoid, one
sees that embedded, and not just immersed, is needed in Theorem
\ref{t:t2}.

\vskip2mm Finally, we recall the chord-arc bound for embedded
minimal disks proven in theorem $0.5$ of  \cite{CM10}:

\begin{Thm}     \label{t:1}
\cite{CM10}. There exists a  constant  $C > 0$  so that if $\Sigma
\subset \RR^3$ is an embedded minimal disk, $\cB_{2R}=\cB_{2R}(0)$
is an intrinsic ball in $\Sigma \setminus
\partial \Sigma$ of radius $2R$, and
$ \sup_{\cB_{r_0}}|A|^2>r_0^{-2}$ where $R>  r_0$,
 then for $x \in \cB_R$ the intrinsic distance is bounded from
 above by the extrinsic distance as follows
 \begin{equation}   \label{e:t1}
    C \, \dist_{\Sigma}(x,0)<   |x| + r_0   \, .
 \end{equation}
\end{Thm}

\part{The singular set $\cS$ and limit lamination $\cL'$}
\label{p:partsi}

\setcounter{equation}{0}

The three main results of this part are the convergence to the
lamination $\cL'$ away from a singular set $\cS$, the description
of a neighborhood of each ULSC singular point, and the description
of the leaves of $\cL'$ whose closure intersects $\cSu$.  We will
explain these in a bit more detail next.

 We start
by defining the singular set $\cS$; roughly speaking, $\cS$ is the
set of points where the curvature blows up (see Definition/Lemma
\ref{l:inftyornot}).  The definition of $\cS$ will immediately
imply that $\cS$ is a closed subset of $\RR^3$. We next show that
in the open subset $\RR^3 \setminus \cS$, a subsequence of the
sequence of embedded minimal surfaces converges to a minimal
lamination $\cL'$ of $\RR^3 \setminus \cS$ (see Lemma
\ref{l:singl}).

The results of \cite{CM3}--\cite{CM6} give a precise description
of a neighborhood of each point in $\cSu$. Namely, for $j$ large,
$\Sigma_j$ must be a double-spiral staircase near each point in
$\cSu$ and the set $\cSu$ must satisfy a local cone property which
gives the regularity of the set.  The description near a singular
point and local cone property are given in Lemma \ref{l:leaf}.  We
also recall in Lemma \ref{l:leaf} that, as $j\to \infty$, this
sequence of double-spiral staircases near a singular point $x$
closes up in the limit to give a leaf of $\cL'$ which extends
smoothly across $x$.  We will say that such a leaf is collapsed;
in a neighborhood of $x$, the leaf can be thought of as a limit of
double-valued graphs where the upper sheet collapses onto the
lower.

Finally, we will show that every collapsed leaf is stable, has at
most two points of $\cSu$ in its closure, and these points are
removable singularities. The key for proving stability  is to use
the separations of the limiting multi-valued graphs to construct a
positive Jacobi field in the limit. The limit Jacobi field is not
a priori well-defined, but is instead well-defined on a covering
space of the collapsed leaf. However, we show in Appendix B that
stability of a covering space implies stability of the surface
itself as long as the covering space has sub-exponential area
growth.  We apply this to show that every collapsed leaf is
stable.  We will also use the fact that the surfaces $\Sigma_j$
have bounded genus, to show that each collapsed leaf has at most
two points of $\cSu$ in its closure.

\vskip2mm  These results on collapsed leaves will be applied first
in the USLC case  in the next part and then later to get the
structure of the ULSC regions of the limit in general, i.e., (C2)
and (D) in Theorem \ref{t:t5.2}.

\section{The singular set $\cS$}

To define the singular set,    recall from \cite{CM6} that
  for any sequence of surfaces (minimal or not) in $\RR^3$, after possibly going
to a subsequence, then there is a well defined notion of points in
$\RR^3$ where the second fundamental form of the sequence blows
up.  The set of such points will below be referred to as the
singular set $\cS$ and is given by an elementary and
straight-forward compactness argument.

\begin{DefLem} \label{l:inftyornot}
(The singular set; Lemma I.1.4 in \cite{CM6}.)
 Let $\Sigma_i\subset B_{R_i}$ with
$\partial \Sigma_i\subset
\partial B_{R_i}$ and $R_i\to \infty$ be a sequence of (smooth)
compact surfaces. After passing to a subsequence, $\Sigma_j$, we
may assume that for each $x\in \RR^3$ either of the two following
properties holds:
\begin{itemize}
\item $\sup_{B_{r}(x)\cap \Sigma_j}|A|^2\to \infty$ for all
$r>0$. (The set of such points $x$ will be denoted by $\cS$.)
\item $\sup_j\sup_{B_r(x)\cap \Sigma_j}|A|^2<\infty$ for some
$r>0$.
\end{itemize}
\end{DefLem}

\subsection{Convergence away from $\cS$}
The first result that we will need is that in the open subset
$\RR^3 \setminus \cS$, a sequence of embedded minimal surfaces has
a subsequence that converges to a minimal lamination $\cL'$ of
$\RR^3 \setminus \cS$.  This is an easy consequence of that the
curvature is bounded on compact sets in the complement of $\cS$
and is proven in the next lemma.

\begin{Lem}     \label{l:singl}
Suppose that $\Sigma_j$ and $\cS$ are as in Lemma
\ref{l:inftyornot}.  If in addition the $\Sigma_j$'s are minimal
and embedded, then there exists a subsequence (still denoted by
$\Sigma_j$) and a lamination $\cL'$ of $\RR^3 \setminus \cS$ so
that the following hold: \begin{itemize} \item $\Sigma_j \to \cL'$
on compact subsets of $\RR^3 \setminus \cS$.
\item
The leaves of $\cL'$ are minimal.
\end{itemize}
\end{Lem}

\begin{proof}
For each compact subset $K$ of $\RR^3 \setminus \cS$, then Lemma
\ref{l:inftyornot} gives an open covering of $K$ by finitely many
balls where the curvatures of the $\Sigma_j$'s are bounded
(independent of $j$) in the concentric double balls.   Both claims
now follow  from proposition B.1 in \cite{CM6} and a diagonal
argument.
\end{proof}

As in \cite{CM6}, convergence to $\cL'$ in the above lemma means
that if we think of the embedded surfaces $\Sigma_j$ themselves as
laminations, then the coordinate charts for these laminations
converge in the $C^{\alpha}$ norm for any $\alpha < 1$ and the
leaves converge as sets. The convergence is actually $C^{\infty}$
tangentially, meaning that if we write a leaf locally as a graph,
then a sequence of local graphs in $\Sigma_j$ converges smoothly
to this leaf.  This tangential regularity follows from the
$C^{\alpha}$ convergence and elliptic estimates.  However, easy
examples show that the convergence in the transversal direction
may only be in the Lipschitz topology; cf. \cite{So}.

\vskip2mm {\it Throughout the rest of this paper, we will assume
that $\Sigma_j\subset B_{R_j}$ with $\partial \Sigma_j\subset
\partial B_{R_j}$ and $R_j\to \infty$ is a sequence of  (smooth)
compact embedded minimal surfaces that converges off of a singular
set $\cS$ to a lamination $\cL'$ of $\RR^3 \setminus \cS$ with
minimal leaves.  The lamination $\cL'$ is given by  Lemma
\ref{l:singl}. In order to obtain additional structure of $\cS$
and $\cL'$, we will need to also make topological assumptions
about the surfaces $\Sigma_j$.  We will always assume that the
$\Sigma_j$'s have bounded genus.  In Part \ref{p:prove1}, we will
assume that the surfaces $\Sigma_j$ are ULSC, i.e, that $\cSt =
\emptyset$;  in Part \ref{p:prove2}, we will consider the other
case where $\cSt \ne \emptyset$.}

\section{The local structure of $\cL'$ near a point of $\cSu$}

We will eventually show that all of the leaves of the lamination
$\cL'$ are flat (see Theorem \ref{t:tab}), but we will need to
first establish some initial structure of $\cL'$.  The first step
will be accomplished in this section where we describe the local
structure of $\cL'$ near a point in $\cSu$.

\begin{figure}[htbp]
    \setlength{\captionindent}{20pt}
    \centering\input{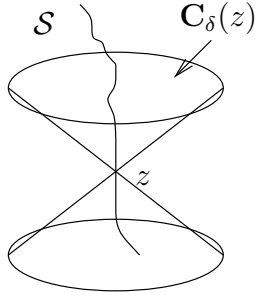}
    \caption{It follows from the one-sided curvature estimate that the ULSC singular
    set $\cSu$ has the local cone property and, as we will see, this gives Lipschitz regularity.}
    \label{f:fcone}
\end{figure}

 The next lemma is going to show that
each ULSC singular point lies in the closure of a leaf of $\cL'$
which extends smoothly across the singular point and, furthermore,
ULSC singular points are leaf-wise isolated and they satisfy a
local cone property. To state this cone property, let
$\cone_{\delta}(z)$ be the (convex) double cone with vertex $z$,
cone angle $(\pi/2 - \arctan \delta)$, and axis parallel to the
$x_3$-axis. That is (see Figure \ref{f:fcone}),
\begin{equation}        \label{e:conedelta}
\cone_{\delta}(z)=\{x\in \RR^3\,|\,(x_3-z_3)^2 \geq
\delta^2\,((x_1-z_1)^2+(x_2-z_2)^2) \} \, .
\end{equation}

The local cone property is now defined as follows. Given $\delta >
0$ and $r_0
> 0$, we will say that a subset $\cSu \subset \RR^3$ has the {\it local
cone property} if $\cSu$ is  nonempty and
\begin{equation}  \label{e:lcpr}
    {\text{if }} z\in \cSu, \, {\text{then }} B_{r_0} (z)
    \cap \cSu \subset \cone_{\delta}(z) \, .
\end{equation}
As in \cite{CM6}, we will see in Section \ref{s:compf} that this
local cone property directly gives Lipschitz regularity of the
subset $\cSu$.

We can now state the lemma which gives the regularity of the
leaves through $\cSu$ and the local cone property for $\cSu$.

\begin{Lem}     \label{l:leaf}
Given a point $x \in  \cSu$,  there exists $r_0 > 0$ so that
$B_{r_0}(x)
 \cap \cL'$ has a component $\Gamma_x$ whose closure
 $\overline{\Gamma_x}$ is a smooth minimal graph containing $x$
and with boundary in $\partial B_{r_0}(x)$ (so $x$ is a removable
singularity for $\Gamma_x$).

 Furthermore, $\overline{\Gamma_x} \cap \cS = \{ x \}$ and, after rotating $\RR^3$ so that
 $\nn_{ \overline{\Gamma_x} } (x) = (0,0,1)$,
 the set $\cSu$ satisfies the local cone property
 \eqr{e:lcpr} for some $\delta >0$ and the above $r_0$.  The
 rotation may vary with $x$, but the dependence is Lipschitz.
\end{Lem}

\begin{proof}
 For simplicity,
translate so that $x = 0$. Since $0 \notin \cSt$, there exists
some $r_0 > 0$  so that the components of $B_{r_0} (0) \cap
\Sigma_j$ are disks for every $j$;
 cf. \eqr{e:ulsc2}.

 The first two properties follow immediately
from theorem 5.8 in \cite{CM4} (this theorem combines the
existence of multi-valued graphs near a blow up point and the
sublinear growth of the separation). Namely, since $0 \in \cS$, we
first get a sequence of points $y_j \in \Sigma_j$ with $|A|^2
(y_j) \to \infty$ and $y_j \to 0$.  Since the component of
$B_{r_0}(0) \cap \Sigma_j$ containing $y_j$ is a disk, theorem 5.8
in \cite{CM4} then gives the following two properties:
\begin{itemize}
\item There is a rotation of $\RR^3$ and a subsequence so that $\Sigma_j$ contains a
 $2$-valued
minimal graph $\Sigma_{d,j} \subset \Sigma_j$ defined over an
annulus $D_{r_0 / C} \setminus D_{r_j}$ where $r_j \to 0$.
\item As $j \to \infty$, the $2$-valued graphs close up in the
limit to converge with multiplicity two to a graph $\Gamma_x$ over
$D_{r_0 /C} \setminus \{ 0 \}$ with $ x \in \overline{\Gamma_x}$.
\end{itemize}
Since any subsequence of a convergent sequence has the same limit,
we conclude that $\Gamma_x$ is contained in a leaf of $\cL'$.
Finally, $x$ is a removable singularity for $\Gamma_x$ by a
standard removable singularity result  for minimal graphs.

The cone property follows easily from corollary I.1.9 in
\cite{CM6} which gives a constant $\delta_0>0$ so that if  $B_{2R}
\cap \Sigma_j$ contains a $2$-valued graph in $\{x_3^2 \leq
\delta_0^2\, (x_1^2+x_2^2)\}$ over $D_{R}\setminus D_{r_j}$ and
with gradient $\leq \delta_0$, then each component of
\begin{equation}    \label{c:conepr}
    B_{R/2}\cap
    \Sigma\setminus (\cone_{\delta_0}(0)\cup B_{2 r_j})
\end{equation}
is a multi-valued graph with gradient $\leq 1$.  After possibly
shrinking the radius above given by theorem 5.8 in \cite{CM4}, we
can assume that $\Gamma_x$ is a graph with small gradient  and
hence corollary I.1.9 in \cite{CM6} applies.  It follows that
\begin{equation}  \label{e:lcprjustx}
     B_{r_0} (x)
    \cap \cSu \subset \cone_{\delta_0}(x) \, .
\end{equation}
Finally, the embeddedness of the $\Sigma_j$'s implies that two
limit minimal graphs through nearby singular points must be
disjoint. It is now easy to see that the map which takes a
singular point $y$ to the tangent plane of the limit minimal graph
through $y$ is Lipschitz, giving the last claim.
\end{proof}

Lemma \ref{l:leaf} shows that each point  $x \in \cSu$ is a
removable singularity for a component $\Gamma_x$ of $B_{r_0}(x)
\cap \cL'$ for some $r_0 > 0$. Furthermore, the local cone
property implies that the intersection of $B_{r_0} (x) \cap
\Sigma_j$ with the complement of (a tubular neighborhood of) a
cone $\cone_{\delta'}(x)$ (for some $\delta'
> 0$) consists of two multi-valued graphs for $j$ large (the fact
that there are exactly two is established in proposition II.1.3 in
\cite{CM6}).
 However, it is worth noting that these two properties alone {\underline{do
 not}}
imply that $x$ is a removable singularity for the lamination
$\cL'$, but rather
  there are two possibilities:
\begin{enumerate}
\item[(P)]
The multi-valued graphs in the complement of the cone
$\cone_{\delta'}(x)$ close up in the limit.
\item[(N-P)]
These multi-valued graphs converge to a collection of graphs (such
as $\Gamma_x$) and {\it{at least one}} multi-valued graph that
spirals infinitely on one side of $\Gamma_x$.
\end{enumerate}
In the first case (P) (we will call this ``properness'' below),
the sequence converges to a foliation in a neighborhood of $x$.
The  second case  (N-P) (``not proper'') is illustrated in
\cite{CM15} by a sequence of embedded minimal disks $\Sigma_i$ in
the unit ball $B_1$ with $\partial \Sigma_i \subset
\partial B_1$ where
 the curvatures blow up only at $0$ and
\begin{equation}
    \Sigma_i  \setminus \{ x_3 = 0 \}
\end{equation}
 converges to two
  embedded minimal disks
 \begin{align}
    \Sigma^- &\subset \{ x_3 < 0 \}  \\
 \Sigma^+ &\subset \{ x_3 > 0 \} \, ,
 \end{align}
  each of which spirals
 into $\{ x_3 =
0 \}$ and thus is not proper. Thus, in the example from
\cite{CM15}, $0$ is the first, last, and only point in $\cSu$ and
the limit lamination consists of three leaves: $\Sigma^+$,
$\Sigma^-$, and the punctured unit disk $B_1 \cap \{ x_3 = 0 \}
\setminus \{ 0 \}$.  In this example of (N-P), the limit
lamination cannot be extended smoothly to any neighborhood of $0$.

 {\underline{To summarize}}, \cite{CM15} shows that (N-P) can occur for a sequence of
{\underline{disks}} $\Sigma_i \subset B_{R_i}$ with $\partial
\Sigma_i \subset
\partial B_{R_i}$; however, \cite{CM6} shows that (N-P) cannot  occur
 for disks if the radii $R_i$  go to infinity.

\subsection{Collapsed leaves of $\cL'$}

One of the difficulties is that the leaves of the lamination
$\cL'$ may not be complete; this occurs at points of $\cS$.  We
will begin by analyzing a particular type of incomplete leaf that
we will call {\it collapsed}.

To define this, note that Lemma \ref{l:leaf} shows that each point
$x \in \cSu$ is a removable singularity for a component $\Gamma_x$
of $B_{r_0}(x) \cap \cL'$. We will say that the leaf $\Gamma$ of
$\cL'$ containing $\Gamma_x$ is collapsed:

\begin{Def}     \label{d:coll}
 A leaf
$\Gamma$ of $\cL'$ is {\it collapsed} if there exists some $x \in
\cSu$ so that $\Gamma$ contains the local leaf $\Gamma_x$ given by
Lemma \ref{l:leaf}.
\end{Def}

For a sequence of rescaled helicoids converging to a foliation by
parallel planes away from an axis, every leaf is collapsed.  We
will eventually show that every leaf of $\cL'$ whose closure
contains a point of $\cSu$ is collapsed.  However, it is worth
pointing out that this is not obvious.  For example, in case (N-P)
of the previous section, we get leaves of $\cL'$ that spiral
infinitely into the collapsed leaf but are not themselves
collapsed (we will eventually rule out this possibility using that
the sequence of outer radii is going to infinity).

 We will describe the structure of the collapsed leaves in the rest of this part.
 It is useful to first
define the closure $\barga$ of a leaf $\Gamma$ of $\cL'$ to be the
union of the closures of all bounded (intrinsic) geodesic balls in
$\Gamma$; that is, we fix a point $x_{\Gamma} \in \Gamma$ and set
\begin{equation}    \label{e:closureaa}
    \barga = \bigcup_{r} \overline{ \cB_r (x_{\Gamma}) } \,
    ,
\end{equation}
where $\overline{ \cB_r (x_{\Gamma}) }$ is the closure of $\cB_r
(x_{\Gamma})$ as a subset of $\RR^3$.

Clearly, a leaf $\Gamma$ is complete if and only if $\barga =
\Gamma$ and we always have that
\begin{equation}
    \barga \setminus \Gamma \subset \cS \, .
\end{equation}
The incomplete leaves of $\Gamma$ can be divided into several
types, depending on how $\barga$ intersects $\cS$:
\begin{itemize}
\item
Collapsed leaves, defined in Definition \ref{d:coll}, where
$\barga \cap \cSu$ contains a removable singularity for $\Gamma$.
\item
Leaves $\Gamma$ with $\barga \cap \cSu \ne \emptyset$, but where
$\Gamma$ does not have a removable singularity.  This would occur,
for example, if $\Gamma$ spirals infinitely into the collapsed
leaf through $\barga \cap \cSu$ as in (N-P). (We will eventually
show that this does not occur.)
\item
Leaves $\Gamma$ where $\barga \setminus \Gamma \subset \cSt$;
these won't be considered until Part \ref{p:prove2}.
\end{itemize}

\section{The structure of the collapsed leaves of $\cL'$}       \label{s:tpt}

In the rest of this part, we will describe the structure of the
collapsed leaves of $\cL'$ defined in Definition \ref{d:coll}.
    The most important properties of a collapsed
leaf $\Gamma$ are given in Proposition \ref{p:cole0} below that
describes a neighborhood of the points of $\cSu$ in $\Gamma$.  The
proposition shows that such a $\Gamma$ is stable and that the
closure of $\Gamma$ intersects  $\cSu$ in at most two points.
These results apply without additional assumptions on the sequence
$\Sigma_j$; we will see in the next part that $\Gamma$ has more
structure when we assume in addition that the sequence is ULSC.

\vskip2mm The next proposition establishes the key properties of a
collapsed leaf in the general case:

\begin{Pro}     \label{p:cole0}
Each collapsed leaf $\Gamma$ of $\cL'$ has the following
 properties:
\begin{enumerate} \item  Given {\underline{any}} $y \in \barga \cap
\cSu$, there exists $r_0 > 0$ so that the closure (in $\RR^3$) of
each component of $B_{r_0}(y) \cap \Gamma$ is a compact embedded
 disk  with boundary in $\partial B_{r_0}(y)$.

Furthermore, $B_{r_0}(y) \cap \Gamma$ must contain the component
$\Gamma_y$ given by Lemma \ref{l:leaf} and $\Gamma_y$ is the only
component of $B_{r_0}(y) \cap \Gamma$ with $y$  in its closure.
\item If $\Gamma$ is oriented, then it is  stable. (Otherwise, its
oriented double cover is stable.)
\item
$\barga$ intersects $\cSu$ in at most two points.  If $\barga \cap
\cSu$ contains two points, then the multi-valued graphs in the
$\Sigma_j$'s spiral in opposite directions around the  two
corresponding axes (see Figure \ref{f:7}).
\end{enumerate}
\end{Pro}

\begin{figure}[htbp]
    \setlength{\captionindent}{20pt}
    %\begin{minipage}[t]{0.5\textwidth}
    \centering\input{gen7aa.pstex_t}
    \caption{The  multi-valued graph converging to
        $\Gamma$ in Proposition \ref{p:cole0}.}
    \label{f:7}
    %\end{minipage}
\end{figure}

Properties (2)  and (3) in Proposition \ref{p:cole0} are
self-explanatory.  However, to appreciate property (1), it may be
useful to observe one implication of (1) and to also see an
example of what it rules out. First, (1) implies that $\barga \cap
\cSu$ consists of a discrete set of points and each of these
points is a removable singularity. Second,  recall from (N-P) -
``not proper'' -  that a priori there may be multi-valued graphs
in $B_{r_0}(y) \cap \cL'$ that spiral infinitely into $\Gamma_y$;
(1) above says that these ``infinite spirals'' are
{\underline{not}} contained in any collapsed leaf.

\vskip2mm Throughout this  section $\Gamma$ will be a collapsed
leaf of $\cL'$.  By definition, a leaf $\Gamma$ is a connected
open surface, but may not be complete (and, in fact, collapsed
leaves are incomplete by definition).  We will let
 $K \subset \Gamma$ denote a connected open subset with compact
closure in $\Gamma$.  Finally, $T(K, \epsilon) $ is the
$\epsilon$-tubular normal neighborhood of $K$, i.e.,
 \begin{equation}
    T(K,
\epsilon)  = \{ x  + s \, \nn_{\Gamma} (x) \, | \, x \in K , \,
|s| <
    \epsilon \} \, .
 \end{equation}

 \subsection{Proving property (1) of Proposition \ref{p:cole0}: Isolated removable singularities}  \label{ss:231}
 To prove (1) of Proposition \ref{p:cole0}, we will show the following claim:
 \begin{quote}
 {\underline{Claim}}: If  $x \in \cSu$ is a
singular point in the closure $\barga$ of  a collapsed leaf
$\Gamma$, then the component of $B_{r_0}(x) \cap
\overline{\Gamma}$ containing $x$ is the one from Lemma
\ref{l:leaf}. \end{quote}
  Property (1) then follows from the
following two properties of the component $\Gamma_x$ from Lemma
\ref{l:leaf}:
\begin{itemize}
\item $\Gamma_x \cup \{ x \}$ is a smooth minimal surface.
\item $\overline{\Gamma_x} \cap \cS = x$.
\end{itemize}
Thus, we will have at the same time shown that each ULSC singular
point in the closure of $\Gamma$ is a removable singularity and
has a neighborhood {\underline{in $\Gamma$}} where there are no
other singular points, as desired.

\begin{figure}[htbp]
    \setlength{\captionindent}{20pt}
    %%\begin{minipage}[t]{0.5\textwidth}
    \centering\input{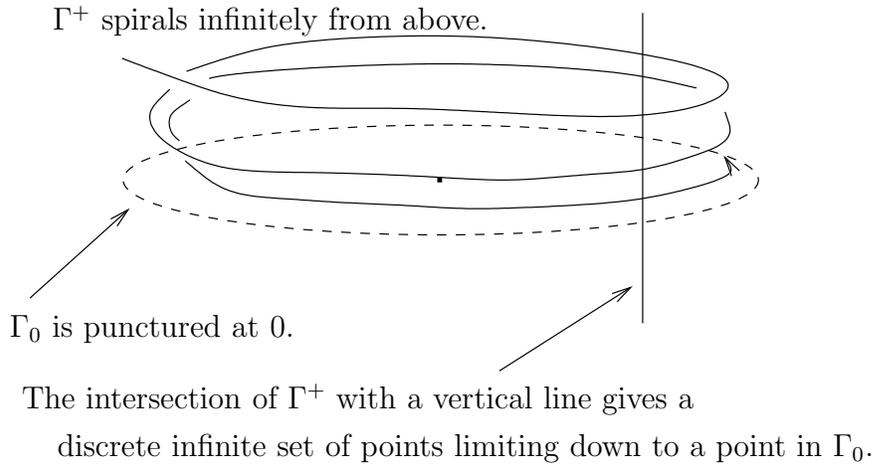}
    \caption{The ``not proper'' example (N-P): The sequence of disks $\Sigma_j$ converges
    in $B_1 \setminus \{0 \}$ to
    a lamination with three leaves: the
punctured disk $\Gamma_0 = D_1 \setminus \{ 0 \}$ (dotted),
$\Gamma^+$ spiralling into $\Gamma_0$ infinitely from above, and
$\Gamma^-$ spiralling into $\Gamma_0$ infinitely from below (not
pictured). The collapsed leaf $\Gamma_0$ is not discrete, but
$\Gamma^+$ and $\Gamma^-$ are.}
    \label{f:fnp}
    %%\end{minipage}
\end{figure}

Lemma \ref{l:onestar} below establishes the above claim about the
singular points in the closure of a collapsed leaf.   The lemma is
best illustrated using the ``not proper'' example in (N-P). In
(N-P), a sequence of embedded minimal disks converges in $B_1
\setminus \{ 0 \}$ to a lamination with three leaves: the
punctured disk $\Gamma_0 = D_1 \setminus \{ 0 \}$, $\Gamma^+$
spiralling into $\Gamma_0$ infinitely from above, and $\Gamma^-$
spiralling into $\Gamma_0$ infinitely from below; see Figure
\ref{f:fnp}.  Notice that all three leaves contain $0$ in their
closure.  The leaf $\Gamma_0$ is collapsed at $0$ (so $0$ is a
removable singularity for $\Gamma_0$), but $\Gamma^+$ and
$\Gamma^-$ cannot extend past the singularity $0$.  The conclusion
of Lemma \ref{l:onestar} is that $\Gamma^+$ and $\Gamma^-$ cannot
be contained in any collapsed leaf of $\cL'$.

The above example from (N-P) also serves to illustrate the idea of
the proof of Lemma \ref{l:onestar}.  Namely, a key distinction
between the collapsed leaf $\Gamma_0$ versus $\Gamma^+$ and
$\Gamma^-$ is that $\Gamma^+$ and $\Gamma^-$ are discrete in the
following sense:
\begin{quote}
Given any point $y$ in $\Gamma^+$ or $\Gamma^-$, there exists
$s>0$ so that $B_s (y) \cap \cL'$ has only one connected component
(i.e., the one containing $y$).
\end{quote}
On the other hand, since $\Gamma^+$ and $\Gamma^-$ spiral
infinitely into $\Gamma_0$, the leaf $\Gamma_0$ is not discrete in
this sense.  Likewise, the description of a neighborhood of a
point in $\cSu$ shows that a collapsed leaf is never discrete.

\begin{Lem}     \label{l:onestar}
Suppose that $x \in \cSu$ and $\Gamma'$ is a component of
$B_{r_0}(x) \cap \cL'$ with $x$ in its closure
$\overline{\Gamma'}$. If $\Gamma'$ is contained in a collapsed
leaf of $\cL'$, then $\Gamma'$ must be the component $\Gamma_x$
given by Lemma \ref{l:leaf}.
\end{Lem}

\begin{proof}
Since the component $\Gamma'$ of $B_{r_0}(x) \cap \cL'$ contains
the point $x \in \cSu$ in its closure,  embeddedness and the cone
property implies that $\Gamma'$ has one of the following two
properties:
\begin{enumerate}
\item[(L1)] $\Gamma'$ is the component $\Gamma_x$
given by Lemma \ref{l:leaf} and hence extends smoothly across $x$.
\item[(L2)] $\Gamma'$ is {\underline{not}} the  component $\Gamma_x$
given by Lemma \ref{l:leaf}.
\end{enumerate}
In Lemma \ref{l:discre} below we will prove that the leaves
satisfying (L2) are discrete in the following sense:
\begin{quote}
Given any point $y$ in a leaf of $\cL'$ satisfying (L2), there
exists  $s>0$ so that $B_s (y) \cap \cL'$ has only one connected
component (i.e., the one containing $y$).
\end{quote}

\vskip2mm {\underline{Completing the proof assuming
discreteness}}:  Suppose now that $\Gamma$ is collapsed, $y \in
\Gamma$, and the ball $B_s (y)$ is disjoint from $\cS$. Let
$\Gamma_{y,s}$ be the component of $B_s(y) \cap \Gamma$ containing
$y$. It follows from the Harnack inequality (since the curvature
is locally bounded on $\Gamma$) that $\Gamma_{y,s}$ is the limit
of distinct leaves of $B_s(y) \cap \cL'$. In particular, $\Gamma$
is {\underline{not}} discrete and hence does not contain any
leaves of $B_{r_0}(x) \cap \cL'$ that satisfy (L2). This completes
the proof of  the lemma modulo  Lemma \ref{l:discre} below.
\end{proof}

The next lemma shows that we always get discreteness for leaves of
$\cL'$ that have a point of $\cSu$ in their closure but are not
collapsed at this point (cf., the picture for $\Gamma^+$ and
$\Gamma^-$ in Figure \ref{f:fnp}).

\begin{Lem} \label{l:discre}
Given any point $y$ in a leaf of $\cL'$ satisfying (L2), there
exists   $s>0$ so that $B_s (y) \cap \cL'$ has only one connected
component (i.e., the one containing $y$).
\end{Lem}

\begin{proof}
Suppose that a component $\Gamma$ of $B_{r_0}(x) \cap \cL'$
contains the point $x \in \cSu$ in its closure but is
{\underline{not}} equal to $\Gamma_x$.   It suffices to find one
point in $\Gamma$ where the leaf  is locally discrete (since the
leaf is connected, the Harnack inequality then implies that every
point is discrete). We will next outline the argument to find this
discrete point.  The key will be to find a sequence of curves
$\gamma_j \subset \Sigma_j$ with uniformly bounded length, where
one sequence of endpoints converges to a point in $\Gamma$, the
$\Sigma_j$'s are {\underline{uniformly}} discrete at the second
endpoint of $\gamma_j$, and the $\gamma_j$'s stay away from the
singular set $\cS$.  These properties are made precise in
(G1)--(G4) below; see Figure \ref{f:fgj}.    Since the
$\gamma_j$'s stay away from $\cS$ and have bounded length, a
subsequence of the $\gamma_j$'s will converge to a curve $\gamma$
in some leaf of $\cL'$.  However, one sequence of endpoints
converges to a point in $\Gamma$ and so the whole curve $\gamma$
is in $\Gamma$. Finally, the second endpoint of $\gamma$ will give
the desired discrete point in $\Gamma$.

\begin{figure}[htbp]
    \setlength{\captionindent}{20pt}
    %%\begin{minipage}[t]{0.5\textwidth}
    \centering\input{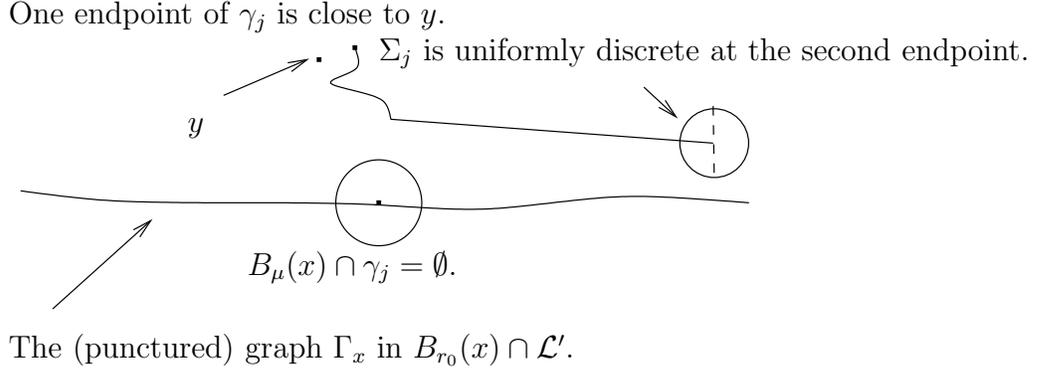}
    \caption{The curves   $\gamma_j$ in $\Sigma_j$.}
    \label{f:fgj}
    %%\end{minipage}
\end{figure}

Before making this precise, we need a few simple preliminaries.
 First, since $\Gamma_x \cup \{ x \}$ separates the ball
$B_{r_0}(x)$ and $\Gamma \subset B_{r_0}(x) \setminus \left(
\Gamma_x \cup \{ x \} \right)$ is connected, we may assume that
$\Gamma$ is contained in the component $B_{r_0}^+(x)$ of
$B_{r_0}(x) \setminus \left( \Gamma_x \cup \{ x \} \right)$ that
is above $\Gamma_x$. Since each point in $\cSu$ comes with a leaf
through it with the same separating properties as $\Gamma_x$, the
distance between $x$ and
 $B_{r_0}^+(x) \cap \cS$  must be positive. Therefore, after shrinking $r_0$, we
may as well assume that
\begin{equation}    \label{e:noothercs}
     B_{r_0}^+(x) \cap \cS = \emptyset \, .
\end{equation}

As mentioned, the key point is to find a sequence of curves
$\gamma_j$ parameterized by arclength
\begin{equation}
    \gamma_j : [0, \ell_j ] \to B_{r_0}^+(x) \cap \Sigma_j
\end{equation}
 with the following properties (see Figure
\ref{f:fgj}):
\begin{enumerate}
\item[(G1)]  The endpoints $\gamma_j (0)$ converge to a point $y
\in \Gamma$.
\item[(G2)]  The lengths $\ell_j$ are uniformly bounded, i.e., $\ell_j \leq \ell$ for every
$j$.
\item[(G3)]  The minimal distance between $\gamma_j$ and $\cS$ is at
least $\mu > 0$.
\item[(G4)]  The $\Sigma_j$'s are ``uniformly discrete'' at the endpoint $\gamma_j (\ell_j)$; precisely,
there exists $\delta > 0$ so that  $B_{\delta}(\gamma_j(\ell_j))
\cap \Sigma_j$ is a (connected) graph over its tangent plane at
$\gamma_j(\ell_j)$ with gradient bounded by one.
\end{enumerate}
The discreteness follows immediately from (G1)--(G4).  Namely,
  (G2) and (G3) imply  that a subsequence of the curves $\gamma_j$ converges to a
curve $\gamma$ contained in a leaf of $\cL'$.  Since the endpoints
$\gamma_j (0)$ converge to the point $y$ in the leaf $\Gamma$, the
entire curve $\gamma$ must be contained in $\Gamma$.  Finally,
(G4) implies that $\Gamma$ is discrete at the second endpoint of
$\gamma$ and, hence, discrete everywhere by the Harnack
inequality.

Before establishing (G1)--(G4), we need to recall
  the following two
additional facts:
\begin{enumerate}
\item
{\underline{Existence of nearby points of large curvature}}:
Given any constants $C_1$ and $C_4$, there exists $\epsilon
> 0$ so that for any $s>0$ and every $j$ sufficiently large
(depending also on $s$) there is a point \begin{equation} q_j \in
B_{s}(x) \cap \Sigma_j \setminus B_{\epsilon \,s}(x)
\end{equation} so that $q_j$ is above $\Gamma_x \cup \{ x\}$ and
$q_j$ satisfies
\begin{equation}    \label{e:largeA}
    |A|^2 (q_j) \geq C_4 \, C_1 \, |x-q_j|^{-2} \, .
\end{equation}
\item
{\underline{Curvature bound away from $x$}}: Given any $\mu > 0$,
there exists a constant $C_2$ so that if $y$ is any point in
$B_{r_0}^+(x) \cap \Sigma_j \setminus B_{\mu} (x)$, then
\begin{equation}    \label{e:smallA}
        |A|^2 (y) \leq C_2 \, .
\end{equation}
\end{enumerate}
Property (1)  was proven in corollary III.3.5 in \cite{CM5}.
Property (2) follows easily since the singular set $\cS$ does not
intersect $B_{r_0}^+(x)$ by \eqr{e:noothercs} (the proof of (2)
can be made precise using  Lemma \ref{l:inftyornot} and a covering
argument).

To complete the proof of discreteness, it suffices to establish
(G1)--(G4).  We will do this next.  First, fix a point $y \in
\Gamma$ in a small ball $B_{h} (x)$ about $x$ ($h$ will need to be
sufficiently small relative to $r_0$ but otherwise does not
matter).  Since $y \in \Gamma$, we can choose a sequence of points
$y_j \in \Sigma_j$ that converge to $y$.  Now choose a constant $s
> 0$ with $s$ much smaller than $|y-x|$.

Observe that property (1) gives points $q_j \in B_{s}(x) \setminus
B_{\epsilon \,s}(x)$ in $\Sigma_j$ satisfying \eqr{e:largeA}.  A
 simple blow up argument (e.g., lemma $5.1$ in \cite{CM4}) then
gives   points  $p_j \in \Sigma_j$ near $q_j$ and radii $r_j$ so
that
\begin{equation}  \label{e:buppj}
    \sup_{B_{r_j}(p_j) \cap \Sigma_j} |A|^2 \leq 4 \, |A|^2 (p_j)
    = 4 \, C_1 \, r_j^{-2} \, ,
\end{equation}
and \begin{equation}
    B_{r_j} (p_j) \subset  B_{ \frac{2\, |x-q_j| }{  \sqrt{C_4} } } (q_j)
      \, . \end{equation}
 In particular, by taking
$C_4$ large in property (1), we can assume that the ratio
\begin{equation}    \label{e:ratsml}
    \frac{r_j}{|p_j-x|}
\end{equation}
is as small as we want and, hence, also that
\begin{equation}
    B_{r_j} (p_j) \subset
    B_{2s}(x) \setminus
B_{\epsilon \,s/2}(x) \, .
\end{equation}
We called the pair $(p_j , r_j)$ a {\it blow up pair} in
\cite{CM6}.  The point about such a pair is that theorem $0.7$ in
\cite{CM6} gives multi-valued graphs
\begin{equation}
\Sigma_j^g \subset
\Sigma_j
\end{equation}
defined outside of a disk of radius $r_j$ centered at $p_j$ and
whose initial separation is proportional to $r_j$.
 On the other hand, since $p_j \notin B_{\epsilon
\,s/2}(x)$, property (2) implies that there is a uniform upper
bound for $|A|^2 (p_j)$ -- and, thus, a uniform lower bound for
the initial scale $r_j$.

We will also need a positive lower bound for the minimum distance
between $\Sigma_j^g$ and $x$. The argument for this is very
similar to an argument in section III.2 of \cite{CM6}.  We will
sketch the argument next.
 The lower bound
follows easily once we have a lower bound for the distance from
$B_{r_j}(p_j)$ to $\Gamma_x$.  Since $p_j \notin B_{\epsilon \,
s/2}(x)$, the one-sided curvature estimate gives a lower bound for
the distance from $p_j$ to $\Gamma_x$ (otherwise $p_j$ would lie
in a narrow cone about $\Gamma_x$ and the one-sided curvature
estimate would contradict \eqr{e:buppj}). Using this and the fact
that $r_j$ is small relative to $|p_j - x|$ (see \eqr{e:ratsml})
then gives the desired lower bound for the distance from
$B_{r_j}(p_j)$ to $\Gamma_x$.  We leave the details to the reader.

To summarize, we have established a positive lower bound for the
distance from $\Sigma_j^g$ to $x$ and for the initial scale $r_j$.
This lower bound on the initial scale also implies a lower bound
for the separation between the sheets of
$\Sigma_j^g$.{\footnote{The existence of some lower bound is easy
and almost obvious; a fairly sharp lower bound is proven in lemma
III.1.6 in \cite{CM6}.}} Moreover, proposition II.1.3 in
\cite{CM6} says that $\Sigma_j$ contains exactly two (oppositely
oriented) multi-valued graphs in this region; the uniform
curvature upper bound given by (2) then also implies a uniform
lower bound for the distance between these two multi-valued
graphs.

As a consequence of these uniform bounds,  a  (sub) sequence of
the two (oppositely oriented) multi-valued graphs is guaranteed to
converge (with multiplicity one) to two multi-valued graphs in
$B_{r_0}^+(x) \cap \cL'$ and these limit multi-valued graphs will
satisfy the same lower bounds.
 Fix
a point $z$ in one of the limit multi-valued graphs.  We will now
find the desired curves $\gamma_j$ from $y_j$ to $z_j$ ((G1) and
(G4) will then automatically be satisfied).  We use two facts to
find these curves.  First, the chord-arc bound of Theorem
\ref{t:1} allows us to connect $y_j$ to the multi-valued graph
$\Sigma_j^g$ by a curve $\gamma_j^+ \subset \Sigma_j$ with length
at most $C_3 \, h$ and, furthermore, we can assume that
$\gamma_j^+$ is {\underline{above}} $\Sigma_j^g$.{\footnote{More
precisely, the curve does not go below the union of $\Sigma_j^g$
and the extrinsic ball $B_{r_j}(p_j)$.}} Now that $\gamma_j^+$
connects $y_j$ to the multi-valued graph, we can use a curve
$\gamma_j^g$ in $\Sigma_j^g$ to connect the endpoint of
$\gamma_j^+$ to points $z_j \in \Sigma_j^g$ converging to $z$; the
$\gamma_j^g$'s automatically have uniformly bounded length and
also stay uniformly away from $x$.
  This completes
the proof of (G1)--(G4) and, consequently, also completes the
proof of discreteness.
\end{proof}

\subsection{Each leaf is a limit of multi-valued graphs in the $\Sigma_j$'s}  \label{ss:232}

Recall that, throughout this section, $\Gamma$ is a leaf of $\cL'$
and $K \subset \Gamma$ is a connected open subset that has compact
closure in $\Gamma$.

 We will first show in Lemma \ref{l:oldclaim} that the
$\Sigma_j$'s are locally graphical over $\Gamma$ in a tubular
neighborhood of $K$.   Corollary \ref{c:leafstuff1} uses the local
description of Lemma \ref{l:oldclaim} to construct multi-valued
graphs $\Sigma_j^g \subset \Sigma_j$ converging to $K$.  Both
Lemma \ref{l:oldclaim} and Corollary \ref{c:leafstuff1} apply to
{\underline{any}} leaf $\Gamma$ and do not require $\Gamma$ to be
collapsed.

 \vskip2mm
The next lemma shows that $\Sigma_j$ is locally graphical over
$\Gamma$ in a small tubular neighborhood of $K$.

\begin{Lem}     \label{l:oldclaim}
  Given any $\delta > 0$, there exist $\epsilon  > 0$ and $J $ so that
 if $j> J $ and $x \in T(K,\epsilon ) \cap \Sigma_j$,
 then  $\cB_{\epsilon }(x) \subset \Sigma_j$ is  a
 graph
 over (a subset of) $\Gamma$ with gradient bounded by $\delta$.
\end{Lem}

\begin{proof}
Since $\Gamma$ is a leaf of $\cL'$, it is disjoint from the
singular set $\cS$.  Therefore, for each point $y \in \Gamma$, the
convergence of the $\Sigma_j$'s to the lamination $\cL'$ away from
$\cS$ gives a ball $B_{\epsilon_y}(y)$ and a $J_y$ so that if $x
\in B_{\epsilon_y}(y) \cap \Sigma_j$ for $j > J_y$, then
$\cB_{\epsilon_y }(x) \subset \Sigma_j$ is  a
 graph
 over (a subset of) $\Gamma$ with gradient bounded by $\delta$.

However, the closure $\bar{K}$ of $K$ in $\Gamma$ is compact, so
it can be covered by a finite subcollection of the half--balls,
i.e.,
\begin{equation}
    \bar{K} \subset \cup_{i=1}^m \, B_{ \frac{\epsilon_{y_i}}{ 2}}(y_i) \,
    .
\end{equation}
It is then easy to see that this implies the lemma with
\begin{equation}
    \epsilon = 1/2 \, \min_{i} \, \epsilon_{y_i} \, .
\end{equation}
\end{proof}

The next corollary uses  Lemma \ref{l:oldclaim} to get
multi-valued graphs $\Sigma_j^g \subset \Sigma_j$ over $K$; see
(A) below.  Furthermore, (C) below shows that $\Sigma_j^g$
contains a point $p_j$ far
 from   the boundary $\partial \Sigma_j^g$ of the multi-valued
graph.  More precisely,   $\partial \Sigma_j^g$ divides naturally
into two parts, depending on whether or not it projects to
$\partial K$; (C) shows that the point $p_j$ is far from the part
of $\partial \Sigma_j^g$ that is not over $\partial K$. We will
later use (C) to get multi-valued graphs with many sheets
converging to a collapsed leaf.

\begin{Cor}         \label{c:leafstuff1}
 Fix a point $p_0 \in K$.
  Given any (small) constant $\delta > 0$ and a (large) constant $N$,
  there   exist $\epsilon  > 0$ and $J $ so that
 for each $j> J $ we get the following:
 \begin{enumerate}
 \item[(A)]
    There is a connected open subset $\Sigma_j^g \subset T(K,\epsilon) \cap
    \Sigma_j$ so that, for each $x \in \Sigma_j^g$,
   the intrinsic ball $\cB_{\epsilon }(x)$ is  a graph
 over (a subset of) $\Gamma$ with gradient bounded by $\delta$.
\item[(B)]
 The normal exponential
map from $K \times (-\epsilon , \epsilon)$ gives a diffeomorphism
to $T(K, 2\,\epsilon)$. Let $\Pi : T(K,\epsilon) \to K$   denote
the projection to $K$ and
  $\Pi_j$  the restriction of $\Pi$ to $\Sigma_j^g$.
\item[(C)]
    There is a point $p_j \in \Sigma_j^g$ with $\Pi_j (p_j) = p_0$ satisfying
    \begin{equation}    \label{e:disttoin}
        \dist_{\Sigma_j^g} ( p_j , \, \partial \Sigma_j^g
        \setminus  \Pi_j^{-1} (\partial K )  ) > N
        \, .
    \end{equation}
\end{enumerate}
\end{Cor}

\begin{proof}
Lemma \ref{l:oldclaim} gives $\epsilon > 0$ (depending only on
$\delta$) so  that for {\underline{every}} point $x$ in
$T(K,\epsilon) \cap
    \Sigma_j$  the intrinsic ball $\cB_{\epsilon }(x)$ is  a graph
 over (a subset of) $\Gamma$ with gradient bounded by $\delta$.

Since $K$ has compact closure in the (open) surface $\Gamma$, we
can shrink $\epsilon > 0$ so that the normal exponential map from
$K \times (-\epsilon , \epsilon)$ gives a diffeomorphism to $T(K,
2\,\epsilon)$.

     Furthermore, since the $\Sigma_j$'s converge to $\cL'$ in a
    neighborhood of the point $p_0 \in \Gamma$, there is a
    sequence of points $p_j \in \Sigma_j$ converging to $p_0$ (in
    fact, there are many such sequences; just pick one).  This
    determines the sequence of multi-valued graphs $\Sigma_j^g \subset
    \Sigma_j$.

     It remains to prove that (C) holds
    for $J$ sufficiently large.  We will do this by contradiction,
    so suppose that no such $J$ exists for some fixed $N$.  In
    particular, we get infinitely many $j$'s where there exist
    curves $\gamma_j \subset \Sigma_j$ with the following properties:
    \begin{itemize}
    \item
    $\gamma_j$ starts at $p_j$ and ends at a point in $\partial T(K,\epsilon)$
    that is distance $\epsilon$
    from $K$.
    \item
    The length of $\gamma_j$ is at most $N$.
    \item
    $\gamma_j$ is contained in $T(K,\epsilon)$.
    \end{itemize}
    After passing to a subsequence, the $\gamma_j$'s must converge
    to a curve
    $\gamma \subset \overline{T(K,\epsilon)}$ that is contained in some leaf of $\cL'$
    (we are using here that $\gamma_j$ stays away from $\cS$).
    Since the $p_j$'s converge to $p_0$, the
    curve $\gamma$ starts at $p_0 \in K$ and, hence, we have
    \begin{equation}
        \gamma \subset  \overline{K} \, .
    \end{equation}
    However, this is impossible since the second endpoints of
    $\gamma_j$ are all distance $\epsilon$ from $K$ and, thus,
    could not have converged to a point in $K$.
This contradiction completes the proof.
\end{proof}

\subsection{Property (2) of Proposition \ref{p:cole0}: Each collapsed leaf is stable} \label{ss:233}

The main result of this subsection is that each oriented collapsed
leaf of $\cL'$ is stable.  The proof of stability has the
following three main steps:
\begin{itemize}
\item
Corollary \ref{c:leafstuff1}  gives multi-valued graphs
$\Sigma_j^g \subset \Sigma_j$ converging to $K$ with large
multiplicity.  The $\Sigma_j^g$'s can be thought of as
single-valued graphs over a (subset of a) covering space $K_j$
over $K$.
\item
  Corollary
\ref{c:leafstuff2} describes the covering spaces $K_j$ by
analyzing the ``holonomy'' action of $\pi_1 (K)$ on the fibers
(the holonomy is defined below).
\item
Lemma \ref{l:seeapp} then shows that a subsequence of the $K_j$'s
satisfies (G1) and (G2) in Appendix B, so we can apply Corollary
\ref{c:collstab} to see that $K$ is stable.
\end{itemize}
 Since this applies
for any such $K$, and $\Gamma$ can be exhausted by such $K$'s by
Lemma \ref{l:simpap} in Appendix \ref{s:a}, we conclude that
$\Gamma$ itself is stable.

\begin{figure}[htbp]
    \setlength{\captionindent}{20pt}
    \centering\input{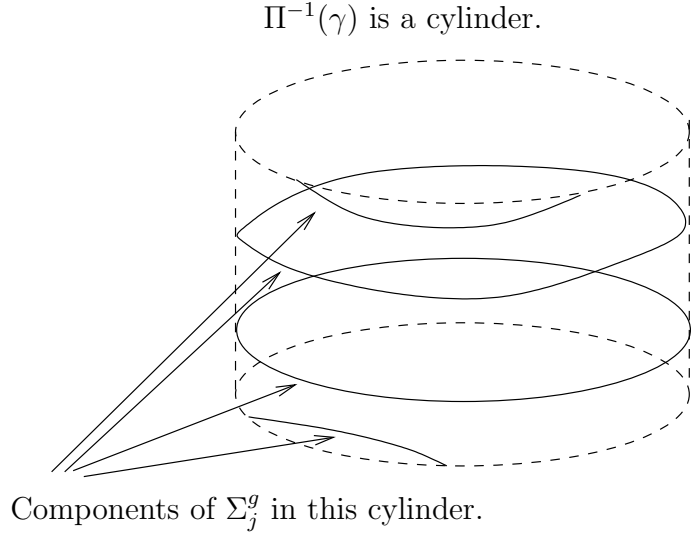}
    \caption{Each component of $\Pi_j^{-1} (\gamma)$ is locally a graph over $\gamma$.}
    \label{f:f1112}
\end{figure}

\vskip2mm The next corollary describes what the multi-valued
graphs $\Sigma_j^g$ look like as we follow them around a simple
closed curve $\gamma$ in $K$.  Obviously, the pre-image
$\Pi_j^{-1} (\gamma)$ consists of a disjoint union of connected
simple curves in the topological annulus $\Pi^{-1} (\gamma)$; see
Figure \ref{f:f1112}.

Some of the  components of $\Pi_j^{-1} (\gamma)$ are more
important than others. To distinguish the components, we will say
that one of these components is ``short'' if it has two endpoints
contained in the same boundary circle of $\Pi^{-1} (\gamma)$;
otherwise, we will say the component is ``long'' (so a long
component either has no boundary, or it has endpoints in distinct
boundary circles of $\Pi^{-1} (\gamma)$). The corollary describes
these long components:

\begin{figure}[htbp]
 \center{Schematic picture of Corollary \ref{c:leafstuff2} (the two boundary circles of
 $\Pi^{-1} (\gamma)$ are dotted):}
        \vskip2mm
    \setlength{\captionindent}{20pt}
    \begin{minipage}[t]{0.5\textwidth}
    \centering\input{gen45.pstex_t}
    \caption{Case (1A): The long components of $\Pi_j^{-1}(\gamma)$ are graphs.}
    \label{f:cleaf1}
    \end{minipage}\begin{minipage}[t]{0.5\textwidth}
    \centering\input{gen46.pstex_t}
    \caption{Case (1B): The long components of $\Pi_j^{-1}(\gamma)$ are multi-valued graphs spiralling together from
    one boundary circle of $\Pi^{-1} (\gamma)$ to the other.}
    \label{f:cleaf2}
    \end{minipage}
\end{figure}

\begin{Cor}         \label{c:leafstuff2}
 Suppose that  $\Pi_j : \Sigma_j^g \to K$ is as in
 Corollary \ref{c:leafstuff1}.

 If
 $\gamma \subset K$ is a simple closed curve,
then either (1A) or (1B) holds: \begin{enumerate} \item[(1A)] Each
long component of $\Pi_j^{-1} (\gamma)$  is closed and is a graph
over $\gamma$; see Figure \ref{f:cleaf1}.
\item[(1B)] The long components of $\Pi_j^{-1} (\gamma)$ are disjoint
simple curves spiralling together from one boundary circle of
$\Pi^{-1} (\gamma)$ to the other; see Figure \ref{f:cleaf2}.
\end{enumerate}
If, in addition,   $K$ contains a simple closed curve $\sigma$
that circles $p \in \barga \cap \cSu$ but is contractible in
$\Gamma \cup \{ p \}$, then
\begin{enumerate}
\item[(2)]
 $\Pi_j^{-1} (\sigma)$ has a single long component{\footnote{Note that even though
 $\Pi_j^{-1} (\sigma)$ is all of where $\Pi^{-1} (\sigma)$
 intersects the multi-valued graph,
$\Pi_j^{-1} (\sigma)$ is
 {\underline{not}}
 all of $\Pi^{-1}(\sigma) \cap \Sigma_j$.  At the least, there must be
 another oppositely oriented
 component of
$\Pi^{-1}(\sigma) \cap \Sigma_j$ that spirals together; cf. the
example of rescaled helicoids.}}; this long component spirals
  from one boundary circle of the
topological annulus  $\Pi^{-1} (\sigma)$ to the other.
\end{enumerate}
\end{Cor}

\begin{proof}
Since we are working in the {\underline{compact}} embedded surface
$\Sigma_j$ and not in the limit, each component of $\Pi_j^{-1}
(\gamma)$ is a simple curve with compact closure.  In particular,
these curves cannot spiral infinitely.  Moreover, since
$\Pi_j^{-1} (\gamma)$ is contained in the multi-valued graph, each
component of $\Pi_j^{-1} (\gamma)$ is also locally a graph over
$\gamma$.

 If any long component is closed
(and hence a graph over $\gamma$), then it separates the two
boundary components of the topological annulus $\Pi^{-1} (\gamma)$
and, by embeddedness, every long component must be a closed graph
over $\gamma$; this is case (1A). Suppose, on the other hand, that
one (and, hence, every) long component connects the two boundary
components of $\Pi^{-1} (\gamma)$.  In this case, the embeddedness
of $\Pi_j^{-1} (\gamma)$ forces all of these curves to spiral
together; this is case (1B).

Suppose now that   a simple closed curve $\sigma \subset K$
circles $p \in \barga \cap \cSu$ but is contractible in $\Gamma
\cup \{ p \}$ and, in particular, does not circle any other points
in $\barga \cap \cSu$. It follows from Lemma \ref{l:leaf} (and its
proof) that the long components of $\Pi_j^{-1} (\sigma)$ do not
close up and, hence, we are in case (1B).  It remains to see that
there is just one long component. This follows immediately from
proposition II.1.3 in \cite{CM6} which shows that $\Pi^{-1}
(\sigma) \cap \Sigma_j$ consists of exactly two oppositely
oriented double spiral staircases.{\footnote{Technically, this
description applies only when $\sigma$ is in a neighborhood of
$p$.  This is sufficient for us since our $\sigma$ is homotopic to
a curve in a neighborhood of $p$ and $\Sigma_j^g$ is locally
graphical over $K$.}} Since $\Sigma_j^g$ is a multi-valued graph
over the connected set $K$, and hence can achieve only one of
these orientations, it can contain only one of these.
\end{proof}

We will say that $K$ is {\it sufficiently large} when it contains
a simple closed curve $\sigma$ that circles exactly one point $p$
in $\barga \cap \cSu$ but is contractible in $\Gamma \cap \{ p
\}$, i.e., when (2) applies in Corollary \ref{c:leafstuff2}.  We
will assume in the rest of this section that $K$ is sufficiently
large.

\begin{Lem}     \label{l:seeapp}
If the $\Sigma_j$'s are planar domains, $K$ is sufficiently large,
and $\gamma \subset K$ is a simple closed curve, then  there can
be only one long curve in (1B) of Corollary \ref{c:leafstuff2} for
$j$ sufficiently large.

More generally, when the $\Sigma_j$'s have bounded genus, then we
get a bound for the number of distinct curves in (1B).
\end{Lem}

\begin{figure}[htbp]
 \center{The contradiction for the proof of Lemma \ref{l:seeapp}:}
        \vskip3mm
    \setlength{\captionindent}{20pt}
    %\begin{minipage}[t]{0.5\textwidth}
    \centering\input{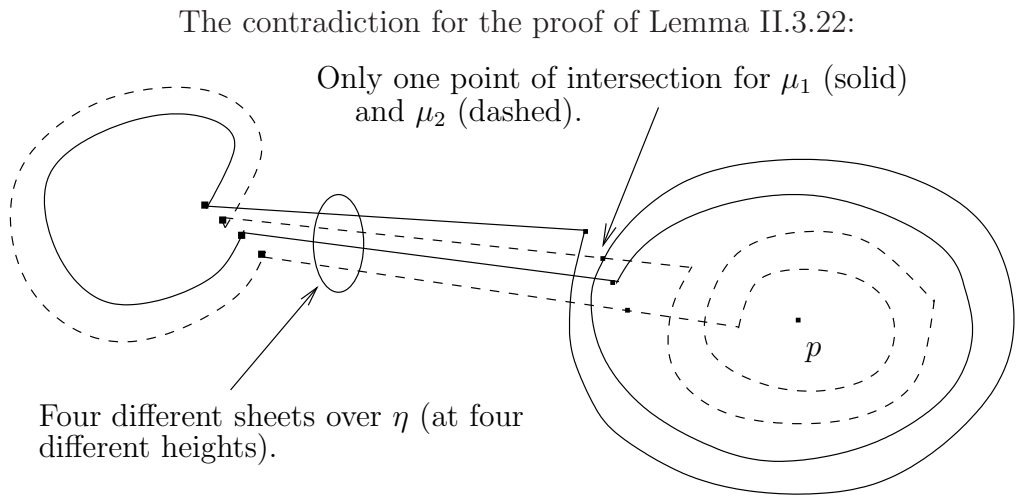}
    \caption{Repeating the construction with $\gamma_1$ in place of $\gamma_2$ gives a second simple closed
    curve $\mu_1$. Perturbing the $1$-valued graph over $\sigma$ slightly outside of
    $\sigma$,  $\mu_1$ and $\mu_2$ intersect in exactly one point and do so
transversely.   Hence, $\mu_1$ and $\mu_2$ have linking number
one, which is impossible in the planar domain $\Sigma_j$.}
    \label{f:sigetagam3}
    %\end{minipage}
\end{figure}

\begin{proof}
We will give the proof for genus zero, i.e., when the $\Sigma_j$'s
are planar domains; the  easy modifications needed for the general
case are left to the reader.

Let $\sigma \subset K$ be a simple closed curve circling $p \in
\barga \cap \cSu$ and so $\sigma$ is contained in a small
neighborhood of $p$ (this exists since $K$ was assumed to be
sufficiently large). Let $\gamma \subset K$ be a second simple
closed curve. After possibly perturbing $\sigma$ slightly, we can
assume that it is disjoint from $\gamma$.  Fix points $x \in
\sigma$ and $y \in \gamma$ and let $\eta \subset K$ be a simple
curve from $x$ to $y$.  Again, after perturbing things, we can
assume that $\eta$ intersects $\sigma$ and $\gamma$ only at its
endpoints $x$ and $y$; see Figure \ref{f:sigetagam}.

\begin{figure}[htbp]
    \setlength{\captionindent}{20pt}
    %%\begin{minipage}[t]{0.5\textwidth}
    \centering\input{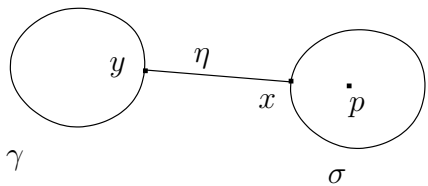}
    \caption{The proof of Lemma \ref{l:seeapp}: The curves $\sigma$, $\eta$, and $\gamma$ in $K$.}
    \label{f:sigetagam}
    %%\end{minipage}
\end{figure}

Suppose now that $\Sigma_j^g \subset \Sigma_j$ contains two
distinct long curves, $\gamma_1$ and $\gamma_2$, in $\Pi_j^{-1}
(\gamma)$ that spiral together; see Figure \ref{f:sigetagam1}.  We
will show that this leads to a contradiction by constructing two
simple closed curves, $\mu_1$ and $\mu_2$, in $\Sigma_j$ that have
linking number one in $\Sigma_j$. This is impossible for a planar
domain (it implies that the genus is at least one).

We will first construct the curve $\mu_2 \subset \Sigma_j$ out of
four parts; see Figure \ref{f:sigetagam2}.  The first part of
$\mu_2$ is a $1$-valued graph over $\gamma$ that is contained  in
$\gamma_2$ and has both of its endpoints over $y$.  These two
endpoints are distinct since they are at different heights over
$y$.  The next two parts of $\mu_2$ are graphs over $\eta$ that
connect these two endpoints to two distinct points over $x \in
\sigma$.  Finally, we close the curve up by connecting the two
points over $x$ by a multi-valued graph over $\sigma$.  Here, we
have used that $\Pi_j^{-1} (\sigma)$ has {\underline{exactly}} one
long component to show that these endpoints can be connected and
to see that the curve connecting them  is at least $2$-valued.
Furthermore, we have also implicitly used that $j$ is large to
ensure that we can find the graphs over $\eta$ and to ensure that
the endpoints of these over $x$ lie in a long component of
$\Pi_j^{-1} (\sigma)$.  (We will use that $j$ is  large in the
same way later in the paper, usually without mentioning that we
are doing so.)

\begin{figure}[htbp]
 \center{The proof of Lemma \ref{l:seeapp}:}
        \vskip3mm
    \setlength{\captionindent}{20pt}
\begin{minipage}[t]{0.5\textwidth}
    \centering\input{gen42.pstex_t}
    \caption{Two curves $\gamma_1$ and $\gamma_2$ in $\Sigma_j$ spiral together over $\gamma$, but only one curve
    spirals over $\sigma$.}
    \label{f:sigetagam1}
    \end{minipage}\begin{minipage}[t]{0.5\textwidth}
    \centering\input{gen43.pstex_t}
    \caption{The (dashed) simple closed curve $\mu_2$ in $\Sigma_j$ has four parts:
        \newline A $1$-valued graph over $\gamma$ in $\gamma_2$.
        \newline A multi-valued graph over $\sigma$.
        \newline Two graphs  over $\eta$.}
    \label{f:sigetagam2}
    \end{minipage}
\end{figure}

The curve $\mu_1 \subset \Sigma_j$ is constructed similarly, with
two notable differences; see Figure \ref{f:sigetagam3}.  First,
the $1$-valued graph over $\gamma$ is chosen to be in $\gamma_1$
this time, as opposed to $\gamma_2$ before.  Consequently, the
one-valued graphs over $\gamma$ in $\mu_1$ and $\mu_2$ are
disjoint and, furthermore, the graphs over $\eta$ are at four
distinct heights.  Second, instead of closing $\mu_1$ up with a
multi-valued graph over $\sigma$, do it over a slight outward
perturbation of $\sigma$ (see Figure \ref{f:sigetagam3}).  This
makes the two ``closing up'' curves for $\mu_1$ and $\mu_2$
disjoint.  However, since there is just one long component over
$\sigma$ (and also over its slight outward perturbation), we see
that the ``closing up'' curve for $\mu_1$ must cross one of the
graphs over $\eta$ in $\mu_2$.  Moreover, this intersection is
transverse and the curves are otherwise disjoint.  This implies
that $\mu_1$ and $\mu_2$ have linking number one in $\Sigma_j$,
which gives the desired contradiction.
\end{proof}

\vskip2mm
 Each $\Sigma_j^g$ is a multi-valued graph over $K$, but
can be thought of as a single-valued graph over a domain $K_j$ in
some covering space of $K$.  However, this covering space may
depend on $j$.  Therefore, in order to apply the results of
Appendix B, we need to pass to a subsequence so that:
\begin{itemize} \item
The $K_j$'s all lie in the same covering space $\hat{K}$
(independent of $j$). \item  The $K_j$'s exhaust $\hat{K}$. \item
The holonomy group of the covering space $\hat{K}$ is $\ZZ$ (the
definition of the holonomy group is recalled below).
\end{itemize}

In order to achieve these three points, we need a few elementary
facts about covering spaces.  First, recall that a covering space
$\hat{\Pi}: \hat{K} \to K$ with base point $x \in K$ is uniquely
determined by the {\it holonomy homomorphism} $\Hol$ from $\pi_1
(K)$ to the automorphisms of the fiber $\hat{\Pi}^{-1} (x)$.  To
define this homomorphism, suppose that
\begin{equation}
\gamma:[0,1] \to K
\end{equation}
 is a curve with $\gamma (0) =
\gamma (1) =
  x$ and    $\hat{x}$ is a point in $\hat{\Pi}^{-1} (x)$.  The
  lifting property for covering spaces gives a unique
  lift{\footnote{Recall that $\gamma_{\hat{x}}$ is said to be a
  {\it lift} of $\gamma$ if $\gamma = \hat{\Pi} \circ
  \gamma_{\hat{x}}$.  The lifting property for covering spaces
  says that we get a unique lift of $\gamma$ for each choice of
 point $\hat{x}$ with $\hat{\Pi} (\hat{x}) = x$.}}
  \begin{equation}
    \gamma_{\hat{x}}:[0,1] \to \hat{K}
  \end{equation}
   of $\gamma$ with
  $\gamma_{\hat{x}}(0) = \hat{x}$.  We define $\Hol (\gamma)
  (\hat{x}) $ to be the endpoint  $\gamma_{\hat{x}}(1)$.  Finally,
   define the holonomy group to be the image $\Hol (\pi_1 (K))$.

\vskip2mm We are now ready to prove that each oriented collapsed
leaf is stable:

\begin{proof}
(of (2) in Proposition \ref{p:cole0}). We will show that any
connected open subset $K \subset \Gamma$ that has compact closure
in $\Gamma$ and is sufficiently large must be stable.  Since
$\Gamma$ can be exhausted by such $K$'s by Lemma \ref{l:simpap} in
Appendix \ref{s:a}, we will conclude that $\Gamma$ itself is
stable.

Fix a point $x \in K$.  By repeatedly applying Corollary
\ref{c:leafstuff1} with $\delta = 1/j$ and passing to a
subsequence, we get a sequence of connected multi-valued graphs
$\Sigma_j^g$ over $K$, covering spaces $\Pi_j : \hat{K}_j \to K$,
domains $K_j \subset \hat{K}_j$, and functions $u_j: K_j \to \RR$
with
 \begin{equation}  \label{e:hh1aa}
    |u_j| + |\nabla u_j| \leq 1/j \, ,
 \end{equation}
 so that there is a bijection  from $K_j$ to $\Sigma_j^g$
 given by
 \begin{equation}
     x   \, \to  \,  \Pi_j (x) + u_j (x) \, \nn_{\Gamma} (\Pi_j (x))
      \, .    \label{e:hh2aa}
 \end{equation}
Furthermore, \eqr{e:disttoin} gives   a point $x_j \in \Sigma_j^g$
with $\Pi_j (x_j) = x$ satisfying
    \begin{equation}    \label{e:disttoinaaa}
        \dist_{\Sigma_j^g} ( x_j , \, \partial \Sigma_j^g
        \setminus  \Pi_j^{-1} (\partial K )  ) > j
        \, .
    \end{equation}

We must do two things in order to apply Corollary \ref{c:collstab}
in Appendix B.  Namely, we must pass to a subsequence so that the
$K_j$'s all sit in the same covering space $\hat{K}$ and we must
show that the holonomy group of $\hat{K}$ is $\ZZ$.  Once we have
done these, \eqr{e:disttoinaaa} will imply that the $K_j$'s
exhaust $\hat{K}$.

We will deal with the second one first, i.e., we will show that
the holonomy group is always $\ZZ$.  This follows immediately from
Corollary \ref{c:leafstuff2}.  Namely,  (2) in Corollary
\ref{c:leafstuff2}
   implies that the fiber over $x$ in each $K_j$ can be identified with $\ZZ$
   and the holonomy from circling the point $p \in \cSu$ is just $n \to (n+1)$ or  $n \to
   (n-1)$, depending on whether the multi-valued graph spirals up or down.
   Suppose now that $\gamma$ is a simple closed curve through $x$
   representing a homotopy class $[\gamma]$ in $\pi_1 (K)$.
   Furthermore, (1A) and (1B) in Corollary \ref{c:leafstuff2}
   imply that  either:
   \begin{itemize}
   \item If (1A) holds, then $\Hol ([\gamma])$ is the identity map, i.e.,  $n \to
   n$.
   \item
   If (1B) holds, then $\Hol ([\gamma])$ maps to $n \to (n\pm k)$, where $k$ is the number of
   disjoint curves spiralling together in (1B).
   \end{itemize}
    In particular, the  image of the holonomy is always in $\ZZ$
    in either case.

   Finally, we will use Lemma \ref{l:seeapp} to prove that only a finite
   set of distinct covering spaces arise as one of the $\hat{K}_j$'s and, consequently,
   one of the $\hat{K}_j$'s occurs infinitely many times.
   We have already established that each holonomy group is $\ZZ$,
   but the covering space is determined by the holonomy homomorphism (and not just
   the group).  Each holonomy homomorphism
   \begin{equation}
        \Hol_j : \pi_1 (K) \to \ZZ
   \end{equation}
   is determined by the image of a fixed (finite) set of
   generators $\gamma_1 , \dots , \gamma_m$ of $\pi_1 (K)$, so we
   need only to show a uniform bound for $\Hol_j (\gamma_n)$ for
   every $j$ and $n$.  However,  (1B) in Corollary
   \ref{c:leafstuff2} implies that $\Hol_j (\gamma_n)$ is just the
   number of disjoint (long) curves  spiralling together in
   $\Pi_j^{-1} (\gamma_n)$ and Lemma \ref{l:seeapp} bounds this
   uniformly, completing the proof.
\end{proof}

\begin{Rem}     \label{r:orient}
We have assumed throughout this subsection that the leaf $\Gamma$
is oriented.  When this is not the case, the same argument applies
to show that the oriented double cover is stable.
\end{Rem}

 \subsection{Property (3) of Proposition \ref{p:cole0}: Opposite orientations at distinct points of $\barga \cap \cSu$}
  \label{ss:nothree}

\begin{proof}
(of (3) in Proposition \ref{p:cole0}).  We must show that if $p$
and $q$ are distinct points in $\barga \cap \cSu$, then the
    multi-valued graphs in $\Sigma_j$ near $p$ spiral in the opposite direction as the ones near
    $q$. Since there are only two possible directions, this
    implies that $\barga \cap \cSu$ contains at most two points (if there were three such points,
    then two would have to be oriented the same way which we will show is impossible).

We will argue by contradiction, so suppose that the
    multi-valued graphs near $p$ and $q$ have the same
    orientation.  In this case, we can choose a  closed  ``figure eight''
    curve $\gamma_j$ in $\Sigma_j$ with the following properties (see Figure \ref{f:3}):
    \begin{itemize}
    \item $\gamma_j$ is a  graph over a fixed (immersed) figure
        eight curve $\gamma$ in $\Gamma$ which circles $p$ and $q$ in opposite
        directions.  Let $r \in \gamma$ be the double point where
        $\gamma$ is not embedded.
    \item The two points in $\gamma_j$ above the double point $r$ are in
        distinct sheets of $\Sigma_j$; hence $\gamma_j$ is embedded.
    \end{itemize}

\begin{figure}[htbp]
    \setlength{\captionindent}{20pt}
    \begin{minipage}[t]{0.5\textwidth}
    \centering\input{gen3.pstex_t}
    \caption{The figure eight curves $\gamma^j$ in $\Sigma_j$.}
    \label{f:3}
    \end{minipage}\begin{minipage}[t]{0.5\textwidth}
    \centering\input{gen4.pstex_t}
    \caption{The stable surface $\Gamma_j$ would be forced to cross an axis.}
    \label{f:4}
    \end{minipage}
\end{figure}

The second condition has a very useful consequence.  Namely, the
unit normal to $\Sigma_j$ is always either upward or downward
pointing along $\gamma_j$ since
 $\Sigma_j$ is graphical along $\gamma_j$; therefore, elementary topology implies that:
 \begin{itemize}
 \item
  The two
 points in $\gamma_j$
  above $r$ are separated by an oppositely--oriented
     sheet of $\Sigma_j$.
  \end{itemize}
  We will now use these properties of the $\gamma_j$'s to
  find stable minimal surfaces $\Gamma_j$ disjoint from the
  $\Sigma_j$'s which contain a graph near either
  $p$ or $q$, contradicting that these points are in $\cSu$.
    Since $\Sigma_j$ has genus zero, the curve $\gamma_j$ separates in $\Sigma_j$; let $\Sigma_j^+$ be one of the two
     components of
    $\Sigma_j \setminus \gamma_j$.    Since $B_{R_j} \setminus \Sigma_j$ is mean convex in the sense of Meeks-Yau,
    the existence theory of \cite{MeYa2} gives a stable
    embedded minimal planar domain
\begin{equation}
    \Gamma_j^+ \subset B_{R_j} \setminus \Sigma_j
    {\text{ with }}
    \partial \Gamma_j^+ = \partial \Sigma_j^+ \, .
\end{equation}
      Let
    $\Gamma_j$    be
 the component of  $\Gamma_j^+$ with $\gamma_j \subset \partial
 \Gamma_j$.
 Using estimates for stable surfaces (\cite{Sc1}, cf. \cite{CM2})
  and the fact that
    $\gamma_j$ is a figure eight, it is now not hard to see that
    $\Gamma_j$ must contain a graph near either
  $p$ or $q$ (see Figure \ref{f:4}).  This can  be seen as
  follows:
    \begin{enumerate}
    \item[1.]
    After leaving the upper portion of $\gamma_j$ over $r$,  the stable surface $\Gamma_j$
      is separated
    from the lower portion of $\gamma_j$ by an oppositely oriented
     sheet of $\Sigma_j$ and, hence, $\Gamma_j$ has an a priori curvature
     bound there
    by \cite{Sc1}, cf. \cite{CM2}.

    To see that $\Gamma_j$ does
    indeed leave the upper portion of $\gamma_j$ over $r$,
    intersect $\Gamma_j$ with a large transverse{\footnote{The application of transversality uses the
regularity of $\tilde \Gamma_j$ up to the interior of $\tilde
\gamma$.  Local boundary regularity was established for
two-dimensional minimal surfaces in \cite{Hi}.}} ball $B_R(z)$ to
    get a collection of closed curves and one segment $\sigma_j$
    where $\sigma_j$ connects the upper and lower portions of
    $\gamma_j$ (see Figure \ref{f:4a}).
    Since these are separated near $r$ by an oppositely--oriented
     sheet of $\Sigma_j$, the segment $\sigma_j$ moves away from
     $\gamma_j$ as desired.
    \item[2.]  Away from the singular points $p$ and $q$, the surface $\Gamma_j$ is
    locally pinched between sheets of $\Sigma_j$.  Combining this pinching with
    the curvature bound from $1.$  implies that
     $\Gamma_j \to \Gamma$ away from $p$, $q$, and
     $\gamma_j$. (Here ``away'' is with respect to distance along
     paths in  $B_{R_j} \setminus \Sigma_j$.)
    \item[3.]    Combining the a priori bound of $1.$ with the flatness given by $2.$, unique
    continuation forces $\Gamma_j \to \barga$ even as it approaches $p$
    or $q$ (this unique continuation argument is spelled out in lemma II.1.38 in \cite{CM5}).
    However, the one-sided curvature estimate, i.e., Theorem \ref{t:t2},
    would then apply to the $\Sigma_j$'s near $p$ or $q$,
      contradicting that $|A| \to \infty$ near $p$ and $q$.
    \end{enumerate}
   This contradiction shows that the
    multi-valued graphs near $p$ and $q$ are oppositely--oriented,
    completing the proof of (3).
    \end{proof}

\begin{figure}[htbp]
    \setlength{\captionindent}{20pt}
    \centering\input{gen4a.pstex_t}
    \caption{The stable surface moves away from its boundary near $r$.}
    \label{f:4a}
\end{figure}

This completes the proof of Proposition \ref{p:cole0}.

\vskip2mm
\begin{Rem}     \label{r:genpoints}
The genus bound on the $\Sigma_j$'s can be used to directly see
that $\barga \cap \cSu$ cannot contain three points.  To see this,
suppose that $p,q$, and $r$ are three distinct points in $\barga
\cap \cSu$ and $\gamma_{pq}$ is a geodesic in $\Gamma$ from $p$ to
$q$.   For $j$ large, Theorem \ref{t:1} allows us to find simple
closed curves $\gamma^j_{pq} \subset \Sigma_j$ with the following
properties:
\begin{itemize}
\item $\gamma^j_{pq}$ is contained in the $\epsilon$-tubular
neighborhood of $\gamma_{pq}$.
\item $\gamma^j_{pq} \setminus \left( B_{\epsilon}(p) \cup B_{\epsilon}(q) \right)$
 consists of two graphs over $\gamma_{pq}$ which
are in distinct sheets of $\Sigma_j$.
\end{itemize}
Since $\Sigma_j$ has genus zero, the curve $\gamma^j_{pq}$ must
separate $\Sigma_j$ into two distinct components.  However, it is
easy to see that this is impossible by using  the local connecting
property near the third point $r$. Namely, we can take two points
near $p$ on opposite sides of $\gamma^j_{pq}$ and connect each of
them to $B_{\epsilon}(r)$ by curves in $\Sigma_j$ which do not
intersect $\gamma^j_{pq}$. These two curves can then be connected
to each other in $B_{\epsilon}(r) \cap \Sigma_j$, giving the
desired contradiction.
\end{Rem}

\part{When the surfaces are ULSC: The proof of Theorem \ref{t:t5.1}}
\label{p:prove1}

\setcounter{equation}{0}

In this part, we will prove Theorem \ref{t:t5.1}, i.e., the main
structure theorem for ULSC sequences where $\cSt = \emptyset$. The
key will be to analyze the ULSC singular set $\cSu$ and, in
particular, the collapsed leaves of $\cL'$.  Although the emphasis
will be on the ULSC case, many of the arguments will actually
apply to a neighborhood of a collapsed leaf whose closure does not
intersect $\cSt$. This will be used later when we analyze the
general case.

In the previous section, we showed that a collapsed leaf $\Gamma$
of $\cL'$ is a stable, incomplete minimal surface with isolated
removable singularities at points in $\cSu$.  In general, $\Gamma$
may have worse singularities at points of $\barga \cap \cSt$, but
we will  assume that $\barga \cap \cSt = \emptyset$ in this part.

In addition to what we have shown in the previous section,  we
need to establish two facts  to complete the proof of Theorem
\ref{t:t5.1}. First, we must show that every collapsed leaf is a
plane (it then follows easily from embeddedness that all of these
planes are parallel). Since we have shown in Proposition
\ref{p:cole0} that the collapsed leaves are stable with isolated
removable singularities at each point of $\cSu$, this follows
easily from the Bernstein theorem for complete stable surfaces in
$\RR^3$. The second additional fact that must be established is
the ``properness'' of the limit in the sense of \cite{CM7}.
Roughly speaking, the local cone property already implies that the
closed set $\cS$ is {\underline{contained}} {\underline{in}} two
Lipschitz curves each of which is transverse to the limit planes.
The properness consists of showing that $\cS$ actually fills out
these curves completely, i.e., there cannot be a first or last
point in $\cS$. See ($\star$) in Section \ref{s:star} for the
precise statement. As in \cite{CM7}, we will prove properness by
showing that  the vertical flux  of a potential non-proper limit
would have to be positive, which is impossible by Stokes' theorem.

All of this will show   the following:
\begin{itemize}
\item
The ULSC sequence of surfaces converges to the foliation by
parallel planes
\begin{equation}
    \cF = \{ x_3 = t \}_t
\end{equation}
     away from
the singular set $\cS$.
\item
($C_{ulsc}$) from Theorem \ref{t:t5.1} holds.
\item
($D_{ulsc}'$): $\cS$ consists of two disjoint Lipschitz graphs
$\cS_1: \RR \to \RR^3$ and  $\cS_2: \RR \to \RR^3$ over the
$x_3$--axis.
\end{itemize}
From this, it follows immediately from the main theorem of
\cite{Me1}  that $\cS_1$ and $\cS_2$ are in fact straight lines
{\it{orthogonal}} to  the leaves of the foliation, giving
($D_{ulsc}$) from Theorem \ref{t:t5.1} and completing the proof of
Theorem  \ref{t:t5.1}.

\vskip2mm Recall that collapsed leaves are the leaves of $\cL'$
that ``go through'' a point of $\cSu$, i.e., that   contain the
local leaf $\Gamma_x$ given by Lemma \ref{l:leaf} for some $x \in
\cSu$; see Definition \ref{d:coll}. The next proposition
establishes the key properties of a collapsed leaf in the ULSC
case:

\begin{Pro}     \label{p:cole}
Suppose that $\Gamma$ is a collapsed leaf of $\cL'$.  If $\barga
\cap \cSt = \emptyset$, then
\begin{enumerate}
\item $\barga$  is  a
plane.
\end{enumerate}
If, in addition, $\cSt = \emptyset$ (i.e., the sequence is ULSC),
then
\begin{enumerate}
\item[(2)]  $\barga$ intersects $\cSu$ in exactly two points and the multi-valued
graphs in the $\Sigma_j$'s spiral in opposite directions around
the  two corresponding axes (see Figure \ref{f:7}).
\end{enumerate}
\end{Pro}

This proposition will be proven over the rest of this section.

  \subsection{Property (1) in Proposition
\ref{p:cole}: Collapsed leaves are planar}
  \label{ss:cflat}
 To prove that $\Gamma$ is flat,  we first use property (3)
 in Proposition \ref{p:cole0} to see that $\barga$ is the union of
 $\Gamma$ together with at most two points in $\cSu$ since we are
 assuming that $\barga
\cap \cSt = \emptyset$.  In particular, since each point in
$\barga \cap \cSu$ is a removable singularity by (1)
 in Proposition \ref{p:cole0}, we conclude that
$\barga$ is a smooth complete surface without boundary.

 Assuming first that $\Gamma$ is oriented, (2)
 in Proposition \ref{p:cole0} implies that $\Gamma$ is  stable.
 We can then
 use a
 standard logarithmic cutoff argument at each point in $\barga \setminus
 \Gamma$ to conclude that $\barga$
 is itself stable.  The Bernstein theorem for
 stable complete minimal surfaces, \cite{FiSc}, \cite{DoPe}, then implies
that $\barga$ is a plane, as desired.  When $\Gamma$ is not
oriented, the preceding discussion applies to show that its
oriented double cover is flat -- and hence so is $\Gamma$. The
obvious details are left to the reader.

 \subsection{Property (2) in Proposition
\ref{p:cole}: Ruling out just one point of $\cSu$ in a
 leaf} \label{ss:notone}
 In contrast to property (1), we will need to use that the sequence
is ULSC in order to prove (2).
 We will later see that this
assumption can be removed. However, the argument we will give to
prove (2) in general will use the ULSC case that we are proving
now (this is why we are not proving the general case directly).

 We have
shown in property (3) of Proposition \ref{p:cole0} that the
closure of a collapsed leaf contains at most two ULSC singular
points and that the $\Sigma_j$'s spiral in opposite directions
around two such points.  Hence, to prove (2) in Proposition
\ref{p:cole}, we must show that $\barga$ cannot intersect $\cS$ in
just one point.

\begin{figure}[htbp]
    \setlength{\captionindent}{20pt}
    \centering\input{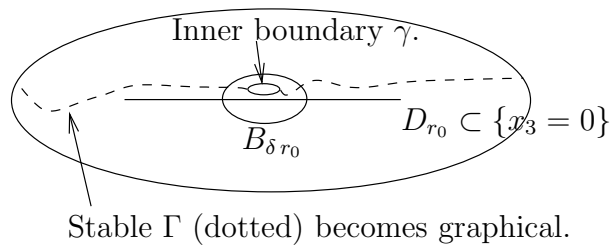}
    \caption{Lemma \ref{l:stab}: The stable surface $\Gamma$ starts off close to - but above - a disk and
    is forced to become graphical.}
    \label{f:stab}
\end{figure}

\vskip2mm
 Before proving (2), we need to recall a useful
property of stable minimal surfaces.  Namely, the following lemma
shows that a stable surface that starts out on one side of a plane
where the interior boundary is in a small ball is graphical away
from its boundary (see Figure \ref{f:stab}):

\begin{Lem} \label{l:stab}
There exists a small constant $0 < \delta < 1$ so that if $r_0 <
\delta \, R_0$ and $\Gamma \subset B_{R_0}$ is a connected
embedded stable minimal planar domain with non-empty inner
boundary $\gamma =
\partial \Gamma \setminus
\partial B_{R_0}$ contained in the small ball
$B_{\delta \, r_0}$, outer boundary $\partial \Gamma \setminus
B_{\delta \, r_0}$ non-empty, and
\begin{equation}
    B_{r_0} \cap \Gamma \cap \{ x_3 = 0 \}  = \emptyset \, ,
\end{equation}
 then $\Gamma$ contains a graph over the annulus $D_{\delta \, R_0} \setminus D_{r_0} \subset
 \{ x_3 = 0 \}$.  Moreover, this graph can be connected to the
 inner boundary
 $\gamma$ by a curve in $B_{2r_0} \cap \Gamma$.
\end{Lem}

\begin{proof}
The proof has two steps.  Namely, we first show that $\Gamma$
contains an initial graph over a small annulus on the scale of
$\delta \, r_0$.  The next step uses the initial graph to apply
the ``stable graph proposition'' - Proposition \ref{l:lext0} in
Appendix \ref{s:aD} - to get the desired graph over the large
annulus $D_{\delta \, R_0} \setminus D_{r_0}$.

{\underline{Producing the initial graph}}: Observe first that the
a priori curvature estimates for stable surfaces of \cite{Sc1},
\cite{CM2} and the gradient estimate imply that each point in
$\Gamma$ close to $D_{r_0}$ - but outside of $B_{2 \,\delta \,
r_0}$ - is graphical (see, e.g., lemma I.0.9 in \cite{CM3} for a
precise statement and proof). Moreover, a standard catenoid
barrier argument (lemma $3.3$ in \cite{CM8}) guarantees the
existence of such ``low points'' that can be connected to the
 inner boundary
 $\gamma$ by a curve in $B_{2 \, C_1 \, \delta \, r_0} \cap \Gamma$.
Starting from one such low point, we can then build up a graph
$\Gamma^g$ over the annulus $D_{2 \, C_1 \, \delta \, r_0}
\setminus D_{C_1 \, \delta \, r_0}$ by
 applying the Harnack inequality a fixed number of times, thereby
 staying low as we go around (we get a
graph and not a multi-valued graph since the surface $\Gamma$ is
embedded and does not spiral infinitely).{\footnote{The constant
$C_1$ depends on the constants from the Harnack inequality and
lemma $3.3$ in \cite{CM8} but is independent of $\delta$.}}

{\underline{Applying Proposition \ref{l:lext0} to get the graph
from $r_0$ to $\delta \, R_0$}}: Let $\hat{\gamma}$   be the graph
in $\Gamma^g$ defined over the inner boundary $\partial D_{C_1 \,
\delta \, r_0}$.   The simple closed curve $\hat{\gamma}$
separates the planar domain $\Gamma$ into two components; let
$\barga$ be the ``outer'' one, so that $\gamma$ is not in
$\partial \barga$. Since the tubular neighborhood $\an_{C_1 \,
\delta \, r_0} (\hat{\gamma})$ of radius $C_1 \, \delta \, r_0$ of
$\hat{\gamma}$ in $\barga$ is contained in the graph $\Gamma^g$,
we get uniform bounds for the area and total curvature of
$\an_{C_1 \, \delta \, r_0} (\hat{\gamma})$.  We can now apply
Proposition \ref{l:lext0} to $\hat{\gamma}$ to get a graph in
$\Gamma$ defined over the annulus $D_{R_0/\omega} \setminus
D_{\omega \, C_1 \, \delta \, r_0}$. The lemma follows by taking
$\delta > 0$ sufficiently small.
\end{proof}

We note next that the local cone property has two important
consequences for the singular set $\cSu$, and in particular, for
how the singular set changes as we move from one collapsed leaf to
the next:
\begin{enumerate}
\item[(S1)]
$\cSu$ cannot run off to infinity.
\item[(S2)]
Distinct points of $\cSu$ in a collapsed leaf cannot combine in
another collapsed leaf.
\end{enumerate}

\begin{proof}
(of  property (2) in Proposition \ref{p:cole}).
 To prove  that $\barga \cap \cSu$ contains
exactly two points, suppose for a moment that there was only one
singular point (there is always at least one by definition). In
particular, after a translation and rotation of $\RR^3$, we may
suppose that
\begin{equation}    \label{e:1point}
 \barga \cap \cSu = \{ 0 \} \, .
\end{equation}
and the collapsed leaf through $0$ is the punctured horizontal
plane $\{ x_ 3 = 0 \} \setminus \{ 0 \}$.
 We will show that \eqr{e:1point}   implies that
{\underline{every}} leaf of $\cL'$ is a plane with one point
removed, these planes foliate $\RR^3$, and as a consequence the
intersection of any fixed ball with the surfaces $\Sigma_j$ is
simply connected for $j$ sufficiently large. However, this is
impossible since we have assumed an upper bound for the
injectivity radii of the $\Sigma_j$'s in Theorem \ref{t:t5.1}, so
we conclude that \eqr{e:1point} cannot hold.

{\underline{Properness}}:  We will show next that every open
neighborhood of $\{ x_3 = 0\}$ contains points of $\cSu$ both
above and below $\{ x_3 = 0\}$; the proof will use only that
$\barga \cap \cSt = \emptyset$.  We called this properness in
\cite{CM6} and the argument is essentially the same, with one
caveat: \cite{CM6} argues for embedded minimal disks, whereas
presently we only know that the $\Sigma_j$'s are ULSC near $0$.
The disk hypothesis was used for two things in \cite{CM6}:
\begin{enumerate}
\item[(D1)]  The $\Sigma_j$'s are multi-valued graphs in the cone $\{ |x_3| <
\mu \, |x| \}$ for some $\mu > 0$.
\item[(D2)]  The portions of these two multi-valued graphs in a fixed ball
combine to be part of a single embedded minimal disk in this small
ball (this disk property was used in \cite{CM6} to apply Stokes'
theorem).
\end{enumerate}
The second fact (D2) holds in this case since $0 \in \cSu$.  The
first fact (D1) will follow immediately from the one-sided
curvature estimate once we establish the following
{\underline{scale invariant}} ULSC property:
\begin{enumerate}
\item[(D)] There exists $\tau > 0$ so that, for $z \in \{ x_3 = 0 \}$ and $j$ large, each
component of $B_{\tau \, |z|} (z) \cap \Sigma_j$ that connects to
the multi-valued graph in $\Sigma_j$ is a disk.
\end{enumerate}
We will next prove (D) by contradiction, using a variation of the
``between the sheets'' estimate of \cite{CM3}.  To do this, assume
that $\tau > 0$ is small and $z \in \{ x_3 = 0 \}$ is the first
time that (D) fails (i.e., $|z|$ is minimal); obviously, we must
have $ |z| > r_0$ since the $0 \in \cSu$. Fix a sequence of simple
closed non-contractible curves $\gamma_j \subset B_{\tau \, |z|}
(z) \cap \Sigma_j$. We will see that this leads to a
contradiction:
\begin{enumerate}
\item
See Figure \ref{f:d1}.  The existence results of Meeks-Yau,
\cite{MeYa2}, gives stable embedded connected minimal surfaces
$\Gamma_j \subset B_{R_j} \setminus \Sigma_j$ with $\partial
\Gamma_j \setminus \partial B_{R_j} = \gamma_j$ and $\partial
\Gamma_j \cap \partial B_{R_j} \ne \emptyset$.{\footnote{This is a
standard application of \cite{MeYa2}; see, e.g., lemma $3.1$ in
\cite{CM8} for details.}}
 Since these stable surfaces start out on
one side of, but close to, the multi-valued graphs in $\Sigma_j$
converging to $\{ x_3 = 0 \} \setminus \{ 0 \}$, Lemma
\ref{l:stab} implies that the $\Gamma_j$'s quickly become
graphical.
\item
See Figure \ref{f:d2}.  The portion $\Sigma_j^+$ of $\Sigma_j$
between $\partial B_{|z|/2} \cap \{ x_3 = 0 \}$ and the graph in
$\Gamma_j$ must be simply connected in extrinsic balls of radius
$\tau \, |z|/2$  since it is trapped ``between the sheets'' of the
multi-valued graph in $\Sigma_j$ and the graph in $\Gamma_j$ (and
these two can be connected).  Namely, if $\Sigma_j^+$ contained a
non-contractible curve   in   an extrinsic ball of radius $\tau \,
|z|/2$  centered there, then we could apply \cite{MeYa2} to get a
second stable surface $\tilde{\Gamma_j}$ disjoint from both
$\Sigma_j$ and $\Gamma_j$. The surface $\tilde{\Gamma_j}$ would be
forced to become graphical (again by Lemma \ref{l:stab}) but would
then have to cross the curve in $\Sigma_j \cup \Gamma_j$ which
connects the two graphical regions. (Compare the proof of theorem
I.0.8 in \cite{CM3}.)
\item
Since each ball of radius $\tau \, |z|/2$ centered on $\Sigma_j^+$
is simply connected by (2), the set $\Sigma_j^+$ is locally
graphical by the one-sided curvature estimate and, since it
contains a multi-valued graph but cannot pass through $\Gamma_j$,
it spirals infinitely. This is impossible since each $\Sigma_j$ is
compact, giving the contradiction needed to establish (D).
\end{enumerate}
\begin{figure}[htbp]
    \setlength{\captionindent}{20pt}
    \begin{minipage}[t]{0.5\textwidth}
    \centering\input{gen1d.pstex_t}
    \caption{The stable surface $\Gamma_j$ is dotted.}
    \label{f:d1}
    \end{minipage}\begin{minipage}[t]{0.5\textwidth}
    \centering\input{gen2d.pstex_t}
    \caption{$\Sigma_j$ is simply connected between the
    multi-valued graph in $\Sigma_j$ and
    the graph in $\Gamma_j$.}
    \label{f:d2}
    \end{minipage}
\end{figure}

Now that we have established (D), we can argue precisely as in
\cite{CM6} using \cite{CM7} (see lemma I.1.10 there) to prove that
every open neighborhood of this plane contains points of $\cSu$
both above and below this plane.

It follows easily from this properness - as in \cite{CM6} - that
an entire slab $\{ - \epsilon < x_3 < \epsilon \}$ must be
foliated by (the closures of) planar leaves of $\cL'$.  In order
to extend this foliated structure to all of $\RR^3$, we will need
to use the ULSC assumption next.

{\underline{Using the ULSC hypothesis to repeat the argument}}:
Since the set $\cSu$ is automatically closed and transverse to
these planes (by the one-sided curvature estimate), we see that a
 neighborhood of this
plane is foliated by parallel planes and, in this neighborhood,
 $\cSu$ is a single
Lipschitz curve.  However, since the constant $\tau > 0$ was
uniform and did not depend on the particular singular point, we
can now repeat the above argument to extend the set of foliated
planes to the whole of $\RR^3$.  Once we have the foliation by
parallel planes, the ULSC condition and one-sided curvature
estimate imply that $\cSu$ is a discrete collection of transverse
Lipschitz curves.  The transversality implies that the each curve
hits every leaf and hence there is only one curve. In sum, the
sequence is converging to a foliation of $\RR^3$ by parallel
planes away from a single Lipschitz curve $\cSu$ transverse to the
planes.  (This was exactly the result of \cite{CM3}--\cite{CM6}
for sequences of disks.)  This has two consequences:
\begin{itemize}
\item
Near $\cSu$, the sequence looks like a double spiral staircase.
\item
Away from $\cSu$, the sequence is locally converging (with bounded
curvature) to a foliation by parallel planes.
\end{itemize}
The second fact allows us to extend the double spiral staircase
structure away from the singular curve $\cSu$, so that
  we get a sequence $R_j' \to \infty$ where the component of $B_{R_j'} \cap \Sigma_j$ intersecting
$B_{R_j'/C}$ is a double spiral staircase.  In particular, this
component is also   a disk. Since we have assumed that no such
sequence of expanding disks in $\Sigma_j$ exists, we rule out
\eqr{e:1point} as promised.
\end{proof}

\section{Properness and the limit foliation}    \label{s:star}

We have now shown that each collapsed leaf is a plane that is
transverse to $\cSu$ (with a definite lower bound on the angle of
intersection). As in \cite{CM6}, we must show that nearby leaves
are also planes; we call this {\it properness} of the limit
foliation.  Since each singular point in $\cSu$ has a plane
through it, this properness will follow from showing that there
cannot be a first or last such singular point. In fact, it is not
hard to see that these properties are equivalent.  Namely, a
planar leaf nearby a collapsed leaf must also contain singular
points since otherwise the one-sided curvature estimate would give
a curvature bound at the singular point in the collapsed leaf.
 In \cite{CM6}, a similar
properness for disks (where each plane had only one puncture as
opposed to the current situation of two) was proven using
\cite{CM7}.

Before giving the precise statement of properness, observe that we
can rotate $\RR^3$ so that  the closure of each collapsed leaf of
$\cL'$ is a horizontal plane, i.e., is given by $\{ x_3 = t \}$
for some $t \in \RR$.   This is because the closure of each
collapsed leaf is some plane by Proposition \ref{p:cole} and these
planes must all be parallel since the surfaces $\Sigma_j$ are
embedded. With this normalization, the precise statement of
properness is:
\begin{enumerate}
\item[($\star$)]
If $t\in x_3(\cSu)$ and $\epsilon>0$, then $\cS \cap
\{t<x_3<t+\epsilon\}\ne \emptyset$ and $\cS \cap
\{t-\epsilon<x_3<t\}\ne \emptyset$.
\end{enumerate}

We should point out that when the sequence is ULSC, ($\star$)
automatically implies that $\{ x_3 = s \} \cap \cSu$  contains at
least {\underline{two}} points (cf. (S1) and (S2)). Namely, once
$\{ x_3 = s \} \cap \cSu$ contains one point $p'$, then
Proposition \ref{p:cole} implies that   there is a second point
$q' \in \{ x_3 = s \} \cap \cSu$ so that
\begin{equation}
    \{ x_3 =
s \} \setminus \{ p'
   , \, q' \}
\end{equation}
 is a collapsed leaf of $\cL'$.

As in \cite{CM6} and \cite{CM7}, the key to proving ($\star$) is a
careful analysis of the vertical flux of the multi-valued graphs.
Recall that if $\Sigma$ is a minimal surface and $\sigma \subset
\Sigma$ is a simple closed curve, then the vertical flux across
$\sigma$ is
\begin{equation}    \label{e:flux}
    \int_{\sigma} \frac{\partial x_3}{\partial n} \, ,
\end{equation}
where $\frac{\partial x_3}{\partial n}$ is the derivative of $x_3$
in the direction normal to $\sigma$ but tangent to $\Sigma$.  By
Stokes' theorem,  the integral \eqr{e:flux} depends only on the
homology class of $\sigma$ since the coordinate function $x_3$ is
harmonic on a minimal surface.

We will prove ($\star$) by contradiction as we now outline: If
($\star$) does not hold for $t=0$,  then we can assume, after
possibly reflecting across $\{ x_3 = 0 \}$, that
\begin{equation}
    \cS \cap \{0<x_3< \epsilon\} = \emptyset \, .
\end{equation}

\begin{figure}[htbp]
    \setlength{\captionindent}{20pt}
    %\begin{minipage}[t]{0.5\textwidth}
    \centering\input{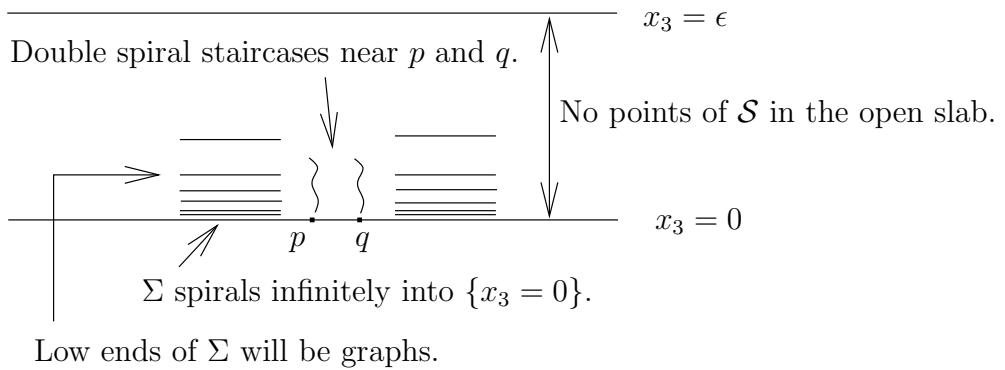}
    \caption{The limit $\Sigma$ when properness fails.}
    \label{f:8a}
    %\end{minipage}
\end{figure}

We will use this to show that there is a unique leaf $\Sigma$ of
$\cL'$ which spirals into the plane $\{ x_3 = 0 \}$ from above and
this leaf is a multiplicity one limit of the $\Sigma_j$'s.
Moreover, near the singular points in $\{ x_3 = 0 \}$, $\Sigma$
will be a double spiral staircase which spirals infinitely into
the plane; see Figure \ref{f:8a}. The ends of $\Sigma$ coming from
circling both double spiral staircases will be  graphs lying above
the plane $\{ x_3 = 0 \}$; we will see that this implies that each
such end has non-negative vertical flux. This structure of the
ends also allows us to cutoff $\Sigma$ below a carefully chosen
horizontal plane $\{ x_3 = \epsilon \}$; it will be almost
automatic that the boundary curve   produced has positive vertical
flux. Using the fact that the plane $\{ x_3 = 0 \}$ contains two
points of $\cSu$, we will find a sequence of separating curves in
$\Sigma$ whose vertical flux goes to zero. This gives a sequence
of compact domains in $\Sigma$ bounded at the top by a curve in
the plane $\{ x_3 = \epsilon \}$ with positive flux, bounded at
the bottom by curves with flux going to $0$, and with boundary
curves on the sides   with non-negative flux.  Combining all of
this will give the desired contradiction since, by Stokes'
theorem, the total flux of any compact domain must sum to zero.

There are in general two ways to show that the vertical flux in
\eqr{e:flux} is small; one can either show that the curve $\sigma$
is short (since $|\nabla x_3| \leq 1$) or show that $|\nabla x_3|$
is small and the length of $\sigma$ is bounded. Since  the length
of any closed non-contractible curve near $\{ x_3 = 0 \}$ is
bounded away from zero,  we must take the second approach here to
bound the flux of the bottom boundary curves (we will use the
first approach in the next part near points in $\cSt$). In our
application, the harmonic function $x_3$ will be positive on the
surface $\Sigma$ and the estimate on $|\nabla x_3|$ will follow
from the gradient estimate.

\subsection{Establishing properness: The proof of ($\star$)} \label{ss:p}
The next lemma shows that ($\star$) holds as long as we have
properties (1) and (2) in Proposition \ref{p:cole}.  In
particular, since these properties always hold when the sequence
is ULSC, we get ($\star$) in the ULSC case.

\begin{Lem}     \label{t:setup}
If (1) and (2) in Proposition \ref{p:cole} hold, then ($\star$)
holds. That is, if the horizontal plane  $\{ x_3 = t \}$ is the
closure of a collapsed leaf satisfying (2)  in Proposition
\ref{p:cole}, and $\epsilon>0$, then
\begin{equation} \cS \cap \{t<x_3<t+\epsilon\}\ne \emptyset {\text{
and }} \cS \cap \{t-\epsilon<x_3<t\}\ne \emptyset \, .
\end{equation}
\end{Lem}

\begin{proof}
For simplicity, we will assume that the $\Sigma_j$'s have genus
zero. The general case follows with easy modifications.

 We will
argue by contradiction, so suppose that $\cSu \cap \{ x_3 = 0 \} =
\{ p , q \}$, $\epsilon > 0$, and
\begin{equation}    \label{e:star}
    \cS \cap
\{0<x_3<\epsilon\} = \emptyset \, .
\end{equation}
Fix a radius $R> 0$ so that the disk $D_R \subset \{ x_3 = 0 \}$
contains both $p$ and $q$.

We will first  record four consequences of \eqr{e:star} that will
be proven below  and
 then use these properties to rule out the
possibility of  such a non-proper limit (see Figure \ref{f:8b}):
\begin{enumerate}
\item[(P1)]
There is exactly one leaf $\Sigma$ of $\cL'$ in $\{ x_3 > 0 \}$
whose closure intersects the plane $\{ x_3 = 0 \}$.  This leaf
$\Sigma$ is a multiplicity one limit of the $\Sigma_j$'s.
Furthermore, after possibly reducing $\epsilon > 0$, the leaf
$\Sigma$ is proper in
 compact subsets of
$\{0<x_3<\epsilon\}$.
\item[(P2)]
Each ``low'' end of $\Sigma$ is an asymptotic graph with
non-negative vertical flux.  Here ``low'' will be made precise
below but roughly means starting off close to $\{ x_3 = 0 \}$ over
the disk $D_R$ containing $p$ and $q$.
\item[(P3)]
Intersecting $\Sigma$ with a carefully chosen horizontal plane
where $x_3$ is constant will give a vertically separating curve
$\gamma_+ \subset \Sigma$
 with positive vertical flux. Here {\it vertically
separating} means that
 if a curve in $\overline{\Sigma}$ is over $D_R$ and intersects both of
\begin{equation}
    \{ x_3 = 0 \}  {\text{ and }}
    \{ x_3 = \epsilon \} \, ,
\end{equation}
then the curve also intersects $\gamma_+$.
\item[(P4)]
There is a sequence of vertically separating curves $\gamma_i
\subset B_1 \cap \Sigma$ with $x_3(\gamma_i) \to 0$ and vertical
flux going to zero.
\end{enumerate}
\begin{figure}[htbp]
    \setlength{\captionindent}{20pt}
    %\begin{minipage}[t]{0.5\textwidth}
    \centering\input{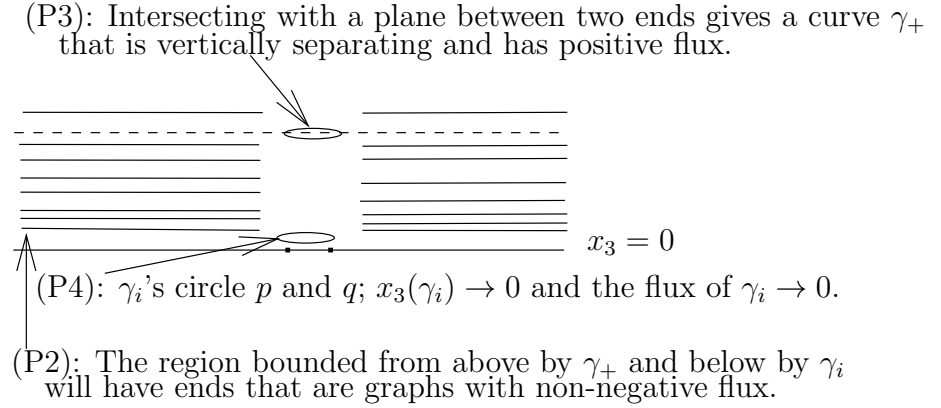}
    \caption{Properties (P2)--(P4) of $\Sigma$ when properness fails.}
    \label{f:8b}
    %\end{minipage}
\end{figure}

{\underline{The proof of (P1)}}.
  Since $p$ and $q$ are locally the
last points in $\cSu $, we are in case (L2) in the proof of Lemma
\ref{l:onestar} near $p$ and $q$ (cf. the ``not proper'' example
(N-P)). Consequently,   Lemma \ref{l:discre} implies that near $p$
and $q$, but above $\{ x_3 = 0 \}$, the $\Sigma_j$'s converge with
multiplicity one to double spiral staircases which spiral
infinitely into $\{ x_3 = 0 \}$.{\footnote{This is  proven within
the proof of Lemma \ref{l:discre}.}} In particular, only one leaf
$\Sigma$ of $\cL$ (namely, the one containing these double spiral
staircases) can spiral into $D_R \subset \{ x_3 = 0 \}$ from
above.  It remains to show that if $\Sigma' \subset \{ x_3 > 0 \}$
is a leaf of $\cL$ that contains a point $q' \in \{ x_3 = 0 \}$ in
its closure $\overline{\Sigma'}$, then $\Sigma' = \Sigma$. This
has already been established when $q'$ is near $p$ or $q$ since
there is only leaf there (cf. the chord arc property of Theorem
\ref{t:1}). However, this easily gives the general case. Namely,
first fix a compact disk
\begin{equation}
    \overline{D_S} = \{ x_1^2 + x_2^2 \leq S^2 \} \cap \{ x_3 = 0
    \} \, .
\end{equation}
containing $p$, $q$, and $q'$. Since, by assumption (see (1) in
Proposition \ref{p:cole}), we have that $\overline{D_S} \cap \cSt
= \emptyset $, an easy compactness argument and the convex hull
property give some $r_1>0$ so that
\begin{equation}
    {\text{each component of }} B_{r_1} (y) \cap \Sigma_j {\text{
    is a disk for every }} y \in \overline{D_S} \, .
\end{equation}
Hence, by the one-sided curvature estimate, the $\Sigma_j$'s are
locally graphical in a cylindrical slab about $\{ x_3 = 0 \}$ and
over $D_S$, as long as we stay away from $p$ and $q$. If we now
choose a point $\bar{q}$ in $\Sigma'$ with $|\bar{q} - q'|$
sufficiently small (depending on both $|p-q'|$ and   $r_1$), then
we can repeatedly apply the Harnack inequality to connect
$\bar{q}$ by a curve in $\Sigma'$ back to a small neighborhood of
$p$. Therefore, since $\Sigma$ is the only leaf in $\{x_3 > 0 \}$
that intersects  a sufficiently small neighborhood of $p$, we
conclude that $\Sigma' = \Sigma$.

Finally, we show that there exists some $\epsilon_0 > 0$ so that
$\Sigma$ is proper in compact subsets of
\begin{equation}
    \{ 0 <
x_3 < \epsilon_0 \} \, .
\end{equation}
This properness will follow by combining the two following facts:
\begin{enumerate}
\item[(Fact 1)]  There exists some $\epsilon_1 > 0$ so that $\cL'$
does not contain any horizontal planes in the slab $\{ 0 < x_3 <
\epsilon_1 \}$. \item[(Fact 2)] There exists $\epsilon_2 > 0$ so
that if $y \in \{ 0 < x_3 < \epsilon_2 \}$ is an ``accumulation
point'' of $\Sigma$, then the leaf of $\cL'$ containing $y$ is  a
horizontal plane.  A point $y$ is said to be ``accumulation
point'' of $\Sigma$ if
  there exists a sequence of points $y_j \in \Sigma$ so
that \begin{align}
    \lim_{j \to \infty} \dist_{\RR^3} (y , y_j) &= 0 \, ,
    \\
     \lim_{j \to \infty} \dist_{\Sigma} (y_1 , y_j) &= \infty \, .
     \end{align}
\end{enumerate}
Obviously, the  two facts together easily imply the properness of
$\Sigma$ in the slab between $0$ and the minimum of $\epsilon_1$
and $\epsilon_2$.

To prove (Fact 1), note that $\inf_{\Sigma} x_3 = 0$ and
$\sup_{\Sigma} x_3 > 0$ (since $\Sigma$ is not flat).
Consequently, since the leaves of $\cL'$ are disjoint and the leaf
$\Sigma$ is connected, $\cL'$ cannot contain any horizontal planes
in the open slab $\{ 0 < x_3 <  \sup_{\Sigma} x_3 \}$.  This gives
(Fact 1).

(Fact 2) follows easily from the proof of lemma $1.3$ in
\cite{MeRo} (in fact, one can even take $\epsilon_2 = \epsilon$
but we will not need this).  We will give the proof next for
completeness.  Suppose therefore that $\Sigma$ accumulates at a
point $y \in \{ 0 < x_3 < \epsilon_2 \}$  in a leaf $\hat{\Sigma}
\in \cL'$ ($y$ is in a leaf since the union of the leaves of
$\cL'$ is a closed subset of $\RR^3 \setminus \cS
 $ and $\cS$ does not intersect the open slab $\{ 0 < x_3 < \epsilon \}$).
  Since $\Sigma$ is discrete by Lemma
 \ref{l:discre}, we must have $ \hat{\Sigma} \ne \Sigma$ and,
 hence, also
 \begin{equation}  \label{e:signothat}
    \overline{\hat{\Sigma}} \cap \{ x_3 = 0 \} = \emptyset \, .
 \end{equation}
 First, we can use the Harnack inequality to extend the local sequence of
graphs converging to $\hat{\Sigma}$ to graphs over a sequence of
expanding subdomains of $\hat{\Sigma}$, eventually obtaining a
positive Jacobi field on the universal cover of $\hat{\Sigma}$
(cf. lemma 2.1 in \cite{CM4}).  In particular, the existence of
such a positive Jacobi field implies stability of the universal
cover of $\hat{\Sigma}$ (see, e.g., proposition 1.26 in
\cite{CM1}).  Combining the local curvature estimate for stable
surfaces, \cite{Sc1}, with the gradient estimate and
\eqr{e:signothat} gives $\epsilon_2 > 0$ (depending only on
$\epsilon$) so that
\begin{equation}        \label{e:loot}
    \{ x_3 \leq \epsilon_2 \} \cap \hat{\Sigma} {\text{ is
    {\underline{locally}}
    graphical over }} \{ x_3 = 0 \} \, .
\end{equation}
(See, e.g., lemma I.0.9 in \cite{CM3} for a detailed proof of
\eqr{e:loot}.) However, the theory of covering spaces (see lemma
$1.4$ in \cite{MeRo}) now implies that each component of $\{ x_3
\leq \epsilon_2 \} \cap \hat{\Sigma}$ is {\underline{globally}} a
graph of a function $u$ with $0 < u \leq \epsilon_2$ over a domain
$\Omega \subset \{ x_3 = 0 \}$ with boundary values $u\, |_{
\partial \Omega} = \epsilon_2$.
However, the proof of the strong halfspace theorem of
Hoffman-Meeks, \cite{HoMe}, implies that such a $u$ must be
identically equal to its boundary values and, by unique
continuation, we see that $\hat{\Sigma}$ is a horizontal plane as
claimed.  This completes the proof of (Fact 2) and, hence, also of
(P1).

{\underline{The proof of (P2)}}. The proof of the asymptotic graph
structure in (P2) will be similar to the proof of property (D1) in
subsection \ref{ss:notone}. First,  we fix a large constant
$\Omega
> 1$ and some disk $D_R \subset \{ x_3 = 0 \}$ containing both $p$
and $q$.  Since $\{ x_3 =0 \} \setminus \{ p , q \}$ is a leaf of
$\cL'$ (and, in particular, disjoint from $\cS$),  an easy
covering argument gives  constants $\mu_A > 0$ and $C_A$ so that
\begin{equation}    \label{e:mua}
    \sup_{ \{ 0 < x_3 < \mu_A \} \cap \{ x_1^2 + x_2^2 \leq
    4 \, \Omega^2  \, R^2 \} \setminus (B_{\epsilon}(p) \cup
    B_{\epsilon}(q)} |A|^2 \leq C_A \, .
\end{equation}
The gradient estimate and \eqr{e:mua} then give a constant $\mu
>0$ (and less than $\mu_A$) so that:
\begin{enumerate}
\item[] See Figure \ref{f:dp0}.
Each point $y \in \{ 0 < x_3 < \mu \} \cap \Sigma$ over $\partial
D_R$ is contained in a graph $\Sigma_y \subset \Sigma$  over the
annulus $D_{\Omega R}\setminus D_R$ with $\Sigma_y \subset \{ 0 <
x_3 < \epsilon \}$.  Furthermore, the graph $\Sigma_y$ extends
over (a large part of) $D_R$ as a graph, connecting to the two
double spiral staircase structures near $p$ and $q$.
\end{enumerate}
\begin{figure}[htbp]
    \setlength{\captionindent}{20pt}
    \centering\input{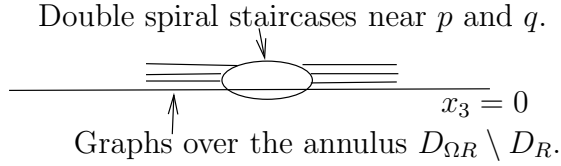}
    \caption{The graphs $\Gamma_y$ in $\Sigma$.}
    \label{f:dp0}
\end{figure}
This allows us to make the notion of a ``low'' end precise.
Namely, a low end is the component of $\Sigma \setminus \Sigma_y$
-- with $\Sigma_y$ as above -- whose closure does not contain the
boundary graph over $\partial D_R$.

It remains to show that each of the graphs $\Sigma_y$ extends as a
graph indefinitely, i.e., past $\partial D_{\Omega R}$.  This is
where we argue as the proof of (D) in subsection \ref{ss:notone},
proving the following {\underline{scale invariant}} ULSC property:
\begin{enumerate}
\item[(D')] There exists $\tau > 0$ so that, for $z \in \{ x_3 = 0 \} \setminus D_R$ and $j$ large, each
component of $B_{\tau \, |z|} (z) \cap \Sigma_j$ that connects to
the multi-valued graph in $\Sigma_j$ is a disk.
\end{enumerate}
Once we have shown (D'), then the one-sided curvature estimate and
gradient estimate will allow us to extend $\Sigma_y$  as a graph
indefinitely.  We will leave this easy extension argument to the
reader.

We will next prove (D') by contradiction; suppose therefore that
some $\Sigma_y$ connects outside $D_R$ to a component of $B_{\tau
\, |z|} (z) \cap \Sigma_j$ which is not a disk. (Here and below we
are identifying $\Sigma_y$ with the portion of the $\Sigma_j$'s
converging to it.  This identification should not lead to
confusion since  we have already shown that the convergence is
multiplicity one.) We observe that this ``non-disk component''
leads to a contradiction, roughly following the proof of (D) in
subsection \ref{ss:notone}
 as follows:
\begin{itemize}
\item  Figure \ref{f:dp1}.  A curve $\gamma_j$ which is non-contractible
in $B_{\tau \, |z|} (z) \cap \Sigma_j$ is also non-contractible in
$\Sigma_j$ by the convex hull property.  We can therefore apply
\cite{MeYa2} to find a stable surface $\Gamma_j \subset B_{R_j}
\setminus \Sigma_j$ with interior boundary $\gamma_j$ as in (1) in
subsection \ref{ss:notone}.  As before,  using the plane $\{ x_3 =
0\}$ allows us to conclude by Lemma \ref{l:stab} that these stable
surfaces quickly become graphs as we move away from the interior
boundary $\gamma_j$.
\item
 Figure \ref{f:dp2}.
Fix a  graph $\Sigma_{y'}$ {\underline{above}} $\Sigma_y$ and let
$\sigma_{y'}$ be the component of $\partial \Sigma_{y'}$ over
$\partial D_R$.  As in (2) in subsection \ref{ss:notone}, we could
then apply \cite{MeYa2} to put in a second stable surface
$\tilde{\Gamma_j} \subset B_{R_j} \setminus (\Gamma_j \cup
\Sigma_j)$ with interior boundary equal to $\sigma_{y'}$.
Furthermore, this surface also quickly becomes graphical by Lemma
\ref{l:stab}.  Finally, the graph in $\tilde{\Gamma_j}$ must start
off between the two graphs $\Sigma_y$ and $\Gamma_j$ because (1)
its interior boundary $\sigma_{y'}$ was chosen to be above
$\Sigma_y$ and (2) every point near $D_R$ in $\Sigma$ which
connects  back to the multi-valued graph in $\Sigma$ must be below
$\Gamma_j$.
\newline This easily gives the desired contradiction:
The construction of $\Gamma_j$ guarantees that there is a curve in
$\Sigma_j \cup \Gamma_j$ connecting $\Sigma_y$ to the graph in
$\Gamma_j$ and moreover is a graph over the $x_1$-axis except for
in a small neighborhood of $z$.   The graph   in
$\tilde{\Gamma_j}$ is  consequently forced to intersect this
curve, giving the desired contradiction.
\end{itemize}
\begin{figure}[htbp]
 \setlength{\captionindent}{20pt}
    \begin{minipage}[t]{0.5\textwidth}
    \centering\input{gen1dp.pstex_t}
    \caption{The proof of (D'): Constructing $\Gamma_j$.}
    \label{f:dp1}
    \end{minipage}\begin{minipage}[t]{0.5\textwidth}
    \centering\input{gen2dp.pstex_t}
    \caption{The proof of (D'): $\tilde{\Gamma_j}$ must intersect $\Gamma_j \cup
\Sigma_j$.}
    \label{f:dp2}
    \end{minipage}
\end{figure}

The final part of (P2), i.e., the non-negativity of the vertical
flux of the low ends, follows immediately since the ends are
asymptotic to planes (vertical flux zero) or
{\underline{upper}}--halves of catenoids (vertical flux positive).
This is because the only other possibility would be an end
asymptotic to the  {\underline{lower}} half of a catenoid which is
impossible since one of these would eventually go below the plane
$\{ x_3 = 0 \}$.  (Recall that any embedded minimal end with
finite total curvature is asymptotic to either a plane or half of
a catenoid by  proposition 1 in \cite{Sc2}.)

{\underline{The proof of (P3)}}.
 To prove (P3), recall first that there is a positive distance
between consecutive ends by the maximum principle at infinity of
\cite{LaRo}.  This positive distance allows us to intersect
$\Sigma$ with a horizontal plane which intersects $\Sigma$
transversely between the heights of two consecutive ends, giving a
finite collection $\gamma_+$ of disjoint simple closed curves
separating these ends; see Figure \ref{f:8}. This finiteness
follows from the compactness of the level set which in turn used
the properness of $\Sigma$.   Since we will be considering the
part of $\Sigma$ {\underline{below}} this plane, the outward
normal derivative of $x_3$ is non-negative  at every point along
$\gamma_+$. However, this plane was chosen to be transverse to the
surface, so this derivative must in fact be pointwise positive
along $\gamma_+$.
 That is, the flux integrand is pointwise positive along
$\gamma_+$ so the vertical flux across $\gamma_+$ is clearly
positive.

Finally, we will sketch briefly why the fact that $\Sigma$ has
genus zero implies that we can choose a single component of
$\gamma^+$ that is vertically separating. We must  show that only
one component of $\Sigma \setminus \gamma^+$ connects $\gamma^+$
to the ``the ceiling'' $\{ x_3 = \epsilon \}$. The  point  is that
if there were two such components of $\Sigma \setminus \gamma^+$,
then we could solve a sequence of Plateau problems to get a stable
surface $\Gamma^+$ between them with the following properties:
\begin{itemize}
\item $\partial \Gamma^+ \subset \gamma^+$ is a finite collection
of disjoint simple closed curves in a plane. \item $\Gamma^+$ does
not go below the plane containing $\gamma^+$. \item $\Gamma^+$ is
above one of the components of $\Sigma \setminus \gamma^+$ and
below the other.
\end{itemize}
It is not hard to see that this is impossible.
 The connectedness
of $\gamma^+$ is not actually necessary for the proof, so we will
leave the details for the reader.

\begin{figure}[htbp]
    \setlength{\captionindent}{20pt}
    %\begin{minipage}[t]{0.5\textwidth}
    \centering\input{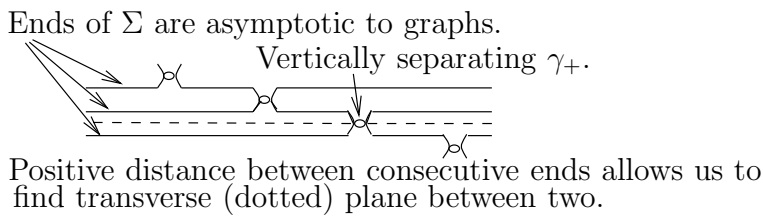}
    \caption{The vertically separating  curve $\gamma_+$.}
    \label{f:8}
    %\end{minipage}
\end{figure}

{\underline{The proof of (P4)}}.
 The last claim (P4) essentially follows from the description of
 $\Sigma$ near the plane
 $\{ x_3 = 0 \}$ and the gradient estimate.
 To see this, let
 \begin{equation}
    \gamma_{pq} \subset \{ x_3 = 0 \}
 \end{equation}
be the line segment from $p$ to $q$.  Using the description of
$\Sigma$ near $\{ x_3 = 0 \}$ and the chord arc bound of Theorem
\ref{t:1}, we can find simple closed curves $\gamma_i \subset
\Sigma$ with the following properties (see Figure \ref{f:9}):
\begin{itemize}
\item $\gamma_i$ is contained in the $\epsilon_i$-tubular
neighborhood of $\gamma_{pq}$ where $\epsilon_i \to 0$.
\item $\gamma_i \setminus \left( B_{\epsilon_i}(p) \cup B_{\epsilon_i}(q) \right)$
 consists of two graphs over $\gamma_{pq}$ which
are in {\underline{oppositely}}--oriented sheets of $\Sigma$.
\item
The length of $\left( B_{\epsilon_i}(p) \cup B_{\epsilon_i}(q)
\right) \cap \gamma_i$ is at most $C\,\epsilon_i$ for a fixed
constant $C$.
\end{itemize}
\begin{figure}[htbp]
    \setlength{\captionindent}{20pt}
    %\begin{minipage}[t]{0.5\textwidth}
    \centering\input{gen9.pstex_t}
    \caption{The  separating curves $\gamma_i$.}
    \label{f:9}
    %\end{minipage}
\end{figure}
It is easy to see that the $\gamma_i$'s are vertically separating
since the sheets containing
\begin{equation}
    \gamma_i \setminus
    \left( B_{\epsilon_i}(p) \cup B_{\epsilon_i}(q) \right)
\end{equation}
 are oppositely--oriented (meaning that the unit normals point in
 nearly
 opposite directions).  Namely, if we fix a spiralling curve in
 either of the two multi-valued graphs near either singular curve,
 then the third component to the
 unit normal to $\Sigma$ does not change sign along this curve. Consequently,
 such a spiralling curve intersects exactly one of the $\gamma_i$'s and does so exactly once.
 It follows that the $\gamma_i$'s are vertically separating as
 claimed.

 We can now use the gradient estimate to bound $|\nabla x_3|$ along $\gamma_i$ away from $p$ and $q$ to see
that the vertical flux on
\begin{equation}
    \gamma_i \setminus \left(
B_{\delta}(p) \cup B_{\delta}(q) \right)
\end{equation}
 goes
to zero for any fixed $\delta
> 0$.   To get this bound, note first that the height function $x_3$ is
positive and harmonic on the multi-valued graphs in $\Sigma$ that
spiral infinitely into $\{ x_3 = 0 \}$ and these graphs have
bounded curvature away from $p$ and $q$. Therefore, the gradient
estimate (for positive harmonic functions) implies that
\begin{equation}
    \sup_{ \gamma_i \setminus \left(
B_{\delta}(p) \cup B_{\delta}(q) \right)} \, \, |\nabla x_3| \to 0
{\text{ uniformly as }} i \to \infty \, .
\end{equation}
Combining this gradient bound with the bound on the length of
\begin{equation}
    \left( B_{\delta}(p) \cup B_{\delta}(q) \right) \cap
    \gamma_i
\end{equation}
gives the last claim.  This completes the proof of the properties
(P1)--(P4) of $\Sigma$.

{\underline{Using (P1)--(P4) to deduce a contradiction}}. We now
see that these properties are contradictory.  Namely, by Stokes'
theorem, the total flux across $\gamma_i$, $\gamma_+$, and the
``ends'' of $\Sigma$ between $\gamma_i$ and $\gamma_+$ must sum to
zero. However, the flux across $\gamma_+$ is positive and every
other flux is either non-negative or approaches zero. This
contradiction shows that \eqr{e:star} could not have held, proving
the lemma.
  \end{proof}

\section{Completing the proof of Theorem \ref{t:t5.1}}
\label{s:compf}

We will now use the properties of the singular set $\cSu$ and
lamination $\cL'$ to show that $\cL'$ is a foliation by parallel
planes with two Lipschitz curves removed, thereby completing the
proof of Theorem \ref{t:t5.1}.   The two main steps are:
\begin{itemize}
\item Using properness (Lemma \ref{t:setup}) to see that
the (collapsed) planar leaves of $\cL'$ intersect every height.
\item
Using the local cone property to get regularity of $\cSu$.
\end{itemize}

\begin{proof}
(of Theorem \ref{t:t5.1}.) Lemma \ref{l:singl} gives a subsequence
$\Sigma_j$, singular set $\cS$, and lamination $\cL'$ of $\RR^3
\setminus \cS$ with minimal leaves.  The set $\cS = \cSu$ is
nonempty by assumption.

For each point $x$ in $\cSu$, properties (1) and (2) of
Proposition \ref{p:cole}
 give  a (collapsed) leaf of $\cL'$ which is a plane with two points removed
 ($x$ is one of the two points).  It follows easily from the convergence to these planes
 and the embeddedness of $\Sigma_j$ that all of the limit planes are
 parallel so, after a rotation of $\RR^3$, we can assume that these planes
 are horizontal, i.e., given as level sets $\{ x_3 = t \}$.
   Furthermore, since
$\cSu \subset \RR^3$ is a nonempty closed set, the local cone
property implies that $x_3(\cSu)\subset \RR$ is also closed (and
nonempty).

 We will show first that the collapsed leaves  (or, rather, their closures) foliate $\RR^3$, more precisely,
 that
 \begin{equation}   \label{e:allofr}
    x_3(\cSu) = \RR \, .
 \end{equation}
  To prove \eqr{e:allofr},
we assume that $\{x_3=t_0\} \cap \cSu  =\emptyset$ for some $t_0
\in \RR$ and will see that this
 leads to a contradiction.    Namely, since $x_3(\cSu)$ is closed,
 there exists $t_s\in x_3(\cSu)$ which is a
closest point in $x_3(\cSu)$ to $t_0$. The desired contradiction
now easily follows from Lemma \ref{t:setup} since either $\{t_s
<x_3<t_0\} \cap \cSu$ or $\{t_0<x_3<t_s \} \cap \cSu$ is
 empty. We conclude therefore that $x_3(\cSu) = \RR$.

Finally, the Lipschitz regularity of the curves now follows as in
lemma I.1.2 of \cite{CM6}; the same argument applies with obvious
minor modifications to deal with the fact that each horizontal
plane now contains two singular points as opposed to just one in
\cite{CM6}.
\end{proof}

\section{Sequences with fixed genus}    \label{s:fixedg}

Most of the arguments in the preceding sections have assumed that
the surfaces $\Sigma_j$ have genus zero as opposed to just some
fixed finite genus.   The arguments for the genus zero case are
slightly simpler, however, the modifications needed for the
general case are straightforward.  The key point is that the
infinite multiplicity of the multi-valued graphs converging to  a
collapsed leaf means that there is arbitrarily large number of
disjoint curves to choose from that ``circle both axes'' and thus
have the desired properties for the preceding arguments.
   In the general case of finite (but nonzero) genus, we can therefore
follow the preceding argument using
   the following
    lemma:

     \begin{Lem}     \label{l:gn}
If $\Sigma$ is oriented with genus $g$ and $\sigma_1 , \dots ,
\sigma_{g+1} \subset \Sigma$ are disjoint simple closed curves,
then $\Sigma \setminus \cup_i \sigma_i$ is disconnected.
\end{Lem}

\begin{proof}
The first integral homology group of $\Sigma$ is  $2g$-dimensional
and the intersection form
 is a bilinear form of full rank
(cf. lemma I.0.9 of \cite{CM5}).  Therefore the maximal subspaces
on which the intersection form vanishes have dimension $g$.
Consequently,
 there is a nontrivial linear (integral) relation between the
$\sigma_i$'s and the lemma follows easily.
\end{proof}

\section{An application: A one-sided property for ULSC surfaces}
\label{s:onesided}

The compactness theorem for ULSC sequences,  Theorem \ref{t:t5.1},
can be used to prove estimates for embedded minimal surfaces that
have a lower bound on their injectivity radius.  We will prove
several such estimates in this paper, including Lemma \ref{l:os}
in this section. This lemma proves a one-sided property for
non-simply connected surfaces on the smallest scale of non-trivial
topology, showing that an intrinsic ball in such a surface cannot
lie on one side a plane and have its center close to the plane on
this scale. This result requires that we work on this scale since,
after all, large balls in the catenoid can be rescaled to lie
above a plane and yet come arbitrarily close to the plane.

  The proof of the lemma   divides
naturally into two extreme cases, depending on whether the
(inverse of the) curvature is comparable to the injectivity radius
or is much larger.  In the first case,  the surface looks more
like a catenoid while in the second it looks like a pair of
oppositely oriented helicoids joined together.  In the first case,
the lemma essentially follows from the logarithmic growth of the
ends of the catenoid; the second case  follows from that these
double-helicoids  converge to a foliation of all of $\RR^3$ by the
compactness theorem for ULSC sequences,  Theorem \ref{t:t5.1}.

\begin{Lem}     \label{l:os}
Let $\Sigma$ be an embedded minimal planar domain with $0 \in
\Sigma$ and so that $\cB_{4r_1}(0) \subset \Sigma$ is
{\underline{not}} a disk.
 Given any $H> 0$, there exists $C_1 > H$ so that if  $\cB_{C_1 \, r_1} (0)
 \cap \partial \Sigma = \emptyset$ and
$\cB_{r_1}(x)$ is a disk for each $x \in \cB_{C_1 \, r_1}(0)$,
then
\begin{equation}    \label{e:os}
     \sup_{ \cB_{C_1 \, r_1} (0) } x_3 >  H \, r_1 \, .
\end{equation}
\end{Lem}

\begin{proof}
After rescaling, we can assume that $r_1 = 1$.  We will argue by
contradiction, so suppose that $\Sigma_j$ is a sequence of
embedded minimal planar domains containing $0$ with
\begin{align}   \label{e:cos1}
    \cB_j (0) &\subset \left( \Sigma_j \setminus \partial \Sigma_j \right) \cap \{ x_3 \leq  H \} \, , \\
    \cB_4(0) &{\text{ is not a disk}} \, , \label{e:cos2} \\
    \cB_1 (x) &{\text{ is a disk for each }} x \in \cB_j (0) \, . \label{e:cos3}
\end{align}

 Observe first that Lemma
\ref{l:cyulsc} in Appendix \ref{s:aC} gives a sequence $R_j \to
\infty$ so that the component $\Sigma_{0,R_j}$ of $B_{R_j} \cap
\Sigma_j$ containing $0$ is compact and has boundary in $\partial
B_{R_j}$.{\footnote{This will be needed later to apply the
compactness results for ULSC sequences.}} Replacing the sequence
$\Sigma_j$ by $\Sigma_{0,R_j}$ gives a ULSC sequence -- still
denoted $\Sigma_j$ -- of embedded minimal planar domains in
extrinsic balls whose radius goes to infinity. We will now divide
the proof into two cases depending on whether or not the
curvatures of the sequence blows up.

\vskip2mm \noindent {\underline{Case 1}}:
 Suppose first that $|A|^2
\to \infty$ in some fixed ball of $\RR^3$ (for some subsequence).
We can then apply the compactness theorem for ULSC sequences,
Theorem \ref{t:t5.1}, to get a subsequence of the $\Sigma_j$'s
that converges to a foliation by parallel planes away from two
lines orthogonal to the leaves of the foliation.   Since the
foliation is of all of $\RR^3$, this contradicts the upper bound
for $x_3$ in \eqr{e:cos1}.

\vskip2mm \noindent {\underline{Case 2}}: Suppose now that $|A|^2$
is uniformly bounded on compact subsets of $\RR^3$ for every
$\Sigma_j$.  In this case, a subsequence of the $\Sigma_j$'s
converges smoothly to a minimal lamination $\cL$ of $\RR^3$ by
proposition B.1 in \cite{CM6}.  We will first see that
\eqr{e:cos1}--\eqr{e:cos3} imply that there is a non-flat leaf
$\Gamma$ of $\cL$ satisfying
\begin{align}   \label{e:cos1a}
    \Gamma &\subset \{ x_3 \leq  H \} \, , \\
    \cB_1 (x) &{\text{ is a disk for each }} x \in \Gamma \, . \label{e:cos3a}
\end{align}
To see this, first note that \eqr{e:cos2} implies that each
$\cB_4(0) \subset \Sigma_j$ cannot be written as a graph over any
plane and hence contains a point $y_j$ with $|A|^2 (y_j) >
\delta_0$ for some $\delta_0 > 0$. A subsequence of these points
converges to a point $y$ in some leaf - call it $\Gamma$ - of
$\cL$ with $|A|^2 (y) \geq \delta_0$. Thus $\Gamma$ is not flat.
Equation \eqr{e:cos1a} follows immediately from \eqr{e:cos1}.

Observe next that the leaf $\Gamma$: \begin{itemize} \item is a
multiplicity one
 limit of
the $\Sigma_j$'s. \item  is locally isolated in $\cL$ in the sense
that  each point $y \in \Gamma$ has a neighborhood $B_{r}(y)$  so
that $B_{r}(y) \cap \cL$ consists only of the component of $B_r(y)
\cap \Gamma$ containing $y$. \end{itemize} If either of these was
not the case, then the universal cover of $\Gamma$ would be stable
and, hence, flat; cf. the proof of Corollary \ref{c:collstab} for
more details. As a consequence, if $K$ is any compact subset of
$\Gamma$, then for $j$ sufficiently large (depending on $K$)
$\Sigma_j$ contains a normal graph $K_j$ over $K$; as $j \to
\infty$, the $K_j$'s converge smoothly to $K$. Using this
convergence and the convex hull property, it is easy to see that
\eqr{e:cos3}  implies \eqr{e:cos3a}.  Moreover, this convergence
and the fact that the $\Sigma_j$'s are planar domains implies that
$\Gamma$ is also a planar domain.  To see this, suppose instead
that $\Gamma$ contains a pair of non-separating curves $\gamma$
and  $\tilde \gamma$ that have linking number one (i.e., so they
are transverse and intersect at exactly one point). Then for $j$
large, we would get a similar pair of curves in $\Sigma_j$; since
this is impossible for planar domains, we conclude that $\Gamma$
is also a planar domain.

We have shown that $\Gamma$ is a planar domain satisfying
   \eqr{e:cos1a} and \eqr{e:cos3a}.  We can assume that
   $H= \sup_{\Gamma} x_3$.  Recall that
\cite{MeRo} gives
\begin{equation}     \label{e:plcu}
    \sup_{ \Gamma \cap \{ H - 1 < x_3 < H
   \} } \, \, |A|^2 = \infty \, .
   \end{equation}
Namely, by the first paragraph of the proof of lemma $1.5$ in
\cite{MeRo}, if instead we had
\begin{equation}     \label{e:plcunot}
    \sup_{ \Gamma \cap \{ H - 1 < x_3 < H
   \} } \, \, |A|^2 <  \infty  \, ,
   \end{equation}
then $\Gamma = \{ x_3 = H \}$.  However, $\Gamma$ is not flat, so
we conclude that \eqr{e:plcu} must hold.

We will now use \eqr{e:plcu} to define a new sequence of planar
domains where we can argue to a contradiction as in Case $1$.
Namely, \eqr{e:plcu} gives a sequence of points $p_n \in \Gamma
\cap \{ H - 1 < x_3 < H
   \}$ with
\begin{equation}    \label{e:dpnto}
    |A|^2 (p_n) \to \infty \, .
\end{equation}
By \eqr{e:cos3a}, we can apply Lemma \ref{l:cyulsc} to conclude
that the component $\Gamma_{p_n , n}$ of $B_n ( p_n) \cap \Gamma$
containing $p_n$ is compact and has boundary in $\partial B_n (
p_n)$.  Translate  $\Gamma_{p_n , n}$ by moving $p_n$ to the
origin to get a ULSC sequence
\begin{equation}
    \Gamma_n = \Gamma_{p_n , n} - p_n
\end{equation}
 of compact embedded minimal planar domains
with $\Gamma_n \subset
 B_n$, $\partial \Gamma_n
\subset
\partial B_n$, and $|A|^2(0) \to \infty$.  As in Case $1$, a
subsequence converges to  a foliation of {\underline{all of
$\RR^3$}} by parallel planes away from two lines orthogonal to the
leaves of the foliation.  However, by \eqr{e:cos1a}, the
translated surfaces $\Gamma_n$ are in the half--space $\{ x_3 < 1
\}$.  This contradiction completes the proof of the lemma in Case
$2$, completing the proof of the lemma.
\end{proof}

\begin{Rem}
Using this one-sided property, we can go back and prove a stronger
version of Lemma \ref{l:cyulsc}.  This stronger result gives that
ULSC surfaces are proper -- as opposed to just knowing that each
component in a ball is proper.
\end{Rem}

\part{When the surfaces are not ULSC: The proof of Theorem
\ref{t:t5.2}}  \label{p:prove2}

\setcounter{equation}{0}

We will now turn to the case where the sequence is not ULSC and
there is consequently no longer a lower bound for the injectivity
radius of the $\Sigma_j$'s in compact subsets of $\RR^3$.  As we
did in the ULSC case, we will initially argue for the genus zero
case and then explain the easy modifications needed for the
general case of fixed finite genus.

We have already defined the singular set $\cS$ in Definition/Lemma
\ref{l:inftyornot} to be the set of points where the curvature
blows up. Furthermore, Lemma \ref{l:singl} gives a subsequence
$\Sigma_j$ that converges to a minimal lamination $\cL'$ of $\RR^3
\setminus \cS$.  This gives (A) and (B) in Theorem \ref{t:t5.2}.

In the ULSC case, every singular point was essentially the same;
namely, in a neighborhood of each singular point, the surfaces
were double-spiral staircases.  However, we now have the
possibility that the injectivity radius of the $\Sigma_j$'s is
going to zero at the singular point. This occurs, for example, by
taking a sequence of rescalings of a catenoid or one of the
Riemann examples.  Recall that the Riemann examples are
singly-periodic embedded minimal planar domains which are
topologically - and conformally -  equivalent to an infinite
cylinder with a one-dimensional lattice of punctures.

For the sequence of rescaled catenoids, the singular set $\cS$
consists of just the origin  and we get $C^{\infty}$ convergence
to a single plane with multiplicity two away from the origin.
Rescaling one of the Riemann examples gives a line of singular
points and convergence to a foliation by parallel planes away from
this line. By choosing different sequences of rescaled   Riemann
examples, we can get different singular sets -- but we always get
a foliation by parallel planes.

The local behavior of the surfaces near a singular point is quite
different, depending on whether or not the injectivity radius is
going to zero there.
 To
account for this, we define   the subset $\cSt \subset \cS$ to be
the set of points where the injectivity radius goes to zero.
Proposition I.0.19 of \cite{CM5} implies that, after passing to a
further subsequence, each point $y \in \cS \setminus \cSt$ has a
radius $r_y > 0$ so that each component of $B_{r_y}(y) \cap
\Sigma_j$ is a disk for every $j$.  In other words, this
proposition implies that $\cS$ is given as the disjoint union
\begin{equation}
    \cS = \cSt \cup \cSu \, .
\end{equation}
Recall that  $\cSu$ was defined in the introduction to be the set
of points where the curvature blows up but where the sequence is
locally ULSC.

\vskip2mm \noindent {\bf{An overview of this part}}:  In Section
\ref{s:ceno}, we prove the main structure result for the non-ULSC
part of the limit lamination $\cL'$, i.e, (C1) in Theorem
\ref{t:t5.2}. Namely, we show that for each point $y$ in $\cSt$,
we get a sequence of {\underline{graphs}} in the $\Sigma_j$'s that
converges to a plane through $y$.  These graphs will be defined
over a sequence of expanding domains and the convergence will be
smooth away from $y$ and possibly one other point in the limit
plane.  In the process of proving this, we will also establish (P)
in Theorem \ref{t:t5.2}.

In Section \ref{s:ctwod}, we prove the main structure results for
the ULSC part of the lamination, i.e., (C2) and (D) in Theorem
\ref{t:t5.2}.  Namely, we show that this part of the lamination is
actually a foliation by an open set of parallel planes in $\RR^3$
and the ULSC singular set $\cSu$ is a collection of Lipschitz
curves transverse to these planes.

In Section \ref{s:puttog}, we combine all of this and complete the
proof of Theorem \ref{t:t5.2}.

\begin{Rem}
Theorem \ref{t:t5.2}  gives a flat leaf of $\cL'$ through every
singular point in $\cS$ but does not show that all of the leaves
of $\cL'$ are flat.  This will be proven in Part \ref{p:prove3}.
\end{Rem}

\section{Proving (C1) in Theorem \ref{t:t5.2}: A plane
through each point of $\cSt$}   \label{s:ceno}

For each  point $y$ in $\cSt$, we will prove in this section that
there is a sequence of graphs in $\Sigma_j$ converging to a plane
through $y$.  The graphs converge smoothly to the plane away from
$y$ and possibly one other point (the other point is also in
$\cSt$).

\begin{figure}[htbp]
    \setlength{\captionindent}{20pt}
    %\begin{minipage}[t]{0.5\textwidth}
    \centering\input{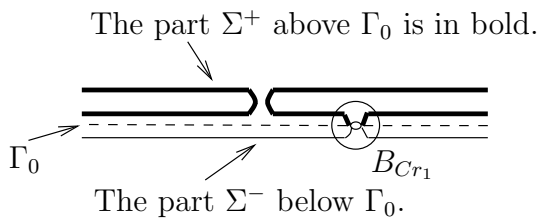}
    \caption{Proposition \ref{l:getsetgo}: A stable graph $\Gamma_0$ separates $\Sigma$ into parts
    above and below the graph.}
    \label{f:getsetgo}
    %\end{minipage}
\end{figure}

The key tool in this section is Proposition \ref{l:getsetgo} that
allows us to
 decompose  an embedded minimal planar domain $\Sigma \subset B_{r_0}$ with $\partial \Sigma \subset \partial
 B_{r_0}$ by ``cutting it'' inside a small ball $B_{r_1}$ whenever some component of $B_{r_1} \cap \Sigma$
 is not a disk.  Moreover, the proposition uses a barrier
 construction to find a stable graph disjoint from $\Sigma$ so
 that the pieces of $\Sigma$ are on opposite sides of this graph; see
Figure \ref{f:getsetgo}.
 The basic example to keep in mind is the
 catenoid: Cutting the catenoid along the unit circle in the $\{
 x_3 = 0 \}$ plane gives two pieces; these pieces are on opposite
 sides of the stable graph $\{ x_3 = 0 \} \cap \{ x_1^2 + x_2^2 >
 1 \}$.

\begin{Pro} \label{l:getsetgo}
There exists a constant $C > 1$ so that the following holds:

Let $\Sigma \subset B_{r_0}$ be an embedded minimal planar domain
with $\partial \Sigma \subset
\partial B_{r_0}$ and $0\in \Sigma$.  If
$\cB_{r_1} \subset \Sigma$ is not a topological disk for some $r_1
< r_0 / C^2$, then there exists a stable embedded minimal surface
$\Gamma \subset B_{r_0} \setminus \Sigma$ with $\partial \Gamma
\subset
\partial B_{r_0} \cup \cB_{r_1}$ and satisfying the following
properties:
\begin{enumerate}
\item[(A)] A component $\Gamma_0$ of $B_{r_0/C} \cap \Gamma
\setminus B_{C r_1}$   is a graph with gradient bounded by one and
so that $\partial \Gamma_0$ intersects both $\partial B_{r_0/C}$
and $\partial B_{C r_1}$. \item[(B)] There are distinct components
$H^{+}$  and $H^-$ of $B_{r_0/C} \setminus (\Gamma_0 \cup B_{C \,
r_1})$, a separating curve $\tilde \sigma \subset B_{C \, r_1}
\cap \Sigma$, and distinct components $\Sigma^+$ and $\Sigma^-$ of
$B_{r_0/C} \cap \Sigma \setminus \tilde  \sigma$ so that
$\Sigma^{\pm} \subset H^{\pm} \cup B_{C \, r_1}$ and $\tilde
\sigma =
\partial \Sigma^+ \cap \partial \Sigma^-$.
\end{enumerate}
\end{Pro}

As mentioned above, Proposition \ref{l:getsetgo} will be the key
tool for getting the limit plane through each point of $\cSt$ that
was promised in (C1) in Theorem \ref{t:t5.2}.  To see why, we will
first use Proposition \ref{l:getsetgo} to get a sequence of stable
graphs that are disjoint from $\Sigma_j$ and converge, away from
$y$, to a plane through $y$.
 Since the  outer
radii $R_j$ go to infinity, applying Proposition \ref{l:getsetgo}
to the sequence $\Sigma_j$ will give a sequence of stable graphs
that are disjoint from $\Sigma_j$ and defined over larger and
larger annuli centered at $y$.  As $j \to \infty$, the inner radii
of these annuli go to zero and the outer radii go to infinity.
Consequently,  the stable graphs
 will converge (subsequentially) to a minimal graph over a plane punctured
 at $y$ and this graph will have $y$ in its closure.  By a standard removable singularity theorem,
 this limit graph  extends smoothly across $y$ to an entire
 minimal graph and, hence, is flat by the Bernstein theorem.
  The easy details will be left
 to the reader.

\begin{figure}[htbp]
    \setlength{\captionindent}{20pt}
    %\begin{minipage}[t]{0.5\textwidth}
    \centering\input{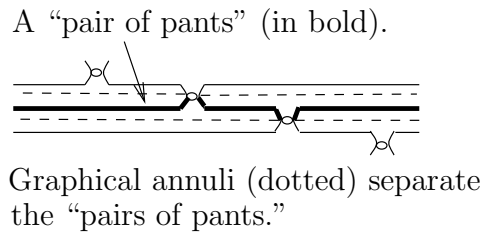}
    \caption{A pair of pants decomposition near a point
    where the injectivity radius goes to zero.}
    \label{f:pants}
    %\end{minipage}
\end{figure}

Now that we have this stable limit plane through $y \in \cSt$, the
proof of (C1) in Theorem \ref{t:t5.2} will consist of two main
steps.  We  will sketch these two steps next:
\begin{enumerate}
\item
{\bf{Decomposing $\Sigma_j$ into ULSC pieces}}: Let $\Sigma_j^+$
denote the portion of $\Sigma_j$ above the stable graph.  There
are now two possibilities:
 \begin{itemize} \item
$\Sigma_j^+$ is scale-invariant ULSC away from $y$.
 More
precisely, there exists a constant $C'$ and a sequence $r_j \to 0$
so that if $x \notin
 B_{ C' \, r_j} (y)$, then each
component  of $B_{|x|/ C'}(x) \cap \tilde{\Sigma}_j^+$ is a
disk.{\footnote{This intersection is empty when  $|x|$ is larger
than $C' \, R_j$.}}
 \item
Otherwise, we can apply Proposition \ref{l:getsetgo} to cut along
a second non-contractible curve; see Figure \ref{f:pants}.
\end{itemize}
In the second case, we replace $\Sigma_j^+$ by the portion
$\Sigma_j^{+-}$ of $\Sigma_j^+$  that is {\underline{below}} the
second stable graph.  After repeating this a finite number of
times, we will eventually get down to a scale-invariant ULSC
subset of $\Sigma_j$ with two interior boundary components (one
component for the cut near $y$ and one component for the last cut
that we make). \item {\bf{The ULSC pieces contain graphs}}:  In
either case, Lemma $3.3$ in \cite{CM8} will then give low points
in $\Sigma_j$ on either side of these stable graphs (see Figure
\ref{f:low}). Here ``low points'' roughly means points close to
the stable graph but away from its boundary. The one-sided
curvature estimate{\footnote{The one-sided curvature estimate is
recalled in this paper in Theorem \ref{t:t2}.}} from \cite{CM6}
and the gradient estimate will imply that the low points in the
resulting ULSC subsets of $\Sigma_j$ are graphical. Piecing this
together will easily give the desired global graphs.
\end{enumerate}

In the first case in (2) above, the graphs in the $\Sigma_j$'s
will be defined over annuli; in the second case, the graphs will
be over pairs of pants, i.e., over disks with two subdisks
removed. We will refer to the second case as
 a ``pair of pants''
decomposition; see Figure \ref{f:pants}.

The steps (1)--(3) above  are modelled on similar arguments for
  topological annuli in \cite{CM8}.  Some new complications will
  arise here because of the more complicated topological types of the
  surfaces, especially in the second case in step (2).

\subsection{The proof of Proposition \ref{l:getsetgo}: A decomposition near each point of $\cSt$}
\label{ss:41}

  The  next
lemma will first give stable surfaces disjoint from $\Sigma$ and
with ``interior boundary'' contained in a small ball.  In order to
prove Proposition \ref{l:getsetgo}, we will later show that these
stable surfaces contain the desired graphs.  More precisely, the
next lemma assumes that a component of a minimal planar domain
$\Sigma$ in a small ball is not a disk so that it must contain a
simple closed curve $\tilde \gamma$ separating two components
$\sigma_1$ and $\sigma_2$ of $\partial \Sigma$.  The lemma then
uses the separating curve $\tilde \gamma$ as ``interior boundary''
for a Plateau problem to get a stable minimal surface ``between''
$\sigma_1$ and $\sigma_2$; see (C) in Lemma \ref{l:getset} below.

\begin{figure}[htbp]
    \setlength{\captionindent}{20pt}
    %\begin{minipage}[t]{0.5\textwidth}
    \centering\input{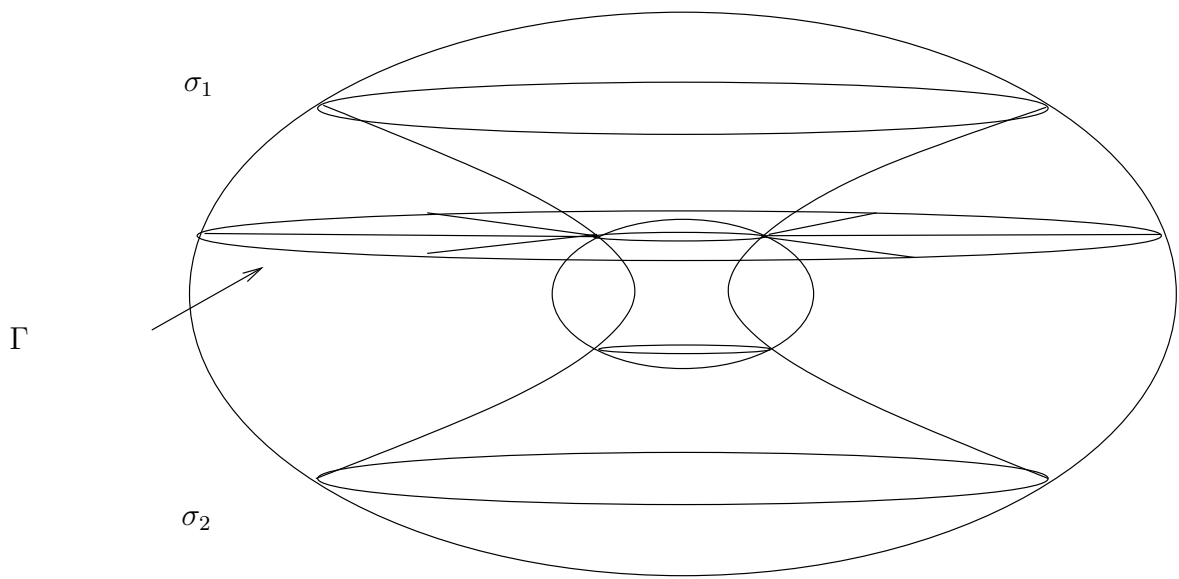}
    \caption{Lemma \ref{l:getset}:  If the planar domain $\Sigma$ contains a closed non-contractible curve in the
    small ball $B_{r_1}$, then $\Sigma$ has distinct boundary components $\sigma_1$ and $\sigma_2$.  Moreover,
    there is a stable surface $\Gamma$ that is disjoint from $\Sigma$ and separates $\sigma_1$ and $\sigma_2$
    in a component
    $\Omega$ of $B_{r_0} \setminus \Sigma$.  The boundary of $\Gamma$ has two parts, an outer boundary in
     $\partial B_{r_0}$ and an inner boundary curve $\tilde \gamma \subset \partial B_{r_1} \cap \Sigma$.}
    \label{f:gen63}
    %\end{minipage}
\end{figure}

\begin{Lem} \label{l:getset}
Let $\Sigma \subset B_{r_0}$ be an embedded minimal planar domain
with $\partial \Sigma \subset
\partial B_{r_0}$.  If $r_1 < r_0$ and a component $\Sigma_{r_1}$ of
$ B_{r_1}\cap \Sigma$ is not a topological disk, then the
following four properties hold (see Figure \ref{f:gen63}):
\begin{enumerate}
\item[(A)]
$\Sigma \setminus \Sigma_{r_1}$ has at least two connected
components; each of these components has at least one component of
$\partial \Sigma$ in its closure.
\item[(B)] If $\sigma_1$ and $\sigma_2$ are components of
$\partial \Sigma$ that are separated by $\Sigma_{r_1}$ (i.e.,
$\sigma_1$ and $\sigma_2$ are in the closure of distinct
components of $\Sigma \setminus \Sigma_{r_1}$), then we can choose
a simple closed curve \begin{equation}
    \tilde \gamma \subset \partial \Sigma_{r_1} \subset \partial
    B_{r_1}
    \end{equation}
that  separates $\sigma_1$ and $\sigma_2$ in $\Sigma$.
\item[(C)]
There is a component $\Omega$ of $B_{r_0} \setminus \Sigma$ and an
embedded {\underline{stable}} minimal surface $\Gamma \subset
\Omega$ with interior boundary $
\partial \Gamma \setminus \partial B_{r_0}$ equal to $\tilde \gamma$, and so that $\Gamma$
separates $\sigma_1$ and $\sigma_2$ in $\Omega$.{\footnote{One
must be careful interpreting this ``separation''  since $\partial
\Gamma$ may intersect $\sigma_1$ or $\sigma_2$.  In this case, we
mean that $\Gamma$ separates points in the interior of $\Sigma$
that are arbitrarily close to $\sigma_1$ and $\sigma_2$.}}
\item[(D)]  $\Gamma$ is area-minimizing amongst surfaces
 in $\Omega$ with boundary equal to $\partial \Gamma$.
\end{enumerate}
\end{Lem}

\begin{proof}
The first two claims, i.e., (A) and (B), follow from that $\Sigma$
is a planar domain and, by the maximum principle, any
homologically nontrivial curve in $\Sigma_{r_1}$ is also
homologically nontrivial in $\Sigma$.

We will next solve the Plateau problem to get the desired stable
surface in the complement of $\Sigma$.  To do this, we need to
choose the boundary of the stable surface and decide which of the
two components of $B_{r_0} \setminus \Sigma$ that we will solve
in.

 To get the boundary, simply let   $\Sigma_0$ be the component
  of $\Sigma \setminus \tilde \gamma$ that contains
$\sigma_1$ in its boundary; we will minimize area among surfaces
with boundary equal to
\begin{equation}
      \partial \Sigma_0 \, .
\end{equation}
Note that $\partial \Sigma_0$ has ``interior boundary'' equal to
$\tilde \gamma$ and $\sigma_2$ is not in $\partial \Sigma_0$.

We will use a simple linking argument to choose the domain
$\Omega$ to solve in.  First,
 fix a smooth curve $\eta \subset \Sigma$ from $\sigma_1$
to $\sigma_2$ that intersects $\tilde \gamma$ exactly once and
does so transversely (such a curve exists since $\tilde \gamma$
separates $\sigma_1$ and $\sigma_2$ in $\Sigma$).  Since $\Sigma$
is compact and embedded, we can ``push $\eta$ off of $\Sigma$'' -
on either side of $\Sigma$ - to get curves $\eta^+$ and $\eta^-$
that are
 disjoint from $\Sigma$ and in distinct components of $B_{r_0} \setminus \Sigma$.
 It follows that the (mod $2$)
linking numbers of $\eta^+$ and $\eta^-$, respectively, with
$\tilde \gamma$ differ by one.{\footnote{Recall that if $\eta
\subset B_r$ is a curve with endpoints in $\partial B_r$ and
$\gamma \subset B_r$ is a closed curve, then their linking number
is defined to be the number of times (mod $2$) that $\eta$
intersects a surface $\Gamma \subset B_r$ with $\partial \Gamma =
\gamma$. As usual, we assume that $\Gamma$ and $\eta$ intersect
transversely when counting intersections.  The point is that this
number does not depend on the particular choice of bounding
surface $\Gamma$.}}
 In particular, one of these - say $\eta^-$  -
has linking number $1$ (mod $2$) with $\tilde \gamma$.  Let
$\Omega$ be the component of $B_{r_0} \setminus \Sigma$ that
contains the {\it other} curve $\eta^+$.

 It follows that we have the
following three properties:
\begin{itemize}
\item The domain $\Omega$ is mean convex in the sense of
\cite{MeYa2}.
 \item  $\partial \Sigma_0$ is contained in $\partial
\Omega$ and bounds the planar domain $\Sigma_0$  in $\partial
\Omega$. \item $\tilde \gamma = \partial \Sigma_0 \setminus
\partial B_{r_0}$ has linking number $1$ (mod $2$) with the curve
$\eta^-$ {\it that is not in} $\Omega$.
\end{itemize}
 Using the first two properties, a result of Hardt-Simon, \cite{HSi},
 gives an embedded minimal
surface $\tilde \Gamma \subset \Omega$ with $\partial \tilde
\Gamma =
\partial \Sigma_0$ and so that $\tilde \Gamma$ minimizes area
amongst surfaces in $\Gamma$ with the same boundary.{\footnote{We
could of course have applied a result of Meeks-Yau to get a stable
planar domain. However, this planar domain would have minimized
area only amongst planar domains.  We will later use this
minimizing property to bound the area of $\Gamma$ by constructing
comparison surfaces.  It will be convenient not to have to
restrict the topological type of the comparison surfaces.}} In
particular, $\tilde \Gamma$ must be stable.

The surface $\tilde \Gamma$ may have several components.  We will
use the third property above to show that the component $\Gamma$
containing $\tilde \gamma$ in its boundary separates $\sigma_1$
and $\sigma_2$ in $\Omega$ and, hence, satisfies (C).  First of
all, $\tilde \gamma$ alone cannot be the entire boundary of
$\partial \Gamma$; indeed, any surface $\Gamma_{\tilde \gamma}
\subset B_{r_0}$ with $\partial \Gamma_{\tilde \gamma} = \tilde
\Gamma$ would be forced to intersect the curve  $\eta^{-}$  that
is not in $\Omega$. Therefore, we must have that
\begin{equation}
    \partial B_{r_0} \cap \partial \Gamma \ne \emptyset \, .
\end{equation}
A similar use of the linking condition implies that $\Gamma$
separates $\sigma_1$ and $\sigma_2$ in $\Omega$, giving (C).
\end{proof}

The minimizing property given in (D) will be used to give an upper
bound for the area of $\Gamma$ by constructing comparison surfaces
with the same boundary.  To carry this out, we will need two
elementary lemmas. The first is a simple topological lemma showing
that any collection of disjoint simple closed curves in a planar
domain is homologous to a collection of {\underline{distinct}}
boundary curves. Moreover, together the initial curves and the
boundary curves bound a subdomain of the planar domain.  The
second lemma uses this to construct comparison surfaces and hence,
using (D),  deduce an area bound for the surface $\Gamma$ above.

\begin{Lem}     \label{l:pdom}
Let $P$ be a (possibly disconnected) compact planar domain with
boundary $\partial P$.  Given any collection $\sigma_1 , \dots ,
\sigma_n$   of disjoint simple closed curves in $P$, then there is
a subdomain $P_0 \subset P$ with
\begin{equation}        \label{e:fouroneseven}
    \partial P_0 = \left( \cup_{i=1}^n \sigma_i \right) \cup   \left( \cup_{i=1}^m \eta_i
    \right) \, ,
\end{equation}
where $\eta_1 , \dots , \eta_m$ are distinct components of
$\partial P$.
\end{Lem}

\begin{Rem}
Before giving the proof of Lemma \ref{l:pdom}, it may be helpful
to make two remarks.  First,  it is possible that $m = 0$ in
\ref{e:fouroneseven}, i.e., that $\partial P_0 = \cup_{i=1}^n
\sigma_i$.  Second, notice that all of the above curves - both the
$\sigma_i$'s and $\eta_i$'s - are thought of as
{\underline{unoriented}} curves.
\end{Rem}

\begin{proof}
(of Lemma \ref{l:pdom}.) Since we can consider each connected
component of $P$ separately, we may as well assume that $P \subset
\RR^2$ is connected.  The set $P_0$ will be given as the level set
$f^{-1} (1)$ of a map
\begin{equation}
    f: P \setminus
    \left( \cup_{i=1}^n \sigma_i \right) \to \{ -1 , + 1 \} \, .
\end{equation}
To define $f$, first fix a point $p_0 \in P \setminus \left(
\cup_{i=1}^n \sigma_i \right)$.  For each point $p \in P \setminus
\left( \cup_{i=1}^n \sigma_i \right)$, choose a curve $\gamma_p
\subset P$ from  $p_0$ to $p$ that is transverse to the
$\sigma_i$'s and let $n(p)$ be the number of times that $\gamma_p$
crosses $\left( \cup_{i=1}^n \sigma_i \right)$; see Figure
\ref{f:gen60}. It follows from elementary topology that $n(p)$
(mod two) does not depend on the choice of the curve $\gamma_p$
and, hence, we can define the function $f$ by
\begin{equation}
    f(p) = (-1)^{n(p)} \, .
\end{equation}
Define $P_0$ by
\begin{equation}
    P_0 = \{ p \, | \, f(p) = + 1 \} \, .
\end{equation}
It follows easily that $f$ changes sign as we cross each
$\sigma_i$ and, therefore, each $\sigma_i \subset \partial P_0$ as
desired.
\end{proof}

\begin{figure}[htbp]
    \setlength{\captionindent}{20pt}
    %\begin{minipage}[t]{0.5\textwidth}
    \centering\input{gen60.pstex_t}
    \caption{The proof of Lemma \ref{l:pdom}: The curves $\sigma_i$ are dashed while
    $\partial P$ is dotted. After fixing a basepoint $p_0$, we define a function $f(p)$
    to be $+1$ or $-1$, depending on whether you have to cross an even or odd number of the $\sigma_i$'s to
    get from $p_0$ to $p$.  The domain $P_0$ is then defined to be the level set where $f$ is $+1$.
    In the example pictured, $P_0$ has two components.}
    \label{f:gen60}
    %\end{minipage}
\end{figure}

\begin{Lem}     \label{l:comparison}
Given a constant $C_1$, there exists  $C_2 > C_1$ so that if
$\Sigma \subset B_{r_0}$, $\tilde \gamma$, and $\Gamma \subset
\Omega \subset B_{r_0} \setminus \Sigma$ are as in Lemma
\ref{l:getset} and for some $r$ between $r_0$ and $r_1$ we have
\begin{equation}    \label{e:thebs}
    r^{-1} \, \Length
    (\partial B_r \cap \Sigma) \leq C_1 \, ,
\end{equation}
then we get an area bound for $B_r \cap \Gamma$
\begin{equation}    \label{e:tdabc2}
    r^{-2} \, \Area \, (B_r \cap \Gamma) \leq C_2 \, .
\end{equation}
\end{Lem}

\begin{proof}
Note first that \eqr{e:thebs} and Stokes' theorem (using that $\dv
(\nabla |x|^2) = 4$ on a minimal surface) gives
\begin{equation}    \label{e:thebs2}
     r^{-2} \, \Area \, (B_r \cap \Sigma) \leq 1/2 \,  r^{-1} \, \Length
    (\partial B_r \cap \Sigma) \leq C_1 / 2 \, .
\end{equation}
The minimizing property (D) in Lemma \ref{l:getset} implies that
$B_r \cap \Gamma$ is itself area minimizing among surfaces in its
homology class in $\Omega$.  It will suffice therefore to
construct a comparison surface in $\overline{\Omega}$ with bounded
area. This follows from the following steps:
\begin{enumerate}

 \item
  The outer boundary of $B_r \cap \Gamma$ -- i.e., $(\partial
B_r) \cap \Gamma$ -- sits inside the (possibly disconnected)
planar domain $\Omega \cap
\partial B_r$.
Consequently, Lemma \ref{l:pdom} gives a subset $P_0$ of the
planar domain $ \Omega \cap
\partial B_r$ with
\begin{equation}
    (\partial B_r \cap \Gamma )\subset \partial P_0 {\text{ and }}
        \partial P_0 \setminus (\partial B_r \cap \Gamma) \subset
(\partial B_r) \cap \Sigma \, .
\end{equation}
 Note that $P_0$ has
bounded area since it is contained in $\partial B_r$.
\item  Let $\Sigma_0$ denote the component of $B_r \cap \Sigma$
containing $\tilde \gamma$ and let $\Sigma_0^+$ be one of the
components of $\Sigma_0 \setminus \tilde \gamma$.  Note that
$\Sigma_0^+$ has bounded area by \eqr{e:thebs2}.
\item  By the previous two steps, the boundary of $P_0 \cup \Sigma_0^+$
contains all of $\partial (B_r \cap \Gamma)$ (both $\tilde \gamma$
and the outer boundary).  Ideally, we would have that $\partial
(P_0 \cup \Sigma_0^+) = \partial (B_r \cap \Gamma)$ so that $P_0
\cup \Sigma_0^+$ would be a valid comparison surface.  However,
this does not have to be the case since $\partial (P_0 \cup
\Sigma_0^+)$ may have some additional components. We can therefore
assume that
\begin{equation}
\partial
(P_0 \cup \Sigma_0^+) \setminus \partial (B_r \cap \Gamma) \ne
\emptyset \, .
\end{equation}
\item
Observe that $\partial (P_0 \cup \Sigma_0^+) \setminus
\partial (B_r \cap \Gamma)$ is itself the boundary of a surface in $\overline
\Omega$, namely of the surface $P_0 \cup \Sigma_0^+ \cup (B_r \cap
\Gamma)$.  We can therefore solve the Plateau problem for a
surface $\tilde \Gamma$ in $B_r \cap \Omega$ with
\begin{equation}
    \partial \tilde \Gamma = \partial (P_0 \cup \Sigma_0^+) \setminus
\partial (B_r \cap \Gamma)  \subset (\partial B_r) \cap \Sigma \,
.
\end{equation}
The length bound on $(\partial B_r) \cap \Sigma$ and the
isoperimetric inequality for minimal surfaces then give an area
bound for $\tilde \Gamma$.
\item
Finally, it follows that
\begin{equation}
    \partial ( P_0 \cup \Sigma_0^+ \cup \tilde \Gamma ) = \partial (B_r \cap
    \Gamma) \, ,
\end{equation}
so that $P_0 \cup \Sigma_0^+ \cup \tilde \Gamma$ is the desired
comparison surface.{\footnote{This surface is not embedded and it
may not even have the same topological type, but it is nonetheless
a valid comparison surface.}}  This gives \eqr{e:tdabc2} since
each of the three pieces of the comparison surface has the desired
area bound.
\end{enumerate}
\end{proof}

The next lemma gives an area estimate for the components of the
 stable surface constructed in Lemma \ref{l:getset} on the largest scale $r_1$ where $\Sigma$ is ULSC.
 Recall that $\Sigma$ contains a non-contractible curve $\tilde
 \gamma$ in $\partial B_{r_1}$, so $\Sigma$ is not ULSC on scales
 larger than $r_1$.  On the other hand, the assumption
 \eqr{e:ulscin} below gives that $\Sigma$ is ULSC on this scale.

\begin{Lem}     \label{l:lareab}
Given a constant $C_1$, there exists  $C_2 > C_1$ so that if
$\Sigma \subset B_{r_0}$, $\tilde \gamma$, and $\Gamma \subset
\Omega \subset B_{r_0} \setminus \Sigma$ are as in Lemma
\ref{l:getset} with $r_1 < r_0/C_2$ and, in addition,
\begin{equation}    \label{e:ulscin}
    \cB_{r_1/4} (x) \subset \Sigma {\text{ is a disk for every }} x \in \cB_{C_2\,r_1} \,
    ,
\end{equation}
then each component $\Gamma'$  of $B_{C_1 \, r_1}  \cap \Gamma
\setminus B_{8 \,
 r_1} $ which can be connected to $\tilde \gamma$ by a curve in $B_{17 \,
 r_1}  \cap \Gamma$ satisfies
\begin{equation}     \label{e:areab}
    \Area \, (\Gamma ' ) \leq C_2 \, r_1^2 \, .
\end{equation}
\end{Lem}

\begin{proof}
 Observe first that the chord-arc bound
for ULSC surfaces, Lemma \ref{l:cyulsc} in Appendix \ref{s:aC},
shows that the ULSC hypothesis \eqr{e:ulscin} also holds for $x$
in the component of $\Sigma$ containing $0$ in an extrinsic ball
$B_{C_2' \, r_1}$ where $C_2'$ goes to infinity as $C_2$ does. In
particular, after replacing $\Sigma$ by this component, we may as
well assume that
\begin{equation}    \label{e:ulscout}
    \cB_{r_1/4} (x) {\text{ is a disk for every }} x \in \Sigma \,
    .
\end{equation}

The proof of  \eqr{e:areab} will be by contradiction, using a
compactness argument.  Suppose therefore that $\Sigma_j$,
$\Gamma_j$ is a sequence of counter-examples where \eqr{e:areab}
fails with $C_2 = j \to \infty$. After translating and rescaling,
we may assume that $r_1 = 4$.

We will consider two cases depending on whether $|A|^2 \to \infty$
on a compact set for the sequence $\Sigma_j$.

\vskip2mm \noindent {\underline{Case 1}}:
 Suppose first that $|A|^2
\to \infty$ in some fixed ball of $\RR^3$ (for some subsequence of
the $\Sigma_j$'s).
 We can then apply the
compactness theorem for non-simply connected ULSC sequences,
Theorem \ref{t:t5.1}, to get a subsequence of the $\Sigma_j$'s
that converges to a foliation by parallel planes away from two
lines orthogonal to the leaves of the foliation.

Observe that both of these orthogonal ``singular'' lines intersect
the leaf through $0$ inside the ball $B_{4}$ since $\Sigma_j$ is
assumed to be non-simply connected inside that ball. It follows
easily from the description of the convergence near the lines (as
oppositely oriented double spiral staircases) that any such
component $\Gamma_j'$ is sandwiched between the almost planar
leaves of the foliation (for $j$ sufficiently large).  This
sandwiching, together with interior curvature estimates for stable
surfaces, implies that $\Gamma_j'$ is itself a graph and hence has
bounded area as desired.

\vskip2mm \noindent {\underline{Case 2}}: Suppose now that $|A|^2$
is uniformly bounded on compact subsets of $\RR^3$ for every
$\Sigma_j$.  In this case, we will get uniform area bounds for the
surfaces $\Sigma_j$ in the ball $B_{4\, C_1}$.  Once we have these
area bounds for the $\Sigma_j$'s, then the comparison argument in
Lemma \ref{l:comparison} will give a uniform bound for area of the
$\Gamma_j$'s in the same ball.  However, we assumed that there was
no such area bound for $\Gamma_j$'s.  This  contradiction will
complete the proof of the lemma.

Therefore, to complete the proof of the lemma, it suffices to
bound the area of $\Sigma_j$ in $B_{4\, C_1}$.  This area bound
follows immediately from combining the following two facts:
\begin{itemize}
\item
By Lemma \ref{l:ca} in Appendix \ref{s:aC}, the uniform curvature
bound implies uniform area bounds for each component of $\Sigma_j$
in extrinsic balls (the bound depends on the ball but not on $j$).
More precisely, if $\Sigma_{j,R}$ is a component of $B_R \cap
\Sigma_j$, then Lemma \ref{l:ca} implies that
\begin{equation}    \label{e:fromlca0}
    \Area \, (\Sigma_{j,R}) \leq C_c \, R^2 \, ,
\end{equation}
where the constant $C_c$ depends only on the supremum of $|A|^2$
on $B_{C_0 \, R} \cap \Sigma_j$.  The constant $C_0$ here is
universal and does not depend on the upper bound for the
curvature. \item  Even though each intrinsic ball of radius one in
$\Sigma_j$ is a disk, there is a component of $B_4 \cap \Sigma_j$
that is {\underline{not}} a disk.
 Therefore, it follows easily from a barrier argument and the one-sided
lemma for non-simply connected surfaces, Lemma \ref{l:os}, that
there exists $R > 4 \, C_1$ so that only one component of $B_R
\cap \Sigma_j$ intersects $B_{4\, C_1}$  for $j$
large.{\footnote{This follows exactly as does the analogous result
for disks given in  corollary $0.4$ in    \cite{CM6}.  Namely, if
there were two such components, then we could put a stable surface
between them.  Interior estimates for stable surfaces then imply
that each of the original components lies on one side of a plane
that comes close to the center of the ball.  However, this would
contradict  the one-sided  lemma for non-simply connected
surfaces, Lemma \ref{l:os}, so we conclude that there could not
have been two such components.}}
\end{itemize}
\end{proof}

The last result that we will need to recall before proving
Proposition \ref{l:getsetgo} is the following elementary property
of connected planar domains:

\begin{Lem}     \label{l:espdq}
Let $ \Sigma$ be a connected planar domain and $\sigma_1 , \dots ,
\sigma_n$ the components of $\partial \Sigma$.  Given $k < n$ and
a collection $\{ \sigma_{i_1} , \dots , \sigma_{i_k} \}$, there is
a simple closed curve $\tilde \sigma \subset \Sigma$ which
separates $\cup_{j\leq k} \sigma_{i_j}$ from $\partial  \Sigma
\setminus \cup_{j\leq k} \sigma_{i_j}$.
\end{Lem}

\vskip2mm

\begin{proof}
(of Proposition \ref{l:getsetgo}).
 We will first use a rescaling
argument to locate the smallest scale of non-trivial topology,
choose a non-contractible curve $\gamma$ on this scale, and then
solve the Plateau problem with $\gamma$ as interior boundary.  We
will then obtain an area bound for the components of $\Gamma$ on
this scale. This area bound will allow us to apply the ``stable
graph proposition'', Proposition \ref{l:lext0}, to get the graph
$\Gamma_0 \subset \Gamma$ and thus prove (A).  Finally, in the
last step of the proof, we will find the separating curve $\tilde
\sigma$ and prove (B).

 {\underline{Blowing up   on the smallest scale of non-trivial topology}}.
Fix a large constant $C_1 > 1$ to be chosen.  Applying the blow up
lemma, Lemma
            \ref{l:l5.1} in Appendix \ref{s:aE}, at $0$ gives an intrinsic
            ball
\begin{equation}   \label{e:bu1}
                \cB_{C_1 s_1} (y_1) \subset \cB_{5 \, C_1 \, r_1}   \, ,
\end{equation}
so that $\cB_{4\, s_1}(y_1)$ is not a disk but
                $\cB_{s_1} (y)$ is a disk for each $y \in \cB_{C_1 s_1}
                (y_1)$.    We can now use this topologically
                non-trivial region $\Sigma$ to solve a Plateau
                problem.  Namely,
                applying  Lemma
                \ref{l:getset} to the component of $B_{4 \, s_1}(y_1) \cap \Sigma$
                containing $\cB_{4\, s_1}(y_1)$ gives  a simple closed
                non-contractible{\footnote{Since $\Sigma$ has non-positive curvature and $\gamma$ is non-contractible in
                the intrinsic ball $\cB_{4\, s_1}(y_1)$, it is also non-contractible in $\Sigma$.}}
                curve $\gamma \subset B_{4 \, s_1}(y_1) \cap \Sigma$, a mean convex domain $\Omega \subset
                B_{r_0} \setminus \Sigma$, and
                a stable embedded minimal surface
\begin{equation}
    \Gamma \subset \Omega \, ,
\end{equation}
with interior boundary $\partial \Gamma \setminus \partial
B_{r_0}$ equal to   $\gamma$. Moreover, there are distinct
components $\sigma_1$ and $\sigma_2$ of $\partial \Sigma \subset
\partial B_{r_0}$ that are separated in $\Sigma$ by $\gamma$ and
separated in $\Omega$ by $\Gamma$.  Finally, $\Gamma$ is
area-minimizing amongst surfaces in $\Omega$ with boundary equal
to $\partial \Gamma$.

{\underline{An area estimate on the smallest scale of non-trivial
topology}}.
    Suppose that $\Gamma'$ is a component of $B_{24 \, s_1}(y_1) \cap \Gamma \setminus B_{8 \,
 s_1}(y_1)$ which can be connected to $\gamma$ by a curve in $B_{17 \,
 s_1}(y_1) \cap \Gamma$.  If
the constant $C_1$ from the previous step is sufficiently large
(independent of $\Sigma$ and $\Gamma$), then Lemma \ref{l:lareab}
gives a constant $C_2$ so that
\begin{equation}     \label{e:areabn}
    \Area \, (\Gamma ' ) \leq C_2 \, s_1^2 \, .
\end{equation}

{\underline{Finding the graph in $\Gamma$}}. Using the area bound
 \eqr{e:areabn}, we can  apply Proposition \ref{l:lext0} to get that each component
of $B_{r_0/C} \cap  \Gamma \setminus B_{Cr_1}$ is a graph.
 A linking
argument as in Lemma \ref{l:getset} then implies that one of these
components $\Gamma_0$ has the property that $\Gamma_0 \cup B_{C
r_1}$ separates $B_{r_0/C}$ into components $H^+$ above and $H^-$
below $\Gamma_0$ where $\sigma_1 \subset H^+$ and $\sigma_2
\subset H^-$.{\footnote{Technically, this is not quite right since
$\partial \Sigma$ is contained in the boundary of the larger ball
$B_{r_0}$.  Rather, the linking argument gives two components --
call them $\tilde \sigma_1$ and $\tilde \sigma_2$ -- of $\partial
B_{r_0/C} \cap \Sigma$ that are separated by $\Gamma_0$.}}

{\underline{Finding the separating curve}}.
To complete the proof,
we need only
 find the separating curve $\tilde \sigma \subset \Sigma$ and prove (B).  In doing
this, we will increase the constant $C$ several times below.

The key to this step is to prove that there is a constant $C_3 >
1$ so that \begin{equation}   \label{e:onlyonecomp}
    {\text{only one component
 of $B_{C_3 r_1} \cap \Sigma$ intersects both $H^+$
and $H^-$.}}
\end{equation}
Before proving \eqr{e:onlyonecomp}, it may be helpful to make a
few remarks. First, it is not hard to see that \eqr{e:onlyonecomp}
is necessary to establish (B).  Namely, if there were two distinct
  components of $B_{C_3 r_1} \cap \Sigma$ that each connected $H^+$
and $H^-$, then it would be impossible to find a single connected
curve in $B_{C_3 r_1} \cap \Sigma$ that separates $H^+$ and $H^-$.
Second, it is easy to see that there must be at least one
component of $B_{C_3 r_1} \cap \Sigma$ that intersects both $H^+$
and $H^-$. This is because $\Sigma$  has boundary components
 $\sigma_1 \subset H^+$ and $\sigma_2 \subset H^-$ and the only
 way to connect these without crossing the annular graph $\Gamma_0$ is to go through
 the ``hole'' in the middle.  Finally,  the basic idea behind
 \eqr{e:onlyonecomp}  is that if there were two components passing
 through the ``hole'' in $\Gamma_0$, then a barrier argument would
 also give a stable surface between the two components
 that also passes through the
hole.  However, such a stable surface would have to be very flat
if $C_3$ is large, so it cannot pass through this hole.  This is
essentially the argument that we will give below, but it will take
a little work to make it precise.

\begin{figure}[htbp]
 \setlength{\captionindent}{20pt}
    \begin{minipage}[t]{0.5\textwidth}
    \centering\input{genp32a.pstex_t}
    \caption{The key step in finding the separating curve: Ruling out that
    two components of $B_{C_3 r_1} \cap \Sigma$ both intersect both
    $H^+$ and $H^-$.}
    \label{f:32a}
    \end{minipage}\begin{minipage}[t]{0.5\textwidth}
    \centering\input{genp32b.pstex_t}
    \caption{The contradiction: We cannot have a second stable graph that is
    between the two components of $\Sigma$
    and on one side of $\Gamma_0$.}
    \label{f:32b}
    \end{minipage}
\end{figure}

  Once we establish \eqr{e:onlyonecomp}, the rest of the proof of the proposition will follow
  easily.  Namely,  if let $\hat \Sigma$  be the component of $B_{C r_1} \cap \Sigma$ intersecting both
$H^+$ and $H^-$, then Lemma \ref{l:espdq} gives a simple closed
curve $\tilde \sigma \subset \hat \Sigma$ separating $H^+ \cap
\partial \hat \Sigma$ from $H^- \cap \partial \hat \Sigma$.  In
particular, the components $\Sigma^{\pm}$ of $\Sigma \setminus
\tilde \sigma$ satisfy $\Sigma^{\pm} \subset H^{\pm} \cup B_{C \,
r_1}$.

\vskip2mm  Finally, to complete the proof of the proposition, it
remains only to prove \eqr{e:onlyonecomp}.  We will do this by
contradiction, so suppose  that $\hat \Sigma_1$ and $\hat
\Sigma_2$ are distinct components of $B_{C_3 r_1}\cap \Sigma$ each
of which intersects both $H^+$ and $H^-$.  Since the only ``hole''
in the graph $\Gamma_0$ is in $B_{Cr_1}$, there must be components
$\tilde \Sigma_i \subset \hat \Sigma_i$  of $B_{ C r_1} \cap
\Sigma$ intersecting both $H^+$ and $H^-$; see Figure \ref{f:32a}.
Label these components so $\gamma \cap \hat \Sigma_2 = \emptyset$.
To get the contradiction, we will solve a Plateau problem to get a
second stable graph that is between the $\tilde \Sigma_i$'s and
also disjoint from the graph $\Gamma_0$; such a graph would be
forced to sit on one side of $\Gamma_0$ and hence would not allow
both of the $\tilde \Sigma_i$'s to intersect both $H^+$ and $H^-$;
see Figure \ref{f:32b}.

To set this up, note first that we can assume that $\Gamma_0$ is a
graph with arbitrarily small gradient -- say at most $\delta >0$
-- after possibly increasing $C$.  This follows from  estimates
for minimal graphs; see  proposition $1.12$ in \cite{CM9}. After a
rotation of $\RR^3$, we can assume that $\Gamma_0$ is a graph over
the horizontal plane $\{ x_3 = 0 \}$.

Fix a point $y_1$ in $B_{C r_1} \cap \tilde \Sigma_1$ and choose a
point $y_2$ in $B_{2 C r_1} \cap \tilde \Sigma_2$ so that the
segment $\gamma_{y_1,y_2}$ from $y_1$ to $y_2$ is ``almost
vertical''; see Figure \ref{f:32a}.  More precisely, applying
lemma A.8 of \cite{CM3} (as in (I.0.20) of \cite{CM3}) gives $y_2
\in B_{2 C r_1} \cap \tilde \Sigma_2$ with
\begin{equation}
  \label{e:aveq} |\Pi (y_2 - y_1)| \leq |y_2 - y_1| \, \cos
  \theta_0 \, ,
 \end{equation}
 where $\Pi$ is orthogonal projection to the horizontal plane and the constant
 $\theta_0$ is defined in the
appendix of \cite{CM3}.   It now follows that there is a component
$\tilde \Omega$ of $B_{C_3 r_1}  \setminus (\Gamma \cup \Sigma)$
so (some subsegment of) $\gamma_{y_1,y_2}$ is linked with
$\partial \hat \Sigma_2$ in $\tilde \Omega$. Note that $\tilde
\Omega$ is mean convex in the sense of \cite{MeYa2}.  A result of
\cite{MeYa1}--\cite{MeYa2} gives a stable embedded minimal surface
\begin{equation}
    \hat
\Gamma^0 \subset \tilde \Omega
\end{equation}
 with $\partial \hat
\Gamma^0 =
\partial \hat \Sigma_2$. Since $\partial \hat \Sigma_2$ and
$\gamma_{y_1,y_2}$ are linked in $\tilde \Omega$, a component
$\hat \Gamma$ of $B_{C_4 C r_1} \cap \hat \Gamma^0$ intersects
$\gamma_{y_1,y_2}$ at least once.  However, combining curvature
estimates \cite{Sc1}, \cite{CM2} for stable surfaces with the fact
that $\hat \Gamma^0$ is disjoint from the graph $\Gamma_0$ with
small gradient (but comes close to this graph) implies that $\hat
\Gamma^0$ is also a graph with small gradient over the horizontal
plane.  Choosing the constants appropriately so the gradient of
these graphs is sufficiently small, we see that $\hat \Gamma^0$
can only intersect the ``almost vertical'' segment $
\gamma_{y_1,y_2}$ exactly once (see \eqr{e:aveq}).   In
particular, $\hat \Gamma^0$ separates $\tilde \Sigma_1$ and
$\tilde \Sigma_2$, forcing one of these to lie on the same side of
$\Gamma_0$ as does $\hat \Gamma^0$. This gives the desired
contradiction; see Figure \ref{f:32b}. Consequently, we conclude
that $B_{C_3 r_1}(x_1) \cap \Sigma$ contains only one component
$\hat \Sigma$ which intersects both $H_1^+$ and $H_1^-$, i.e.,
\eqr{e:onlyonecomp} holds.
\end{proof}

\subsection{Step (1): Decomposing $\Sigma_j$ into ULSC pieces}	\label{ss:step1}

 Suppose
now that $0 \in \cSt$, so that Proposition \ref{l:getsetgo} gives
\begin{enumerate}
\item A sequence of stable graphs $\Gamma_j$ that are disjoint
from $\Sigma_j$ and
 that converge to a punctured plane through $0$; after rotating $\RR^3$, we can assume that
 the stable graphs converge to $\{ x_3 = 0 \} \setminus \{ 0 \}$. \item
    A sequence of closed curves
$\tilde \sigma_j \subset B_{r_j} \cap \Sigma_j$ with $r_j \to 0$
and so that $\tilde \sigma_j$ divides $\Sigma_j$ into a component
$\Sigma_j^+$ above $\{ x_3 = 0 \}$ and a component $\Sigma_j^-$
below.{\footnote{More precisely, there are shrinking extrinsic
balls $B_{r_j}$ so that $\Sigma_j^+ \setminus B_{r_j}$ is above
$\Gamma_j$ and similarly for $\Sigma_j^-$.}}
\end{enumerate}

We will show next that each $\Sigma_j^+$ contains a large
scale-invariant ULSC piece $\Sigma_j^{ulsc}$. Before stating this
precisely, it may be helpful to recall two simple examples:
\begin{itemize}
\item If we consider  a sequence of shrinking catenoids,
then $0$ is the only point in $\cSt$ and each half of the
catenoids is easily seen to be scale-invariant ULSC (in fact,
given any point $x \ne 0$ in one of the catenoids, the ball
$B_{|x|}(x)$ has two simply connected components).
\item Consider now a sequence of rescalings of one of the Riemann
examples.  In this case, $\cSt$ is a line through the origin and
  $\Sigma_j^+$ is {\it{not}} scale-invariant ULSC.  However, if we
  cut $\Sigma_j^+$ along a second short curve (the ``neck'' immediately above the first curve),
  then the resulting ``pair of pants'' is scale-invariant ULSC with respect to the distance
  to the closer of the two necks.
\end{itemize}

The precise statement of the decomposition into ULSC pieces is
given in the next lemma.  For simplicity, we will suppose that $0
\in \cSt$, $\Sigma_j^+ \subset \Sigma_j$ are as above, and the
constant $C$ is given by Proposition \ref{l:getsetgo}.

\begin{Lem}     \label{l:ulscpie}
Let $\tilde{\Sigma}_j^+$ denote the connected component of
$B_{R/(2C)} \cap \Sigma_j^+$ with $\tilde \sigma_j$ in its
boundary and fix a constant $\alpha
> 1$. For each $j$ sufficiently large, one of the following two cases
holds:
 \begin{enumerate} \item
$\tilde{\Sigma}_j^+$ is scale-invariant ULSC: Given any $x \notin
 B_{ \alpha \, C \, r_j}$, then each
component  of $B_{|x|/(\alpha C)}(x) \cap \tilde{\Sigma}_j^+$ is a
disk.
\item
$\tilde{\Sigma}_j^+$ contains a non-contractible curve $\tilde
\sigma_j^+$ in a ball $B_{s_j}(y_j)$ with
\begin{equation}
    |y_j| > \alpha \, C
\, r_j {\text{ and }} s_j < |y_j|/(C\alpha) \, ,
\end{equation}
 so
that the component $\Sigma_j^{ulsc}$ of $\tilde{\Sigma}_j^+
\setminus \tilde \sigma_j^+$ with $\tilde \sigma_j$ in its
boundary is scale-invariant ULSC:  Given any $x \notin \left( B_{
\alpha \, C \, r_j} \cup B_{  \alpha \, C \, s_j}(y_j) \right)$,
then each component of
\begin{equation}
    B_{\frac{\min \{ |x| , \, |x-y_j| \} }{C \alpha} }(x) \cap
\Sigma_j^{ulsc}
\end{equation}
 is a disk.
\end{enumerate}
\end{Lem}

\begin{figure}[htbp]
    \setlength{\captionindent}{20pt}
    %\begin{minipage}[t]{0.5\textwidth}
    \centering\input{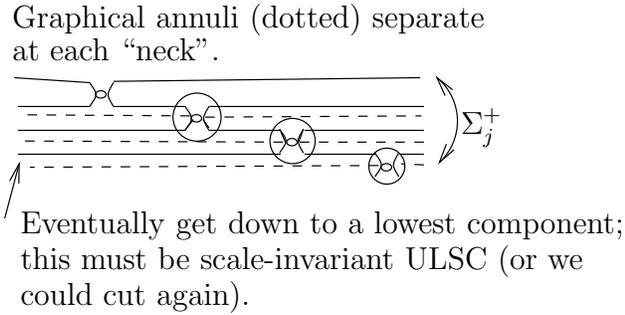}
    \caption{Cutting repeatedly to get the pair of pants decomposition.}
    \label{f:pp}
    %\end{minipage}
\end{figure}

\begin{proof}
 The key for establishing this lemma is that the decomposition into
 a $\Sigma_j^+$ and a $\Sigma_j^-$ can be repeated anywhere that
 the topology is concentrating.  Namely, suppose that (1) does not hold and, hence,
 there exists some
$z_1$ in $B_{R_j/(2C)} \setminus B_{ \alpha \, C \, r_j}$ so that
some component  of
\begin{equation}
    B_{|z_1|/(\alpha C)}(z_1) \cap  \tilde{\Sigma}_j^+
\end{equation}
 is not a disk.
We can repeat the argument of Proposition \ref{l:getsetgo} to get
a second stable graph $\Gamma_j'$, separating curve $\tilde
\sigma_j'$, and components $\Sigma_j^{++}$ and  $\Sigma_j^{+-}$ of
$\Sigma_j^+ \setminus \tilde \sigma_j'$ that are above and below,
respectively, the graph $\Gamma_j'$; see Figure
\ref{f:pp}.{\footnote{Proposition \ref{l:getsetgo} directly gives
the second stable graph disjoint from $\Sigma_j$.  However, the
proposition does not explicitly give that the two stable surfaces
can be chosen to be disjoint.  This is easy to achieve since the
components of $B_{R_j} \setminus (\Sigma_j \cup \Gamma_j)$ are
also mean convex in the sense of Meeks-Yau.}}
 Observe that the ``middle
component'' $\Sigma_j^{+-}$ is between the two stable graphs and
has only two components in its interior boundary.  If
$\Sigma_j^{+-}$  satisfies (2), then we are done. Otherwise, there
is a third non-contractible simple closed curve.  We can repeat
the argument to cut $\Sigma_j^{+-}$ to get an even lower component
$\Sigma_j^{+--}$.  The key point is that this new surface
$\Sigma_j^{+--}$ also has only two components in its interior
boundary.  Since $\Sigma_j$ is compact, this process must
eventually terminate to give a lowest component
$\Sigma_j^{+-\dots-}$ that satisfies (2).
\end{proof}

\begin{Rem}
The reader may find the constant $\alpha$ in Lemma \ref{l:ulscpie}
somewhat mysterious.  The point is that taking $\alpha$ large
forces the two interior boundary components in case (2) to be
relatively far apart.  This will be used to guarantee that the
ULSC piece is sufficiently large, i.e., goes all the way out to
the outer boundary in $\partial B_{R_j/(2C)}$.
\end{Rem}

\subsection{Step (2): The ULSC pieces of $\Sigma_j$ contain graphs} \label{ss:step2}

We will next find the graphs in $\Sigma_j$ converging to the plane
$\{ x_3 = 0 \}$ away from $0$ and possibly one other point.  The
argument for this is slightly simpler in case (1) where
$\tilde{\Sigma}_j^+$ is itself scale-invariant ULSC and we do not
need to cut along a second curve, but this simpler case already
illustrates the key ideas.  The argument follows a similar one in
\cite{CM14}.

\begin{figure}[htbp]
    \setlength{\captionindent}{20pt}
    %\begin{minipage}[t]{0.5\textwidth}
    \centering\input{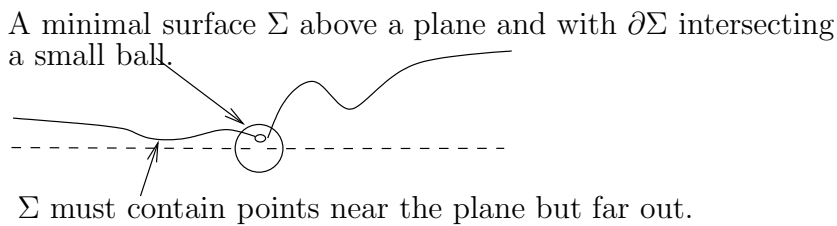}
    \caption{The existence of low points near a plane.}
    \label{f:low}
    %\end{minipage}
\end{figure}

Suppose now that case (1) in Lemma \ref{l:ulscpie} holds for (a
subsequence of) the $\Sigma_j$'s.  The existence of the graphs in
the $\tilde{\Sigma}_j^+$'s converging to $\{ x_3 = 0 \} \setminus
\{ 0 \}$ follows immediately from combining three facts:
\begin{itemize}
\item As $j \to \infty$, the minimum distance between $\partial
B_1 \cap \tilde{\Sigma}_j^+$ and $\{ x_3 = 0 \}$ goes to zero.
This  was actually proven in lemma $3.3$ of \cite{CM8} that gave
the existence of low points in a {\underline{connected}} minimal
surface contained on one side of a plane and with interior
boundary close to this plane.{\footnote{The argument for this was
by contradiction. Namely, if there were no low points, then we
would get a contradiction from the strong maximum principle by
first sliding a catenoid up under the surface and then sliding the
catenoid horizontally away, eventually separating two boundary
components of the surface. Here the strong maximum principle is
used to keep the sliding catenoids and the surface disjoint.  See,
for instance, corollary $1.18$ in \cite{CM1} for a precise
statement of the strong maximum principle.}} We will recall this
lemma from \cite{CM8} next:
\begin{Lem} \label{l:annbd1}
Lemma $3.3$ in \cite{CM8}; see Figure \ref{f:low}.  If $0<
\epsilon < 4 r_0 / 5$ and $\Sigma\subset B_{r_0}$ is a connected
immersed minimal surface with $B_{\epsilon} \cap \Sigma \ne
\emptyset$,   $\Sigma \setminus B_{\epsilon} \ne \emptyset$, and
\begin{equation}
    \partial \Sigma
\subset B_{\epsilon} \cup (\partial B_{r_0} \cap \{ x_3 > - 3 r_0
/ 5 \}) \, ,
\end{equation}
then
\begin{equation}    \label{e:lpts}
    \min_{\Sigma \cap \{ x_1^2 + x_2^2 \geq (4r_0/5)^2 \} }
     x_3 \leq
    4 \, \epsilon \, \cosh^{-1} (3 r_0 / \epsilon) <
    4 \,\epsilon \, \log (6 r_0 / \epsilon)
    \, .
\end{equation}
\end{Lem}
 \item The one-sided curvature estimate
and the scale-invariant ULSC property give a scale-invariant
curvature estimate for the $\tilde{\Sigma}_j^+$'s in a narrow cone
about the plane $\{ x_3 = 0 \}$.   Here we have used that the
$\tilde{\Sigma}_j^+$'s stay on one side of the graphs $\Gamma_j$
converging to $\{ x_3 = 0 \} \setminus \{ 0 \}$.  Similarly, this
curvature estimate and the barrier limit plane imply that the
$\tilde{\Sigma}_j^+$'s are locally graphical in a slightly
narrower cone about $\{ x_3 = 0 \}$. \item The first step gives a
sequence of points in the $\tilde{\Sigma}_j^+$'s converging to a
point in $\partial B_1 \cap \{ x_3 = 0 \}$.  The second step
allows us to apply the Harnack inequality to build this out into
expanding, locally graphical, subsets of the
$\tilde{\Sigma}_j^+$'s that are converging to the plane.

These locally graphical regions piece together to give graphs over
expanding annuli; the other possibility would be to form a
multi-valued graph, but this is impossible since such a
multi-valued graph would be forced to spiral infinitely (since it
cannot cross itself and also cannot cross the stable graph
$\Gamma_j$).
\end{itemize}

Finally, we will briefly describe the modifications needed for
case (2) in Lemma \ref{l:ulscpie} when the  $\Sigma_j^{ulsc}$'s
have two interior boundary components.  The complication arises in
the second step. Namely, we can no longer locally extend the graph
over $\{ x_3 = 0 \} \setminus \{ 0 \}$; this is because
$\tilde{\Sigma}_j^+$ is not
 scale-invariant ULSC  in the second ball $B_{s_j}(y_j)$.
To deal with this, we will consider several different cases.

 The
two simplest possibilities are when the points $y_j$ go to either
zero or infinity.  When $y_j \to 0$, then we can replace the radii
$r_j$ by another sequence $r_j' > \max \{ r_j , |y_j| \}$ where
$r_j' \to 0$; with the new choice of $r_j'$, the new
$\tilde{\Sigma}_j^+$'s are ULSC and we can proceed as in case (1).
On the other hand, when $|y_j| \to \infty$, we can replace the
outer radii $R_j$ by $|y_j|$ and the new sequence of
$\tilde{\Sigma}_j^+$'s  will again be scale-invariant ULSC.

Suppose therefore that the points $y_j$ converge to a finite point
$y \ne 0$.  We will consider two separate subcases here (we can
reduce to these after taking subsequences):
\begin{itemize}
\item
Suppose first that $s_j$ goes to $0$.  In this case,    the
one-sided curvature estimate gives estimates for the
$\Sigma_j^{ulsc}$'s as long as we stay away from the points $0$
and $y$.  We can then argue as in (1) to get the get the desired
graphs -- these graphs converge to $\{ x_3 = 0 \} \setminus \{ 0,
\, y \}$.
\item Suppose now that $\liminf s_j = s_{\infty} > 0$.  In this
case, the sequence is ULSC away from $0$ but  not scale-invariant
ULSC (i.e., the injectivity radius stays away from zero, but it
does not necessarily grow as we go away from $0$). To make this
precise, we will need an additional property of the balls
$B_{s_j}(y_j)$ that was not recorded in Lemma \ref{l:ulscpie} but
follows easily from its proof.{\footnote{See ``{\underline{Blowing
up   on the smallest scale of non-trivial topology}}.'' in the
proof of Proposition \ref{l:getsetgo}.}} Namely, we can assume
that $s_j$ is the ``smallest scale of non-trivial topology.'' More
precisely, we can assume that the component of $B_{s_j}(y_j) \cap
\Sigma_j$ containing the second interior boundary curve has
injectivity radius at least $\beta \, s_j$ for some fixed constant
$\beta > 0$.  In particular, since $\liminf s_j  > 0$, the
one-sided curvature estimate gives uniform estimates on these
components of $B_{s_j}(y_j) \cap \Sigma_j$.  We can now argue as
in (1) to get the get the desired graphs; this time the graphs
converge to $\{ x_3 = 0 \} \setminus \{ 0 \}$.
\end{itemize}

This completes the proof of (C1) in Theorem \ref{t:t5.2}.

\section{The ULSC regions of the lamination: (C2) and (D) in Theorem \ref{t:t5.2}}
\label{s:ctwod}

In this section, we will prove that the ULSC regions of the
lamination have the same structure as in the globally ULSC case of
Theorem \ref{t:t5.1}.  Namely, we will prove that:
\begin{itemize}
\item The
leaves intersecting the ULSC part of the singular set $\cSu$ are
parallel planes.  Each plane intersects $\cSu$ at two points.
\item
$\cSu$ is a union of Lipschitz curves transverse to the  leaves.
The leaves intersecting $\cSu$ foliate an open subset of $\RR^3$
that does not intersect $\cSt$.
\end{itemize}

The key for the proof of these two properties will be to show that
each collapsed leaf has a neighborhood that is ULSC; this will be
done in Proposition \ref{l:ptw}.  Recall that a leaf $\Gamma$ of
$\cL'$ is said to be collapsed if its closure $\barga$ contains a
point in $\cSu$ and this point is a removable singularity for
$\Gamma$; see Definition \ref{d:coll}.  We have already
established a great deal of structure for collapsed leaves in
Proposition \ref{p:cole0} and much of this will be used below.

\vskip2mm  Here, and elsewhere in this section, the closure
$\barga$ of a leaf $\Gamma$ is defined to be the union of the
closures of all bounded geodesic balls in $\Gamma$; that is, we
fix a point $x_{\Gamma} \in \Gamma$ and set
\begin{equation}    \label{e:closure}
    \barga = \bigcup_{r} \overline{ \cB_r (x_{\Gamma}) } \,
    ,
\end{equation}
where $\overline{ \cB_r (x_{\Gamma}) }$ is the closure of $\cB_r
(x_{\Gamma})$ as a subset of $\RR^3$.  Eventually we will show
that $\barga$ is a flat plane and hence, in particular, $\barga =
\bar{\Gamma}$.  However, a priori   $\Gamma$  may not be proper,
and thus the two notions could a priori differ.

The main result of this section is Proposition \ref{l:ptw} below
showing that $\barga$ does not intersect $\cSt$.    Since $\cSt$
is a closed subset of $\RR^3$, it follows that every compact
subset of $\barga$ has a neighborhood in $\RR^3$ that does not
intersect $\cSt$.

\begin{Pro}     \label{l:ptw}
If $\Gamma$ is a  collapsed leaf   of $\cL'$, then $\barga \cap
\cSt = \emptyset$.
\end{Pro}

  The rough idea of the proof is to
first   show that $\barga \setminus \Gamma$ consists of exactly
two points (this is analogous to each leaf having at most two
singular points in the ULSC case); see Corollary
\ref{c:missingtwo} below. Consequently, the union of $\Gamma$ and
the given point in $\barga \cap \cSu$ will give a stable surface
that is ``complete away from a point''  and, hence flat by Lemma
\ref{l:punstab}. Finally, once we know that $\Gamma$ is flat, it
will be easy to check that $\barga \cap \cSt = \emptyset$.

Before we can get into the proof just outlined, we will need to
recall a little of the structure  that has already been proven. We
will do this in the next two subsections.  The next subsection
establishes a key property of the    stable limit planes that we
get through each point of $\cSt$.   The second subsection below
reviews the properties of a general collapsed leaf of $\cL'$.

\subsection{The leaf $\Gamma$ cannot cross the limit planes}

The structure result ($C_{neck}$) from Theorem \ref{t:t5.2a}
  gives graphs
$\tilde \Sigma_j^+$ and $\tilde \Sigma_j^-$ in $\Sigma_j$ that
converge to a plane $P_z$ through each point $z \in \cSt$; see
Figure \ref{f:gexc}.  The graphs $\tilde \Sigma_j^+$ converge
smoothly away from $z$ and possibly one other point (call it
$z^+$); this second point must also be in $\cSt$.  Similarly, the
$\tilde \Sigma_j^+$ converge away from $z$ and possibly a point
$z^- \in \cSt$.  Furthermore, $\tilde \Sigma_j^+$ and $\tilde
\Sigma_j^-$ are separated in $\Sigma_j$ by the curve
$\sigma_j$.{\footnote{Property ($C_{neck}$) from Theorem
\ref{t:t5.2a} holds at $z$ by (C1) in Theorem \ref{t:t5.2}; this
was proven in Section \ref{s:ceno}. The last ``separation'' claim
is not explicit in ($C_{neck}$) but follows immediately from
Proposition \ref{l:getsetgo}; using the notation from that
proposition, we have that $\tilde \Sigma_j^+ \subset \Sigma_j^+$
and $\tilde \Sigma_j^- \subset \Sigma_j^-$.}} One expects that the
limit plane $P_z$ should be the closure of a leaf of $\cL'$, but
this is not a priori clear; for instance, $\cS$ might even be
dense in $P_z$.

\begin{figure}[htbp]
 \setlength{\captionindent}{20pt}
    \centering\input{gen16c.pstex_t}
    \caption{The structure near $z$: The two graphs
    $\tilde \Sigma_j^+$ and $\tilde \Sigma_j^-$ are separated by curves
$\sigma_j$ shrinking to $z$.  The multi-valued graph $\Sigma^g_j$
comes near $z$.}
    \label{f:gexc}
\end{figure}

Using this structure,  the next lemma proves that the leaves of
$\cL'$ do not cross any of these limit planes.  This is almost
obvious since the stable surfaces converging to the limit plane
are disjoint from the $\Sigma_j$'s that are converging to the
leaves of $\cL'$. The only possible difficulty comes from that the
convergence is only away from the singular set $\cS$, but this
will be easy to handle.  The lemma applies to an arbitrary leaf
$\Gamma$ of $\cL'$, i.e., we do not need $\Gamma$ to be collapsed.

\begin{Lem} \label{l:sm1}
Suppose that  $\Gamma$ is an arbitrary leaf of $\cL'$, collapsed
or not.   If  $z$ is any point in $\cSt$
 and $P_z$ is the
corresponding limit plane through $z$, then $\Gamma$  does not
cross $P_z$.
\end{Lem}

\begin{proof}
Fix an
  open connected set $K \subset \Gamma$ with compact closure in $\Gamma$ and
 recall that
 the $\Sigma_j$'s contain:
 \begin{itemize}
 \item Graphs $\tilde
\Sigma_j^+$ and $\tilde \Sigma_j^-$ that both converge to $P_z$
away from a finite set  $\cP$ of points; see Figure \ref{f:gexc}.
Moreover, $\tilde \Sigma_j^+$ and $\tilde \Sigma_j^-$ are
separated in $\Sigma_j$ by curves $\sigma_j$ shrinking to $z$.
\item Connected subsets $\Sigma^g_j \subset \Sigma_j$ given by
 Lemma \ref{l:oldclaim} that are locally graphical over $K$ and
 that converge with
multiplicity to $K$.  These locally graphical subsets might
globally be graphs or multi-valued graphs over $K$.
\end{itemize}

We will show first that the $\Sigma^g_j$'s cannot intersect both
$\tilde \Sigma_j^+$ and $\tilde \Sigma_j^-$. First, using that $z$
is not in $\Gamma$ and $K \subset \Gamma$ has compact closure, we
can fix a ball $B_{s} (z)$ so that
\begin{equation}
    B_s (z) \cap \Sigma^g_j = \emptyset
\end{equation}
for all $j$ sufficiently large.  On the other hand, the curves
$\sigma_j$ separating $\tilde \Sigma_j^+$ and $\tilde \Sigma_j^-$
are shrinking to $z$.  Therefore,  the curves $\sigma_j$ don't
intersect $\Sigma^g_j$ when $j$ is large and, hence, the connected
set $\Sigma^g_j$ cannot intersect both $\tilde \Sigma_j^+$ and
$\tilde \Sigma_j^-$.  Without loss of generality, we can assume
that
\begin{equation}    \label{e:saythat}
    \Sigma^g_j \cap \tilde \Sigma_j^+ = \emptyset \, .
\end{equation}

We will next use \eqr{e:saythat} to show that the two smooth open
surfaces $P_z$ and $K$ do not have any points of transverse
intersection. Namely, if $P_z$ and $K$ have points of transverse
intersection, then, since the singular set $\cP$ for the
convergence to $P_z$ is finite, $P_z\setminus \cP$ and $K$ would
also have points of transverse intersection. However, this
contradicts \eqr{e:saythat} since $\tilde \Sigma_j^+ \to P_z$
smoothly away from $\cP$ and $\Sigma_j^g\to K$.

Finally, recall that if a connected minimal surface intersects
both sides of a plane, then  the surface and plane must have a
point of transverse intersection; this follows from the local
structure of the nodal set of a harmonic function, see, e.g.,
lemma $4.28$ in \cite{CM1}. Therefore, since $K$    is connected
and does not intersect $P_z$ transversely at any point,  we see
that $K$ must be on one    side of $P_z$.  Since this holds for
every such $K$ and these    exhaust $\Gamma$ by Lemma
\ref{l:simpap} in Appendix \ref{s:a}, we see that $\Gamma$ also
lies on one
    side of $P_z$.
\end{proof}

\vskip2mm {\it{Throughout the rest of this section, $\Gamma$ will
be a collapsed leaf of $\cL'$.}}

\subsection{The properties of a collapsed leaf $\Gamma$}

Before getting  into the   proof, it may be useful to recall the
properties of the leaf $\Gamma$. Eventually, we will use these
properties to show that $\barga$ is a plane.
\begin{itemize}
\item
    $\Gamma$ is by definition an injective immersion of a
{\underline{connected}} surface without boundary, but
{\underline{not}} necessarily complete. Furthermore, the immersion
  is not necessarily proper.
\item
Since $\Gamma$ is a leaf of $\cL'$, it follows that $\Gamma$ does
not intersect $\cS$ -- and hence, since $\cS$ is closed (as a
subset of $\RR^3$), each point in $\Gamma$ has a neighborhood
where the curvatures of the $\Sigma_j$'s are uniformly bounded.
\item The following   local structure of $\Gamma$ near a point of
$\barga \cap \cSu$  was established in (1) in Proposition
\ref{p:cole0}:
\begin{enumerate}
  \item[(Loc)]
  Given {\underline{any}} $y \in \barga \cap
\cSu$, there exists $r_0 > 0$ so that the closure (in $\RR^3$) of
each component of $B_{r_0}(y) \cap \Gamma$ is a compact embedded
 disk  with boundary in $\partial B_{r_0}(y)$.

Furthermore, $B_{r_0}(y) \cap \Gamma$ must contain the component
$\Gamma_y$ given by Lemma \ref{l:leaf} and $\Gamma_y$ is the only
component of $B_{r_0}(y) \cap \Gamma$ with $y$  in its closure.
\end{enumerate}
\item  $\Gamma$ (or its oriented double cover) must be stable by (2)
in Proposition \ref{p:cole0}. \item
 $\barga$ intersects $\cSu$ in
at most two points by (3) in Proposition \ref{p:cole0}.
\end{itemize}

\begin{figure}[htbp]
 \setlength{\captionindent}{20pt}
    \centering\input{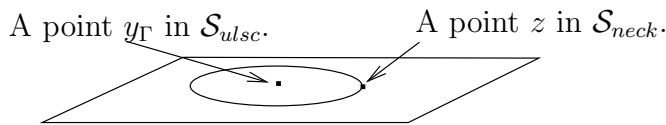}
    \caption{A priori $\Gamma$ could be a punctured disk in a plane.}
    \label{f:gex}
\end{figure}

These  properties will be essential for proving Proposition
\ref{l:ptw}.  The main difficulty will be that $\Gamma$ is not
complete. This occurs where   $\barga$ intersects
  $\cS$; see Figure
\ref{f:gex} for such an example.  By (Loc), the points in $\barga
\cap \cSu$ are isolated removable singularities of $\Gamma$ and
thus are easily dealt with.  Consequently, the first step will be
to control the number of points   of $\barga \cap \cSt$   when
$\Gamma$ is collapsed.  This will be done in the next subsection.

\subsection{$\barga \setminus \Gamma$ consists of exactly two
points}

The next corollary is the first step needed for the proof of
Proposition \ref{l:ptw} that gives $\barga \cap \cSt = \emptyset$.
This corollary  shows that if $\barga \cap \cSt \ne \emptyset$,
then $\barga \setminus \Gamma$ consists of exactly two points.

\begin{Cor}     \label{c:missingtwo}
Let $\Gamma$ be a collapsed leaf of $\cL'$.  If $\barga \cap \cSt
\ne \emptyset$, then $\barga \setminus \Gamma$ consists of exactly
two points with one each in $\cSu$ and $\cSt$.
\end{Cor}

 This will be an easy corollary of Lemma \ref{l:neckulsc} below
 that shows
that the sheets of the multi-valued graphs over $\Gamma$ connect
in a small neighborhood of any singular point.  Previously, we
used the one-sided curvature estimate to establish a similar
connecting property near a point of $\cSu$.

Before making this connecting property precise, we need to set up
some notation.
 Recall  that if $K$ is a  ``sufficiently large'' open connected
subdomain of $\Gamma$ with compact closure in $\Gamma$, then
Corollary \ref{c:leafstuff1} gives a sequence of multi-valued
graphs $\Sigma_{j}^g  \subset \Sigma_j$ that converges to $K$ with
infinite multiplicity. Since   $\Gamma$ can be exhausted by a
nested sequence $K_j$ of such $K$'s  by Lemma \ref{l:simpap},
    we can assume that the following holds  (after passing to a
 subsequence):
 \begin{enumerate}
 \item[(Graph)] $\Sigma_j$
 contains a $j$-valued graph $\Sigma^g_j$ over $K_j$ of a function whose values are bounded by
$1/j$ and whose gradient is bounded by $1/j$.  Here $K_j \subset
\Gamma$ is a nested sequence of connected open sets with compact
closure in $\Gamma$ with $\Gamma = \cup_{j} K_j$.
 \end{enumerate}
We actually know a good deal more about these multi-valued graphs,
but this additional structure will not be needed until the proof
of Corollary \ref{c:missingtwo}.

\begin{figure}[htbp]
 \setlength{\captionindent}{20pt}
    \centering\input{gen82.pstex_t}
    \caption{Lemma \ref{l:neckulsc}:
    The sheets of the multi-valued graph $\Sigma_j^g$ must connect near $z
\in \barga \cap \cSt$.}
    \label{f:neckulsc}
\end{figure}

 \begin{Lem}    \label{l:neckulsc}
The sheets of the multi-valued graph $\Sigma_j^g$ connect near $z
\in \barga \cap \cSt$; see Figure \ref{f:neckulsc}. Precisely,
given any $r>0$, there exist $\delta
> 0$ and $J$ so that if $x \in B_{\delta}(z) \cap \Gamma$, $j >
J$, and
\begin{equation}
    z_j^+ {\text{ and }} z_j^- {\text{ are points in the multi-valued graph }}
    \Sigma_j^g {\text{ over }} x \, ,
\end{equation}
then $z_j^+$ and $z_j^-$ are in the same connected component of
$B_r(z) \cap \Sigma_j$.
\end{Lem}

\begin{proof}
 We will argue by contradiction, so suppose
there exists   some $r> 0$ so that  for every $\delta > 0$ there
exists $x \in B_{\delta} (z) \cap \Gamma$ and infinitely $j$'s so
that $B_r(z) \cap \Sigma_j$ has (at least) two distinct components
that both contain points in $\Sigma^g_j$ over $x$. After passing
to a subsequence, we can assume that there is a sequence of points
$z_j^+$ and $z_j^-$ in $\Sigma^g_j$ with
\begin{equation}
    z_j^+ {\text{ and }} z_j^- {\text{ converging to }} z \, ,
\end{equation}
so that $z_j^+$ and $z_j^-$ are in distinct components of $B_r(z)
\cap \Sigma_j$.

We will use the $\Sigma_j$'s as barriers for a Plateau problem to
construct stable surfaces $\tilde \Gamma_j$ between these distinct
components. Before constructing the stable surfaces $\tilde
\Gamma_j$, recall the following useful consequence of the interior
curvature estimates for stable surfaces of \cite{Sc1}, \cite{CM2}
(see, e.g., lemma $2.2$ in \cite{CM1}):
\begin{enumerate}
\item[(Stab)] There exists a positive constant $\alpha < 1$ so that
if $\Gamma_s$ is a stable embedded minimal surface with $\partial
\Gamma_s \subset
\partial B_R$, then each component of $B_{\alpha \, R} \cap
\Gamma_s$ is a graph  over some plane with  gradient bounded by
one.
\end{enumerate}
Set
\begin{equation}  \label{e:rprime}
    r' = \alpha \, r \, .
\end{equation}

\vskip2mm \noindent {\bf{The properties of the stable graphs
$\tilde \Gamma_j$}}:  We will below find (stable) graphs $\tilde
\Gamma_j$ between these distinct components so that the following
three properties hold:
\begin{align}   \label{e:neep1}
    \tilde \Gamma_j  \subset  B_{r'}(z) \setminus \Sigma_j &{\text{ with }} \partial \tilde \Gamma_j \subset \partial
    B_{r'}(z) \, . \\
    B_{\epsilon_j}(z) \cap \tilde \Gamma_j \ne \emptyset &{\text{ where }} \epsilon_j \to 0  \, .
    \label{e:neep} \\
    {\text{The multi-valued graph }} \Sigma^g_j {\text{ in }} &\Sigma_j
{\text{ intersects both sides of }} \tilde \Gamma_j  \, .
    \label{e:neep2}
\end{align}
Recall that a properly embedded (connected) surface in $\RR^3$
will automatically have two sides. Properties \eqr{e:neep} and
\eqr{e:neep2} will follow from a standard linking argument.

\begin{figure}[htbp]
    \setlength{\captionindent}{20pt}
    %\begin{minipage}[t]{0.5\textwidth}
    \centering\input{gen17.pstex_t}
    \caption{We argue by contradiction to show  that the sheets connect near $z$. Assuming that
    they don't, we first construct  stable surfaces $\tilde \Gamma_j$.}
    \label{f:2a}
    %\end{minipage}
\end{figure}

\vskip2mm \noindent {\bf{Constructing the stable graphs $\tilde
\Gamma_j$}}: We will construct $\tilde \Gamma_j$ in two steps,
first finding stable surfaces $\tilde \Gamma_j^0$ in the larger
ball $B_{r}(z)$ and then letting $\tilde \Gamma_j$ be an
appropriate component of
  $B_{r'}(z) \cap \tilde \Gamma_j^0$ where $r'$ is given by
  \eqr{e:rprime}.

 To construct $\tilde \Gamma_j^0$, first   choose two distinct
components $\Sigma_{r,j}^+$ and $\Sigma_{r,j}^-$ of $B_r(z) \cap
\Sigma_j$ that both contain points in $\Sigma^g_j$ over $z_j$;
these exist by assumption.  Let $\ell_j$ be a line segment
connecting the two points over $x_j$ in $\Sigma_{r,j}^+$ and
$\Sigma_{r,j}^-$;  see Figure \ref{f:2a}. Fix a component
$\ell_j^0$ of $\ell_j \setminus (\Sigma_{r,j}^+ \cup
\Sigma_{r,j}^-)$ that also connects $\Sigma_{r,j}^+$ and
$\Sigma_{r,j}^-$ but
  intersects
$\Sigma_{r,j}^+$ and $\Sigma_{r,j}^-$ only at its endpoints
$\partial \ell_j^0$.  The surface  $\Sigma_{r,j}^+$ sits in the
boundary of two components of $B_r(z) \setminus \Sigma_j$ -- one
component on each side of $\Sigma_{r,j}^+$.{\footnote{In the
simplest case where $\Sigma_{r,j}^+$ and $\Sigma_{r,j}^-$ are the
only components of $B_r(z) \cap \Sigma_j$, we would choose
$\Omega_j$ to be the component of $B_r(z) \setminus \Sigma_j$
between them, i.e., the component containing the interior of
$\ell_j^0$.  In general, there are other components of $B_r(z)
\cap \Sigma_j$ intersecting $\ell_j^0$ so we cannot do this. This
slightly complicates the choice of $\Omega_j$.}} Let
 $\Omega_j$ be the component  that $\ell_j^0$
 points into as it leaves $\Sigma_{r,j}^+$;
 i.e., let $\Omega_j$ be the component of $B_r(z) \setminus
 \Sigma_j$ with
\begin{equation}
 \Sigma_{r,j}^+ \subset \partial \Omega_j {\text{ and }} \Omega_j
 \cap (\ell_j^0 \setminus \partial \ell_j^0)  \ne \emptyset  \, .
\end{equation}
The point of choosing $\Omega_j$ in this way is that the curve
$\ell_j^0$ has linking number one with $\partial \Sigma_{r,j}^+$
in $\Omega_j$.
 The domain $\Omega_j$ is mean
convex and, hence, \cite{MeYa1}--\cite{MeYa2} gives a stable
embedded minimal surface $\tilde \Gamma_j^0 \subset \Omega_j$ with
$\partial \tilde \Gamma_j^0 = \partial \Sigma_{r,j}^+$.

Since $\partial \tilde \Gamma_j^0$ has linking number one with
  $\ell_j^0$ in $\Omega_j$, the endpoints of $\ell_j^0$
are separated in $B_r(z)$ by $\tilde \Gamma_j^0$.  However, the
endpoints of $\ell_j^0$ connect to the endpoints of $\ell_j$ by
curves in $\Sigma_{r,j}^+$ and $\Sigma_{r,j}^-$; these curves do
not cross $\tilde \Gamma_j^0$ and, consequently,
 the endpoints of $\ell_j$ are also separated in
$B_r(z)$ by $\tilde \Gamma_j^0$.  We can therefore choose
  a component $\tilde \Gamma_j$ of $B_{r'}(z) \cap
\tilde \Gamma_j$ that separates the endpoints of $\ell_j$; so
\begin{equation}    \label{e:link2}
    \tilde \Gamma_j \cap \ell_j \ne \emptyset \, .
\end{equation}
Since $z_j^+$ and $z_j^-$  go to $z$, it follows that each
$\ell_j$ is contained in a ball $B_{\epsilon_j}(z)$ where
$\epsilon_j \to 0$. In particular, \eqr{e:link2} gives
\eqr{e:neep}. Since the endpoints of $\ell_j$ are both in the
multi-valued graph, we also get \eqr{e:neep2}.

\vskip2mm \noindent {\bf{Using the stable graphs $\tilde \Gamma_j$
to show that $ \Gamma  \subset P_z$}}: Now that we have
constructed the $\tilde \Gamma_j$'s, we are ready to return to the
proof of the lemma. The first step will be to show that a
subsequence of the $\tilde \Gamma_j$'s converges to a subset of
$P_z$.   First, by (Stab), the surface $\tilde \Gamma_j$ is a
graph with gradient bounded by one. After passing to a
subsequence, we can therefore assume that $\tilde \Gamma_j$
converges to a minimal graph $\tilde \Gamma$. Since $\epsilon_j
\to 0$,  $\tilde \Gamma$ contains $z$. On the other hand, $\tilde
\Sigma_j^+ \subset \Sigma_j$ does not intersect the graph $\tilde
\Gamma_j$; since $\tilde \Sigma_j^+$ converges to $P_z$ away from
a finite set, we conclude that $\tilde \Gamma$ must be on one side
of $P_z$. Since $\tilde \Gamma$ is on one side of $P_z$ and
intersects it at $z$, the strong maximum principle implies that
$\tilde \Gamma$ is contained in $P_z$, as desired.

We will next show that
\begin{equation} \label{e:gipz}
    \Gamma  \subset P_z \, .
\end{equation}
  First, by
\eqr{e:neep2}, there is at least one  sheet of the multi-valued
graph $\Sigma^g_j$ on each side of the graph $\tilde \Gamma_j$ for
every $j$. Since both of these sheets  converge to $\Gamma$ by
(Graph), we can fix a point $\tilde y$ in $B_{r'} \cap \Gamma$
that is both a limit of points $\tilde y_j^+$ above $\tilde
\Gamma_j$ and a limit of points $\tilde y_j^-$ below $\tilde
\Gamma_j$.  For each $j$, the line segment connecting $\tilde
y_j^+$  to $\tilde y_j^-$ must intersect $\tilde \Gamma_j$ at a
point $\tilde y_j$. The sequence of points $\tilde y_j \in \tilde
\Gamma_j$ must also converge to the common limit $\tilde y$ of
$\tilde y_j^+$  and $\tilde y_j^-$. In particular, we conclude
that $\tilde y \in P_z \cap \Gamma$, so that $\Gamma \subset P_z$
by the strong maximum principle.

\vskip2mm \noindent {\bf{The contradiction: We cannot have $
\Gamma \subset P_z$}}: To complete the proof of the lemma, we
explain next how \eqr{e:gipz} leads to a contradiction. Since we
will need the same argument later, it will be useful to isolate it
out as a claim: \vskip2mm \noindent {\bf{Claim}}: \eqr{e:gipz}
cannot hold, i.e., we cannot have $\Gamma  \subset P_z$.

\vskip2mm \noindent {\bf{Proof of Claim}}: We will argue by
contradiction, so suppose that $\Gamma \subset P_z$.  Since
$\barga$ is the closure of $\Gamma$, it follows that $\barga
\subset P_z$ as well.  Since $\Gamma$ is collapsed, $\barga$
contains a point
\begin{equation}
    y_{\Gamma} \subset \cSu \cap
\barga \subset P_z \, .
\end{equation}
Both of the graphs $\tilde \Sigma_j^+$ and $\tilde \Sigma_j^-$ are
converging to
 $P_z$ away from punctures in $\cSt$, so we get
sequences of points $y_j^+ \in \tilde \Sigma_j^+$ and $y_j^- \in
\tilde \Sigma_j^-$
 with
 \begin{equation}   \label{e:yjpm}
 y_j^+ \to y_{\Gamma} {\text{ and }} y_j^- \to y_{\Gamma} \, .
 \end{equation}

We will next use the one-sided curvature estimate to prove that
$y_j^+$ and $y_j^-$ can be connected in $\Sigma_j$ in any small
neighborhood  of $y_{\Gamma}$ as $j\to \infty$.
 To see this,  note first that
$y_{\Gamma}$ is in $\cSu$ and hence each component of
$B_{r}(y_{\Gamma}) \cap \Sigma_j$ is a disk for some $r > 0$;
after possibly choosing $r$ smaller, we can assume that
$|y_{\Gamma} - z| > r$.    If there were at least two of these
disks in  $B_{r}(y_{\Gamma}) \cap \Sigma_j$ intersecting the
concentric sub--ball $B_{C'' \, r}(y_{\Gamma})$ where $C'' > 0$ is
a sufficiently small constant, then the one-sided curvature
estimate would give a uniform curvature bound for each component
of $\Sigma_j$ in this sub--ball.{\footnote{We actually use a
corollary of the one-sided curvature estimate recorded in
corollary $0.4$ in \cite{CM6}. This corollary states that if there
are two disjoint surfaces in a ball in $\RR^3$, both intersect a
sufficiently small ball around the center, and one is a disk, then
we get an interior curvature estimate for the disk-type
component.\label{fn:repeat1}}} Since no such curvature bound holds
near a point of $\cS$ by definition, we conclude that   the points
$y_j^+$ and $y_j^-$ must be
  in the same connected component of
 $B_{r}(y_{\Gamma}) \cap \Sigma_j$ for  all $j$
 sufficiently large.

 This local connectedness near $y_{\Gamma}$ will easily lead to a contradiction.
This is because $y_j^+$ and
 $y_j^-$ were proven to be separated in $\Sigma_j$ by the curve
 $\sigma_j$ and the $\sigma_j$'s are shrinking to the point $z \ne
 y_{\Gamma}$.  This contradiction completes the proof of the Claim and hence of
 the lemma.
\end{proof}

\begin{figure}[htbp]
    \setlength{\captionindent}{20pt}
    %\begin{minipage}[t]{0.5\textwidth}
    \centering\input{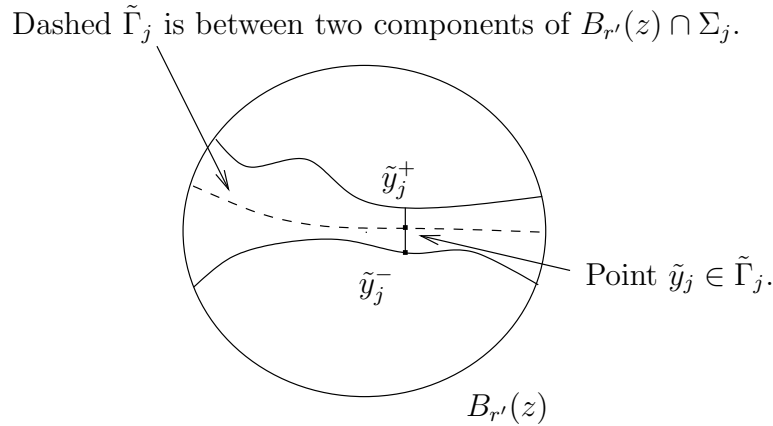}
    \caption{The points $\tilde y_j^+$ and $\tilde y_j^-$ converge to a point $\tilde y \in \Gamma$ from
    opposite sides of the graph $\tilde \Gamma_j$.  Thus, the points $\tilde y_j \in \tilde \Gamma_j$
    between them also converge to $\tilde y \in \Gamma$, giving Lemma \ref{l:neckulsc}.}
    \label{f:2a2}
    %\end{minipage}
\end{figure}

\begin{proof}
(of Corollary \ref{c:missingtwo}.) By assumption, $\barga \cap
\cS$ contains at least one point in each of $\barga \cap \cSu$ and
$\barga \cap \cSt$.  We will argue by contradiction to prove that
$\barga \cap \cS$ cannot contain a third point.  For simplicity,
we will assume that the $\Sigma_j$'s are planar domains; the
general finite genus case follows with easy modifications.

The proof follows the argument given in Remark \ref{r:genpoints}
and the basic idea is simple:
\begin{quote}
Suppose that $p$, $q$, and $r$ are distinct points in $\barga \cap
\cS$. The local connecting property near each point of $\cS$
allows us to construct  closed non-contractible curves in the
$\Sigma_j$'s that  converge with multiplicity two to a curve in
$\Gamma$ connecting $p$ and $q$.  These curves must separate in
the planar domain $\Sigma_j$.  However, if we connect  points on
opposite sides of these curves to  the third point $r$, the local
connecting property near $r$ gives a contradiction.
\end{quote}
The only difficulty in carrying out this argument will be that the
surface $\Gamma$ is not complete.

\vskip2mm \noindent {\bf{Step 1: Choosing the singular points and
curves in $\Gamma$.}}  We will first choose the points $p$, $q$,
and $r$ in $\barga \cap \cS$. Let $p$ be the given point in
$\barga \cap \cSu$ and then let $q$ be a closest point in $\barga
\cap \cS$ to $p$ (a priori there may be many possible choices).
Since we are arguing by contradiction, there is a third distinct
point $r \in \barga \cap \cS$.

By our choice of $q$, there must be a minimizing geodesic
$\gamma_{pq}:[0,L] \to \barga$ parameterized by arclength and with
the following properties:
\begin{itemize}
\item $\gamma_{pq} (0) = p$ and $\gamma_{pq} (L) = q$.
\item The
interior of $\gamma_{pq}$ is contained in $\Gamma$.
\end{itemize}
Since the closed geodesic $\gamma_{pq}$ is compact, we must have
\begin{equation}    \label{e:gapqd}
    \dist_{\RR^3} (\gamma_{pq} , r ) > 0 \, .
\end{equation}
Since $p \in \cSu$ and $\Gamma$ is collapsed, property (1) in
Proposition \ref{p:cole0} gives a ball $B_{\delta} (p)$ and a
component $\Gamma_p$ of $B_{\delta} (p) \cap \Gamma$ so that
$\Gamma_p \cup \{ p \}$ is a smooth minimal graph.  Since
$\gamma_{pq}$ is minimizing, it is not hard to see that $\partial
\Gamma_p$ intersects $\gamma_{pq}$ in a single point $p'$.

Fix a constant $\epsilon > 0$ that is much smaller than the
distance from $r$ to $\gamma_{pq}$.  By the definition of
$\barga$, there must be a point $r'$ in $\Gamma$ that is distance
$\epsilon$ from $r$.
 Since $\Gamma$ is connected, we can  choose a compact curve
 $\tilde \gamma_{p'r'}$ that is contained in $\Gamma$ and
 connects $p'$ to $r'$.  The curve $\tilde \gamma_{p'r'}$ may
 intersect $\gamma_{pq}$ many times, so we replace it with the
 component of $\tilde \gamma_{p'r'} \setminus \gamma_{pq}$ with
 $r'$ in its boundary.  This gives a curve in $\Gamma$ from
 $\gamma_{pq}$ to $r'$ and whose interior does not intersect
 $\gamma_{pq}$.  After adding a subsegment of $\gamma_{pq}$ and
 perturbing the resulting curve slightly off of $\gamma_{pq}$, we
 get a compact curve $\gamma_{p'r'} \subset \Gamma$ from $p'$ to
 $r'$ and whose interior does not intersect $\gamma_{pq}$;
 see Figure \ref{f:gampprp}.

The point about the curve $\gamma_{p'r'}$ is that it will give a
way to connect points near $p'$ to $r'$ in $\Gamma \setminus
\gamma_{pq}$; see Figure \ref{f:gampprp}.  This will be especially
useful since the curve $\partial \Gamma_p$ allows us to connect
points near $p'$ that are on the opposite sides of $\gamma_{pq}$.

\begin{figure}[htbp]
    \setlength{\captionindent}{20pt}
    %\begin{minipage}[t]{0.5\textwidth}
    \centering\input{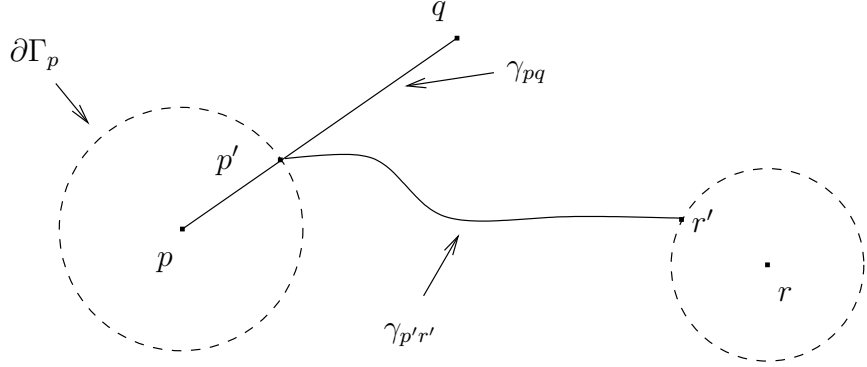}
    \caption{The curves $\gamma_{pq}$ and $\gamma_{p'r'}$.}
    \label{f:gampprp}
    %\end{minipage}
\end{figure}

\vskip2mm \noindent {\bf{Step 2: Choosing the curves in the
$\Sigma_j$'s.}} We can now argue as in Remark \ref{r:genpoints}.
 The key point is that Theorem \ref{t:1} and Lemma \ref{l:neckulsc}
imply that the $\Sigma_j$'s are locally connected in a small
neighborhood of {\underline{any}} of the singular points in
$\barga \cap \cS$. These connecting properties allow us to find
simple closed curves $\gamma^j_{pq} \subset \Sigma_j$ with the
following properties:
\begin{itemize}
\item $\gamma^j_{pq}$ is contained in the $\epsilon$-tubular
neighborhood of $\gamma_{pq}$. \item $\gamma^j_{pq} \setminus
\left( B_{\epsilon}(p) \cup B_{\epsilon}(q) \right)$
 consists of two graphs over $\gamma_{pq}$ which
are in distinct sheets of $\Sigma_j$.
\end{itemize}

\vskip2mm \noindent {\bf{Step 3: The contradiction.}} Since
$\Sigma_j$ has genus zero, the curve $\gamma^j_{pq}$ must separate
$\Sigma_j$ into two distinct components.    However, it is easy to
see that this is impossible by using  the local connecting
property near the third point $r$. Namely, we can take two points
near $p$ on opposite sides of $\gamma^j_{pq}$ and connect each of
them to $r'$ by curves in $\Sigma_j$ which do not intersect
$\gamma^j_{pq}$.   One of these connecting curves will be a graph
over $\gamma_{p'r'}$ while the other is a graph over $\partial
\Gamma_p \cup \gamma_{p'r'}$. These two connecting curves can then
be connected to each other in $B_{C \, \epsilon}(r) \cap
\Sigma_j$, giving the desired contradiction.
\end{proof}

\subsection{The proof of Proposition \ref{l:ptw}}

We can now use Corollary \ref{c:missingtwo} and the properties of
a collapsed leaf to prove Proposition \ref{l:ptw}.  Recall that
this proposition claims that $\barga \cap \cSt = \emptyset$
whenever $\Gamma$ is a collapsed leaf of $\cL'$.

\begin{proof}
(of Proposition \ref{l:ptw}.) We will argue by contradiction, so
suppose that $\Gamma$ is a collapsed leaf of $\cL'$ and $\barga
\cap \cSt \ne \emptyset$.   By Corollary \ref{c:missingtwo}, we
know that $\barga \setminus \Gamma$ consists of one point $y$ in
$\cSu$ and one point $z$ in $\cSt$. Furthermore, property (1) in
Proposition \ref{p:cole0} implies that the point $y$ is a
removable singularity for $\Gamma$ so that $\Gamma \cup \{ y \}$
is smooth and complete away from the point $z$.

We will show next that $\Gamma$ is flat.  The starting point for
this is that $\Gamma$ (or its oriented double cover) is stable by
property (2) in Proposition \ref{p:cole0}.  A standard
``logarithmic cutoff function'' argument then implies that $\Gamma
\cup \{ y \}$ (or its oriented double cover) is also stable; we
leave the simple argument to the reader.  If $\Gamma \cup \{ y \}$
had been complete, then the Bernstein theorem for stable surfaces
would have implied that it was flat.  However, even in this case
where $\Gamma \cup \{ y \}$ is complete away from the single point
$z$, Lemma \ref{l:punstab} in Appendix \ref{s:apB} implies that
$\Gamma$ is flat.

We have now established that $\Gamma$ is a plane with two points
removed and one of these points (namely $z$) is in $\cSt$.  Since
$\Gamma$ cannot cross the limit plane $P_z$ through $z$ by Lemma
\ref{l:sm1}, it follows that
\begin{equation}    \label{e:gamsubpz}
    \Gamma \subset P_z \, .
\end{equation}
 However, we already saw in the Claim at the end of the proof of Lemma
\ref{l:neckulsc} that \eqr{e:gamsubpz} is impossible. This
contradiction completes the proof of the proposition.
\end{proof}

\subsection{The proof of (C2) and (D) in Theorem \ref{t:t5.2}}

We can now argue as in the ULSC case of Part \ref{p:prove1} to
prove (C2) and (D) in Theorem \ref{t:t5.2}.  As we saw in Part
\ref{p:prove1}, the argument requires that we establish the
following four properties of an arbitrary collapsed leaf $\Gamma$
of $\cL'$:
\begin{enumerate}
\item[(0)] $\barga \cap \cSt = \emptyset$.
 \item[(1)]  $\barga$ is a plane
and $\barga \setminus \Gamma$ contains at most two points.
\item[(2)] $\barga \setminus \Gamma$ contains exactly two points.
\item[($\star$)] If $t\in x_3(\cSu)$ and $\epsilon>0$, then
\begin{equation}
    \cSu \cap \{t<x_3<t+\epsilon\}\ne \emptyset {\text{ and }} \cSu \cap
\{t-\epsilon<x_3<t\}\ne \emptyset \, .
\end{equation}
\end{enumerate}
Once we show that (0), (1), (2), and ($\star$) hold, then (C2) and
(D) in Theorem \ref{t:t5.2} follow exactly as in Part
\ref{p:prove1}; we will not repeat the argument here.

It suffices therefore to check that (0), (1), (2), and ($\star$)
hold in this setting.  The first two are quite easy:  (0) is
exactly the conclusion of Proposition \ref{l:ptw} and (1) follows
from (0) together with (1) in Proposition \ref{p:cole}.

We will prove (2) by contradiction, so suppose that   $\barga \cap
\cSu = \{ 0 \}$.  It follows immediately that the $\Sigma_j$'s
contain multi-valued graphs over (subsets of) the punctured plane
$\Gamma = \{ x_3 = 0 \} \setminus \{ 0 \}$ that converge to
$\Gamma$ with infinite multiplicity. Moreover,   (D) in the proof
of property (2) in Proposition \ref{p:cole} gives the following
{\underline{scale invariant}} ULSC property:
\begin{enumerate}
\item[(D)] There exists $\tau > 0$ so that, for $z \in \{ x_3 = 0
\}$ and $j$ large, each component of $B_{\tau \, |z|} (z) \cap
\Sigma_j$ that connects to the multi-valued graph in $\Sigma_j$ is
a disk.
\end{enumerate}
Note that the proof of (D) did not use that the sequence was ULSC.
It follows from (D) and the one-sided curvature estimate that the
multi-valued graphs in $\Sigma_j$ converging to $\Gamma$ spiral
through an entire cone about $\Gamma$. Note that  the $\Sigma_j$'s
are assumed to be uniformly non-simply connected.  Therefore,
Proposition \ref{l:getsetgo} gives stable graphs $\Gamma_j$ that
are disjoint from   $\Sigma_j$.  The $\Gamma_j$'s are graphs with
bounded gradient that start out in a fixed ball and are defined
over annuli with a fixed inner radius and with outer radii going
to infinity.  Standard results for exterior graphs then imply that
the $\Gamma_j$'s grow sublinearly.  Consequently, the $\Gamma_j$'s
are eventually contained in the narrow cone that the multi-valued
graphs in the $\Sigma_j$'s spiral through.  However, this is
impossible since the two are disjoint.  This contradiction
completes the sketch of the proof of (2); we leave the details to
the reader.

\begin{Rem}
 The argument that we   gave here is actually simpler
than in the ULSC case; cf.  (2) in Proposition \ref{p:cole}.
  However, we could not yet use this
argument for the ULSC case since Proposition \ref{l:getsetgo}
relies on the ULSC case.
\end{Rem}

Finally, ($\star$) follows from (1) and (2) together with
  Lemma \ref{t:setup} that  proved
  ($\star$) in the ULSC case.  However,   Lemma \ref{t:setup}
  did not actually require the sequence to be ULSC, but rather
  requires only that (1) and (2) above hold.
  This completes the sketch of (C2) and (D) in Theorem
\ref{t:t5.2}.

\section{Putting it all together: The proof of Theorem
\ref{t:t5.2}}   \label{s:puttog}

We have now completed the proof of all six of the claims in
Theorem \ref{t:t5.2} over the course of this part (the six claims
are (A), (B), (C1), (C2), (D), and (P) in Theorem \ref{t:t5.2}).
For the reader's convenience, we will review next where each was
proven:

\begin{proof}
(of Theorem \ref{t:t5.2}).  The singular set $\cS$ is defined in
Definition/Lemma \ref{l:inftyornot}, where we also prove property
(B).   Lemma \ref{l:singl} gives a subsequence $\Sigma_j$ that
converges to a minimal lamination $\cL'$ of $\RR^3 \setminus \cS$,
thus giving (A). Property (C1) that describes the points in $\cSt$
is  established in Section \ref{s:ceno}. The properties (C2) and
(D) that describe the points in $\cSu$ are established in Section
\ref{s:ctwod}. Finally, property (P) that shows that the leaves of
$\cL'$ don't cross the limit planes given by (C1) is proven in
Lemma \ref{l:sm1}.
\end{proof}

\part{The no mixing theorem, Theorem \ref{c:main}}   \label{p:nomix}

\setcounter{equation}{0}

This part is devoted to the proof of the no mixing theorem, i.e.,
Theorem \ref{c:main}.  Recall that this theorem asserts that the
singular set $\cS$ consists of either exclusively helicoid points
or exclusively catenoid points, i.e., either $\cSt = \emptyset$ or
$\cSu = \emptyset$.  We have already shown in (D) of Theorem
\ref{t:t5.2} that the leaves intersecting $\cSu$ foliate an open
subset of $\RR^3$ that does not intersect $\cSt$.    Using in part
that $\cS$ is closed, we will show that the closure of this
foliated region will also intersect $\cS$; therefore, the boundary
of the foliated region must intersect $\cSt = \cS \setminus \cSu$.
We will prove the no mixing theorem by showing that also the
closure of this foliated region does not intersect $\cSt$ and,
hence, the foliated region is either empty or all of $\RR^3$.  The
argument for this will be very similar to the one that we used
earlier to show that the leaves intersecting $\cSu$ do not
intersect $\cSt$.

\begin{figure}[htbp]
 \setlength{\captionindent}{20pt}
    \begin{minipage}[t]{0.5\textwidth}
    \centering\input{genp4a.pstex_t}
    \caption{A limit lamination where the ULSC region has non-empty
    boundary; we will rule this out.}
    \label{f:p4a}
    \end{minipage}\begin{minipage}[t]{0.5\textwidth}
    \centering\input{genp4b.pstex_t}
    \caption{Properties of a sequence $\Sigma_j$ that converges to a lamination
    where the ULSC region has non-empty boundary.}
    \label{f:p4b}
    \end{minipage}
\end{figure}

\begin{proof}
(of Theorem \ref{c:main}).  Suppose that $\cSu \ne \emptyset$; we
will show that $\cSt = \emptyset$. By (D) of Theorem \ref{t:t5.2},
the set $\cSu$ is a union of Lipschitz curves transverse to the
leaves of the lamination and
  the leaves intersecting $\cSu$ foliate an open subset of
$\RR^3$ that does not intersect $\cSt$; we will call this foliated
region the ``ULSC region''.
  Moreover,   each connected component of the ULSC region contains
  exactly two curves in $\cSu$ and
  each of these two curves intersects each leaf exactly once.
 We will prove the theorem by showing
that these ULSC leaves foliate all of $\RR^3$.  It is easy to see
that this is equivalent to showing that a curve in $\cSu$ cannot
just stop.

Note first that each curve in $\cSu$ is in fact a line segment
{\it{orthogonal}} to  the leaves of the lamination.  This follows
from the main theorem of \cite{Me1} since we have already proven
here that  the ULSC regions are foliated.

Suppose now that one of the line segments in $\cSu$ does stop,
i.e., has an endpoint $z$.  Since the set $\cS$ is closed and
$\cSu$ is open in $\cS$, the endpoint must be in $\cSt$.
 To complete the proof, we will use the
properties of the sequence $\Sigma_j$ to show that having   $z$ in
$\cSt$ leads to a contradiction. This contradiction follows from
the following four steps:
\begin{enumerate}
\item Since $z \in \cSt$, Proposition \ref{l:getsetgo} gives a sequence
of stable graphs $\Gamma_j \subset B_{R_j} \setminus \Sigma_j$
converging to a horizontal plane through $z$; this convergence is
smooth away from the point $z$. Proposition \ref{l:getsetgo} also
gives
 separating  curves $\gamma_j \subset \Sigma_j$ with $\gamma_j \to
 z$ where $\gamma_j$ divides $\Sigma_j$ into  a component
 $\Sigma_j^+$ above $\Gamma_j$ and a component $\Sigma_j^-$
 below.{\footnote{When we say that ``$\Sigma_j^+$ is above $\Gamma_j$'', we
 have to be a little bit careful since each $\Gamma_j$ is defined
 only over an annulus.  The precise statement is given in
 Proposition
 \ref{l:getsetgo}: There are shrinking balls $B_{r_j}(z)$ with
 $r_j \to 0$ so that $\Gamma_j \cup B_{r_j}(z)$ divides
 $B_{R_j/C}$ into components $H^+_j$ above $\Gamma_j$ and $H^-_j$
 below $\Gamma_j$; we then have that $\Sigma_j^+$ is contained in
 $H^+_j\cup B_{r_j}(z)$.}}
 After possibly
reflecting about the horizontal plane through $z$, we can assume
that the segment in $\cSu$ lies above the plane.
\item  Since $z$ is in the closure of $\cSu$, (C2) of Theorem
\ref{t:t5.2} gives   double spiral staircases in $\Sigma_j^+$
above the horizontal plane through $z$.  More precisely,
 fix a ball $B_{2s}(z)$ that intersects only one
of the two vertical line segments in $\cSu$ approaching the plane
through $z$ from above.   After possibly shrinking $s$, we can
also assume that $B_{2s}(z) \cap \{ x_3
> x_3 (z) \}$  is contained in one connected
component of the ULSC foliated region.  Then, by (C2) of Theorem
\ref{t:t5.2}, in each compact subset of $B_{2s}(z) \cap \{ x_3
> x_3 (z) \}$, $\Sigma_j^+$  will consist of a double spiral
staircase for all $j$ sufficiently large.
 This will be used in (4) to pull the double spiral staircases
into an appropriate region near the plane.
\item
We will show next that $\Sigma_j^+$ must be ULSC away from $z$.
Namely, we will show that:
\begin{itemize} \item There exists some
$\epsilon
> 0$ so that if $y \in (B_{2s}(z) \setminus B_{s}(z)) \cap \{ x_3
= x_3 (z) \}$, then each component of $B_{\epsilon \, s}(y) \cap
\Sigma_j$ is a disk for $j$ sufficiently large.
\end{itemize}
We will prove this by contradiction, so suppose that there is a
sequence of non-contractible curves $\tilde \gamma_j$ in
$B_{\epsilon \, s} (y) \cap \Sigma_j^+$; the constant $\epsilon$
will be given by  Proposition \ref{l:getsetgo}.   Applying
Proposition \ref{l:getsetgo} to the $\tilde \gamma_j$'s will lead
to the desired contradiction. Namely, Proposition \ref{l:getsetgo}
gives a second stable graph $\tilde{\Gamma}_j$ that is disjoint
from both $\Sigma_j$ and $\Gamma_j$.  Since $\tilde{\Gamma}_j$
starts off from $\tilde \gamma_j$, we see that $\tilde{\Gamma}_j$
is above $\Gamma_j$.  However, it follows from (1) that  the axis
of the double spiral staircase in $\Sigma_j^+$ can be connected to
$\Gamma_j$ by a short curve $\sigma_j$ in $\Sigma_j$; here short
means that the length of $\sigma_j$ goes to zero as $j \to
\infty$.  In particular, the short curve $\sigma_j$ does not pass
through the ``hole'' in the annulus $\tilde{\Gamma}_j$. Therefore,
the graph $\tilde{\Gamma}_j$ must intersect $\Sigma_j \cup
\Gamma_j$ which is a contradiction.
\item  Finally, (3) will allow us to apply the one-sided curvature
estimate{\footnote{The one-sided curvature estimate from
\cite{CM6} is recalled in Theorem \ref{t:t2}.}} to show that the
$\Sigma_j^+$'s continue to spiral as graphs below the plane $\{
x_3 = x_3 (z) \}$, contradicting (1). To do this, suppose that $y$
is any given point in $(B_{2s}(z) \setminus B_{s}(z)) \cap \{ x_3
= x_3 (z) \}$. Observe that (1) gives the sequence of stable
$\Gamma_j$'s disjoint from $\Sigma_j^+$ converging in $B_s(y)$ to
the horizontal disk $B_{s}(y) \cap \{ x_3 = x_3(z) \}$.  It
follows from this and (3) that we can apply the one-sided
curvature estimate to get that each component of $B_{\epsilon' \,
s}(y) \cap \Sigma_j$ is a graph for $j$ large; here $\epsilon' >
0$ depends on $\epsilon$ and the constant from the one-sided
curvature estimate.  Since these components are graphical for
every such point $y$ and start out as part of a multi-valued
graph, there are now two possibilities:
\begin{itemize}
\item
The multi-valued graph can be continued down to  $x_3 = x_3(z) -
\epsilon' \, s$. \item The multi-valued graphs spiral infinitely
into some horizontal plane above $x_3 = x_3(z) - \epsilon' \, s$.
\end{itemize}
The latter is impossible since each $\Sigma_j$ is a compact
surface.  This completes the proof of the fourth step and, hence,
gives the promised contradiction to (1).
\end{enumerate}
\end{proof}

\part{Completing the proofs of Theorem \ref{t:tab} and Theorem  \ref{t:t5.2a}}
\label{p:prove3}

\setcounter{equation}{0}

The only thing that remains to be proven is that
{\underline{every}} leaf of the lamination $\cL'${\footnote{Recall
that $\cL'$ is a lamination of $\RR^3 \setminus \cS$ given in
Lemma \ref{l:singl}.}} is contained in a plane.  This is the
remaining claim in Theorem \ref{t:tab} (the planar lamination
convergence theorem).  We have already proven that the leaves of
$\cL'$ are planes when the sequence is ULSC; thus, by the no
mixing theorem, the only remaining case is when $\cS = \cSt \ne
\emptyset$.

Recall that each point in $\cSt$ comes with a plane through it
that is a limit of stable graphs in the complement of the sequence
$\Sigma_j$.    Since $\cSt \ne \emptyset$ by assumption,  there is
at least one such plane and, hence, every leaf of $\cL'$ is
contained in a half-space (by Lemma \ref{l:sm1}).  We will divide
the proof that the leaves of $\cL'$ are contained in planes into
two cases, depending on whether or not the leaf is complete.
Recall that a leaf $\Gamma$ is complete when $\barga = \Gamma$,
where the closure $\barga$ is defined by fixing a point
$x_{\Gamma} \in \Gamma$ and setting
\begin{equation}
    \barga = \bigcup_{r} \overline{ \cB_r (x_{\Gamma}) } \,
    ;
\end{equation}
see \eqr{e:closureaa}.  We will prove that complete leaves of
$\cL'$ are planes in Lemma \ref{l:smoothleaves}; the in--complete
leaves will be shown to be planes in Lemma
\ref{l:nonsmoothleaves}.

It may be useful to give an example of the kind of thing that we
need to rule out and a rough idea of why it cannot happen. Suppose
therefore that a leaf $\Gamma$ of $\cL'$ contains infinitely many
necks, one on top of the next, and that these necks ``shrink'' to
a point $p \in \cSt$. It follows that we have a limit plane
through $p$ and that $\Gamma$ is contained on one side of this
plane.  We will use a flux argument to rule out such an example.
Roughly speaking, we will find a ``top'' curve with positive flux
and then find a sequence of ``bottom'' curves shrinking to $p$
whose flux goes to zero.  We will then show that all of the ends
of $\Gamma$ between these curves are asymptotic to planes or
upward sloping catenoids - and hence make a non-negative
contribution to the total flux. This will give the desired
contradiction since Stokes' theorem implies that the total flux is
zero.

\section{Blow up results for ULSC surfaces}

We will later need to analyze the structure of the sequence
$\Sigma_j$  in  Theorem \ref{t:tab} (the planar lamination
convergence theorem) near points where the topology is
concentrating.  In doing so, it will often be useful to work on
the smallest scale of non-trivial topology (this is similar to
blowing up on the scale of the curvature in the ULSC case; cf. the
notion of blow up pairs in \cite{CM6}). This can be achieved using
a simple rescaling argument given in Lemma \ref{l:l5.1} in
Appendix \ref{s:aE}.

The advantage of working on the smallest scale of non-trivial
topology is that we can use the compactness theorem for ULSC
sequences -- Theorem \ref{t:t5.1} -- to prove a great deal of
structure for the surfaces on this scale.  We will use two such
structure results below:
\begin{itemize}
\item   Lemma \ref{l:os} proves a one-sided property for non-simply connected
    surfaces on the smallest scale of non-trivial
topology.  This shows that an intrinsic ball in such a surface
cannot lie on one side a plane and have its center close to the
plane on this scale; in the extreme case as the radius of the
intrinsic ball goes to infinity, the surface would be forced to
grow out of any half-space.
\item   Lemma \ref{l:shortc} finds short curves on the smallest scale of non-trivial
topology separating the ends.  These short curves will be used in
the flux argument for the main results in this part.
\end{itemize}

\subsection{Finding short separating curves}

 The next lemma finds short separating curves and stable graphs near points with small
 injectivity radius.  These curves will separate the surface into two parts: a part $\Sigma^+$ above
the graph and a part $\Sigma^-$ below the graph; cf. Figure
\ref{f:getsetgo}.
 Earlier, in Proposition \ref{l:getsetgo}, we found separating curves contained in small extrinsic
 balls; in fact (1) and (2) below are proven in Proposition \ref{l:getsetgo}.
 The new point here is the bound on the length of the curves in (3) below.
 To prove this length bound, we will work on the smallest scale of non-trivial
topology.

\begin{Lem}     \label{l:shortc}
Let $\Sigma \subset B_{r_0}$ be an embedded minimal planar domain
with $\partial \Sigma \subset
\partial B_{r_0}$.  Suppose also that $\cB_{r_1} \subset \Sigma$
is {\underline{not}} a disk and $\Gamma \subset B_{r_0} \setminus
\Sigma$ is the stable surface given by Lemma
\ref{l:getset}.{\footnote{$\Gamma$ is given by solving a Plateau
problem using a non-contractible curve in $\cB_{r_1} \subset
\Sigma$ as ``interior'' boundary.}}

The following three properties hold (the first two are just
Proposition \ref{l:getsetgo}):
\begin{enumerate}
\item
Given $\tau
> 0$, there exists $C \geq 1$ so  a component $\Gamma_0$ of
$B_{r_0/C} \cap \Gamma \setminus B_{C r_1}$ is a graph of a
function $v$ with $|\nabla v| \leq \tau$ and $\partial \Gamma_0$
intersects both $\partial B_{r_0/C}$ and $\partial B_{C r_1}$.
\item
There are distinct components $H^{+}$  and $H^-$ of $B_{r_0/C}
\setminus (\Gamma_0 \cup B_{C \, r_1})$, a separating curve
$\tilde \sigma \subset B_{C \, r_1} \cap \Sigma$, and components
$\Sigma^{\pm}$ of $B_{r_0/C} \cap \Sigma \setminus \tilde  \sigma$
so that $\Sigma^{\pm} \subset H^{\pm} \cup B_{C \, r_1}$ and
$\tilde \sigma \subset
\partial \Sigma^{\pm}$.
\item
There exists $C_1$ so that if $r_0 \geq C_1 \, r_1$ and
$\cB_{r_1/4}(x)$ is a disk for each $x \in \cB_{C_1 \, r_1}$, then
$\tilde \sigma$ is homologous to a {\underline{collection}} of
curves whose total length   is at most $C_1 \, r_1$.
\end{enumerate}
\end{Lem}

\begin{Rem}
We will call this the ``short curve lemma'' and call the
separating curve $\tilde \sigma$ the ``short curve.''  Of course,
$\tilde \sigma$ itself may not be short; rather it is homologous
to a collection of curves whose total length is at most $C_1 \,
r_1$.
\end{Rem}

\begin{proof}
Claims (1) and (2) are proven in Proposition \ref{l:getsetgo}.  We
will prove the last claim by contradiction, so suppose that (3)
fails with $C_1 = j$ for a sequence $\Sigma_j'$.  After rescaling,
we can assume that $r_1 = 4$.

Observe first that Lemma \ref{l:cyulsc} in Appendix \ref{s:aC}
gives a sequence $R_j \to \infty$ so that the component
$\Sigma_{j}$ of $B_{R_j} \cap \Sigma_j'$ containing $0$ is
contained in the intrinsic ball $\cB_{4 \, j}$ and, hence, each
intrinsic ball of radius one in $\Sigma_j$ is a disk.  The
  $\Sigma_{j}$'s therefore give a ULSC sequence   of embedded minimal planar domains in
extrinsic balls whose radii go to infinity. As in the proof of
Lemma \ref{l:os}, we will now divide into two cases depending on
whether or not the curvatures of the sequence blows up.

\vskip2mm \noindent {\underline{Case 1}}:
 Suppose first that there exists some $R$ so that
 \begin{equation}
    \limsup_{j \to \infty} \, \sup_{B_R \cap \Sigma_j} |A|^2
=  \infty \, .
\end{equation}
 We can then apply the
compactness theorem for ULSC sequences,   Theorem \ref{t:t5.1}, to
get a subsequence of the $\Sigma_j$'s that converges to a
foliation by parallel planes away from two lines orthogonal to the
leaves of the foliation.  It follows easily from the description
of the convergence near the lines (as double spiral staircases)
that we get a uniform length bound for the separating curve. This
length bound is proven in ``{\underline{The proof of (P4)}}''
within the proof of Lemma \ref{t:setup} and will not be repeated
here.

\vskip2mm \noindent {\underline{Case 2}}: Suppose now that $|A|^2$
is uniformly bounded on each compact subset of $\RR^3$ for the
sequence $\Sigma_j$. The length bound in this case follows by
combining three facts:
\begin{itemize}
\item  By Lemma \ref{l:ca} in Appendix \ref{s:aC},
the uniform curvature bounds implies uniform area bounds for each
component of $\Sigma_j$ in extrinsic balls (the bound depends on
the ball but not on $j$). More precisely, if  $\Sigma_{j,R}$
denotes the component of $B_R \cap \Sigma_j$ containing $0$, then
Lemma \ref{l:ca} implies that
\begin{equation}    \label{e:fromlca}
    \Area \, (\Sigma_{j,R}) \leq C_c \, R^2 \, ,
\end{equation}
where the constant $C_c$ depends only on
\begin{equation}
    \sup_{B_{C_0 \,
    R} \cap \Sigma_j} \, |A|^2 \, .
\end{equation}
  The constant $C_0$
here is universal and does not depend on the upper bounds for the
curvature.
\item  The area bound \eqr{e:fromlca} and the co-area formula give uniform length
bounds for the boundary $\partial \Sigma_{j,R}$ for most values of
$R$.  Precisely, at least one-half of the $R$'s between $R_0/2$
and $R_0$ must satisfy
\begin{equation}    \label{e:fromlca2}
    \Length \, (\partial \Sigma_{j,R}) \leq 2 \, C_c \, R_0 \leq 4 \, C_c \, R \,
    .
\end{equation}
\item  In the proof of Proposition \ref{l:getsetgo}, the components of
$\partial \Sigma_{j,R}$ were divided into two groups, depending on
whether they connected to $\Sigma_j^+$ or $\Sigma_j^-$;  the
separating curve $\tilde \sigma_j$  was then chosen to separate
these two groups.{\footnote{This  is described in more detail in
the proof of  Proposition \ref{l:getsetgo}.}}   However,  if we do
not ask for a {\underline{single}} connected separating curve,
then we can instead use either of the two  groups to separate.
Finally, \eqr{e:fromlca2} gives a uniform bound for the total
length.
\end{itemize}

Therefore, in either case, we get uniform length bounds, hence
proving the lemma.
\end{proof}

\section{Complete leaves of $\cL'$}

We will show next that any complete leaf $\Gamma$ of the
lamination $\cL'$ is a plane ($\cL'$ is the lamination of $\RR^3
\setminus \cS$ given in Lemma \ref{l:singl}).  Such a leaf
$\Gamma$ is a complete embedded minimal surface in $\RR^3$, but is
not {\it a priori} known to be proper.

\begin{Lem}     \label{l:smoothleaves}
Suppose that $\Gamma$ is a complete leaf of the lamination $\cL'$,
i.e, suppose that $\barga = \Gamma$.  Then $\Gamma$ must be a
plane.
\end{Lem}

\begin{proof}
We will assume that $\Gamma$ is not a plane and show that this
leads to a contradiction.

The ULSC case was already completed in Theorem \ref{t:t5.1}, so we
can assume that
\begin{equation}
    \cS = \cSt \ne \emptyset
\end{equation}
     by the
no mixing theorem (Theorem \ref{c:main}).  Recall that, by
Proposition \ref{l:getsetgo}, each point in $\cSt$ comes with a
plane through it that is a limit of stable graphs in the
complement of the sequence $\Sigma_j$.{\footnote{More precisely,
repeatedly applying Proposition \ref{l:getsetgo} gave a sequence
of stable graphs defined over larger and larger annuli and this
sequence converges to a limit graph over a punctured plane. The
limit graph is bounded at the puncture and extends smoothly across
the puncture to an entire graph; consequently, the limit graph is
a plane by the Bernstein theorem.}}
 Furthermore, by   Lemma \ref{l:sm1}, the leaves of $\cL'$ do not cross any
of these planes.  Since there is at least one such plane, $\Gamma$
is contained in a half-space. After a translation and a rotation,
we may
    assume that $\Gamma \subset \{ x_3 \geq 0 \}$ and
    \begin{equation}    \label{e:mayassu}
        \inf_{\Gamma} \, x_3 = 0 \, .
    \end{equation}

\vskip2mm \noindent {\underline{Claim}}:  If $\Gamma$ satisfies
\eqr{e:mayassu} and is not a plane, then there is a sequence of
points $p_n \in \Gamma$ satisfying:
\begin{align}
     \label{e:pnti3}
                \cali (p_n) &\to 0 \, , \\
\label{e:pnti3a}
                x_3 (p_n) &\to 0 \, .
\end{align}
Here $\cali (p_n)$ is the
        injectivity radius  of $\Gamma$ at $p_n$.
Since $\Gamma$ is a complete smooth surface,  \eqr{e:pnti3}
immediately implies that
\begin{equation}
         \label{e:pnti2}
     \dist_{\Gamma}(p_1 ,  p_n) \to
\infty \, .
\end{equation}

\vskip2mm \noindent {\underline{Proof of Claim}}: Since $\Gamma$
satisfies \eqr{e:mayassu} but is not a plane, \cite{MeRo} (see the
first paragraph of the proof of lemma $1.5$ there) implies that
for any $\epsilon > 0$ we have
\begin{equation}    \label{e:merogives}
    \sup_{\Gamma \cap \{ 0 < x_3 < \epsilon \} } \, |A|^2 = \infty
    \, .
\end{equation}
Therefore, we get a sequence of points $p_n$ in $\Gamma$
satisfying \eqr{e:pnti3a} and with
\begin{equation}    \label{e:asqu}
    |A|^2 (p_n) \to \infty \, .
\end{equation}
Equation \eqr{e:pnti3} then follows immediately from this and the
one-sided curvature estimate. {\underline{QED of Claim}}

\vskip2mm \noindent {\underline{Return to the proof of the
lemma}}:
  We will use the Claim above to deduce a flux contradiction (similar to
    the proof of ($\star$) in the ULSC case given in Subsection
    \ref{ss:p}) as follows:
    \begin{enumerate}
        \item[(a)]
            {\underline{The leaf $\Gamma$ must be a multiplicity one limit of the
            $\Sigma_j$'s.}}   To see this, observe that if this was not the case,
            then the universal cover of $\Gamma$ would be stable and, hence,
            flat; cf. the proof of Corollary \ref{c:collstab} for more details.
        \item[(b)]
            {\underline{Blowing up at $p_1$ to get a separating
            curve.}}
         Fix a large constant $C_1 > 1$ (it will be chosen depending on both Lemma
            \ref{l:shortc} - the ``short curve lemma'' - and Lemma \ref{l:os} -
            the ``one-sided lemma'' for non-simply
            connected surfaces).
         Applying the blow up lemma,  Lemma
            \ref{l:l5.1}, at $p_1$ gives an intrinsic
            ball \begin{equation}   \label{e:gets1p1}
                \cB_{C_1 s_1} (y_1) \subset \cB_{5 \, C_1 \, \cali
                (p_1)} (p_1) \, ,
                \end{equation}
                so that $\cB_{4\, s_1}(y_1)$ is not a disk but
                $\cB_{s_1} (y)$ is a disk for each $y \in \cB_{C_1 s_1}
                (y_1)$.

            Taking $C_1$ sufficiently large, it follows that $\cB_{C_1 s_1}
                (y_1)$ satisfies the  hypotheses of both Lemma
            \ref{l:shortc} and Lemma \ref{l:os} (where the constant $H$ in Lemma \ref{l:os} is set equal to
            a large constant $C_2 > 1$).
            Hence, the short curve lemma, Lemma
            \ref{l:shortc}, gives an initial short separating
            curve
            \begin{equation}
                \gamma_{1} \subset B_{4\, s_1}(y_1) \cap \Gamma
            \end{equation}
                 and a stable graph
            \begin{equation}
                \Gamma_0 \subset \RR^3 \setminus \Gamma \, .
            \end{equation}
            Since the surface $\Gamma$ is not known to be proper in all of $\RR^3$,
            the graph $\Gamma_0$ would at first appear to be defined
            only over a bounded annulus.  However, the multiplicity one convergence of (a) implies that
            the short curve $\gamma_{1} \subset \Gamma$ is
            actually a smooth limit of  curves $\gamma_{1,j}$
            contained in the proper surfaces $\Sigma_j$.  We can
            therefore apply the barrier construction to these
            curves in $\Sigma_j$ and take the limit of the
            resulting stable graphs to get the desired $\Gamma_0$
            as a graph defined outside the ball $B_{C\, s_1}(y_1)$.  It
            follows that
            the  graph  $\Gamma_0$ is
             asymptotic to either a plane or an upward sloping
            half--catenoid (the other possibility would be a downward sloping half--catenoid which is clearly
             impossible since
            $\Gamma$ is above $\{ x_3 = 0 \}$).

            Moreover,  since $\Gamma \subset \{ x_3 \geq 0 \}$, the one-sided lemma for non-simply
            connected surfaces, Lemma \ref{l:os}, guarantees
            that
            \begin{equation}    \label{e:fromcos}
                C_2 \, s_1 < x_3 (y_1) \, ,
            \end{equation}
            where $C_2 > 1$ is a large fixed constant (we can make
            $C_2$ as large as we want by increasing $C_1$).  By the same argument,
            the extrinsic ball $B_{C_2 \, s_1} (y_1)$
            does not intersect any of the horizontal planes
            associated to the singular set $\cS$.

            Finally, since $x_3 (p_n)$ and $\cali (p_n)$ both go  to zero, we
            can pass to a subsequence of the $p_n$'s so that
            \begin{equation}    \label{e:dividingpr1}
                \sup_{\Gamma} x_3 > \sup_{ B_{C\, s_1}(y_1) } x_3
                \, ,
            \end{equation}
            and then for $n \geq 1$
            \begin{equation}    \label{e:dividingpr2}
                \inf_{ B_{C\, s_n}(y_n) } x_3 > \sup_{ B_{C\, s_{n+1}}(y_{n+1}) } x_3
                \, .
            \end{equation}

\begin{figure}[htbp]
 \setlength{\captionindent}{20pt}
    \begin{minipage}[t]{0.5\textwidth}
    \centering\input{genp30.pstex_t}
    \caption{(b): Lemma
            \ref{l:shortc} gives a short curve $\gamma_1 \subset \Gamma$ and a stable graph
            $\Gamma_0$.}
    \label{f:fora}
    \end{minipage}\begin{minipage}[t]{0.5\textwidth}
    \centering\input{genp31.pstex_t}
    \caption{(c'): $\RR^3 \setminus (\Gamma_0 \cup B_{C s_1}(y_1))$ has components $H^+$ above
            and $H^-$ below $\Gamma_0 \cup B_{C s_1}(y_1)$.
            $\Gamma$ is the only leaf of $\cL'$
            intersecting both $H^+$ and $H^-$.}
    \label{f:forb}
    \end{minipage}
\end{figure}

 \item[(c)]
         {\underline{$\Gamma$ is the only leaf of $\cL'$ that
            intersects  $\overline{B_{C \, s_1}(y_1)}$.}}
             A barrier argument and the one-sided
lemma for non-simply connected surfaces, Lemma \ref{l:os} give
$\tilde{C} > C$ so that if $C_1 \geq \tilde{C}$, then only one
{\underline{proper}} component of $B_{\tilde{C} \, s_1}(y_1) \cap
\cL'$ intersects $\overline{B_{C \, s_1}(y_1)}$.{\footnote{This
follows exactly as does the analogous result for disks given in
corollary $0.4$ in \cite{CM6}.  Namely, if there were two such
components, then we could put a stable surface between them.
Interior estimates for stable surfaces then imply that each of the
original components lies on one side of a plane that comes close
to the center of the ball.  However, this would contradict  the
one-sided  lemma for non-simply connected surfaces, Lemma
\ref{l:os}, so we conclude that there could not have been two such
components.}}  Here  $\tilde{C}$ depends only on $C$.

To complete the argument for (c), we need to verify that each
component of any leaf of $\cL'$
            in $B_{\tilde{C} \,
            s_1}(y_1)$ is proper. Fortunately, this will follow directly from
            Lemma \ref{l:ca} that
            gives the compactness of each component of
            an embedded minimal
            surface in a ball $B_R$ if there is
            {\underline{some}} curvature bound  in the fixed larger
            ball $B_{C_d \, R}$.{\footnote{Clearly, it is crucial
            here that $C_d$ does not depend on the bound for the
            curvature.}}
            Namely, \eqr{e:fromcos} implies that $B_{2 \, C_d \, \tilde{C} \,
            s_1}(y_1)$ is disjoint from $\cS$ so long as   $C_2$
              is sufficiently large.  We can
            then conclude that every leaf of $\cL'$ has bounded
            curvature in $B_{C_d \, \tilde{C} \,
            s_1}(y_1)$ and hence has compact components in $B_{\tilde{C} \,
            s_1}(y_1)$ by Lemma \ref{l:ca}.

         \item[(c')]
         {\underline{$\Gamma$ is the only leaf of $\cL'$ that
            intersects both sides of            $B_{C \, s_1}(y_1) \cup \Gamma_0$.}}
            Since the graph $\Gamma_0$ is a limit of surfaces that
            are disjoint from the $\Sigma_j$'s, it follows that
            none of the leaves of $\cL'$ can cross $\Gamma_0$.
            However, $\Gamma_0$ is a graph over an annulus, so the
            leaves of $\cL'$ may ``go through the hole'' to get from one side of $\Gamma_0$ to the
            other; this is exactly what $\Gamma$ does.  However,
             by (c), $\Gamma$ is the only leaf that intersects $B_{C \,
             s_1}(y_1)$, so we conclude that $\Gamma$ is the only leaf of $\cL'$ that
            intersects both sides of            $B_{C \, s_1}(y_1) \cup \Gamma_0$.

         \item[(d)]
            {\underline{Repeating (b) at each $p_n$.}}
            Using  \eqr{e:pnti3}, we can argue as in (b) at each point $p_n$ to
                get shrinking curves
                \begin{equation}
                    \gamma_n \subset B_{4 \, s_n} (y_n) {\text{ where }}
                    \cB_{C_1 s_n} (y_n) \subset \cB_{5 \, C_1 \, \cali
                (p_n)} (p_n)   \, ,
                \end{equation}
                as well as stable graphs that are defined outside $B_{C\,s_n}(y_n)$ and are
                disjoint from $\Gamma$.
                Since $\cali
                (p_n) \to 0$, Lemma
            \ref{l:shortc} gives that the flux across $\gamma_n$ also goes
                to zero.  Furthermore,  \eqr{e:pnti2}  guarantees that the shrinking curves are separated; the points
                $p_n$ may be close in $\RR^3$ but they are far apart in
                $\Gamma$.

                Let $\Gamma_n$ denote the connected component of
                $\Gamma \setminus (\gamma_1 \cup \gamma_n)$
                containing both $\gamma_1$ and $\gamma_n$ in its
                boundary.  Note that we used that $\gamma_1$ and
                $\gamma_n$ are separating to guarantee that such a
                component exists.

                Finally, let $E_n$ denote the ``sandwiched'' region in $\RR^3$
                that is between the stable graphs
                associated to $p_1$ and $p_n$ together with the balls $B_{C_1
                    \, s_1} (y_1)$ and $B_{C_1 \, s_n} (y_n)$.
                    We will need the following two properties of $E_n$:
                    \begin{equation}		\label{e:propen}
                      \Gamma_n \subset E_n {\text{ and }}
                      \overline{E_n} \cap \cS =
                     \emptyset \, .
                    \end{equation}
                    The first property follows immediately from (2) in Lemma
            \ref{l:shortc}.  To see the second, note that a point
            of $\cS$ in $\overline{E_n}$ would come with a
            horizontal plane through it that is disjoint from $\Gamma$; this is
            impossible since the connected leaf $\Gamma$
            intersects both above and below $E_n$.

            \item[(e)]  {\underline{$\Gamma_n$ is properly    embedded.}}
    We will prove this by contradiction, so
    suppose that some $\Gamma_n$ is not proper.  In this case,
    we would be able to choose a sequence $y_j \in \Gamma_n$
    with
    \begin{equation}
        \dist_{\Gamma_n} (y_1 , y_j) \to \infty {\text{ and }}
        |y_j - y| \to 0 {\text{ for some }}  y \in \overline{\Gamma_n} \subset \overline{E_n} \,
        .
    \end{equation}
    Since the union of the
    leaves of $\cL'$ is closed in $\RR^3 \setminus \cS$ and $\overline{E_n} \cap \cS =
                     \emptyset$, the point $y$ must be contained in some leaf $\tilde{\Gamma}$ of $\cL'$.
    As we have used several times, this implies that the    universal
                cover of $\tilde{\Gamma}$ must be stable; cf.
                the proof of Corollary \ref{c:collstab} for more details.
                Since $\Gamma$ is not stable (see the proof of
                (a)), it follows that $\tilde \Gamma \ne \Gamma$.

                 We claim
                that
                \begin{equation}   \label{e:tildegammac}
                    \tilde \Gamma {\text{ is complete.}}
                \end{equation}
                {\underline{Proof of \eqr{e:tildegammac}.}}
                We know from (c) that $\Gamma$ is the only leaf of
                $\cL'$ that intersects $\overline{B_{C \,
                s_1}(y_1)}\cup \overline{B_{C \, s_n}(y_n)}$, so
                \begin{equation}
                    y \notin \overline{B_{C \,
                s_1}(y_1)}\cup \overline{B_{C \, s_n}(y_n)} \, .
                \end{equation}
                It is also easy to see that the stable graphs
               that form the top and bottom of the boundary of
               $E_n$ cannot be contained in leaves of $\cL'${\footnote{These stable graphs were
               obtained using limits of solutions to Plateau problems using the $\Sigma_j$'s as barriers.}}, so
               we must have that
               \begin{equation}
                    y \in E_n \setminus ( \overline{B_{C \,
                s_1}(y_1)}\cup \overline{B_{C \, s_n}(y_n)} ) \, .
               \end{equation}
                However, (c') then implies that the entire leaf
                $\tilde \Gamma$ must be trapped inside of $E_n$.
                Since $\overline{E_n} \cap \cS = \emptyset$, it follows that
                $\tilde \Gamma$ must be complete.
                 {\underline{QED of \eqr{e:tildegammac}.}}

                Now that we have established \eqr{e:tildegammac}, the Bernstein theorem for stable surfaces
                implies that $\tilde
                \Gamma$ is a  plane.  Since $\tilde \Gamma$ does not cross   $\{ x_3 = 0 \}$,
                it must be a horizontal plane.  However, this is impossible
                since $\tilde \Gamma \cap \Gamma = \emptyset$ and
                $\Gamma$ intersects both above and below $E_n$.
                Therefore, we conclude that $\Gamma_n$ must be proper.

         \item[(f)]
            {\underline{The ends of $\Gamma_n$ are graphs.}}
                We claim next that for each fixed $n$, there is a constant $r_n$ so
                that
                \begin{equation}
                    \Gamma_n \cap \{ x_1^2 + x_2^2 \geq r_n^2 \}
                \end{equation}
                 consists of a
                finite collection of graphs over $\{ x_3 = 0 , \,  x_1^2 + x_2^2 \geq
                r_n^2 \}$.

                We will show first that $\Gamma_n \cap \{ x_1^2 + x_2^2 \geq r_n^2 \}$ is locally
                graphical.  The starting point is to observe that
                $\Gamma_n$ is contained in the sandwich $E_n$
                and the height of this sandwich grows at most logarithmically.
                Therefore, by the one-sided curvature estimate, it suffices
                to prove that $\Gamma_n$
                is scale-invariant ULSC with respect to the distance to $0$;
                see, for instance, (D) in Subsection \ref{ss:notone}.
                This follows from the ``between the sheets'' argument that we
                have used several times already, so we will just sketch
                the
                proof this time. Namely, since $\Gamma_n$ is
                connected, we can fix a curve $\sigma_n \subset
                \Gamma_n$ that connects $\gamma_1$ to $\gamma_n$;
                we will choose $r_n$ so that
                \begin{equation}
                    \sigma_n \subset \{ x_1^2 + x_2^2 <
                            r_n^2 / 4 \} \, .
                \end{equation}
                  If $\Gamma_n \cap \{ x_1^2 + x_2^2 \geq
                r_n^2 \}$ were to contain a scale-invariant small neck, then a
                barrier argument would give a stable surface $\Gamma_{barrier}$ in the
                complement of $\Gamma_n$ that is also sandwiched
                in $E_n$.
                This sandwiching and the curvature estimates for stable surfaces imply
                that the stable surface $\Gamma_{barrier}$  is
                graphical away from its boundary.  Since the curve
                $\sigma_n$ is away from the boundary of the stable surface and
                connects the top and bottom of the sandwich, the
                stable surface $\tilde \Gamma$ is
                forced to intersect the curve $\sigma_n$, giving
                the desired contradiction.

                After increasing $r_n$, we conclude that $\Gamma_n \cap \{ x_1^2 + x_2^2 \geq r_n^2 \}$
                is locally graphical and hence a union of graphs
                over $\{ x_3 = 0 , \,  x_1^2 + x_2^2 \geq
                r_n^2 \}$.  (The other possibility is that it
                could contain a multi-valued graph; as we have argued before, this is
                impossible since such a multi-valued graph would
                have to spiral through the separating plane.) The properness of
                $\Gamma_n$ proven in (e) implies that
                there can only be finitely
                many such graphs.

           Note
                that, by the isoperimetric inequality, this gives
                area bounds for $\Gamma_n$ in compact subsets of
                $\RR^3$.

    \item[(g)] {\underline{Slicing $\Gamma_n$ with a plane to get the top curve.}}
            Each graphical end of each $\Gamma_n$ is above $\{
            x_3 = 0 \}$ and, consequently, is asymptotic to either a plane or to
            an upward sloping half-catenoid.  Since there are only
            finitely many such planes for each $n$,  we can choose a
            height
            $h$ between $\sup_{ \gamma_2 } x_3$ and $\inf_{\gamma_1}
            x_3$ that misses all of the heights of   the
            planar ends for every $\Gamma_n$ and
            so that the plane $\{ x_3 = h \}$ intersects $\Gamma$
            transversely.  It follows that $\{ x_3 = h \}$
            intersects each $\Gamma_n$ transversely in a
            finite collection of simple closed curves.
            Note that this plane separates $\gamma_1$
            from $\gamma_n$ (and, in particular, does not intersect $\partial
            \Gamma_n$).

            Let $\Gamma_n'$ denote the component of $\{ x_3 < h \} \cap \Gamma_n$
            with $\gamma_n$ in its boundary.

    \item[(h)] {\underline{The flux contradiction.}}
         The boundary of each $\Gamma_n'$ consists of a ``bottom curve'' $\gamma_n$ together with
             a  collection of closed ``top   curves'' in
             the plane
             $\{ x_3 = h \}$.    The collection of top curves is ``increasing'' in the following sense
        \begin{equation}    \label{e:topcus}
            \{ x_3 = h \} \cap \partial \Gamma_n \subset \{ x_3 = h \} \cap \partial
            \Gamma_{n+1} \, .
        \end{equation}
        Generally, one might expect equality in \eqr{e:topcus};
        however, if $\Gamma_{n+1}$ contained a catenoidal end that
        was not in $\Gamma_n$, then we would have a strict
        containment.

             The integrand for the vertical flux is point--wise positive along
             the increasing
             boundary in $\{ x_3 = h \}$ and, hence, the vertical
             flux of $\Gamma_n'$ across $\{ x_3 = h \}$ is positive and non-decreasing as a function of $n$.  On
             the other hand, the flux across the bottom curve
             $\gamma_n$   goes to zero as $n \to \infty$ by
             (d).  We can therefore fix some large $n$ so that (the absolute value of) the flux across $\gamma_n$
             is less than the flux across $\{ x_3 = h \}$.
                Since $\Gamma_n'$ has only finitely many ends  and each of these ends
                has non-negative flux at infinity, the total flux of $\Gamma_n'$ is
                    positive. This gives the desired contradiction since, by Stokes' theorem, the
                    total flux of $\Gamma_n'$ must be zero.
  \end{enumerate}
\end{proof}

\begin{Rem} \label{r:notcomplete}
The above argument did not really need that the leaf $\Gamma$ was
complete in order to conclude that it must be flat.  Rather, we
showed that $\Gamma$ must be flat as long as there exists a
sequence of points $p_n \in \Gamma$ satisfying
\begin{equation}
    x_3 (p_n) \to 0 , \, i(p_n) \to 0 , {\text{ and }}
    \frac{i(p_n)}{ \dist_{\barga} (p_n , \cS) } \leq C_0 \, ,
\end{equation}
where $C_0$ is a fixed constant that does not depend on $\Gamma$.
This will be useful when we consider in-complete leaves in the
next section.
\end{Rem}

\section{Incomplete leaves of $\cL'$}

It remains to show that each in-complete leaf $\Gamma$ of  $\cL'$
also must be flat.  We do this in the next lemma.

\begin{Lem}     \label{l:nonsmoothleaves}
Suppose that $\Gamma$ is an incomplete leaf of the lamination
$\cL'$, i.e, suppose that $\barga \ne \Gamma$.  Then $\Gamma$ is
contained in a plane.
\end{Lem}

As proven earlier in Theorem \ref{t:t5.1}, every leaf of $\cL'$ is
flat when the sequence is ULSC.   Therefore, by the no mixing
theorem,  i.e., Theorem
        \ref{c:main},
        we can assume that $\cSu = \emptyset$ and $\cS = \cSt$.

Before getting into the proof, it is useful to consider an example
of what a possible incomplete non-flat leaf $\Gamma$ of $\cL'$
would have to look like.  By assumption,  $\barga \setminus \Gamma
\ne \emptyset$ and, hence, $\barga \cap \cSt \ne \emptyset$. Since
each point of $\barga \cap \cSt$ comes with a plane through it and
none of the leaves of $\cL'$ can cross these planes, such a
$\Gamma$ would be contained in either
\begin{itemize}
\item an open slab between two singular planes, or
\item an open half-space bounded by a singular plane.
\end{itemize}
Note that, by the strong maximum principle, $\Gamma$ cannot
intersect a singular plane and, hence, we can take the above slab
and half-space to be open. We will see in the next subsection that
$\barga \cap \cS$ consists of only one point in the boundary
plane(s).

The basic idea behind the proof of Lemma \ref{l:nonsmoothleaves}
is again that a potential counterexample would lead to a flux
contradiction.  Much of the argument is very similar to the
complete case: \begin{itemize}  \item $\Gamma$ will be
scale-invariant ULSC away from the singular points.  \item
$\Gamma$ will be proper in an open slab or open half-space.  \item
The ends of $\Gamma$ will   be asymptotic to planes or
upward-sloping catenoids. \item We will slice between two planar
ends to get a ``top curve'' with strictly positive flux. \item We
will find a sequence of ``bottom curves'' where the flux goes to
zero.
\end{itemize}
The   main difficulty lies in finding the sequence of ``bottom
curves'' where the flux goes to zero.   One  expects that the
injectivity radius of $\Gamma$ goes to zero as we approach the
singular points.  However, the rate at which it does so may be
quite slow, so we cannot find large regions in $\Gamma$ ``on the
smallest scale of non-trivial topology'' as the injectivity radius
goes to zero. The key for overcoming this will be to get some
additional control over $\Gamma$ near a singular point; in
particular, we will prove scale-invariant curvature and area
bounds for $\Gamma$ near each singular point. Once we have this,
we can use the co-area formula to find a sequence of ``bottom
curves'' whose length goes to zero.

\subsection{If $\Gamma$ is not flat, then $\barga \cap \cS$ consists of at most two points}

As mentioned, we have already shown that the complete leaves of
$\cL'$ must be flat, so the remaining case is when
\begin{equation}
    \barga \cap \cSt \ne \emptyset \, .
\end{equation}
Each point of $\barga \cap \cSt$ comes with a plane through it and
none of the leaves of $\cL'$ can cross this plane.  Hence, by the
strong maximum principle,   this plane does not intersect any of
the non-flat leaves of $\cL'$.
 The starting point for Lemma \ref{l:nonsmoothleaves} is to
show that this plane contains exactly one point of $\barga \cap
\cSt$; see Lemma \ref{l:strbarcS} below.  It follows immediately
from this that
 $\barga \cap \cS$ consists of at most two points
for any non-flat $\Gamma$.

\begin{Lem}     \label{l:strbarcS}
Suppose that $\Gamma \subset \{ x_3 > 0 \}$ is a non-flat leaf of
$\cL'$ with $0 \in \barga \cap \cSt$ and $\{ x_ 3 = 0 \}$ is the
associated stable limit plane through $0$.  Then we must have
\begin{equation}    \label{e:proofforgamma}
    \barga \cap \{ x_3 = 0 \} = \{ 0 \} \, .
 \end{equation}
 In fact,
 if  $\Gamma'
\subset \{ x_3 > 0 \}$ is {\underline{any}} non-flat leaf of
$\cL'$ with $\bargaprime
 \cap \{ x_3 = 0 \} \ne \emptyset$, then
 \begin{equation}   \label{e:proofforgammaprime}
    \bargaprime \cap \{ x_3 = 0 \} = \{ 0 \} \, .
 \end{equation}
\end{Lem}

\begin{proof}
We will first argue by contradiction to prove
\eqr{e:proofforgamma}.  Suppose therefore
 that there exists $p \ne 0$
with
\begin{equation}
     p \in \barga  \cap \{ x_3 = 0
    \} \, .
\end{equation}

We begin by constructing a curve $\gamma$ in $\Gamma$ that
connects $\Gamma$ to $0$ - or a singular point near $0$ - and
stays away from $p$. Precisely, $\gamma$ will have the following
properties:
\begin{equation}
    \gamma:[0,1) \to  B_{|p|/3} \cap \Gamma \, ,
\end{equation}
\begin{equation}
    \Length (\gamma ) \leq |p|/3  \, ,
\end{equation}
\begin{equation}
    \lim_{t \to 1} \, \gamma (t) \in \barga \cap \{ x_3 = 0 \} \,
    .
\end{equation}
To construct $\gamma$, first use the definition of $\barga$  to
choose a point $y \in \Gamma$ so that the closure of $\cB_{|p|/6}
(y) \subset \Gamma$ contains $0$.  Then choose a sequence of
length minimizing curves in $\Gamma$ that start at $y$ and whose
second endpoints converge to $0$.  The Arzela-Ascoli theorem gives
a subsequence of these curves that converges to a curve $\tilde
\gamma$ that starts at $y$, ends at $0$, and is contained in
$\overline{\cB_{|p|/6} (y)}$.  Finally, let  $\gamma$ be the
component of $\Gamma \cap \tilde \gamma$ that starts at $y$.

Note that the curve $\tilde \gamma$ might hit another point of
$\cS$ before it gets to $0$. However, this point must be close to
$0$ and, hence, far from $p$; this is all that the argument will
use. For simplicity, we will assume that $0$ was the first point
of $\cS$ hit by $\tilde \gamma$ so that $\lim_{t \to 1} \, \gamma
(t) = 0$.

 Since   $\gamma$ is contained in $\Gamma$, we get a
sequence of curves $\gamma_j: [0, t_j] \to  \Sigma_j$ with $t_j
\to 1$ and so that the $\gamma_j$'s converge to $\gamma$.  In
particular, $\gamma_j (t_j) \to 0$.

\vskip2mm \noindent {\bf{Claim}}:  The injectivity radius of
$\Sigma_j$ at $\gamma_j (t_j)$ must go to zero.

\vskip2mm \noindent {\bf{Proof of Claim}}:
  Proposition \ref{l:getsetgo} gives a  stable
  graph disjoint from $\Sigma_j$ for each $j$ and this sequence
  is converging to $\{ x_3 = 0 \} \setminus \{ 0 \}$ as $j$ goes
  to infinity.  Moreover,
  exactly one
component of $\Sigma_j$ in a small ball near $0$ intersects both
sides of the stable graph.  The injectivity radius of this
component (obviously) goes to zero as $j$ goes to infinity. It
follows that every other component sits on one side of this stable
graph; see (B) in Proposition \ref{l:getsetgo}. In particular, if
the component of $B_{\epsilon} \cap \Sigma_j$ containing $\gamma_j
(t_j)$ was a disk for some fixed $\epsilon > 0$ and all
sufficiently large $j$, then the one-sided curvature estimate
would imply that this component was graphical in a neighborhood of
$0$.  Moreover, by the strong maximum principle, this sequence of
graphs would have to converge to a subset of $\{ x_3 = 0 \}$.
However, these graphs contain subsets of $\gamma_j$ that are
converging to (a component of)
\begin{equation}
    B_{\epsilon} \cap \gamma \subset
        \Gamma \, .
\end{equation}
 It follows that $\{ x_3 = 0 \}$ would
have to contain a (smooth) point of the leaf $\Gamma$, violating
the strong maximum principle. {\bf{QED for Claim}}.

We can repeat the construction of $\gamma$ near $p$ to get  curves
$\gamma_j' :[0, t_j'] \to \Sigma_j$ converging to a curve $\gamma'
:[0,1) \to \Gamma$ so that the endpoints $\gamma_j' (t_j')$
converge   to a singular point near $p$.  For simplicity, we will
assume that this second singular point is actually equal to $p$.
 Arguing as in the Claim,
we see that the injectivity radius of $\Sigma_j$ at $\gamma_j'
(t_j')$ also goes to zero.

We can now apply
 Proposition
\ref{l:getsetgo} to shrinking balls centered at $\gamma_j(t_j)$
and $\gamma_j'(t_j')$ to get {\underline{disjoint}} stable graphs
$\Gamma_j$ and $\Gamma_j'$ that are disjoint from $\Sigma_j$ and
so
\begin{equation}
    \Gamma_j \to \{ x_3 = 0 \} \setminus \{ 0 \} {\text{ and }}
    \Gamma_j' \to \{ x_3 = 0 \} \setminus \{ p \} \, .
\end{equation}
 Since $\Gamma_j$ and $\Gamma_j'$ are disjoint, one must be above the
other.  After passing to a subsequence (and possibly switching
$\Gamma$ and $\Gamma'$), we can assume that $\Gamma_j$ is always
above $\Gamma_j'$.    It follows easily from the barrier
construction used for the proof of Proposition \ref{l:getsetgo}
that the curve $\gamma_j'$ must also be below the graph
$\Gamma_j$.{\footnote{Namely, the stable graph is actually a
subset of a stable surface that is disjoint from $\Sigma_j$ and
has  interior boundary lying in $\Sigma_j$; this interior boundary
connects within $\Sigma_j$ to the curve $\gamma_j'$.  This barrier
construction is given in Lemma \ref{l:getset}.}} However, this
forces $\gamma_j'$ to converge to a curve in $\{ x_3 = 0 \}$,
contradicting the strong maximum principle as in the proof of
Claim above.  This completes the proof of \eqr{e:proofforgamma}.

Finally,  when $p \in \bargaprime \cap \{ x_3 = 0 \}$, the same
argument applies with obvious changes.  Hence, we also get
\eqr{e:proofforgammaprime}.
\end{proof}

\subsection{The proof of Lemma \ref{l:nonsmoothleaves}}

As mentioned earlier, we can assume that we are in the case where
$\cSu = \emptyset$ and we will use a flux argument to rule out the
possibility of a non-flat leaf of $\cL'$.

\begin{proof}
(of Lemma \ref{l:nonsmoothleaves}). We will prove the lemma by
contradiction, so suppose that $\Gamma \subset \{ x_3
> 0 \}$ is a non-flat leaf of $\cL'$ with $0 \in \barga \cap \cSt$
and $\{ x_ 3 = 0 \}$ is the associated stable limit plane through
$0$.  By Lemma \ref{l:strbarcS}, there are two possibilities:
\begin{itemize}
\item $\barga \cap \cS = \{ 0 \}$. \item  $\barga \cap \cS = \{ 0
, \, p \}$ for some point $p$ with $x_3 (p) > 0$.
\end{itemize}

\vskip2mm \noindent {\underline{$\Gamma$ is scale-invariant ULSC
near $0$}}.  More precisely,   there exist  $\delta
> 0$ and $r_0 > 0$ so that
\begin{equation}    \label{e:notcase1}
    \cB_{\delta \, |x| } (x) {\text{ is a disk for every }} x \in
      B_{r_0} \cap \{ x_3 > 0 \} \cap \cL'  \, .
\end{equation}

 Recall  that the argument used to prove that complete leaves of $\cL'$ must be
 flat actually gave a stronger statement; see Remark
\ref{r:notcomplete}.  This stronger statement implies that
\eqr{e:notcase1} holds.

\vskip2mm \noindent {\underline{$\Gamma$ has quadratic curvature
blowup near $0$}}: We will next use a compactness argument to
prove that there exist constants $C_d$  and $r_1 > 0$ so that
\begin{equation}    \label{e:qcdcd}
    |A|^2 (x) \leq  C_d \, |x|^{-2}
    {\text{ for every }} x \in
      B_{r_1} \cap \{ x_3 > 0 \} \cap \cL'  \, .
\end{equation}
The constant $C_d$ above might depend on $\cL'$, but it will be
fixed throughout this proof.

 \vskip2mm \noindent {\underline{Proof of \eqr{e:qcdcd}}}:
We will argue by contradiction, so suppose that there is a
sequence
    of points $q_n \in \Gamma$ with $q_n \to 0$ and
    \begin{equation}    \label{e:qnc}
        |q_n|^2 \, |A|^2 (q_n) > n \, .
    \end{equation}
    The idea of the proof is that dilating
    $\cL'$ by the factor $|q_n|^{-1}$ about the point $q_n$  gives a sequence
    of laminations
    \begin{equation}
        \cL_n = |q_n|^{-1} \, \left( \cL' - q_n \right)
    \end{equation}
    with $|A|^2 (0) > n$ and so that $\partial B_1$
    intersects $\cSt (\cL_n)$; here $\cSt (\cL_n)$ is the singular set for the
    rescaled lamination $\cL_n$.   Moreover,
    \eqr{e:notcase1} gives a uniform lower bound for the
    injectivity radius of the leaves of $\{ x_3 > 0 \} \cap \cL_n$ in  $B_{1/2}$; see below for more details.  Consequently, as
    $n$ goes to infinity, a subsequence of the $\cL_n$'s would converge to a lamination
    $\cL_{\infty}$ with
    \begin{equation}
        0 \in \cSu (\cL_{\infty}) {\text{ and }}
        \partial B_1 \cap \cSt (\cL_{\infty}) \ne \emptyset \, .
    \end{equation}
     However, this would contradict the no mixing theorem, so
    we conclude that the sequence $q_n$ could not have existed.

    We need two things to make this outline rigorous.  First, we do not have a compactness
    theorem for sequences of laminations, but rather only for
    sequences of embedded minimal surfaces.  This is easily dealt
    with since the limit $\cL_{\infty}$ can be realized as a limit of a
      diagonal sequence of rescalings of the $\Sigma_j$'s; we will
      omit this standard argument.
      Second, we showed in  \eqr{e:notcase1} above only that the leaves  of
      $\cL'$ were scale-invariant ULSC near $0$; what we need instead is
      that the sequence $\Sigma_j$ is itself scale-invariant ULSC near
      $0$.  More precisely, we must show that there exists some
      $\delta_0 > 0$ so that for each fixed $n$ we have
      \begin{equation}  \label{e:disklag}
            {\text{every component of }} B_{\delta_0 \, |q_n|} (q_n) \cap
            \Sigma_j {\text{ is a disk for $j$ large}}.
      \end{equation}
      Since the $\Sigma_j$'s are
            converging to $\cL'$ away from $\cS$ - and $\Gamma$
            does not intersect $\cS$ - then
      the component of $B_{\delta \, |q_n|} (q_n) \cap
            \Sigma_j$ that is converging to $\Gamma$ is a disk
            with large curvature.
            However, the intrinsic version
            of the one-sided curvature estimate implies that this is
            the only component of this ball intersecting a
            smaller concentric sub-ball about $q_n$.  This gives
            the remaining ingredient needed to make the proof
            rigorous.

\vskip2mm \noindent {\underline{Extending flatness}}:  We claim
that there exist constants $C_{flat} > 0$ and  $r_2 > 0$ so that
if $r < r_2$, $x \in
\partial B_r \cap \Gamma$, and
\begin{equation}    \label{e:cflat}
     \cB_{r/4} (x) {\text{ is a graph with gradient less than }}
     C_{flat} {\text{ over }} \{ x_3 = 0 \} \, ,
\end{equation}
then $x$ is contained in a graph $\Gamma_x \subset \Gamma$ defined
over (at least) the annulus
 \begin{equation}
  \{ x_3 = 0 , \, r^2/4
< x_1^2 + x_2^2
 < r_2 \} \, .
\end{equation}
 In other words, once $\Gamma$ becomes very flat, then it
extends to a very flat graph defined over some annulus of a
definite size surrounding the singular point $0$.  Note that the
outer radius $r_2$ of this annulus is independent of $r$.

It is easy to prove from the gradient estimate and the quadratic
curvature bound \eqr{e:qcdcd} that  $x$ is contained in a very
flat graph $\Gamma_{x,r}$ defined over the annulus
 \begin{equation}
  \{ x_3 = 0 , \, r^2/4
< x_1^2 + x_2^2
 < C \, r^2 \} \, ,
\end{equation}
where the constant $C = C( \alpha)$ can be as large as we want for
$\alpha$ sufficiently small.
 A priori, one might worry that this would give a multi-valued
graph. However, by the usual argument, $\Sigma_j$ cannot contain a
multi-valued graph and, therefore, neither can
    $\Gamma$.    It remains to extend $\Gamma_{x,r}$ as a graph all the way
    out to $\partial D_{r_2}$ for some fixed $r_2$.  As long as
    $C$ is sufficiently large, this can be done using the
    sublinear growth of the height of the graph.  This sublinear
    growth is proven in proposition II.2.12 in \cite{CM3}.
    The details of the proof will be left
to the reader.

\vskip2mm \noindent {\underline{$\Gamma$ cannot be too
``horizontal'' near $0$}}: We will show next that
\begin{equation}    \label{e:doesblowup}
    \limsup_{s \to 0} \,  \inf_{\Gamma \cap (B_{2s} \setminus B_s)}
    | \langle \nn , (0,0,1) \rangle | < 1 \, ,
\end{equation}
where $\nn$ is the unit normal to the surface $\Gamma$.  Note that
$| \langle \nn (x) , (0,0,1) \rangle |$ is equal to one if and
only if the tangent plane at $x$ is horizontal.

 \vskip2mm \noindent {\underline{Proof of \eqr{e:doesblowup}}}:
Suppose first for some $s$ that
\begin{equation}
\label{e:doesblowup2} \inf_{\Gamma \cap (B_{2s} \setminus B_s)}
    | \langle \nn , (0,0,1) \rangle | > 0 \, .
\end{equation}
It follows $\Gamma$ is locally graphical with bounded gradient
    in $B_{2s} \setminus
    B_s$.  By the usual argument, $\Sigma_j$ cannot contain a multi-valued graph and, therefore, neither can
    $\Gamma$.   Hence, \eqr{e:doesblowup2} would imply that $\Gamma \cap (B_{2s} \setminus B_s)$ is a
    collection of graphs.

    Consequently, if we had a uniform lower bound for $| \langle \nn , (0,0,1) \rangle
    |$ in any neighborhood of $0$, then
 standard removable singularity theorems for minimal graphs would imply that $\Gamma$ has a
 removable singularity at $0$.
However, $\Gamma$ would have to be flat by the strong maximum
principle if the singularity at $0$ was removable.  We conclude
therefore that
\begin{equation}    \label{e:doesblowup3}
    \liminf_{s \to 0} \,  \inf_{\Gamma \cap (B_{2s} \setminus B_s)}
    | \langle \nn , (0,0,1) \rangle | = 0 \, .
\end{equation}
Finally, \eqr{e:doesblowup}  follows easily from ``Extending
flatness'' and \eqr{e:doesblowup3}.

\vskip2mm \noindent {\underline{Using \eqr{e:doesblowup} to blow
up $\cL'$}}.  The point about \eqr{e:doesblowup} is that any limit
of rescalings of $\cL'$ about $0$ will have a non-flat leaf. More
precisely, if $s_n$ is any sequence going to zero, then a
subsequence of the rescaled laminations
\begin{equation}    \label{e:rightafterthis}
    \cL_n = \frac{1}{s_n} \, \left( \cL' \right) \, ,
\end{equation}
will converge to a lamination $\cL_{\infty}$ of $\RR^{3} \setminus
\cS (\cL_{\infty})$ with the following properties:
\begin{enumerate}
\item[(P1)]
    The origin $0$ is still in $\cSt (\cL_{\infty})$ and
    $\{ x_3 = 0 \}$ is the corresponding limit plane.
\item[(P2)]
    The leaves of $\cL_{\infty}$ satisfy the quadratic curvature bound
    \eqr{e:qcdcd} in all of $\{ x_3 > 0 \}$ (not just in
    $B_{r_1}$),
    the singular set $\cS
    (\cL_{\infty})$ does not intersect the half-space $\{ x_3 > 0
    \}$, and
    $0$ is the only singular point in $\{ x_3 = 0 \}$ ``reachable'' from $\{ x_3
    > 0 \}$.
\item[(P3)]
    $\cL_{\infty}$ contains a  non-horizontal, and hence non-flat, leaf
    in $\{ x_3 > 0 \}$.
\end{enumerate}
The lamination $\cL_{\infty}$ is given as a limit of a subsequence
of rescalings of the $\Sigma_j$'s; see  the proof of \eqr{e:qcdcd}
for such a diagonal argument.  The first property (P1) follows
immediately from this.  The second property (P2) follows from
immediately from \eqr{e:qcdcd}.  Finally, \eqr{e:doesblowup}
implies that $\cL_{\infty}$ contains a non-horizontal leaf.  This
non-horizontal leaf cannot be flat since it would otherwise
intersect $\{ x_3 = 0 \} \setminus \{ 0 \}$, thus giving (P3).

The key point about the rescaled limit lamination $\cL_{\infty}$
is that it has all of the same properties that $\cL'$ did.
Therefore, we can repeat the construction to get that limits of
rescalings of $\cL_{\infty}$ also satisfy (P1), (P2), and (P3).
This will be important below, so we record it next:
\begin{equation}  \label{e:anysuch} {\text{
{\underline{Any}}  limit of rescalings of $\cL_{\infty}$   will
also satisfy (P1), (P2), and (P3).}}
\end{equation}

\vskip2mm {\it{From now on, we will assume that $\{ x_3 > 0 \}
\cap \cL'$ has quadratic curvature decay and
\begin{equation}
    \label{e:pnotthere}
    \{
x_3
> 0 \} \cap \cS = \emptyset \, ;
\end{equation}
this can be achieved by rescaling as above.}}

\vskip2mm \noindent {\underline{No stable leaves in $\{ x_3 > 0
\}$}}: We will show that a lamination
\begin{equation}    \label{e:nostablel}
{\text{$\cL_{\infty}$ satisfying (P1), (P2), and (P3) cannot have
a stable leaf in $\{ x_3
> 0 \}$.}}
\end{equation}
The same argument also rules out a leaf in $\{ x_3 > 0 \}$ whose
oriented double cover is stable.

\vskip2mm \noindent
 {\underline{Proof of no stable leaves in $\{ x_3 > 0
\}$}}: Suppose instead that $\cL_{\infty}$ did contain a stable
leaf $\tilde \Gamma$ in $\{ x_3 > 0 \}$.  We will show first
 that $\tilde \Gamma$ must be a flat plane $\{ x_3
= t \}$ for some $t > 0$. Namely, if it wasn't flat, then it would
be complete away from $0$ by Lemma \ref{l:strbarcS} (see equation
\ref{e:proofforgammaprime}) and then   Lemma \ref{l:punstab} in
Appendix \ref{s:apB} would give a contradiction.

Since the leaves of $\cL_{\infty}$ are - by definition - disjoint,
it follows that the non-flat leaf $\Gamma$ of $\cL_{\infty}$ must
be contained in the open slab $\{ 0 < x_3 < t \}$.  Set $t_0 =
\sup_{\Gamma} x_3$, so that \begin{equation} \label{e:mero1}
    \Gamma  \subset \{ x_3 < t_0 \} \, ,
    \end{equation}
     and
     \begin{equation} \label{e:mero2}
 {\text{$\Gamma$ intersects every tubular neighborhood
of the plane $\{ x_3 = t_0 \}$.}} \end{equation}
  Moreover, the
quadratic curvature bound \eqr{e:qcdcd} for the leaves of $\cL'$
implies that
\begin{equation}    \label{e:mero3}
    \sup_{ \{ t_0 / 2 < x_3 < t_0 \} \cap \Gamma } |A|^2 \leq 4 \, C_d \, t_0^{-2} < \infty \, .
\end{equation}
 However, the three properties \eqr{e:mero1}, \eqr{e:mero2}, and \eqr{e:mero3}
 are impossible by the first paragraph
of the proof of lemma $1.5$ in \cite{MeRo}.
 We conclude that \eqr{e:nostablel} must
    hold.

\vskip2mm \noindent {\underline{$\Gamma$ is proper}}: The first
application of \eqr{e:nostablel} will be to show that $\Gamma$
must be proper in compact subsets of $\{ x_3
> 0 \}$.

\vskip2mm \noindent
 {\underline{Proof of properness}}:
The starting point is that  $\Gamma$ would otherwise accumulate
into a stable leaf $\tilde \Gamma$; we have used this argument
several times and will omit the details (see, e.g., (e) in the
proof of Lemma \ref{l:smoothleaves}).{\footnote{To be precise,
either $\tilde \Gamma$ is stable or
    its oriented double cover is stable.}}  Clearly, $\tilde \Gamma$
intersects the open half-space $\{ x_3 > 0 \}$ and, hence, $\tilde
\Gamma$ must be contained in $\{ x_3 > 0 \}$ by the strong maximum
principle.   However, this is impossible by \eqr{e:nostablel}.

\vskip2mm \noindent {\underline{Scale-invariant area bounds}}:
Given any $\alpha > 0$, there exists a constant $C_{\alpha}$ so
that
\begin{equation}    \label{e:scabb}
    \Area \, ( (B_{2r} \setminus B_{r}) \cap \{ x_3 > \alpha \, |x| \} \cap
    \Gamma ) \leq C_{\alpha} \, r^2 \, .
\end{equation}

\vskip2mm \noindent
 {\underline{Proof of \eqr{e:scabb}}}:
 This will be pretty much the same argument as in the proof of
 ``$\Gamma$ is
 proper'' combined with a compactness argument.
    We will argue by contradiction, so suppose that \eqr{e:scabb} fails with $r= r_n$
    and $C_{\alpha} =n$ for every integer.
    By a diagonal argument and rescaling, we get a sequence of embedded minimal
planar domains $\tilde \Sigma_j$ with
\begin{equation}    \label{e:atoiny}
        \Area \, \left( (B_{2} \setminus B_{1}) \cap \{ x_3 > \alpha \, |x|
    \} \cap \tilde \Sigma_j  \right) \to \infty  \, .
    \end{equation}
    Recall that we have proven in \eqr{e:anysuch} that a subsequence of the
    $\tilde \Sigma_j$'s converges to a limit lamination
    $\cL_{\infty}$ off of a singular set $\cS (\cL_{\infty})$
    satisfying (P1), (P2), and (P3).

   Next, we will use \eqr{e:atoiny} to show that  $\cL_{\infty}$   contains a stable
   leaf in $\{ x_3 > 0 \}$, contradicting \eqr{e:nostablel}.
   This would be obvious if $\cL_{\infty}$ itself had
   infinite area in $(B_{2} \setminus B_{1}) \cap \{ x_3 > \alpha \, |x|
    \}$.   On the other hand, if $\cL_{\infty}$  had
   finite area in $(B_{2} \setminus B_{1}) \cap \{ x_3 > \alpha \, |x|
    \}$, then the
    $\tilde \Sigma_j$'s must converge with infinite multiplicity to
    some
    leaf $\tilde \Gamma$ of $\cL_{\infty}$ that intersects
    \begin{equation}
        \overline{(B_{2} \setminus B_{1})} \cap \{ x_3 \geq \alpha \, |x|
    \} \, .
    \end{equation}
    Note that we   used that $\overline{(B_{2} \setminus B_{1})} \cap \{ x_3 \geq \alpha \, |x|
    \}$ does not intersect the singular set $\cS_{\infty}$ (by (P2)) to
    guarantee the convergence of the $\tilde \Sigma_j$'s in this
    set.
      However, as we have used several times, this
    convergence with multiplicity
    implies that the leaf $\tilde \Gamma$ is stable; see, e.g.,
    the proof of Corollary \ref{c:collstab} for more
    details.{\footnote{To be precise, this convergence with
    multiplicity implies that either $\tilde \Gamma$ is stable or
    its oriented double cover is stable.}}
    Finally, since $\tilde \Gamma$
intersects the  half-space $\{ x_3 \geq \alpha \}$, it must be
contained in the open half-space $\{ x_3 > 0 \}$ by the strong
maximum principle.

\vskip2mm \noindent {\underline{Low points in $\Gamma$ are
contained in graphs}}:  We will need the following complete
version of ``Extending flatness'': There exists $\alpha > 0$  so
that if  $x \in \Gamma$ is in the ``low cone'' $\{ x_3 < \alpha \,
|x| \}$, then $x$ is contained in a graph $\Gamma_x \subset
\Gamma$ defined over (at least)
 \begin{equation}
  \{ x_3 = 0 , \, r^2/4
< x_1^2 + x_2^2
 < \infty \} \, ,
\end{equation}
where $r = |x|$. Moreover, the graph $\Gamma_x$ must be asymptotic
to a plane or to an upward-sloping half-catenoid.  Finally, there
is a positive lower bound for the height of the graph $\Gamma_x$,
i.e.,
\begin{equation}    \label{e:htlobd}
    \inf_{\Gamma_x} x_3 > 0 \, .
\end{equation}

\vskip2mm \noindent {\underline{Proof that low points in $\Gamma$
are contained in graphs}}: It follows from the gradient estimate
and the quadratic curvature decay of $\Gamma$ that $\Gamma$ is
``very flat'' in a neighborhood of $x$ in the sense of ``Extending
flatness.''  It then follows from the sublinear growth of the
height of the graph that  $\Gamma_x$ can then be extended over $\{
x_3 = 0 , \, r^2/4 < x_1^2 + x_2^2
 < \infty \}$ as long as $\alpha$ is sufficiently small.
    The proof of this extension will be left
to the reader.

Now that we know that $\Gamma_x$ is defined over $\{ x_3 = 0 , \,
r^2/4 < x_1^2 + x_2^2
 < \infty \}$, it follows that $\Gamma_x$ is asymptotic to either
 a plane, an upward-sloping half-catenoid, or a downward-sloping
 half-catenoid.  The last is impossible since $\Gamma_x$ is
 contained in $\{ x_3 > 0 \}$.

 Finally, \eqr{e:htlobd} follows from the maximum principle at
 infinity of \cite{LaRo}.

\vskip2mm \noindent {\underline{The components of $\Gamma
\setminus B_r$ are proper}}: Given any $r>0$, then
\begin{equation}   \label{e:procomps}
    {\text{each component of $\Gamma \setminus B_r$ is
proper.}}
\end{equation}

\vskip2mm \noindent {\underline{Proof of \eqr{e:procomps}}}:
 To prove \eqr{e:procomps}, we must show that any such
component $\Gamma_r$ cannot accumulate into $\{ x_3 = 0 \}$; this
is because we already know that $\Gamma$ itself is proper in $\{
x_3
> 0 \}$. We divide this into two cases.

First, suppose that the boundary $\partial \Gamma_r$ of the
component $\Gamma_r$ intersects the ``low cone'' - i.e., suppose
that
\begin{equation}    \label{e:caseonelow}
    \inf_{\partial \Gamma_r} \, x_3 <   \alpha \, r  \, .
\end{equation}
In this case, it follows that the entire component $\Gamma_r$ is a
graph and, hence, proper.

Suppose now that \eqr{e:caseonelow} does not hold.  In this case,
we will find a low component (for some smaller radius) that
extends as a graph underneath $\Gamma_r$, thus keeping $\Gamma_r$
strictly away from $\{ x_3 = 0 \}$.  To get this barrier
component, note that
  Lemma \ref{l:oc2} implies that $\Gamma$ contains a sequence of points $y_n \to 0$ contained
  in the low cone $\{ x_3 < \alpha \, |x| \}$.  If we choose $y_n$
  close enough to zero, then the resulting graph $\Gamma_{y_n}$
  must pass underneath $\partial \Gamma_r$ in $\partial B_r$.  It
  follows that $\Gamma_r$ sits above $\Gamma_{y_n} \cup B_r$ and,
  hence, cannot accumulate into $\{ x_3 = 0 \}$.
  This completes the proof of \eqr{e:procomps}.

\vskip2mm \noindent {\underline{The flux contradiction}}: We will
show that the non-flat leaf $\Gamma$ must contain a sequence of
{\underline{proper}} subdomains $\Gamma_n$ with the following
properties:
\begin{enumerate} \item[(top)] $\partial \Gamma_n$ contains an increasing
sequence of compact ``top curves'' in a fixed plane $\{ x_3 = h
\}$ for some $h$. Here, increasing means that $\{ x_3 = h \} \cap
\partial \Gamma_n \subset \partial \Gamma_{n+1}$ for every $n$.
\item[(ends)]  $\Gamma_n$ has finitely many ends and each end is
asymptotic to a plane or an upward-sloping half-catenoid.
\item[(bottom)] The rest of $\partial \Gamma_n$ consists of a
finite collection of ``bottom curves'' whose total length goes to
zero as $n$ goes to infinity.
\end{enumerate}
This will give a flux contradiction just as in the last step of
the proof of Lemma \ref{l:smoothleaves}.  Namely, the flux of
$\Gamma_n$ across the  top curves in $\{ x_3 = h \}$ is strictly
positive and non-decreasing in $n$, the ends have non-negative
flux, and the flux across the  bottom curves goes to zero.
However, this is impossible since the total flux for each
$\Gamma_n$ must be zero by Stokes' theorem.  It remains to
construct the $\Gamma_n$'s with these properties.

We will start with the ``top curve'' for $\partial \Gamma_n$.  As
in the proof of \eqr{e:procomps},  Lemma \ref{l:oc2} implies that
$\Gamma$ has infinitely many ``low'' ends that are asymptotic to
either planes or upward-sloping half-catenoids.  For simplicity,
we will assume that these ends are planar; the catenoid case
follows similarly and will be left to the reader.  Since
\cite{LaRo} ensures
 that the planar ends are asymptotic to different planes, we
 can choose some $h> 0$ between two consecutive planar ends
  so that $\{ x_3 = h \}$ intersects
 $\Gamma$ transversely.  Let $\Gamma_h$ be a component of $\{
 x_3 < h \} \cap \Gamma$ containing $0$ in its closure and fix
 some component $\gamma_h$ of $\partial \Gamma_h \subset \{ x_3 = h \}$.

Combining the coarea formula with the area bounds from
\eqr{e:scabb}, we can choose a sequence $r_n \to 0$ so that
\begin{equation}    \label{e:fromareatolen}
    \Length \, ( \partial B_{r_n} \cap \{ x_3 \geq \alpha \, r_n
    \} \cap \Gamma ) \leq C \, r_n \, ,
\end{equation}
for a uniform constant $C$ independent of $n$.  The point here is
that the length of these curves goes to zero as $n$ goes to
infinity.

For each $n$, let $\Gamma_n$ be the component of $\Gamma_h
\setminus B_{r_n}$ with $\gamma_h$ in its boundary.  First, it
follows immediately that (top) holds. Second, \eqr{e:procomps}
implies that $\Gamma_n$ is proper. Next, when $r_n$ is
sufficiently small, then each point in $\partial B_{r_n} \cap \{
x_3 < \alpha \, r_n
    \} \cap \Gamma$ is contained in a graphical (planar) end that
    never intersects $\{ x_3 = h \}$.  In particular, we must have
    that
    \begin{equation}
        \partial B_{r_n} \cap  \partial \Gamma_n  \subset \{ x_3 \geq \alpha \, r_n
    \} \, ,
    \end{equation}
    so that the length bound \eqr{e:fromareatolen} gives (bottom).
     By construction, each $\Gamma_n$ has compact boundary,
     is  contained in the slab
     $\{ 0 < x_3 < h \}$, and has quadratic curvature decay.  Therefore,
     the gradient estimate implies that each component of $\Gamma_n$ outside of a
     cylinder $\{ x_1^2 + x_2^2 \leq R^2 \}$
     must be either an asymptotically planar  graph or a multi-valued graph.
     However, as we have used several times, $\Gamma$ cannot contain such a multi-valued graph, so
     we conclude that each component of
     \begin{equation}
        \Gamma_n \cap \{ x_1^2 + x_2^2 > R^2 \}  {\text{ is an asymptotically planar
        graph.}}
     \end{equation}
      There are
     only finitely many such ends for each $n$ because $\Gamma_n$ is proper.
     This gives (ends) and, hence, completes the proof.
\end{proof}

\section{The proofs of Theorem \ref{t:tab} and Theorem \ref{t:t5.2a}}

We now have all of the necessary ingredients to prove Theorems
\ref{t:tab} and \ref{t:t5.2a}.

\begin{proof}
(of Theorem \ref{t:tab}). We have already established properties
(A) and (B) of Theorem \ref{t:tab} in Lemma \ref{l:singl}
 and Definition/Lemma \ref{l:inftyornot}, respectively.
 Therefore, it remains to show that
 every leaf of the
lamination $\cL'$ is contained in a horizontal plane.  Once we
have shown this, then the lamination $\cL$ is obtained by taking
the union of the horizontal planes in $\cL'$ together with a
horizontal plane through each point in $\cS$.

 We
have already proven that the leaves of $\cL'$ are planar when the
sequence is ULSC in Theorem \ref{t:t5.1}.  Therefore, by the no
mixing theorem, Theorem \ref{c:main}, the only remaining case is
when $\cS = \cSt \ne \emptyset$.
 However,
   Lemma \ref{l:smoothleaves} and Lemma
 \ref{l:nonsmoothleaves} together prove that every leaf of $\cL'$ is flat in this case. This completes
the proof of the theorem.
\end{proof}

Theorem \ref{t:t5.2a} now follows immediately:

\begin{proof}
(of Theorem  \ref{t:t5.2a}).  Now that we have established Theorem
\ref{t:tab},  it only remains to show that  property ($C_{neck}$)
holds. However, property ($C_{neck}$) was proven in (C1) in
Theorem \ref{t:t5.2}.
\end{proof}

\part{Modifications in the positive genus case}    \label{s:highergenus}

As we noted earlier, the main theorems were stated for sequences
of planar domains, i.e., for genus zero.   In this section, we
will give the versions of these theorems for sequences with
bounded genus  and describe the necessary modifications for the
proofs. The main change in the theorems is a change in the
definitions of the singular sets $\cSt$ and $\cSu$.  The new
definitions of $\cSt$ and $\cSu$, as well as an example showing
why a change is necessary, can be found in Subsection
\ref{ss:newdefs}.

   Many
aspects of the proofs in the genus zero case were essentially
local and will, therefore, extend easily once we have the
local structure near $\cSt$ and $\cSu$.  However, there are some
global aspects to the proofs and these will require some work. The
two main global facts are the existence of planes through each
singular point and the flatness of nearby leaves (which we often
call ``properness'').  These are  ``global''  in the sense that they fail to hold 
in the local example constructed in \cite{CM15}.

The definitions of $\cS$ and $\cL'$ are unchanged since Definition/Lemma \ref{l:inftyornot}
(that defines the singular set) and Lemma \ref{l:singl} (that constructs $\cL'$)
 did not assume genus zero.

\section{The definitions and  statements for positive genus}

\subsection{The sets $\cSt$ and $\cSu$ for positive genus}  \label{ss:newdefs}

 We will begin with an example
illustrating why we have to change the definitions of
$\cSt$ and $\cSu$ in the
case of positive genus. Namely, let the sequence
$\Sigma_j$ be   a sequence of rescalings (``blow downs'') of the
genus one helicoid constructed in \cite{HoWeWo}. Since the genus
one helicoid is asymptotic to the standard helicoid, the
$\Sigma_j$'s converge to a foliation by horizontal planes away
from the vertical axis. However, the vertical axis contains both
the origin where the injectivity radius goes to zero   - since the
genus concentrates there - and uniformly locally simply connected
points.   This was impossible in the case of genus zero because of
the no mixing theorem.

This example of rescalings of the genus one helicoid illustrates
that even if the injectivity radius goes to zero at a point, the
point still might not belong in $\cSt$.  It is then reasonable to
ask what it was about the injectivity radius going to zero that
was useful in the genus zero case.  The answer is that this
allowed one to use a barrier argument near a point $y \in \cSt$ to
find stable graphs disjoint from the $\Sigma_j$'s that converge to
a punctured plane through $y$.  This motivates the following
re-definition of $\cSu$ and $\cSt$:
\begin{itemize}
\item A point $y$ in $\cS$ is in $\cSu$ if there exist $r_y > 0$
and a sequence $r_{y,j} \to 0$ so that for any $r \in [r_{y,j} , r_y]$ and  {\underline{any}}
connected component $\Sigma_j'$  of $B_{r}(y) \cap \Sigma_j$ we have:
\begin{align}    
	&\partial \Sigma_j' {\text{ is connected}}. \\
    &  \Sigma_j' {\text{ has the same genus as {\underline{one}} of the components of }}
    B_{r_{y,j}} (y) \cap \Sigma_j'  \, . \label{e:genusstable}
\end{align}
\item  A point $y$ in $\cS$ is in $\cSt$ if there exist $r_y > 0$
and a sequence $r_{y,j} \to 0$ so that {\underline{some}}
component of $B_{r_{y,j} }(y) \cap \Sigma_j$ has dis-connected
boundary and \eqr{e:genusstable} holds.
\end{itemize}
Note that these definitions agree with the earlier ones when the
sequence is uniformly locally genus zero, i.e., when  the genus of
$B_{r_y}(y) \cap \Sigma_j$ is zero for every $j$.  In particular,
these definitions agree with the earlier ones when the sequence
$\Sigma_j$ has genus zero.

In the positive genus case, the  set $\cSt$  is divided into two subsets:
 \begin{itemize}
 \item
  A point $z \in \cSt$ is in $\cSta$ if the locally separating curves in $\Sigma_j$ that are shrinking to $z$ are either
  \begin{itemize}
  \item
       globally separating in $\Sigma_j$ (like in the genus 0 case) or, more generally,
     \item globally separating in $\Sigma_j$ once we combine them with at most $g$ other shrinking curves at other points of $\cSt$.
  \end{itemize}
\item  The set    $\cStb = \cSt \setminus \cSta$ consists of  at most $g$ ``exceptional points''  where this does not happen.
  \end{itemize}

 The sets $\cSt$ and $\cSu$ are
obviously disjoint subsets of $\cS$.    It follows from
  proposition I.0.19 in \cite{CM5} that, after passing to a subsequence, we
can assume that{\footnote{More precisely, this follows from the proof of 
proposition I.0.19 in \cite{CM5}; that proposition was stated 
for the complementary case where the inner radius is fixed and the outer radii go to infinity.}}
\begin{equation}    \label{e:onemoretime}
    \cS=  \cSt \cup \cSu \, .
\end{equation}
The fact that there are at most $g$ ``exceptional points'' follows immediately from Lemma \ref{l:gn}.

\subsection{The statements of the theorems for positive genus} \label{ss:newthms}

We will next run through the changes to the statements of the five theorems -  Theorem \ref{c:main},
Theorem \ref{t:tab}, Theorem \ref{t:t5.1}, Theorem \ref{t:t5.2a}, and Theorem \ref{t:t5.2} -
when the surfaces have positive genus.

The first theorem is the no-mixing theorem, Theorem \ref{c:main}; in the positive genus case, this becomes:

  \begin{Thm}	\label{t:nomixg}
(No-mixing theorem in the positive genus case).
If $\Sigma_{i} \subset B_{R_i} = B_{R_i} (0) \subset \RR^3$ is a sequence of compact embedded minimal surfaces of genus at most $g$ with $\partial \Sigma_{i} \subset \partial B_{R_i}$ where $R_{i} \to \infty$, then there is a subsequence so that $\cStb$ consists of at most $g$ points and either $\cS_{ulsc} = \emptyset$ or $\cSta = \emptyset$.

Moreover, if $\cS_{ulsc} \ne \emptyset$, then the lamination $\cL'$ given by Lemma $II.1.2$ consists of a foliation of (all of)  $\RR^3$ by parallel planes away from a singular set $\cS$ consisting of either one or two lines perpendicular to the planes together with at most $g$ points of $\cStb$.
  \end{Thm}

Theorem
\ref{t:tab}  applies verbatim to the general case of bounded
genus with the  new definitions of $\cSu$ and $\cSt$.

On the other hand,   Theorem \ref{t:t5.1} holds also for sequences
with fixed genus with one minor change in the conclusion and one
in the hypothesis.  The change in the hypothesis is that we do not
assume \eqr{e:notulsc}. The change in the conclusion is that there
might be either one or two singular curves.
 The assumption
\eqr{e:notulsc}, which says that the $\Sigma_j$'s are ``uniformly
not-disks'', was used in the genus zero case to rule out the
possibility of just one singular curve (as occurs both  for
sequences of disks and for rescalings of the genus one helicoid).
However, we cannot rule out the possibility of just one singular
curve in the fixed genus case regardless of whether we assume
\eqr{e:notulsc}.  For this reason, we will not assume
\eqr{e:notulsc} and we will allow for
 the possibility of
just one singular curve.

Similarly,
  Theorem \ref{t:t5.2a} and Theorem
\ref{t:t5.2} require  small changes in ($C_{neck}$). Recall that
for each point $y$ in $\cSt$,  ($C_{neck}$) gives a sequence of
graphs in the $\Sigma_j$'s that converges to a plane through $y$
  away from  at most two punctures.  In the positive genus
  case, there are now two types of points in $\cSt$ and the results are different for each:
  \begin{itemize}
  \item[($\cSta$)]  If $y \in \cSta$, then there is a sequence of
graphs in the $\Sigma_j$'s that converges to a plane through $y$
  away from  at most  $(g+2)$ punctures.
  \item[($\cStb$)]  If $y \in \cStb$, then there is a sequence of graphs or multi-valued graphs in the
  $\Sigma_j$'s that converges to a plane through $y$
  away from  at most  $(g+2)$ punctures.
  \end{itemize}

\vskip2mm {\it{The remainder of this part will be devoted to
sketching the modifications needed to prove the main theorems in
the general case of bounded genus.}}

\section{The local structure near points in $\cSu$ and $\cSt$}

The starting point for understanding the sequence $\Sigma_j$ is to
describe the sequence in a neighborhood of each singular point,
depending on whether the point is in  $\cSu$ or $\cSt$. Roughly
speaking, we will get the same picture as in the case of planar
domains. The precise statements are:
\begin{enumerate} \item[($\alpha$)] Given a point $x$ in $\cSu$, there is a
 ball $B_{r}(x)$ so that:
 \begin{enumerate}
 \item The $\Sigma_j$'s contain
   multi-valued
graphs that ``collapse'' to a punctured graph
in $B_{r}(x)$ with a removable singularity at $x$.
\item The set $\cS$ satisfies the cone property with respect to this graph in $B_r(x)$.
\item For $j$ sufficiently large, the $\Sigma_j$'s are connected near $x$.{\footnote{The precise statement is that there exists $C>1$ so that if $Cs< r$ and $j$ is sufficiently large, then there is only one connected component of $B_{Cs}(y) \cap \Sigma_j$ that intersects $B_s(y)$.}}
\end{enumerate}
\item[($\beta$)] Given a point $x$ in $\cSt$, there is a ball
$B_{r}(x)$ and a sequence of graphs in the $\Sigma_j$'s that
converges (with multiplicity at least two) to a finitely punctured graph in
$B_{r}(x)$ with a removable singularity at $x$.
\end{enumerate}

We will prove ($\alpha$)   first  and then ($\beta$).  Properties (a) and (b) in
($\alpha$) give the same structure that Lemma \ref{l:leaf} gave in the genus zero case.  The proof will follow the same outline as in the genus zero case, with modifications that are standard by this point in the series of papers.

\begin{proof}
(Sketch of proof of ($\alpha$)).
Suppose that  $y \in \cSu$.  We can assume that there is genus concentrating at $y$ (otherwise the genus $0$ argument applies).
Thus,  (by definition) there exist $r> 0$, 
 a sequence $r_{j} \to 0$, and points $y_j \in B_{r_{j}} (y) \cap \Sigma_j$  so that:
\begin{itemize}
\item  For any $s \in [r_j , r]$,  the
component $\Sigma_{s,y_j}$ of $B_{s}(y) \cap \Sigma_j$ containing $y_j$ has positive genus and has connected boundary.
\item If $r_j \leq s_1 < s_2 \leq r$ (and $s_1 , s_2$ are regular values), then $\Sigma_{s_2,y_j} \setminus \Sigma_{s_1,y_j}$ is a topological annulus with one inner boundary component in $\partial B_{s_1}(y)$ and one
    outer boundary component in $\partial B_{s_2}(y)$.
\end{itemize}

\noindent{\bf{The proof of (a)}}:  The first step in (a) is to prove the existence of small multi-valued graphs near $y$ in the $\Sigma_j$'s; when
the $\Sigma_j$'s  had genus zero, this was done in \cite{CM4} by identifying blow up pairs and working on the scale of the maximum of the curvature.  This approach does not work here because we do not have any a priori relationship between the radii $r_j$ and the maximum of the curvature near $y$ on $\Sigma_j$.  This difference is the biggest change in the extension of (a) to  the positive genus case.  Instead, we argue as follows:
\begin{enumerate}
\item The first observation is that, by theorem $1.22$ in \cite{CM4}, the area of intrinsic sectors over the inner boundaries in the annuli $\Sigma_{s_2,y_j} \setminus \Sigma_{s_1,y_j}$ must grow faster than quadratically.
\item Observe next that these annuli are scale-invariant simply-connected.  Namely, if $x \in \Sigma_{s_2,y_j} \setminus \Sigma_{s_1,y_j}$ has intrinsic distance $t> 2s_1$ to the inner boundary (and is also not too close to the outer boundary), then $\cB_{t/4}(x)$ is a topological disk.

To prove this, suppose instead that the exponential map from $x$ is not injective on $\cB_{t/4}(x)$; this would give two geodesics from $x$ with the same endpoint that combine to give a simple closed curve with two break-points.    Using the non-positive curvature and Gauss-Bonnet, we see that this curve cannot bound a disk and, thus, must be homologous to the inner boundary component.   However, using Stokes' theorem (applied to $\Delta |x-y|^2$), this would imply small area growth which is impossible by (1).

\item Using corollary $II.2.10$ in \cite{CM5}, we can now divide the intrinsic tubular neighborhood of the inner boundary into sectors whose sides are minimizing geodesics (in fact, even minimizing back to the entire inner boundary).  The bases of these sectors will be chosen to have a length comparable to a fixed large multiple of $s_1$.

\item By (1), if we choose $s_1$ so that $\frac{s_1}{r_j}$ is large enough, then we can make the number of disjoint sectors in (3) as large as we would like.

\item By (2) and the intrinsic version of the one-sided curvature estimate from \cite{CM9} (recorded here in Theorem \ref{t:t2}), if any two of these sectors are sufficiently (scale-invariant) close extrinsically, then they both satisfy a uniform scale-invariant curvature estimate.

\item Combining (4) and (5) with corollary $2.13$ of \cite{CM4}, we can arrange that at least one of these sectors is $1/2$-stable (with the width and length of the sector fixed, but as large as we wish).

\item Finally, (6) allows us to apply corollary $II.1.45$ of \cite{CM5}  to get  the desired multi-valued graph on a fixed scale.
\end{enumerate}

The extension of this multi-valued graph now follows from (a slight variation of) theorem $II.0.21$ in \cite{CM3} which showed that stable multi-valued graphs extend.  Stability there was used for two different purposes:  To get some a priori scale-invariant bound on curvature and then to come back and get a better (global) estimate leading to almost flatness.  The only notable difference in the current case is that
we only have the $1/2$-stability as long as the sheets stay close together; this is easily overcome by using the   sublinear growth of the separation (i.e., proposition $II.2.12$ in \cite{CM3}) to keep them together.

\noindent{\bf{The proof of (b)}}:
The second property that we need is the local cone property.  This follows immediately as in (5) from the intrinsic version of the one-sided curvature estimate from \cite{CM9} together with (2) above.

\noindent{\bf{The proof of (c)}}:
Finally, we need      the local connecting property. This follows immediately from (b) and a barrier argument. Namely, if there were multiple components, then we could use them as barriers to get a stable (thus very flat) surface between them.  Using the intrinsic one-sided curvature estimate and simple-connectivity of (2), the multi-valued graph forming in $\Sigma_j$ would then be forced to spiral graphically forever.  This is impossible since each surface is proper.
\end{proof}

\begin{proof}
(Sketch of proof of ($\beta$)).

The structure ($\beta$) near a point in $\cSt$ follows
from a local version of the results of Section \ref{s:ceno} for
the genus zero case.  As in the genus zero case, there are three main steps:
\begin{enumerate}
\item Using the (local) topology to put in a sequence of stable barrier surfaces that converge to a graph 
through the singular point; see Proposition \ref{l:getsetgo}.  This goes through as before, except that the
outer radii of the extrinsic
balls remains bounded. Hence, the limiting stable graph is defined
over a disk and not the entire plane.
\item Decomposing $\Sigma_j$ into ULSC pieces by cutting along ``small necks''.  This goes through 
as in Subsection \ref{ss:step1} with only obvious changes.

\item  Showing that these ULSC pieces contain graphs that converge to the limiting stable graph through the singular point.  This  goes through as in Subsection \ref{ss:step2} with obvious changes.

\end{enumerate}
\end{proof}

\subsection{Collapsed leaves}

The key properties (1), (2), and (3)  of collapsed leaves are recorded in
  Proposition \ref{p:cole0}.  We next extend the proofs of these properties to the positive genus case.
  The proof of (1) goes through as in Subsection \ref{ss:231} using the local structure ($\alpha$) above.  The proof of  (2) in
  Subsections \ref{ss:232}  and \ref{ss:233} goes through with minor modifications that are noted there (see, e.g., the second paragraph of  Lemma \ref{l:seeapp}).  The proof (3) in Subsection \ref{ss:nothree} goes through with the following minor changes:
  \begin{itemize}
  \item
  	In the proof of (3), we used that the ``figure eight'' curves $\gamma_j$ were separating in $\Sigma_j$.  In the genus $0$ case, this is automatic since all curves are separating.  When the genus is positive, note the infinite multiplicity of the convergence allows us to choose $g+1$ distinct graphs (all on different sheets)
	$\gamma_j^1 , \dots , \gamma_j^{g+1}$ that are embedded graphs over the curve $\gamma$ in $\Gamma$ and all of these have the same orientation (meaning all are on sheets where the  ``normal points upward'').     If none of these is separating on its own, then (since the genus is $g$ and they are all disjoint)  Lemma \ref{l:gn} gives a collection of them that together separate in $\Sigma_j$.  We will use this collection as the inner boundary in the Plateau problem and follow the rest of the argument.
	
	\item The other small change is that in the barrier construction, we apply \cite{HSi} rather than \cite{MeYa2}.  Thus, we do not get an explicit bound on the   topology of the stable surface, but this bound was never used in the argument.
  \end{itemize}

\section{Part \ref{p:prove1}: When the surfaces are ULSC}

This part completed the proof of
  Theorem \ref{t:t5.1} in the   genus zero case, using the tools already developed along with two new ingredients developed there:
  \begin{itemize}
  \item
Proposition \ref{p:cole} shows that the closure of a collapsed leaf is a plane.
\item
Lemma \ref{t:setup} proves ``properness''.
\end{itemize}

We follow the same approach in the positive genus case, with minor changes.  The first changes are in the statement of Proposition \ref{p:cole}, where:
\begin{enumerate}
\item We no longer assume that $\barga
\cap \cSt = \emptyset$ but instead make the weaker assumption that  $\barga
\cap \cSta = \emptyset$.  (Weakening this assumption is not necessary for the ULSC results in Part III, but will be needed later for the generalization of the no-mixing theorem to the positive genus case.)
\item We omit (2) since we make no assumption in the fixed genus case
to ensure that  there   are two axes.   As a result, we will need to also consider the case of ULSC, one axis, and finite genus; in the genus zero case this follows already from \cite{CM6}.
\end{enumerate}

Once we have these two things, then the modified Theorem \ref{t:t5.1} will follow as in the genus zero case with one last small change.  Namely, we can only apply Meeks' result, \cite{Me1}, at points in the traditional $\cSu$ (where the sequence is locally simply connected).  It follows that the singular set is (one or two) Lipschitz curves in $\cSu$ and these curves are orthogonal to the planar foliation at all but a finite collection of points; this of course implies that they are orthogonal everywhere.

\vskip2mm
{\bf{The proof of (1)}} in the (modified) Proposition \ref{p:cole} goes through with the following changes:
\begin{itemize}
\item Since there are at most $g$ points in $\cStb$, we get that $\barga \setminus \Gamma$ consists of at most two points in $\cSu$ together with at most $g$ of the ``exceptional points''.
\item
The points in $\cSu$ are already known to be removable singularities and (a cover of) $\Gamma$ is already known to be stable.  This stability 
together with ($\beta$) allows us to apply a local version of Lemma \ref{l:punstab} to conclude that the isolated exceptional singular points are also removable.{\footnote{This local version states that: Suppose that $\Gamma$ is a connected embedded minimal surface with trivial normal bundle, $\Gamma$ (or a cover) is stable, and $B_1 \cap \barga \setminus \Gamma = \{ 0 \}$.  Then $\Gamma$ has a removable singularity at $0$.  There are a number of ways to prove this, but perhaps the simplest is to deduce it from Lemma \ref{l:punstab} and a compactness theorem (the curvature must blow up at least quadratically if the singularity is not removable).}}
  The claim now follows from a Bernstein theorem as in the genus zero case.
\end{itemize}

\vskip2mm
{\bf{The proof of properness}} when the genus is zero was given in  Lemma \ref{t:setup} using a global flux argument.  We will describe the necessary modifications next.  Suppose first that the leaf has only one point of $\cSu$ in it.  As in Lemma \ref{t:setup}, we need to rule out the possibility of one leaf that spirals into the plane $\barga$.   We would like to appeal to corollary $0.7$ of \cite{CM7} as in  part $I$ of \cite{CM6} to get a contradiction, but we will need some modifications:
\begin{enumerate}
 \item We can find the ``short curves'' (required in \cite{CM7}) by using the multi-valued graph structure that we have already obtained together with corollary $III.3.5$ of \cite{CM5} to get blow up pairs converging down to the singular point from above and then following the argument in \cite{CM6} (see corollary $IV.0.10$ there).
\item  The leaf is not known to be locally graphical above the plane since we cannot directly apply the one-sided curvature estimate.  In particular, as we extend the sheets of  the multi-valued graph, we may come to an intrinsic ball that is not scale-invariant simply connected.  Because of the closeness to the plane and the sublinear growth proven in \cite{CM3}, we can take the scale-invariant constant to be very small.
However, there are at most $g$ of these ``bad balls''; otherwise, some combination of curves in these balls would be (globally) separating and we could put in a stable barrier that is forced to ``cut the axis'' near the singular point; this is a standard variation on the ``estimates between the sheets'' argument from \cite{CM3} that we have now employed a number of times.{\footnote{The original ``estimate between the sheets'' was proven in theorem $I.0.8$ in \cite{CM3}; the version that we use here is essentially (D) in the proof
of property (2) in Proposition \ref{p:cole}.  The difference   is that the stable surface may have up to $g+1$ inner boundary components and we use the existence theory of \cite{HSi} instead of Meeks-Yau.}}
  In particular, the ``bad balls'' (where it is not simply-connected) can be surrounded by ``good balls" and the sheets can be continued globally (with at most $g$ disks removed).

  \item The last modification is that we may need to start ``lower'' to ensure that we do not hit any of these ``bad balls'' as we extend the sheets of the multi-valued graph.  Since there are at most $g$ of these and the multi-valued graph has infinitely many sheets, this is not a problem.  The argument now goes through as in part $I$ on pages $584$ to $593$ of \cite{CM6}.
\end{enumerate}
When the leaf contains two   points in $\cSu$ (as was the case in Lemma \ref{t:setup} because of (2) in
 Proposition \ref{p:cole}), the modifications are similar.  Namely, the local picture near each singularity is identical and
 the leaf may fail to be locally graphical over the plane, but only at at most $g$ ``bad balls'' as in (2).    We use the argument in (1) to find the short curves and we argue as in (3) to work ``below''  these ``bad balls'' and then follow the proof of Lemma \ref{t:setup}.

\section{Parts  \ref{p:prove2} and \ref{p:nomix}: When the surfaces are not ULSC}

We will next turn to analyzing the structure of non-ULSC singular points, including the proofs - in the positive genus case  -  of
Theorems \ref{t:t5.2} and  \ref{c:main} (the no-mixing theorem).  To do this, we must prove:
\begin{itemize}
\item (C1) in Theorem \ref{t:t5.2}; this will follow from Proposition \ref{l:lplane} below.
\item (C2) and (D) in Theorem \ref{t:t5.2}.
\item Theorem \ref{c:main}.
\end{itemize}

    As in the genus zero case, a key point will be to prove that there
is a limit plane through each point in the singular set $\cS$.
These planes were actually (the closure of) leaves of $\cL'$ when
the sequence was ULSC, but this was not the case in general.
However, these planes were always given as smooth limits of
subsets of the $\Sigma_j$'s; cf. ($C_{neck}$) in Theorem
\ref{t:t5.2a}. This is the motivation for the following
definition:

\begin{Def}
Let $\Sigma_j$ be a sequence of surfaces with limit lamination
$\cL'$ and singular set $\cS$.  We will say that a surface
$\Gamma$ is a {\it{pseudo-leaf}} of $\cL'$ if it is
{\underline{connected}} and there is a sequence of subsets
$\Sigma^{\Gamma}_j \subset \Sigma_j$ that converges smoothly to
$\Gamma$.  We will also require that $\Gamma$ is maximal with
respect to these properties, so  that $\Gamma$ is not  a proper
subset of a connected surface that is also a limit of subsets of
the $\Sigma_j$'s.

Here ``converges'' means that for each open subset $\Gamma_c
\subset \Gamma$ with compact closure in $\Gamma$, then   the
$\Sigma_j$'s contain a sequence of graphs - or multi-valued graphs
- over $\Gamma_c$ and these converge smoothly to $\Gamma_c$.  If
we get multi-valued graphs, then we require that the number of
sheets goes to infinity as $j$ goes to infinity.
\end{Def}

Note that every leaf of $\cL'$ is also a pseudo-leaf. We have
already come across pseudo-leaves that may not be leaves. Namely,
($C_{neck}$) in Theorem \ref{t:t5.2a} implies that, for each point
$x$ in $\cSt$, we get a flat pseudo-leaf whose closure is a plane
through $x$.  This pseudo-leaf is a plane punctured at $x$ and
possibly at one other point.

One useful property of a pseudo-leaf is that none of the leaves of
$\cL'$ can intersect a pseudo-leaf transversely.  It then follows
from the local structure of nodal sets that the leaves of $\cL'$
cannot cross a pseudo-leaf.

 The key
point for generalizing the main results for non-ULSC sequences from genus zero to fixed
genus is to show that:

\begin{Pro} \label{l:lplane}
For each point $x \in \cS$, we get a flat pseudo-leaf whose
closure is a plane through $x$.  This pseudo-leaf is a plane punctured at at most 
  $g+2$   points; each puncture is in $\cS$.
\end{Pro}

We will need one more definition before proving  Proposition
\ref{l:lplane}. Recall that when we studied the leaves of $\cL'$,
we began with the collapsed leaves, i.e., the ones ``through'' a
point in $\cSu$. The collapsed leaves were shown to be stable and
to have removable singularities at points in $\cSt$. With this in
mind, we will say that a pseudo-leaf $\Gamma$ is {\it{pinched}} if
it goes ``through'' a point in $\cS$.  There are two local models
for the $\Sigma_j$'s near a point $x$ in $\cS$, depending on
whether $x \in \cSu$ or $x \in \cSt$.  First, if $x \in \cSu$,
then we know that there is a collapsed leaf of $\cL'$ through $x$;
see ($\alpha$). Second, if $x \in \cSt$, then it follows from ($\beta$) that there
is a pinched pseudo-leaf through $x$.

\subsection{The local structure ($\beta$)}

 We begin by recalling the local structure ($\beta$) near points of $\cSt$:
 \begin{enumerate}
 \item[($\beta$)] Given a point $x$ in $\cSt$, there is a ball
$B_{r}(x)$ and a sequence of graphs in the $\Sigma_j$'s that
converges (with multiplicity at least two) to a finitely punctured graph in
$B_{r}(x)$ with a removable singularity at $x$.
\end{enumerate}

\begin{Rem}	\label{r:globeta}
The structure above is forced to be local because the curves that are shrinking off may not be globally separating in the $\Sigma_j$'s.  However, if $y \in \cSta$, then we can argue as in ($\beta$) to solve a 
sequence of {\underline{global}} Plateau problems using the $\Sigma_j$'s as barriers to get a limiting plane $P_y$ through $y$ so that:
\begin{itemize}
\item $P_y$ is a smooth limit (of stable graphs disjoint from the $\Sigma_j$'s) away from at most $g+1$ points in $\cStb$.   In particular, $P_y$ does not cross any leaves (or pseudo-leaves).
\item Observe that $P_y$ cannot contain any points of $\cSu$. (If it did, then the multi-valued graphs that developed would be forced to spiral forever contradicting properness of the $\Sigma_j$'s.)
\item As in ($\beta$), we can cut the $\Sigma_j$'s along a collection of at most $g+1$ small necks to get graphs in the $\Sigma_j$'s that converge to $P_y$ away from at most $g+2$ points.   This follows as in 
steps (1) and (2) in Subsection \ref{ss:41}.
\item Finally, as in the genus zero case, we can do this both above and below the stable barriers.
\end{itemize}
\end{Rem}

{\it{From now on, we will assume that every pseudo-leaf is
oriented.  This slightly simplifies some of the arguments below
involving stability.  As we have seen several times, the
unoriented case can be dealt with by going to a double cover.  We
will leave the easy modifications needed for this case to the
reader.}}

\section{Planes through $\cSu$ and the proofs of (C2) and (D)}

Suppose now that $x \in \cSu$ and $\Gamma$ is the collapsed leaf through it.   Generalizing
Proposition \ref{l:ptw}, we first show that:
  \begin{enumerate}
 \item[(0)]
$\cSta$ does not intersect the closure of the collapsed leaf $\Gamma$.
 \end{enumerate}
 
 The keys for showing this are the structure result ($\alpha$) and the following lemma (that generalizes \eqr{e:gipz}):
 
 \begin{Lem}	\label{l:py}
 If $\Gamma$ is a pinched pseudo-leaf (e.g., a collapsed leaf), $y$ is a point in $\cSta$,   $P_y$ is the plane through $y$ given by Remark \ref{r:globeta},  and $\barga \cap P_y \ne \emptyset$, then 
   \begin{equation}    \label{e:newold}
      \Gamma \subset P_y \, .
 \end{equation}  
 \end{Lem}
 
 \begin{proof}
 Because of embeddedness of the sequence, it is not hard to see that $\Gamma$ and $P_y$ cannot cross.{\footnote{In the genus zero case, this is Lemma \ref{l:sm1}; the lemma extends easily to the finite genus case.}}  We will argue by contradiction and, thus, assume that $\Gamma$ is above $P_y$.
 
 The key point will be the following claim:\\
 {\bf{Claim}}:   $\barga \cap \cSt$  is a finite collection of points.
 
 \vskip2mm
 {\bf{Proof of claim}}:   (This is proven by a modification of the proof of Lemma \ref{l:strbarcS}.)
  
 Since $\Gamma$ cannot cross any of the separating planes through the points in $\cSta$, it follows that
 $\barga \cap \cSta$ is contained in $P_y$ together with at most one other plane parallel to $P_y$ (and above it) and $\Gamma$ is contained either in the half-space above $P_y$ or in the slab between the two planes.
 We will  show that $P_y$ contains finitely many points in $\barga \cap \cSta$; the claim follows from this (together with a similar argument for the second plane in the case of two planes).    
 
 We already know that $P_y$ is a smooth limit of stable graphs that are disjoint from the
  $\Sigma_j$'s away from at most $g+1$ points and that there are at most $g$ points where the genus is concentrating; let $G$ denote these (at most $2g+1$) ``bad points''.   
  
  We will prove the claim by showing that $P_y \cap \barga \cap \cSta$ cannot contain  $4g+3$ distinct points.  
  Namely, if it did, then 
  $P_y \cap \barga \cap \cSta$ contains two collections 
  \begin{equation}
  	\{ y_1 , \dots , y_{g+1} \}  {\text{ and }} \{ z_1 , \dots , z_{g+1} \}   
\end{equation}
so that all these points are distinct and disjoint from (the at most $2g+1$ points in) $G$.
  
  It follows from the one-sided curvature estimate that the injectivity radius of the $\Sigma_j$'s is going to zero at the points in $\Sigma_j$ that are converging to $\Gamma$ near the $y_i$'s and $z_i$'s.    Thus, since the genus is at most $g$, we can choose subcollections of each collection so that they separate in the $\Sigma_j$'s and we can put in stable barriers (using the local version of Proposition \ref{l:getsetgo}).  This leads to a contradiction as in the end of the proof of Lemma \ref{l:strbarcS}: Namely, we get two distinct stable barriers that separate in space, thus they must be ordered by height, but the limiting surface somehow goes ``through'' both of them.  This contradiction proves the claim.
  
  Once we have shown that $\Gamma$ is complete away from isolated points on the boundary of a half-space and $\Gamma$ is contained in this half-space, then (a local version of)  Lemma \ref{l:punstab} implies that the isolated exceptional singular points are also removable.  The strong maximum principle then gives \eqr{e:newold}.
   \end{proof}

Here is why ($\alpha$) and Lemma \ref{l:py} imply (0):
 \begin{quote}
 
Suppose that (0) fails and $\barga$ contains $y \in \cSta$.  Let $P_y$ be the limiting plane through $y$ given by Remark \ref{r:globeta}, so that Lemma \ref{l:py} implies that $\Gamma \subset P_y$.

  We get the contradiction from using the barrier graphs to separate the sheets (cf. Remark \ref{r:globeta})  which is impossible because of the local connecting property  near $\cSu$ given in part (c)  of  ($\alpha$).  
 \end{quote}

Using (0), we can now apply the modified (1) from Proposition \ref{p:cole} to get that
 the collapsed leaves are all punctured planes and we can apply the modified  Lemma \ref{t:setup}
 to get that a neighborhood of each point in $\cS_{ulsc}$ is foliated by collapsed leaves.  Thus,    
  (the modified) (C2)  and (D) hold.

\section{The remaining cases of Proposition
\ref{l:lplane}}

\begin{proof}
(Sketch of proof of
 Proposition \ref{l:lplane}). 
 Suppose that $x \in \cS$.  We have already dealt with the cases where $x$ is in $\cSu$ or $\cSta$, so we may assume that $x \in \cStb$.  Let $\Gamma$ be the pinched pseudo-leaf through $x$ guaranteed by ($\beta$).
  We will show next that $\Gamma_x$ is flat.  This  follows from
stability when  $\Gamma_x$ is complete  or if it has only isolated removable singularities (by the usual logarithmic cut-off argument).   We will divide into several cases:
\begin{itemize}
\item  Suppose that $\Gamma_x$ contains a point of $\cSu$ in its closure.  Since we have already shown that a neighborhood of each point in $\cSu$ is foliated by flat leaves, we conclude that $\Gamma_x$ is contained in one of these flat leaves and is, thus, itself a punctured plane.  (The bound on the number of punctures has also already been established.)
\item Suppose next that $\Gamma_x$ contains a point $y$ of $\cSta$ in its closure; let $P_y$ be the corresponding plane through $y$.  It follows from Lemma \ref{l:py}   that $\Gamma \subset P_y$, giving the desired flatness.
\item Finally, suppose that $\Gamma_x$ contains only points in $\cStb$ in its closure.  Since there are at most $g$ of these points and each is a removable singularity (by stability, the structure ($\beta$), and
 the usual slight variation on Lemma
 \ref{l:punstab} in Appendix \ref{s:apB}), we  can apply the Bernstein type argument to get flatness.
\end{itemize}
This completes the proof.
\end{proof}

\subsection{Part \ref{p:nomix}: The no-mixing theorem}

We will now combine the previous results to extend the no-mixing theorem, Theorem \ref{c:main}, to the positive genus case; that is, we will prove Theorem \ref{t:nomixg}.  We must show that
\begin{itemize}
\item If $\cSu \ne \emptyset$, then the planar collapsed leaves (through $\cSu$)  foliate all of $\RR^3$.
\end{itemize}
We have already shown that  the foliated region consists  of an open set of planes and   $\cSu$ is either one or two straight lines perpendicular to these planes.  We must rule out that one of these lines has an endpoint.  However, the singular set is closed so this endpoint would have to be in $\cSt$.  Thus ($\beta$) gives graphical stable barriers near the endpoint (for $j$ large) which force the spiralling in $\Sigma_j$ (from the nearby $\cSu$ points) to continue forever.  This contradicts that the $\Sigma_j$'s are proper.

\subsection{Part \ref{p:prove3}: The leaves are all flat}

The other two global flux arguments are used to show that the
leaves are flat in the non-ULSC case, i.e., when $\cS = \cSt$.
This is divided into two cases, depending on whether or not the
leaf $\Gamma$ is complete.  
The complete leaves were shown to be flat in
Lemma \ref{l:smoothleaves} and the in-complete leaves were handled
in Lemma \ref{l:nonsmoothleaves}.   We will next explain how to extend the proofs of these to the positive genus case.

{\bf{$\Gamma$ is complete: Lemma \ref{l:smoothleaves}}}:  The point is that $\Gamma$ must lie in a half-space (since it cannot cross any of the limit planes through $\cS$) and, after a translation and a rotation,
we may
    assume that $\Gamma \subset \{ x_3 \geq 0 \}$ and
    \begin{equation}    \label{e:mayassupg}
        \inf_{\Gamma} \, x_3 = 0 \, .
    \end{equation}
Arguing as in the claim after \eqr{e:mayassu} (using the intrinsic version of the
one-sided curvature estimate), we get a sequence of points
 $p_n \in \Gamma$ satisfying:
\begin{align}
     \label{e:pnti3pg}
                \cali (p_n) &\to 0 \, , \\
\label{e:pnti3apg}
                x_3 (p_n) &\to 0 \, ,
\end{align}
where $\cali (p_n)$ is the
        injectivity radius  of $\Gamma$ at $p_n$.  Thus far, there is no difference in the positive genus case.

        The contradiction comes from   cutting $\Gamma$ along these ``short curves'' to get that the flux of $\Gamma$ is arbitrarily small, which contradicts the strict positivity of the flux ``at the top''  that comes from slicing $\Gamma$ by a plane between two of its ends.  This is carried out in steps (a) through (h) of the proof of 
Lemma \ref{l:smoothleaves}.  In the positive genus case, we need the following modifications:
\begin{enumerate}
\item[(a)] No changes.
\item[(b)]  This is where we find the separating curves; this comes almost for free in the genus zero case just because the injectivity radius is going to zero at the $p_n$'s.  When the genus is at most $g$, then at most $g$ of the balls centered at the $p_n$'s can have positive genus{\footnote{We are using that we can assume that all the balls $\cB_{5C_1 \, i(p_n)}(p_n)$ are disjoint; cf. \eqr{e:pnti2}.}}, so we throw these out; this allows us to apply the one-sided lemma for non-simply
            connected surfaces, Lemma \ref{l:os}, on the remaining balls.   There is still another difficulty; namely, the curves in the balls $\cB_{5C_1 \, i(p_n)}(p_n)$  are locally separating, but they may not be globally separating.  To deal with this, we group the $p_n$'s together with $(g+1)$ of them in each group.  We know that some subcollection of each group must be globally separating.   This requires obvious changes when we introduce the stable barriers (as we have now done many times).
            \item[(c),(c')] No changes.
\item[(d)]  The first part of this is just repeating (b), with the same changes.  The second part is to get the properties in \eqr{e:propen}.  This follows without change because we only work on the balls that are genus zero and, thus, can still apply Lemma 
\ref{l:shortc}.

\item[(e)] No changes.
\item[(f)] Here we  use the one-sided curvature estimate and a decomposition into ULSC pieces to show that ends of $\Gamma$ (above where we cut) are graphical.   This is dealt with exactly as in the decomposition around necks.  Namely, this can only fail on at most $g$ ``bad balls'' each of which connects to a finite number of sheets and each bad ball can be surrounded with graphical pieces.
\item[(g)]  We choose the slicing plane below all of the ``bad sheets'' from (f).
\item[(h)] This is where the flux contradiction comes in.  The only difference is that instead of one ``bottom curve'' there may be $(g+1)$ bottom curves.
\end{enumerate}

{\bf{$\Gamma$ is not complete: Lemma \ref{l:nonsmoothleaves}}}:
 Suppose instead that $\Gamma \subset \{ x_3
> 0 \}$ is a non-flat leaf of $\cL'$ with $0 \in \barga \cap \cSt$
and $\{ x_ 3 = 0 \}$ is the associated stable limit plane through
$0$.   As in the genus zero case, the argument for incomplete leaves uses short curves to get a flux contradiction.  The issue is the construction of the  ``bottom curves''  which required that the injectivity radius was small relative to the distance to the boundary (cf. Remark \ref{r:notcomplete}).
\begin{itemize}
\item
We first modify Lemma \ref{l:strbarcS}  for the positive genus case to get that $\barga \cap \cS$ consists of at most $3g+1$ points and
 there are two possibilities:
\begin{itemize}
\item $\barga \cap \cS  \subset \{ x_3 = 0 \}$. 
\item  $\barga \cap \cS \subset \{ x_3 = 0 \} \cup \{ x_3 = x_3 (p) > 0\}$ for some $p$.
\end{itemize}
(The proof of this modification follows the original proof of Lemma \ref{l:strbarcS}  with the obvious modifications that we throw away the (at most $g$) points where the  genus is concentrating and then we need to work with two collections of $g+1$ points in order to guarantee that they separate globally.)
The important point is that $\Gamma$ fails to be complete only at isolated points.
\item Next note that 
$\Gamma$ is scale-invariant ULSC
near each singular point (cf. \eqr{e:notcase1}); 
this follows as before, except that we may need to throw away $g$ bad balls and work below these.
\item Next we show that the curvature blows up at most quadratically near each singular point; 
see \eqr{e:qcdcd}.  We will argue by contradiction, so suppose instead that
there is a
sequence
    of points $q_n \in \Gamma$ with $q_n \to 0$ and
    \begin{equation}    \label{e:qncfg}
        |q_n|^2 \, |A|^2 (q_n) > n \, .
    \end{equation}
    Thus, the sequence of dilated and translated laminations
    \begin{equation}
        \cL_n = |q_n|^{-1} \, \left( \cL' - q_n \right)
    \end{equation}
    satisfies
   $|A|^2 (0) > n$ and   the point $\frac{-q_n}{|q_n|} \in \partial B_1$
    is in $\cSt (\cL_n)$  (where $\cSt (\cL_n)$ is the non-ULSC singular set for the
    rescaled lamination $\cL_n$).  Since $\cL'$ is a limit lamination, we can apply these rescalings to a subsequence of the original sequence and use a diagonal argument to get the $\cL_n$'s to converge to a limit $\cL_{\infty}$.  It follows that $\cS_{\infty} = \cS (\cL_{\infty})$ has a ULSC singularity at $0$ and 
    the points $\frac{-q_n}{|q_n|}$ converge to a point $q \in \partial B_1$ that is  
    a
    non-ULSC singularity for $\cL_{\infty}$.   
    When the genus is zero, this violates the no-mixing theorem giving the desired contradiction.  In the present case of positive genus, this alone is not enough.   However, observe that the plane $x_3 = 0 $ is the limit of punctured graphs in the $\Sigma_j$'s and {\underline{not}} multi-valued graphs;  this is because $\cSu = \emptyset$.    It follows that the horizontal planes through  $\frac{-q_n}{|q_n|} $ are also   limits of graphs (in the dilated and translated $\Sigma_j$'s).    From this, we conclude that the plane $\{ x_3 = x_3 (q) \}$ is also a limit of graphs and, thus, that $\cSu_{\infty} = \cSu (\cL_{\infty})$ does not intersect this plane.  However, this violates the generalized no-mixing theorem, Theorem \ref{t:nomixg} (which gives that once $\cSu \ne \emptyset$, then $\cSu$ contains a line that intersects every one of the limiting planes), giving the desired contradiction.
    \item Once we have the quadratic curvature bounds, the rest of the proof follows as in the genus zero case.
    \end{itemize}

\appendix

\part{Appendices}

\section{Compact connected exhaustions} \label{s:a}

We recall next that any connected surface $\Gamma$ without
boundary, but {\underline{not}} necessarily complete, may be
exhausted by
 connected open sets with compact closure.  Just to be clear,
 the punctured plane $\{ (x,y) \, | \, x^2 + y^2 > 0 \}$ is an
 example of such a surface.
 This lemma
  is elementary, but  we include a
proof below for completeness.

\begin{Lem}     \label{l:simpap}
Given a  connected surface $\Gamma$ without boundary, there exists
a sequence of connected open sets $K_j$ with compact closure that
exhaust  $\Gamma$. That is, we have $\Gamma = \cup_{j=1}^{\infty}
K_j$  and $K_{j} \subset K_{j+1}$ for every $j$.
\end{Lem}

\begin{proof}
To see this, first fix a point $x \in \Gamma$.  For each $r
> 0$, let $\Gamma_r$ be the set of points $y$ in $\Gamma$ such
that the geodesic ball of radius $1/r$ about $y$ is complete. Next
define a {\underline{compact}} subset $K_r^{\comp} \subset
\Gamma_r$ to be the set of points $z \in \Gamma$ such that there
is a Lipschitz map
\begin{equation}
    \gamma_z : [0,r] \to \Gamma_r {\text{ with }} \gamma_z (0) = x
    \, , \, \gamma_z (r) = z \, , \, {\text{ and }} |\gamma_z'|
    \leq 1 \, .
\end{equation}
It follows immediately that each $K_r^{\comp}$ is connected.
Moreover, since $\Gamma$ is locally compact, the Arzela-Ascoli
theorem implies that each $K_r^{\comp}$ is compact.

We will show that the $K_r^{\comp}$'s exhaust $\Gamma$.  Suppose
that $z \in \Gamma$ is a fixed but arbitrary point; we will show
that $z \in K_r^{\comp}$ for some sufficiently large $r$.  Observe
first that $\Gamma$ is path connected since it is connected and
locally path connected; we therefore get a continuous map
\begin{equation}
    f_z : [0,1] \to \Gamma  {\text{ with }} f_z (0) = x
    {\text{ and }} f_z (1) = z \, .
\end{equation}
Since $[0,1]$ is compact and $f_z$ is continuous, the image $f_z
([0,1]) \subset \Gamma$ is also compact.  In particular, there
exists $r'$ so that
\begin{equation}
    f_z ([0,1]) \subset \Gamma_{r'} \, .
\end{equation}
Finally, we will replace the curve $f_z ([0,1])$ with a broken
geodesic to get a Lipschitz curve in $\Gamma_{2r'}$ from $x$ to
$z$.  Namely, compactness allows us to cover $f_z ([0,1])$ by a
finite collection of balls of radius $1/(4r')$ with centers on the
curve.  Replacing segments in $f_z ([0,1])$ connecting the centers
of overlapping balls by intrinsic geodesics gives a broken
geodesic -- with finitely many breaks  -- connecting $x$ to $z$.
The triangle inequality guarantees that this broken geodesic stays
in $\Gamma_{2r'}$.  This completes the proof that $z$ is in some
$K_r^{\comp}$ since the length of this curve is finite.

Finally, we define the {\underline{open}} sets $K_j$ by
\begin{equation}
    K_j = \cup_{r < j} K_r^{\comp} \, .
\end{equation}
The $K_j$'s are obviously nested and exhaust $\Gamma$ since the
$K_j^{\comp}$'s do.  The set $K_j$ is path connected since it is
the union of nested path connected sets; it is contained in the
compact set $K_j^{\comp}$ and therefore has compact closure.
Finally, it is easy to see that each $K_j$ is open.
\end{proof}

\section{Surfaces with stable covers}  \label{s:apB}

\subsection{Going from stability of a covering space to stability of a surface itself}

If an oriented minimal surface is stable, then any covering space
is also stable.  However, the converse may not always be true. The
next lemma states that the converse {\underline{is}
{\underline{true} if in addition the holonomy group of the
covering space has sub-exponential  growth.

Before showing this, we will need to recall  a   few elementary
properties of groups and covering spaces.

\vskip2mm {\underline{Growth of groups}}. Suppose  that $\Lambda$
is a finitely generated group and fix a set of generators.  Such a
choice of generators induces a natural metric on $\Lambda$ called
the word metric, cf. \cite{Gr}. Let $\Lambda_n$ denote the ball of
radius $n$ about the identity in this metric. The group is said to
have sub-exponential growth if we have for every $\epsilon
> 0$ that
\begin{equation}    \label{e:dds}
   \lim_{n \to \infty} \,  \frac{ |\Lambda_n| }{ \e^{\epsilon \, n} } = 0
\end{equation}
where $|\Lambda_n|$ denotes the number of elements of
$\Lambda_n$.{\footnote{It is not hard to see that having
sub-exponential growth is independent of the choice of
generators.}} Given any fixed integer $k$, it follows, almost
immediately, that sub-exponential growth guarantees that
 there is a sequence $n_j \to
\infty$ with
\begin{equation}    \label{e:segra}
    \frac{ |\Lambda_{n_j+k} \setminus \Lambda_{n_j}| }{ |\Lambda_{n_j}|} \to 0 \,
    .
\end{equation}

\vskip2mm {\underline{Covering spaces}}. Recall that a connected
covering space $\hat{\Pi}: \hat{\Gamma} \to \Gamma$ with base
point $x \in \Gamma$ is uniquely determined by the {\it holonomy
homomorphism} $\Hol$ from $\pi_1 (\Gamma)$ to the automorphisms of
the fiber $\hat{\Pi}^{-1} (x)$.  To define this homomorphism,
suppose that $\gamma:[0,1] \to \Gamma$ is a curve with $\gamma (0)
= \gamma (1) =  x$ and    $\hat{x}$ is a point in $\hat{\Pi}^{-1}
(x)$.  The
  lifting property for covering spaces gives a unique lift
  $\gamma_{\hat{x}}:[0,1] \to \hat{\Gamma}$ of $\gamma$ with
  $\gamma_{\hat{x}}(0) = \hat{x}$.  We define $\Hol (\gamma)
  (\hat{x}) $ to be the endpoint  $\gamma_{\hat{x}}(1)$.

  We call the image $\Hol (\pi_1 (\Gamma))$ the holonomy group of
  the covering space; to keep the notation simple, set $\Lambda =
  \Hol (\pi_1 (\Gamma))$.{\footnote{If $\hat{\Gamma}$ is the
  universal cover, then the holonomy group is exactly the group of
  deck transformations and, hence, isomorphic to $\pi_1 (\Gamma)$.
   However, the deck group   acts transitively on the fiber
   $\Pi^{-1} (x)$ if and only if $\pi_{\hat{\Gamma}}$ is a normal
   subgroup of $\pi_1 (\Gamma)$; when this is not the case, the
   holonomy group is bigger than the deck group.}}
If we fix a point $\hat{x}$ with $\Pi (\hat{x}) = x$, then we can
define a fundamental domain $\Gamma_0$ in $\hat{\Gamma}$ by
\begin{equation}
    \Gamma_0 = \{ y \in \hat{\Gamma} \, | \, \dist_{\hat{\Gamma}}
    (y,\hat{x}) \leq \dist_{\hat{\Gamma}}
    (y, z ) \, {\text{ for all }} z \in \Pi^{-1} (x) \} \, .
\end{equation}
Using this, let $\hat{\Gamma}_n = \cup_{\lambda \in \Lambda_n}
\lambda (\Gamma_0)$ be the covering of $\Gamma$ corresponding to
$\Lambda_n$.

The next property that we will need is a positive lower bound for
the distance between $\hat{\Gamma}_n$ and $\hat{\Gamma} \setminus
\hat{\Gamma}_{n+k_0}$ when $k_0$ is sufficiently large. Precisely,
if $\Gamma$ has compact closure (so, in particular, $\pi_1
(\Gamma)$ is finitely generated and $\diam (\Gamma_0)$ is finite),
then an easy compactness argument gives a constant $k_0$ so that
\begin{equation}    \label{e:kzdz}
    \dist_{\hat{\Gamma}}  ( \hat{\Gamma}_n ,  \hat{\Gamma} \setminus \hat{\Gamma}_{n+k_0}
    ) > 1  \, .
\end{equation}
Here $k_0$   depends on $\hat{\Gamma}$, $\Gamma$, and $\Lambda$
but does not depend on $n$.

The last fact that we will need is that the holonomy group extends
to an action on $\hat{\Gamma}$ when it is abelian; we include a
proof for completeness.

\begin{Lem}     \label{l:simcov}
If $\hat{\Gamma} \to \Gamma$ is a connected covering space with
{\underline{abelian}} holonomy group $\Lambda$, then $\Lambda$
extends to an action on $\hat{\Gamma}$ as the group of deck
transformations as follows:

Suppose that $\gamma:[0,1] \to \Gamma$ is a curve with $\gamma (0)
= \gamma (1) =  x$ (where $x$ is the base point in $\Gamma$). We
have to define the action of $\Hol (\gamma)$ on an arbitrary point
$\hat{y}$ in $\hat{\Gamma}$. To do this, choose a curve $\sigma:
[0,1] \to \Gamma$ from $y$ to $x$ and define $\Hol (\gamma)
(\hat{y})$ to be the second endpoint of the curve starting at
$\hat{y}$ that lifts the curve
\begin{equation}    \label{e:simcov}
    (-\sigma) \circ \gamma \circ \sigma \, ,
\end{equation}
where $(-\sigma)$ denotes the curve $\sigma$ traversed in the
opposite direction.
\end{Lem}

\begin{proof}
The only thing to check is that this definition does not depend on
the choice of the curve $\sigma$.  Suppose therefore that $\mu :
[0,1] \to \Gamma$ is a second curve from $y$ to $x$.  It is then
easy to see that $\sigma$ and $\mu$ give the same endpoint in
\eqr{e:simcov} if and only if
\begin{equation}    \label{e:simcov2}
    \mu \circ (-\sigma) \circ \gamma \circ \sigma \circ (-\mu) \circ (-\gamma)
\end{equation}
lifts to a closed curve in $\hat{\Gamma}$ starting at $\hat{x}$.
However, the second endpoint of the curve in \eqr{e:simcov2} is
nothing more than
\begin{align}
    &\Hol \, ( \mu \circ (-\sigma) \circ \gamma \circ \sigma \circ (-\mu) \circ
    (-\gamma)) (\hat{x})  = \notag \\
    &\quad \quad \quad \quad \Hol \, ( \mu \circ (-\sigma))  \circ \Hol (\gamma )
    \circ \Hol ( \sigma \circ (-\mu))   \circ
     \Hol (-\gamma) (\hat{x})
     = \hat{x} \, ,
\end{align}
since the holonomy group is abelian. This completes the proof.
\end{proof}

\vskip2mm {\underline{Stability of covering spaces}}. We will next
show that if a cover of a minimal surface is stable and its
holonomy group has sub-exponential growth, then the surface itself
is stable. This would be obvious for finite covers; in that case,
any compactly supported function on $\Gamma$ lifts to a compactly
supported function on $\hat{\Gamma}$. When the holonomy group is
infinite, the lift of a compactly supported function on $\Gamma$
no longer has compact support.  To deal with this, we have to
introduce a second cutoff function.

\begin{Lem}     \label{l:punstabIII}
Suppose that $\Gamma \subset \RR^3$ is an oriented minimal surface
with compact closure, possibly with boundary, and $\hat{\Gamma}$
is a covering space of $\Gamma$.  If $\hat{\Gamma}$ is stable and
its holonomy group has sub-exponential growth, then $\Gamma$
itself is stable.
\end{Lem}

\begin{proof}
 We will show that, for each
function $0 \leq \phi \leq 1$ compactly supported on $\Gamma
\setminus \partial \Gamma$, we have the following stability
inequality
\begin{equation}    \label{e:goal}
    \int_{\Gamma} |A|^2 \phi^2 \leq \int_{\Gamma} |\nabla \phi|^2
    \, .
\end{equation}
  Since the holonomy group $\Lambda$   of  the covering space  has sub-exponential growth,
\eqr{e:segra} gives   a sequence $n_j \to \infty$ with
\begin{equation}    \label{e:segr}
    \frac{ |\Lambda_{n_j+k_0} \setminus \Lambda_{n_j}| }{ |\Lambda_{n_j}|} \to 0 \,
    ,
\end{equation}
where $k_0$ is given by \eqr{e:kzdz}.

Define a sequence of functions $\psi_j$ on $\hat{\Gamma}$ by
\begin{equation}
 \psi_j =
\begin{cases}
  1 &\hbox{ on } \hat{\Gamma}_{n_j} \, ,    \\
  1 -\dist_{\hat{\Gamma}} (\hat{\Gamma}_{n_j} , \, \cdot)    &\hbox{ on }
   \{ 0 <   \dist_{\hat{\Gamma}} (\hat{\Gamma}_{n_j} \, , \, \cdot)  < 1   \}    \, ,\\
  0 &\hbox{ otherwise }  \, . \\
\end{cases}
\end{equation}
In particular, $\psi_j$ is one on $\hat{\Gamma}_{n_j}$, zero
outside the $1$--tubular neighborhood of $\hat{\Gamma}_{n_j}$ and
hence zero outside $\hat{\Gamma}_{n_j + k_0}$ by \eqr{e:kzdz}.
Moreover, $\psi_j$ decays linearly in the distance to
 $\hat{\Gamma}_{n_j}$ and hence satisfies
 \begin{equation}
    |\nabla \psi_j| \leq   1 \, .
\end{equation}
Below, we will  identify the functions $\phi$ and $|A|^2$ on
$\Gamma$ with their lifts to the cover $\hat{\Gamma}$.

Although the function $\psi_j$ does not vanish on all of $\partial
\hat{\Gamma}$, the function $\psi_j \, \phi$ does.
  We can therefore use $\psi_j \, \phi$ in the stability inequality for
$\hat{\Gamma}$ to get
\begin{align}
    |\Lambda_{n_j}| \, \int_{\Gamma} |A|^2 \phi^2 &=
    \int_{\hat{\Gamma}_{n_j}} |A|^2 \phi^2  \leq  \int_{\hat{\Gamma}_{n_j +
    k_0}} |A|^2 \, (\psi_j \, \phi)^2 \notag \\
        &\leq \int_{\hat{\Gamma}_{n_j +
    k_0}}   |\nabla (\psi_j \, \phi)|^2
    = \int_{\hat{\Gamma}_{n_j}} |\nabla \phi|^2
    +  \int_{\hat{\Gamma}_{n_j + k_0} \setminus  \hat{\Gamma}_{n_j}
    }   |\nabla (\psi_j \, \phi)|^2  \notag \\
  &= |\Lambda_{n_j}| \, \int_{\Gamma} |\nabla \phi|^2
    +  \int_{\hat{\Gamma}_{n_j + k_0} \setminus  \hat{\Gamma}_{n_j}
    }   |\nabla (\psi_j \, \phi)|^2
     \, .  \label{e:goal2}
\end{align}
  Since $\phi$ is smooth
and has compact support,  there is a constant $C_{\phi}$ so that
$2 \, |\nabla \phi|^2 +2 \leq C_{\phi}$; hence
\begin{equation}    \label{e:thegrb}
     |\nabla (\psi_j \, \phi)|^2 \leq  2 \, ( |\nabla \phi|^2 + |\nabla \psi_j|^2 )
    \leq  2 \,  |\nabla \phi|^2  + 2
     \leq
    C_{\phi}\, .
\end{equation}
We can use this to bound the last term in \eqr{e:goal2} as follows
\begin{equation}    \label{e:thelastt}
\int_{\hat{\Gamma}_{n_j + k_0} \setminus  \hat{\Gamma}_{n_j}
    }   |\nabla (\psi_j \, \phi)|^2 \leq C_{\phi} \, \Area
        ( \hat{\Gamma}_{n_j + k_0} \setminus  \hat{\Gamma}_{n_j} )
        = C_{\phi} \, \Area (\Gamma) \, |\Lambda_{n_j + k_0} \setminus  \Lambda_{n_j
    }| \, .
\end{equation}
Substituting  \eqr{e:thelastt} into \eqr{e:goal2} gives
\begin{equation}    \label{e:goal3}
    \int_{\Gamma} |A|^2 \phi^2 \leq \int_{\Gamma} |\nabla \phi|^2
    + C_{\phi} \, \Area (\Gamma) \, \frac{ |\Lambda_{n_j + k_0} \setminus  \Lambda_{n_j
    }| }{|\Lambda_{n_j }| }
    \, .
\end{equation}
Finally, \eqr{e:segr} implies that \eqr{e:goal3} goes to
\eqr{e:goal} as $j \to \infty$, completing the proof.
\end{proof}

\subsection{A surface and stable cover with cyclic holonomy group where
the previous lemma applies}

We will show, in Corollary \ref{c:collstab} below,   that a
certain minimal surface $\Gamma$ given as a limit of embedded
minimal multi-valued graphs $\Sigma_j$ must be stable. This will
follow from Lemma \ref{l:punstabIII}   once we show that there is
a connected covering space $\hat{\Gamma}$ satisfying the following
two properties:
\begin{itemize} \item The holonomy group $\Lambda$ of the covering
space is cyclic (and, hence, has sub-exponential growth). \item
The cover $\hat{\Gamma}$ is stable.
\end{itemize}
Throughout this subsection,
 $\Gamma \subset \RR^3$ will be an oriented minimal surface with compact closure, possibly
with boundary, and  $\Pi : \hat{\Gamma} \to \Gamma$    a covering
map with holonomy group $\ZZ$ (in fact, abelian is sufficient)
with the following properties:
\begin{enumerate}
\item[(G1)]
 $\Sigma_j$ is a sequence of embedded minimal multi-valued
 (normal exponential) graphs
 over $\Gamma$.
\item[(G2)]  There is a sequence $K_1 \subset K_2 \subset \dots \subset
\hat{\Gamma}$ of open domains exhausting $\hat{\Gamma}$ and
functions $u_j: K_j \to \RR$ with
 \begin{equation}  \label{e:hh1}
    |u_j| + |\nabla u_j| \leq 1/j \, ,
 \end{equation}
 so that there is a bijection from   from $K_j$ to $\Sigma_j$
 given by
 \begin{equation}
     x   \, \to  \,  \Pi (x) + u_j (x) \, \nn_{\Gamma} (\Pi (x))
      \, .    \label{e:hh2}
 \end{equation}
\end{enumerate}
The condition (G2) says that the $\Sigma_j$'s can be thought of as
one to one graphs over the domains $K_j$ in the cover
$\hat{\Gamma}$.

\begin{Cor}     \label{c:collstab}
If $\hat{\Gamma} \to \Gamma$ satisfies (G1) and (G2), then the
surface $\Gamma$ is stable.
\end{Cor}

\begin{proof}
By assumption, the holonomy group $\Lambda$ is cyclic and, thus,
 has sub-exponential growth.  Therefore,
to apply Lemma \ref{l:punstabIII}, we must show that the cover
$\hat{\Gamma}$ is stable.  We will prove the stability of
$\hat{\Gamma}$ by constructing a positive solution $w$ of the
Jacobi equation on $\hat{\Gamma}$.

First,   since the holonomy group $\Lambda$ is abelian, Lemma
\ref{l:simcov} implies that it acts as the deck group of
$\hat{\Gamma}$.

Next, define a sequence of subsets $\tilde{K}_j^{o} \subset
\hat{\Gamma}$ by
\begin{equation}    \label{e:wjwhere}
    \tilde{K}_j^{o} = \{ x \in K_j \, | \, h(1)(x) \in K_j \} =
    K_j \cap h(1)^{-1} (K_j) \, ,
\end{equation}
where $h(1) \in \Lambda$ is the generator of the infinite cyclic
subgroup $\Lambda = \ZZ$.  Fix a point $p \in \tilde{K}_1$ and let
$\tilde{K}_j$ be the connected component of $\tilde{K}_j^{o}$
containing $p$.

We will need below that the $\tilde{K}_j$'s are nested, open,
connected sets that exhaust $\hat{\Gamma}$.  The only point to
check is that they exhaust $\hat{\Gamma}$.  To see this, suppose
that $y \in \hat{\Gamma}$ and choose a path $\sigma : [0,1] \to
\hat{\Gamma}$  from $p$ to $y$.   Since the $K_j$'s are open and
exhaust $\hat{\Gamma}$, the compact set $\sigma  ([0,1]) \cup
h(1)( \sigma ( [0,1]))$ is entirely contained in some $K_{j}$ for
$j$ sufficiently large and, in particular, $\sigma ([0,1]) \subset
\tilde{K}_{j}$.

 Given $x \in \tilde{K}_j$,
both $x$ and $h(1)(x)$ are in $K_j$ and, therefore, we can define
functions $w_j$  on $\tilde{K}_j$ by
\begin{equation}    \label{e:defwja}
    w_j (x) = u_j ( h(1) (x)) - u_j (x) \, .
\end{equation}
Since the bijection   \eqr{e:hh2} takes $x$ and $h(1)(x)$ to
distinct points in the embedded surface
  $\Sigma_j$ and  these distinct points have the same projection to
  $\Gamma$, we conclude that
\begin{equation}
    w_j (x) \ne 0 \, .
\end{equation}
Therefore, we may as well assume that $w_j$ is positive on the
connected set $\tilde{K}_j$. Since $u_j$ and $|\nabla u_j|$ are
going to zero by \eqr{e:hh1}, a standard calculation (cf. lemma
$2.4$ in \cite{CM4}) gives that $u_j$ almost satisfies the Jacobi
equation.{\footnote{Precisely, $\Delta u_j + |A|^2 \, u_j = Q
(u_j)$, where the nonlinear term $Q (u_j)$ is at least quadratic
in $u_j$ and $\nabla u_j$.}} Likewise, the positive function $w_j$
is almost a solution of the Jacobi equation. In particular, if we
define normalized functions
\begin{equation}
    \tilde{w}_j = \frac{w_j}{w_j (p)} \, ,
\end{equation}
then a subsequence of the  $\tilde{w}_j$'s converges
  to a positive solution $w$ of the Jacobi equation on
  $\hat{\Gamma}$
and, thus, $\hat{\Gamma}$ is stable.
\end{proof}

\subsection{A Bernstein theorem for {\underline{incomplete}} surfaces}

The results of the previous subsections will be used show that
certain {\underline{incomplete}} minimal surfaces must be stable.
We will next prove a Bernstein theorem showing that such a stable
surface $\Gamma$ must then be flat, as long as it is ``complete
away from a single point.''  This generalizes the well-known
Bernstein theorem for complete stable surfaces of \cite{FiSc},
\cite{DoPe}.

 More precisely, we will
assume that the closure $\barga$ of $\Gamma$ is equal to the union
of $\Gamma$ and a single point. Recall that the closure $\barga$,
defined in \eqr{e:closureaa}, is given by
\begin{equation}
    \barga = \bigcup_{r} \overline{ \cB_r (x_{\Gamma}) } \,
    .
\end{equation}
  The flatness of such a
$\Gamma$ follows from an argument of Gulliver and Lawson,
\cite{GuLa}; for completeness, we recall this in the next lemma.

\begin{Lem}     \label{l:punstab}
Suppose that $\Gamma \subset \RR^3$ is a  connected stable minimal
surface without boundary and with trivial normal bundle.  If
\begin{equation}
    \barga \setminus \Gamma = \{ 0 \} \, ,
\end{equation}
  then $\Gamma$ is   a (punctured) plane.
\end{Lem}

\begin{proof}
 We will use an argument of Gulliver and
Lawson, \cite{GuLa}, to conformally change the metric $ds^2$ on
$\Gamma$ so that:
\begin{enumerate}
\item The universal cover $\Gamma_{U}$ of $\Gamma$ is complete in the new metric $d \tilde s^2$.
\item The operator $\tilde L =
\tilde \Delta - 2 \tilde K$ is non-negative on $\Gamma_U$; i.e.,
if $\phi$ is any compactly supported function on $\Gamma_U$, then
\begin{equation}
    \int \phi \, \tilde{L} \, \phi \leq 0 \, .
\end{equation}
Note that the sign convention here may be the opposite of what one
would expect.
\end{enumerate}
Once we have done this, it follows from \cite{FiSc} that $(
\Gamma_{U} , \, d \tilde s^2)$ is conformal to $\RR^2$ with the
standard flat metric. Translating back to the original metric
$ds^2$ will then imply that the original $\Gamma$ was flat.

Following \cite{GuLa},  we make the conformal change of metric
\begin{equation}
    d\tilde s^2 = \frac{ds^2}{|x|^2} \, .
\end{equation}
Since the covering map from $\Gamma_U$ to $\Gamma$ is an
immersion, the metric $d\tilde s^2$ on $\Gamma$ pulls back to give
a metric on $\Gamma_U$; we will also use $d\tilde s^2$ to denote
this pull back metric.
 It follows immediately that $\Gamma_{U}$ is complete
in the new metric $d\tilde s^2$.  Set $\tilde L = \tilde \Delta -
2 \tilde K$ where the Laplacian $\tilde \Delta$ and the curvature
$\tilde K$ are computed with respect to the metric $d \tilde s^2$.
Corollary $2.13$ in \cite{GuLa}{\footnote{Note that our operator
$L$ has the opposite sign convention from the operator $L_2$ in
\cite{GuLa}.}} gives that
\begin{equation}    \label{e:gl}
    \tilde L = |x|^2 \, L - 4 (1- |\nabla \, |x||^2) \, .
\end{equation}
Combining \eqr{e:gl} with the stability of $\Gamma_{U}$ gives for
any compactly supported $\phi$ that
\begin{equation}    \label{e:stai}
    \int_{\Gamma_{U}} \phi \, \tilde L \phi \, d\tilde \mu =
    \int_{\Gamma_{U}} \phi \, \left( |x|^2 \, L \phi - 4 (1- |\nabla \, |x||^2) \phi \right) \, |x|^{-2} \, d \mu
    \leq \int_{\Gamma_{U}} \phi \, L \phi \, d \mu \leq 0 \, ;
\end{equation}
that is, the operator $\tilde L$ is non-negative on $\Gamma_{U}$
with the complete metric $d\tilde s^2$.
 However,  theorem $2$ in \cite{FiSc} states that, for
 {\underline{any}} complete surface conformal to the
 disk, the intrinsically defined operator
 $\Delta - 2 K$ must
 be {\underline{negative}}.{\footnote{Note that \cite{FiSc} does not assume a
sign on the curvature $K$.}}  Therefore, since the plane is the
only other possible conformal type,
 we conclude that $(\Gamma_{U} , \, d \tilde s^2)$ -- and hence also $( \Gamma_{U} ,
\, d  s^2)$ -- is conformally equivalent to $\RR^2$. In
particular, there is a sequence of compactly-supported logarithmic
cutoff functions $\phi_j$ defined on $\Gamma_{U}$ with
\begin{align}
     \phi_j  \leq \phi_{j+1}    &{\text{ for every $j$ and }}
    \phi_j (x)  \to 1  {\text{ for every }} x \in
\Gamma_{U} \, . \\
    \lim_{j \to \infty} \, \int_{\Gamma_U} |\tilde \nabla \phi_j |^2 \, d \tilde \mu &=
        \lim_{j \to \infty} \, \int_{\Gamma_U} |\nabla \phi_j |^2 \, d \mu =  0 \, .
\end{align}
Using the functions $\phi_j$ in the stability inequality for $L$
on $\Gamma_U$ gives
\begin{equation}
    -2 \int_{\Gamma_{U}} K \, \phi_j^2 \, d \mu \leq    \int_{
    \Gamma_{U}} |\nabla \phi_j|^2 \, d\mu \to 0 \, .
\end{equation}
Since $K \leq 0$ and the functions $\phi_j$ go to $1$, we conclude
that $\Gamma_{U}$ is flat.  This completes the proof.
\end{proof}

\section{An extension of \cite{CM10}}       \label{s:aC}

In the current paper, we need slight modifications of several
results in \cite{CM10}.  We will give these results in this
appendix and explain whatever modifications are needed for their
proofs.

\subsection{Chord-arc bounds for ULSC surfaces}

The next lemma  extends the chord-arc bounds of \cite{CM10} from
disks to ULSC surfaces.

\begin{Lem}     \label{l:cyulsc}
Given a constant $r$, there exists $R > r$ so that if $\Sigma$ is
an embedded minimal surface with  $\cB_R (x_0) \subset \Sigma
\setminus \partial \Sigma$ and
\begin{equation}    \label{e:1ulsc}
    \cB_{1}(x)   {\text{ is a disk  for each $x \in \cB_R(x_0)$ }}
    \, ,
\end{equation}
then the connected component of $B_{r}(x_0) \cap \cB_R (x_0)$
containing $x_0$ has boundary in  $\partial B_{r}(x_0)$.
\end{Lem}

\begin{proof}
The proof follows the proof of lemma $2.23$ given in appendix B of
\cite{CM10} with one modification (the statement of  lemma $2.23$
from \cite{CM10} is recalled below in Lemma \ref{l:ca}). The
difference is that lemma $2.23$ assumed a curvature bound and used
this to show that two disjoint intrinsic balls whose centers were
close (in $\RR^3$) could be written as graphs over each other.  In
the current case, the required curvature bound is not assumed but
rather comes from the intrinsic version of the one-sided curvature
estimate (corollary $0.8$ in \cite{CM10}).
\end{proof}

\subsection{Chord-arc and area bounds for surfaces with bounded
curvature}    \label{ss:aC}

We also needed following lemma from \cite{CM10} which gives
chord-arc bounds for surfaces with bounded curvature{\footnote{Of
course, any surface with bounded curvature is also ULSC and is
therefore already covered by Lemma \ref{l:cyulsc}.  The usefulness
of Lemma \ref{l:ca} is that it makes the dependence very
precise.}}:

\begin{Lem}  \label{l:ca}  (lemma $2.23$ in \cite{CM10}.)
There exists $C_0 > 1$ so that given a constant $C_a$, we get
another constant $C_b$ such that the following holds:

\noindent If $\Sigma \subset \RR^3$ is an embedded minimal surface
with $0 \in \Sigma \subset B_{C_0 \, R}$ and $\partial \Sigma
\subset
\partial B_{C_0 \, R}$  and in addition
\begin{equation}        \label{e:cacb}
    \sup_{ B_{C_0 \, R} \cap \Sigma } |A|^2 \leq C_a \, R^{-2} \, ,
\end{equation}
then the component $\Sigma_{0, R}$ of $B_R  \cap \Sigma$
containing $0$ satisfies
\begin{equation}    \label{e:first}
    \Sigma_{0, R} \subset \cB_{C_b \, R}(0) \, .
\end{equation}
In particular, we also get a constant $C_c$ (depending only on
$C_a$) so that
\begin{equation}    \label{e:second}
    \Area \, (\Sigma_{0, R}) \leq   C_c \, R^2  \, .
\end{equation}
\end{Lem}

\begin{proof}
The first claim \eqr{e:first} follows precisely from the proof of
lemma $2.23$ in \cite{CM10} that is given in appendix B in
\cite{CM10}.{\footnote{We should point out that we have slightly
modified the statement of  lemma $2.23$ from \cite{CM10}; in
particular, the statement in \cite{CM10} assumes that $\Sigma$ is
a disk. However, this was not used in the proof of the lemma given
in appendix B in \cite{CM10}.}}

  Since
$|A|^2$ is bounded on $\Sigma_{0,R}$ by assumption, \eqr{e:first}
and  standard comparison theorems give the area bound
\eqr{e:second}.
\end{proof}

A key point in Lemma \ref{l:ca} is that the constant $C_0$ does
not depend on the constant $C_a$ in the curvature bound
\eqr{e:cacb}.

\section{Estimates for stable surfaces}  \label{s:aD}

Throughout this section, $\Gamma$ will be a stable surface with
connected ``interior boundary'' $\gamma$.  We will use $\an_r
(\gamma)$ to denote the intrinsic tubular neighborhood of radius
$r$ about a curve $\gamma$, i.e., \begin{equation}
    \an_r
(\gamma) = \{ x \in \Gamma \, | \, \dist_{\Gamma} (x , \gamma) < r
\} \, .
\end{equation}
Similarly,
 we will write
$\an_{s,t}(\gamma)$ for the ``annulus''
$\an_{t}(\gamma) \setminus
\an_{s}(\gamma)$.

The main result of this appendix is the following ``stable graph''
proposition.  This proposition shows that a stable embedded
minimal surface with a single interior boundary curve $\gamma$ and
an area bound near $\gamma$ is graphical away from its boundary.

\begin{Pro} \label{l:lext0}
Given a constant $C$, there exists $\omega > 1$ so that if $\Gamma
\subset B_R$ is a  stable embedded minimal surface whose
``interior boundary'' $\partial \Gamma \setminus
\partial B_R$ is a simple closed curve $\gamma \subset B_4$
satisfying
\begin{align}    \label{e:assume0}
    \Area ( \an_2 (\gamma) )  \leq C \, ,
\end{align}
then each component of $B_{R/\omega} \cap \Gamma \setminus
B_{\omega}$ is a graph with gradient bounded by one.
\end{Pro}

\subsection{The regularity of the distance function to the
interior boundary}

In   proving the proposition, we will need some basic results on
the level sets of the distance function to an interior boundary
curve.  Before stating these results, it will be helpful to recall
the Gauss-Bonnet theorem with corners and set the notation.

The Gauss-Bonnet theorem with corners implies that a surface
$\Sigma$ with piecewise smooth boundary $\partial \Sigma$
satisfies
\begin{equation}        \label{e:gbcorn}
    \int_{\partial \Sigma} k_g + \int_{\Sigma} K_{\Sigma} + \sum \alpha_i
    = 2 \pi \, \chi (\Sigma) \, .
\end{equation}
Here $K_{\Sigma}$ is the Gauss curvature of $\Sigma$, $\chi
(\Sigma)$ is its Euler characteristic, and $k_g$ is the geodesic
curvature of $\partial \Sigma$.  The sign convention of $k_g$ is
such that it is {\underline{positive}} on the boundary of the unit
disk in the plane.  Finally,  the $\alpha_i$'s are the ``jump
angles'' at the corners of $\partial \Sigma$; see Figure
\ref{f:jump}.  By convention, $\alpha_i$ is positive at a corner
where $\Sigma$ is locally convex.  For instance, on each corner of
a square, $\alpha_i$ is $\pi/2$.

\begin{figure}[htbp]
 \setlength{\captionindent}{20pt}
    \begin{minipage}[t]{0.5\textwidth}
    \centering\input{gen71.pstex_t}
    \caption{The jump angle $\alpha_i$ at a corner.}
    \label{f:jump}
    \end{minipage}\begin{minipage}[t]{0.5\textwidth}
    \centering\input{gen72.pstex_t}
    \caption{The level sets $S(t)$ of the distance function to a curve.}
    \label{f:soft}
    \end{minipage}
\end{figure}

\vskip2mm

\vskip2mm The next lemma of  Shiohama and Tanaka contains the main
results that we will need
 (cf. the proof of theorem $1$ in
\cite{Ro}):

\begin{Lem} \label{p:shta}
\cite{ShTa1}, \cite{ShTa2} Suppose that $\Gamma$ is a complete
noncompact oriented surface whose boundary $\partial \Gamma$ is a
smooth simple closed curve. The set
\begin{equation}
    S(t) = \{ x \in \Gamma \, | \, \dist_{\Gamma} ( x , \partial \Gamma
) = t \}
\end{equation}
 satisfies the following properties:
\begin{enumerate}
\item[($\star$1)] For almost every $t$, the set $S(t)$ is a finite
union of piecewise smooth curves with length $\ell (t)$. Moreover,
the ``jump angle'' $\alpha_i (t)$ at each corner is
{\underline{negative}} and always between $-\pi/2$ and $0$; let
$-\theta_i(t)$ denote this negative angle at the $i$-th
corner.{\footnote{This definition of $\theta_i(t)$ is chosen for
consistency with \cite{Fa} and \cite{Ha}. Note that each
$\theta_i(t)$ is positive.}}
 \item[($\star$2)]  For almost
every $t$, the derivative $\ell'(t)$ exists and satisfies
\begin{equation}        \label{e:fvgb}
    \ell' (t)  = \int_{S(t)} k_g - \sum_{i}  \tan
    (\theta_i(t))
    \leq \int_{S(t)} k_g - \sum_{i}  \theta_i(t) \,
    .
\end{equation}
The key for the inequality in \eqr{e:fvgb} is that each
$\theta_i(t)$ is between $0$ and $\pi / 2$ by ($\star$1).

  Notice that the right-hand side of
\eqr{e:fvgb} is exactly the boundary term corresponding to $S(t)$
in the Gauss-Bonnet formula with corners.{\footnote{Unfortunately,
the convention here is that $\alpha_i = - \theta_i (t)$.}}
 \item[($\star$3)] Given any $s
> r \geq 0$, we get
\begin{equation}
    \ell (s) - \ell (r) \leq \int_r^s \ell'(t) \, dt \, .
\end{equation}
\item[($\star$4)] The area of the ``annulus'' $\an_{r,s}(\partial
\Gamma) = \{x \in \Gamma \, | \,  r \leq \dist_{\Gamma} ( x ,
\partial \Gamma ) < s \}$ is
\begin{equation}
    \Area \, ( \an_{r,s}(\partial \Gamma) ) = \int_r^s \ell(t) \, dt \, .
\end{equation}
\end{enumerate}
\end{Lem}

\begin{Rem}
The papers \cite{ShTa1} and \cite{ShTa2}  extend earlier results
of Fiala, \cite{Fa}, for analytic surfaces and Hartman, \cite{Ha},
for simply connected surfaces. Since our surfaces are minimal in
$\RR^3$ and, thus, analytic, the classical results of Fiala could
be applied here. However, it is useful not to require analyticity
so that the results easily generalize to local ones in a
Riemannian $3$-manifold.

 The claim ($\star$1) was
proven in \cite{ShTa1}, while the claims ($\star$2), ($\star$3),
and ($\star$4) appear in \cite{ShTa2}. Note also that ($\star$4)
follows from the coarea formula.  We should note that the formula
 \eqr{e:fvgb} does not appear explicitly in
\cite{ShTa2},  but is implicit there and can also be found in
section $9.6$ of \cite{Fa}.
\end{Rem}

We will need two additional properties of the level sets $S(t)$
that hold if in addition $\Gamma$ is stable:
\begin{enumerate}
\item[($\star$5)] There is a constant $C_g$ so that if $\Gamma$ is
embedded and {\underline{stable}}, then we get the upper bound
\begin{equation}    \label{e:starfive}
    \sup_{S(t)} k_g \leq C_g \, t^{-1} \, .
\end{equation}
Recall that given our sign convention for $k_g$, \eqr{e:starfive}
means that   $S(t)$ cannot be ``too convex'' when it is thought of
as part of the boundary of $\an_t(\gamma)$.
 \item[($\star$6)] There is a constant $\epsilon_g$ so
that if $\Gamma$ is embedded and {\underline{stable}} and $\sigma
\subset S(t)$ is a closed curve with
\begin{equation}
    \Length \, (\sigma) \leq \epsilon_g \, t \, ,
\end{equation}
 then $\sigma$ bounds a disk $\Gamma_{\sigma} \subset \Gamma$ and
 $\Gamma_{\sigma} \subset \an_{3t/4, \, 4t/3}$.
\end{enumerate}

\begin{proof}
(of ($\star$5) and ($\star$6).)
 The upper bound ($\star$5) follows immediately from a barrier argument
 using
standard comparison theorems and the interior curvature estimate
for stable surfaces.  Namely, suppose that $p \in S(t)$ is a
smooth point.  Let
\begin{equation}
    \gamma_p:[0,t] \to \Gamma
\end{equation}
 be a minimizing geodesic connecting $\partial \Gamma$ to $p$.
The triangle inequality then implies that $S(t)$ does not
intersect the interior of the geodesic ball $\cB_s( \gamma_p(t-s)
)$ for any $s$ between $0$ and $t$.   Standard comparison theorems
and the curvature estimate for stable surfaces then give some $a >
0$ and $C_g$ so that $\partial \cB_{at}( \gamma_p(t-at) )$ is a
smooth curve with geodesic curvature at most $C_g \, t^{-1}$.
Since the $p$ is in the boundaries of these balls and is a smooth
point of $S(t)$, we conclude that the geodesic curvature of $S(t)$
at $p$ is also at most $C_g \, t^{-1}$.{\footnote{This argument
also shows that the jump angles at the corners of $S(t)$ are
negative as claimed in ($\star$1).}}
 Note that we
do not claim a lower bound for $k_g$ (in fact, easy examples show
that $k_g$ can go to $-\infty$; see \cite{Ha}).

To see ($\star$6), fix a point $p \in \sigma$ and note that the
entire curve $\sigma$ is contained in the intrinsic geodesic ball
$\cB_{\epsilon_g \, t}(p)$ and this ball stays away from $\partial
\Gamma$.  Taking $\epsilon_g$ small, the interior curvature
estimates for stable surfaces imply that $\cB_{\epsilon_g \,
t}(p)$ is a graph over some plane.  In particular, the curve
$\sigma$ is contractible in $\cB_{\epsilon_g \, t}(p)$, giving the
desired $\Gamma_{\sigma}$.
\end{proof}

\begin{Rem}
We will actually use a very slight generalization of these
results.  Namely, in applications, $\Gamma$ will not be complete,
but will rather be allowed to have other boundary components. This
does not matter since we will always work with level sets $S(t)$
where $t$ is less than the distance to any of the other boundary
components.  It's easy to see that the above results extend to
this case.
\end{Rem}

\subsection{The proof of the ``stable graph'' proposition}

The key point for proving Proposition \ref{l:lext0} will be to
show that $\Gamma$ has quadratic area growth.  This quadratic area
estimate formally follows from the argument in \cite{CM2}, but we
need the results of the previous subsection to deal with technical
difficulties that arise from the lack of regularity of the level
sets of the distance function.

\begin{proof}
(of Proposition \ref{l:lext0}.) The key point is to prove that the
intrinsic annuli $\an_r (\gamma)$ have quadratic area growth:
\begin{equation}    \label{e:qag}
    \Area \, ( \an_r (\gamma) ) \leq C_1 \, r^2 + C \, ,
\end{equation}
where the constant $C_1$ depends on the constant $C$ in
\eqr{e:assume0}.  Once we have \eqr{e:qag}, the lemma follows
easily from the proof of lemma II.1.34 in \cite{CM5}.  For the
reader's convenience, we will sketch the proof of the lemma
assuming \eqr{e:qag} next:
\begin{quote}
First, \eqr{e:qag} allows
us to use a logarithmic cutoff function to get sub-annuli with
small total curvature. Since these sub-annuli have small total
curvature and
 are stable, the mean value inequality gives a   small
 scale-invariant pointwise curvature estimate.
  Here the scale refers roughly to the
 distance to $\gamma$.  In particular, integrating this curvature bound implies that
 each component of a level set of
 the distance to $\gamma$ is itself a graph over (a curve
 in) {\underline{some}} plane.  Moreover, proposition $1.12$ in \cite{CM9} uses the
fact that the Gauss map is conformal to piece these together and
get a graph over one fixed plane, as desired.
\end{quote}

It remains therefore to establish \eqr{e:qag}. Note that
proposition II.1.3 in \cite{CM5} gives \eqr{e:qag} directly under
the additional assumption that $\Gamma$ is a topological annulus.
We will see that the general case follows similarly if we also use
the regularity of the length of level sets of the distance
function from $\gamma$  given by Lemma \ref{p:shta}.

\vskip2mm {\underline{The proof of \eqr{e:qag}.}} There are two
steps in the proof of \eqr{e:qag}:
\begin{enumerate}
\item The stability inequality allows us to bound the total
curvature in terms of the energy of a cutoff and this in turn is
bounded by the area.
\item The area growth is always controlled by the total curvature.
This follows easily from Gauss-Bonnet when the exponential map is
smooth but holds more generally by Lemma \ref{p:shta}.
\end{enumerate}

\vskip2mm \noindent{\underline{Step (1)}}:
 Set   $\dd (\cdot) = \dist_{\Gamma} ( \gamma , \cdot)$ and
define a (radial) cut-off function $\phi$ by
\begin{equation}
 \phi =
\begin{cases}
  \dd &\hbox{ on } \an_{1}(\gamma) \, ,    \\
  (r -\dd) / (r -1) &\hbox{ on } \an_{1,r}(\gamma)\, ,\\
  0 &\hbox{ otherwise }  \, . \\
\end{cases}
\end{equation}
 By the stability inequality applied to $\phi$, we get
\begin{align}   \label{e:intcut}
    \int_{\an_{1,r}(\gamma)}|A|^2\,[(r-\dd)/(r -1)]^2
    &\leq \int |A|^2\,\phi^2 \leq   \int |\nabla \phi|^2
        \notag \\
    &\leq  \Area\,(\an_{1}(\gamma))+(r -1)^{-2}\,
    \Area\,(\an_{1,r}(\gamma))   \, .
\end{align}
If we set $K(s)= \int_{\an_{1,s}(\gamma)}|A|^2$, then  the coarea
formula gives
\begin{equation}    \label{e:coara}
    K(s)= \int_{\an_{1,s}(\gamma)}|A|^2 = \int_1^s \, \left[ \int_{ \{ \dd =
        t \} } |A|^2 \right] \, dt \, .
\end{equation}
In particular,   we can integrate by parts twice to get
\begin{align} \label{e:intcut2}
    &2\,(r-1)^{-2}\int_1^{r} \int_1^t K(s)ds\,dt
    \leq 2/(r-1) \int_1^{r} K(s) (r-s)/(r -1)ds\notag \\
    &\quad \leq \int_{1}^{r}K'(s)  \,((r-s)/(r -1))^2 ds
      \leq  \Area\,(\an_{1}(\gamma))+(r -1)^{-2}\,
    \Area\,(\an_{1,r}(\gamma))    \, ,
\end{align}
where the last inequality is \eqr{e:intcut}.

\vskip2mm \noindent {\underline{Step (2)}}: We will now use Lemma
\ref{p:shta} to estimate the area by the total curvature. Set
$\ell ( t)$ equal to the length of the level set $\{ \dd = t \}$.
The key will be to prove the following estimate for $\ell (t)$ for
$t \geq 1$:
\begin{equation}    \label{e:ellte}
    \ell (t) \leq C_2 \, (1+t) + \frac{1}{2} \,
    \int_1^t \, \int_{\an_{1,s}(\gamma)} |A|^2 \, \, ds
    = C_2 \, (1+t) + \frac{1}{2} \,
    \int_1^t \, K(s) \, ds
    \, \,  ,
\end{equation}
where $C_2$ depends only on the constant $C$ in \eqr{e:assume0}.

\vskip2mm \noindent {\underline{The proof of the proposition
assuming \eqr{e:ellte}}}: Integrating the length bound
\eqr{e:ellte} gives the area bound
\begin{align}   \label{e:intcut3}
    \Area\,(\an_{r}(\gamma)) &\leq    \Area\,(\an_{1}(\gamma))
    + \int_{1}^{r} \ell (t) \, dt \notag \\
    &\leq  C
     + C_2 \, r + C_2 \, r^2 / 2
 +
 \int_{1}^{r} \, \int_1^t \frac{K(s)}{2} \, ds \,dt
      \, .
\end{align}
 Combining \eqr{e:intcut2} and
\eqr{e:intcut3} gives the needed bound \eqr{e:qag}.

\vskip2mm \noindent {\underline{The proof of  \eqr{e:ellte}}}:  We
will prove \eqr{e:ellte} by integrating a bound on $\ell' (t)$.
There will be two steps; namely, bounding $\ell'(t)$ and then
finding some value of $t$ where $\ell (t)$ is bounded (this is
where we will integrate the bound on $\ell'(t)$ from).

\vskip2mm \noindent {\underline{The bound on $\ell'(t)$}}: Roughly
speaking, we will bound $\ell'(t)$ in terms of the total curvature
by using the Gauss-Bonnet theorem in the  ``annulus''
$\an_{t_0,t}$ for a specific choice of $t_0$. Recall that
($\star$2) implies that $\ell'(t)$ is bounded by the Gauss-Bonnet
terms corresponding to $S(t)$.   To get the desired upper bound,
we will need to control the contributions from the geodesic
curvature of the ``inner boundary'' $S(t_0)$ as well as the Euler
characteristic of the ``annulus''. We will do this next.

First, the area bound \eqr{e:assume0} and ($\star$4) imply that
there must exist some  $t_0 \in (1/3 , 2/3)$ with
\begin{equation}    \label{e:elltz}
    \ell (t_0) \leq 3 \, \Area \, (\an_{1/3,2/3}(\gamma)) \leq 3
    \, C \, .
\end{equation}
Moreover, by the regularity property ($\star$1), we may assume
that the level set $S(t_0)$ is a finite union of simple closed
piecewise smooth curves.  We will sort these curves into two
groups, depending on their length. Namely, let $\sigma^{long}_1 ,
\dots , \sigma^{long}_n$ be the components of $S(t_0)$ with
\begin{equation}    \label{e:defsii}
    \Length (\sigma^{long}_i ) \geq \epsilon_g / 3 \, ,
\end{equation}
where $\epsilon_g$ is given by ($\star$6). Let $\sigma^{short}_1 ,
\dots , \sigma^{short}_m$ be the remaining components.  Combining
\eqr{e:defsii} with the upper bound on the total length of
$S(t_0)$ from \eqr{e:elltz} immediately gives the bound
\begin{equation} \label{e:claimnc}
        n \leq n(C) \, ,
   \end{equation}
    where $n(C)$  depends only on the area bound
    \eqr{e:assume0}.{\footnote{We are not claiming a bound on the
    total number $m+n$ of components of $S(t_0)$.}}

We will not actually apply the Gauss-Bonnet theorem to all of
$\an_{t_0,t}(\gamma)$, but rather to the subset $\Gamma_t$  that
``sees the outer boundary'' $S(t)$.  To be precise, define
$\Gamma_t$ to be the union of all connected components of
$\an_{t_0,t}(\gamma)$ whose boundaries intersect $S(t)$; see
Figure \ref{f:defgammat}.

\begin{figure}[htbp]
 \setlength{\captionindent}{20pt}
    \centering\input{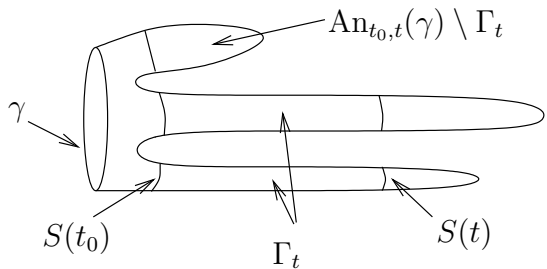}
    \caption{An example illustrating $\Gamma_t$ in a case
    where $\Gamma_t \ne \an_{t_0,t}(\gamma)$.}
    \label{f:defgammat}
\end{figure}
By construction, we have
\begin{equation}     \label{e:propsft}
     S(t) \subset \partial \Gamma_t {\text{ and }} \partial \Gamma_t \setminus
   S(t) \subset S(t_0) \, .
   \end{equation}
Consequently, combining the length bound \eqr{e:elltz} with the
pointwise geodesic curvature bound ($\star$5), we get
 a total (geodesic)
curvature bound for $\partial \Gamma_t \setminus
   S(t)$
   \begin{equation} \label{e:claimkc}
        \int_{\partial \Gamma_t \setminus S(t)} k_g \leq \ell (t_0) \, \sup_{S(t_0)} k_g \leq k(C) \, ,
   \end{equation}
   where $k(C)$  depends only on
   the area bound \eqr{e:assume0}. We should make two remarks about \eqr{e:claimkc}:
   \begin{itemize}
   \item The integration in \eqr{e:claimkc} is over only the smooth
   part of $\partial \Gamma_t \setminus S(t)$.
   \item  The sign convention on $k_g$ in \eqr{e:claimkc} is as part of the boundary of $\an_{t_0} (\gamma)$; this
   is the opposite as it would be as part of the boundary of
   $\Gamma_t$.  This is important later when we apply the
   Gauss-Bonnet theorem.
   \end{itemize}
The last ingredient that we will need to bound $\ell'(t)$ is a
bound on the Euler characteristic $\chi (\Gamma_t)$ that depends
only on  the area bound \eqr{e:assume0}.  This bound  follows
immediately from the bound \eqr{e:claimnc} on the number of long
components of $S(t_0)$ together with the following claim:
\begin{enumerate}
\item[(Claim)] For $t \geq 3/4$, each connected component of
$\Gamma_t$, i.e., each component of $\an_{t_0,t}(\gamma)$ whose
boundary touches $S(t)$, contains at least one long component
$\sigma^{long}_i$ in its boundary.
\end{enumerate}
The point  here is that the short components  of $S(t_0)$ are
contractible near $S(t_0)$, so $S(t)$ never sees them.  More
precisely, ($\star$6) implies that each $\sigma_i^{short}$ bounds
a disk
\begin{equation}    \label{e:gamdiski}
    \Gamma^{disk}_i \subset \an_{1/4 , \, 3/4}(\gamma) \, .
\end{equation}
Therefore, if $p$ is an arbitrary point in $S(t)$, then we know
that $p$ and $\gamma$ are in the same connected component of
\begin{equation}
    \overline{\an_t(\gamma)} \setminus \cup_{i=1}^m \sigma_i^{short}  \, .
\end{equation}
Note that we have used here that $\an_t (\gamma)$ is itself
connected.  Since $S(t_0)$ separates $\gamma = S(0)$ from $S(t)$,
we conclude that it must be $\cup_i \sigma^{long}_i$ that
separates $p$ and $\gamma$.  In particular, the component of
$\an_{t_0,t}$ with $p$ in its boundary also contains at least one
$\sigma^{long}_i$ in its boundary.  This completes the proof of
(Claim).

We can now  bound $\ell'(t)$ for $t \geq 3/4$. Namely, ($\star$2)
implies that $\ell'(t)$ is bounded by the Gauss-Bonnet integrand
along $S(t)$ so the Gauss-Bonnet theorem gives for almost every
$t$ that
\begin{equation}    \label{e:thrownaw}
    \ell '(t) \leq \int_{S(t)} k_g - \sum_{i} \theta_i (t) \leq
    \frac{1}{2} \, \int_{\Gamma_t} |A|^2 + 2 \pi
    \chi (\Gamma_t) +   \int_{\partial \Gamma_t \setminus S(t)} k_g \, .
\end{equation}
We have thrown away the angle contributions at the corners of
$\partial \Gamma_t \setminus S(t)$ in \eqr{e:thrownaw} since these
are all negative by ($\star$1). Since $\Gamma_t \subset \an_{1/3 ,
t}(\gamma)$, we can use interior curvature estimates for stable
surfaces and the area bound on $\an_1 (\gamma)$ to get
\begin{equation}    \label{e:antot}
    \int_{\Gamma_t} |A|^2 \leq \int_{\an_{1,t}(\gamma)} |A|^2 +
    \Area (\an_1 (\gamma)) \, \sup_{ \an_{1/3 , 1}(\gamma) } \, |A|^2
    \leq \int_{\an_{1,t}(\gamma)} |A|^2 + C_3 \, ,
\end{equation}
where $C_3$ depends only on the initial area bound
\eqr{e:assume0}.    Substituting the  above bounds into
\eqr{e:thrownaw}, we get for almost every $t \geq 3/4$ that
\begin{equation}    \label{e:ellpb}
    \ell '(t)  \leq \frac{1}{2} \, \int_{\an_{1,t}(\gamma)} |A|^2 +
C_4 \, ,
\end{equation}
where $C_4$ depends only on the initial area bound
\eqr{e:assume0}.

To complete the proof, use the area bound and ($\star$4) again to
find $t_1$ between $3/4$ and $1$ with $\ell (t_1) \leq 4 \, C$.
Given $t \geq 1$, we can then use
 ($\star$3) to integrate \eqr{e:ellpb}:
\begin{align}
    \ell (t) &\leq \ell (t_1) + \int_{t_1}^t \ell'(s) \, ds \notag
    \\
    &\leq 4 \, C + C_4 (t-t_1) + \frac{1}{2} \,  \int_{t_1}^t \, \int_{\an_{1,s}(\gamma)}
    |A|^2 \, ds \, .
\end{align}
This gives \eqr{e:ellte}, thus completing the proof.

\vskip2mm
 We should point out that we have actually shown only
that the components coming from the tubular neighborhood $\an_r
(\gamma)$ are graphs.  However, the other components are easily
also seen to be graphs by combining the curvature estimate and
embeddedness. Namely, any other component is intrinsically far
from the boundary and hence graphical over some plane.  By
embeddedness, these graphs do not cross, and we can take these
planes to be parallel.
\end{proof}

\section{Blowing up intrinsically on the scale of non--trivial
topology}   \label{s:aE}

The next lemma uses a standard blowup argument  to locate the
smallest scale of non-trivial topology:

\begin{Lem}     \label{l:l5.1}
Suppose that $\Sigma \subset \RR^3$ is a smooth minimal surface,
possibly with boundary $\partial \Sigma$.  If the ball $\cB_{5 \,
C_1 \, r_1}(y_0) \subset \Sigma$ is disjoint from $\partial
\Sigma$ for some $C_1 > 1$ and
\begin{equation}    \label{e:startf}
    \cB_{r_1}(0) {\text{ is not a disk}},
\end{equation}
 then there exists a sub--ball $\cB_{C_1 s}(y_1) \subset \cB_{4 \, C_1 \, r_1}(y_0)$ so
 that
\begin{align}    \label{e:b1}
    \cB_{4s}(y_1) &{\text{ is {\underline{not}} a disk}} \, , \\
    \cB_{s}(y) &{\text{ is  a disk for any }} y \in \cB_{C_1 s}(y_1) \,
    .  \label{e:b2}
\end{align}
\end{Lem}

\begin{proof}
After rescaling, we can assume that $r_1 = 1$. The lemma will
follow from a simple rescaling argument as in Lemma $5.1$ of
\cite{CM4}, except we define $F$ intrinsically on $\cB_{4 \,
C_1}(y_0)$ by
\begin{equation}
    F (x) = d^2 (x) \, \cali^{-2} (x) \, ,
\end{equation}
 where $\cali (x)$ is the injectivity radius of $\Sigma$ at
$x$  and
\begin{equation}
    d(x) = 4 \, C_1   - \dist_{\Sigma} (x,
y_0)
\end{equation}
is the distance to $\partial \cB_{4 \, C_1}(y_0)$.
 It follows that
 $F=0$ on $\partial \cB_{4 \, C_1}(y_0)$ and
 $ F(y_0) \geq 16 \, C_1^2$.  Also, since $\cB_{5 \, C_1}(y_0)$ is smooth, it follows that $F$ is bounded
 from above on $\cB_{4 \, C_1}(y_0)$.
   We can therefore choose a
 point $y_1$ where $F(y_1)$ is at least half of its supremum{\footnote{Note that we are
 not claiming that   $\cali$, or $F$, is continuous, so we do not know that it achieves its maximum.  However,
 since it is bounded, there must be points where it is at least half of its supremum.}}, i.e.,
 \begin{equation}   \label{e:b3}
        F(y_1) > 1/2 \, \sup_{\cB_{4 \,
        C_1}(y_0)} \, F  \, .
 \end{equation}
Set $s^2 = \cali^2 (y_1)/8$.

To see that \eqr{e:b1} holds, first note that $4\,s > \cali
(y_1)$. In particular, there must be two distinct geodesics,
$\gamma_1$ and $\gamma_2$,  contained in $\cB_{4s}(y_1)$ with
\begin{equation}
    \gamma_1 (0) = \gamma_2 (0) = y_1 {\text{ and }} \gamma_1
    (\cali (y_1)) = \gamma_2
    (\cali (y_1)) \, .
\end{equation}
Since $\Gamma$ has non--positive curvature, it follows immediately
from the Gauss-Bonnet theorem with corners   that
 the closed curve $\gamma_1 \cup \gamma_2$ cannot bound  a disk
in $\Gamma$, thus giving \eqr{e:b1}.{\footnote{The Gauss-Bonnet
theorem with corners implies that a disk $D$ has  $\int_{\partial
D} k_g + \int_D K + \sum \alpha_i = 2 \pi$, where $\partial D$ has
jump angles $\alpha_i$ at the corners.  In this case, both
integrals are non-positive and there are only two corners with
each contributing {\underline{less}} than $\pi$, so no such disk
can exist.}}

We will use \eqr{e:b3} twice to prove \eqr{e:b2}.  First, since $d
\geq d(y_1)/2$ on $\cB_{ \frac{d(y_1)}{2} }(y_1)$, \eqr{e:b3}
implies that
\begin{equation}
    \sup_{\cB_{ \frac{d(y_1)}{2}} (y_1)} \, \cali^{-2} \leq \frac{4}{d^2(y_1)} \,
    \sup_{\cB_{\frac{d(y_1)}{2}} (y_1) } \, F < \frac{8 \, F(y_1)}{d^2(y_1)} =  s^{-2}  \, \,
    ,
\end{equation}
so that $\cali > s$ on $\cB_{\frac{d(y_1)}{2}}(y_1)$.  Second,
using $F(y_0)$ as a lower bound for the sup of $F$ in \eqr{e:b3}
implies that
\begin{equation}
         8 \, C_1^2 \leq F(y_1) = \frac{d^2 (y_1)}{8 \, s^2} \, ,
 \end{equation}
so that
 $d(y_1)/2 > C_1 \, s$.
\end{proof}

\section{Minimal surfaces with a quadratic curvature bound in a
half--space}

\begin{figure}[htbp]
    \setlength{\captionindent}{20pt}
    \begin{minipage}[t]{0.5\textwidth}
    \centering\input{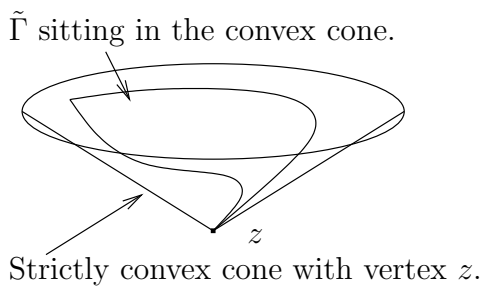}
    \caption{We prove Lemma \ref{l:oc2} by contradiction,
    so suppose that $\tilde \Gamma$ lies in a strictly mean convex cone.}
    \label{f:cvx}
    \end{minipage}
\end{figure}

  The next lemma deals with a minimal
surface  whose curvature blows up at most quadratically at a point
$z$ in its closure. The lemma   shows that the surface must come
arbitrarily scale-invariant close to {\underline{any}} plane
through $z$. Roughly speaking, this means that   the surface does
not lie
 in any strictly mean convex cone through $z$; see Figure \ref{f:cvx}.

To state the lemma precisely,  given a plane through $z$, we will
 define a scale-invariant function
$\beta (s)$ that measures how close in the sphere $\partial
B_{s}(z)$ a surface $\tilde \Gamma \subset \RR^3 \setminus  \{ z
\}$ comes to the  plane.  After a rotation, we can assume that the
plane is the horizontal plane $\{ x_3 = x_3 (z) \}$. Define the
function $\beta (s)$ by setting
\begin{equation}        \label{e:betas2}
    \beta (s) = \frac{ \inf_{ \partial B_{s}(z) \cap \tilde
    \Gamma} \, \, |x_3 - x_3 (z)| }{s}  \, .
\end{equation}
  The next lemma shows that the liminf of  $\beta (s)$ is zero, so
 that $\tilde \Gamma$  comes arbitrarily scale-invariant close to
 the plane as we approach $z$.

\begin{Lem}     \label{l:oc2}
Let $\tilde \Gamma \subset   \RR^3 \setminus  \{ z \}$ be a
minimal surface (embedded or not) with $z$ in its closure. Suppose
that, for each $\epsilon > 0$, each component of $\tilde \Gamma
\setminus B_{\epsilon}(z)$ is complete  and has boundary in
$\partial B_{\epsilon}(z)$.

 If there exist constants $r_0 > 0$ and $C$ so that
for $x \in B_{r_0}(z) \cap \tilde \Gamma$ we have
\begin{equation}    \label{e:ocqb}
    |A|^2 (x) \leq C \, |x-z|^{-2} \, ,
\end{equation}
 then the function $\beta (s)$ defined in \eqr{e:betas2} satisfies
\begin{equation}    \label{e:infbs22}
    \liminf_{s \to 0} \, \beta (s) = 0 \, .
\end{equation}
\end{Lem}

\begin{proof}
  We will prove \eqr{e:infbs22} by contradiction (see Figure \ref{f:cvx}), so suppose
that
\begin{equation}    \label{e:infbs2}
    \liminf_{s \to 0} \, \beta (s) = \beta_0  > 0 \, .
\end{equation}
 In particular, given any $\delta > 0$, equation \eqr{e:infbs2}
 implies that there exists $s_0>0$ so
that  $\beta (s) > \beta_0 -
    \delta$ for every $s < s_0$ and, hence, $B_{s_0}(z) \cap \tilde
    \Gamma$ lies inside a strictly mean convex  (double) cone:
\begin{equation}    \label{e:bigg1}
    B_{s_0}(z) \cap \tilde \Gamma \subset \{ |x_3 - x_3 (z)| > (\beta_0 -
    \delta ) \, |x-z| \} \, .
\end{equation}
On the other hand, \eqr{e:infbs2} also implies that there is an
$s< s_0/2$ with $\beta (s) < \beta_0 + \delta$ and, consequently,
there is a point $y_s \in \partial B_s(z) \cap \tilde \Gamma$
close to the strictly mean convex cone:
\begin{equation}    \label{e:bigg2}
    y_s \in \partial B_s(z) \cap \{ |x_3 - x_3 (z)| <  (\beta_0 + \delta) \, s
    \} \cap \tilde \Gamma \, .
\end{equation}
Note that \eqr{e:bigg1} and \eqr{e:bigg2} imply that the intrinsic
ball $\cB_{s/2} (y_s)$ stays inside, but comes close to, the
{\underline{strictly}} mean convex cone
\begin{equation}
\label{e:scc}
    \{ |x_3 - x_3 (z)| = (\beta_0 -
    \delta ) \, |x-z| \} \, .
\end{equation}
We will assume below that   $\delta < \beta_0$.

  We will see that \eqr{e:bigg1} and \eqr{e:bigg2} lead to a
contradiction for $\delta$ sufficiently small, thus proving
\eqr{e:infbs22}.

First, recall  that the quadratic curvature bound \eqr{e:ocqb}
gives an $\alpha > 0$ so that the component $\tilde \Gamma_{\alpha
\,s}$ of $B_{\alpha \, s}(y_s) \cap \tilde \Gamma$ containing
$y_s$ is a graph with gradient bounded by one  (see, e.g., lemma
$2.2$ in \cite{CM1}).
 After possibly reducing $\alpha$, we can therefore assume that
\begin{equation}
    \tilde \Gamma_{\alpha \,s} \subset \cB_{s/2} (y_s) \, .
\end{equation}
Since $\tilde \Gamma_{\alpha \,s}$ is connected and does not
intersect the (double) cone \eqr{e:scc}, it must be in one of the
two components of $\{ |x_3 - x_3 (z)| > (\beta_0 -
    \delta ) \, |x-z| \}$.  After possibly reflecting, we can
    assume that
\begin{equation}
    \tilde \Gamma_{\alpha \,s} \subset \{ x_3 - x_3 (z) >
(\beta_0 -
    \delta ) \, |x-z| \} \, .
\end{equation}
Define a function $f$ that vanishes on the cone   $\{ x_3 - x_3
(z) = (\beta_0 -
    \delta ) \, |x-z| \}$ by
setting
\begin{equation}
    f (x) = x_3 - x_3 (z) - (\beta_0 - \delta) \, |x - z| \, .
\end{equation}
Note that \eqr{e:bigg1} and \eqr{e:bigg2} imply that
\begin{equation}    \label{e:bigg3}
    0 \leq \inf_{ \tilde \Gamma_{\alpha \,s} } \, f \leq
    f(y_s) < 2 \, \delta \, s \, .
\end{equation}
 Using that
$\tilde \Gamma_{\alpha \,s}$ is minimal, we have that
\begin{equation}
    \Delta f   =   - (\beta_0 - \delta) \, \Delta |x - z|  < - \frac{\beta_0 -
    \delta}{|x-z|}
    \, .
\end{equation}
Define a function $g$ on $\tilde \Gamma_{\alpha \,s}$ by setting
\begin{equation}
    g = f +
    |x-y_s|^2 \, \frac{\beta_0 - \delta}{6s} \, .
\end{equation}
 Using that $|x-z| < 3s/2$ on
$\tilde \Gamma_{\alpha \,s}$, we get that $g$ is superharmonic
since
\begin{equation}
    \Delta g   <     - \frac{\beta_0 -
    \delta}{|x-z|} + 4 \, \frac{\beta_0
- \delta}{6s} < 0    \, .
\end{equation}
Therefore, the minimum of $g$ is achieved on $\partial \tilde
\Gamma_{\alpha \,s}$ and thus
\begin{equation}    \label{e:notbigg3}
    \min_{\partial \tilde
\Gamma_{\alpha \,s}} \left[ f + (\alpha \, s)^2 \, \frac{\beta_0 -
\delta}{6s} \right] = \min_{\tilde \Gamma_{\alpha \,s}} \, g <
g(y_s) = f(y_s) < 2 \, \delta \, s \, ,
\end{equation}
where the last inequality is from \eqr{e:bigg3}. Combining the
first inequality from \eqr{e:bigg3} and \eqr{e:notbigg3}   gives
\begin{equation}    \label{e:notbigg33}
    0 \leq \min_{\partial \tilde
\Gamma_{\alpha \,s}}   f  < 2 \, \delta \, s - \alpha^2 \, s  \,
(\beta_0 - \delta) / 6
  \, .
\end{equation}
This gives the desired contradiction for $\delta$ sufficiently
small.
\end{proof}

\end{document}